\documentclass[11pt, reqno]{amsart}

\usepackage[margin=0.9in]{geometry}
\usepackage[utf8]{inputenc}
\usepackage[T1]{fontenc}

\usepackage{subfiles}
\usepackage{comment}
\usepackage{float}
\usepackage{marginnote}
\usepackage{euscript}
\usepackage[dvipsnames]{xcolor}   		             		
\usepackage{graphicx}			
\usepackage{amssymb}
\usepackage{mathrsfs}
\usepackage{amsthm}
\usepackage{amsmath, amsrefs}
\usepackage[foot]{amsaddr}
\usepackage{mathtools}
\usepackage[cal=cm]{mathalfa}
\usepackage{stmaryrd}
\usepackage{upgreek}
\usepackage{bbm}
\usepackage[hypertexnames=false, hyperfootnotes=false, colorlinks=true,linktocpage=true, allcolors=highlight, bookmarksdepth = 2]{hyperref}
\usepackage{cleveref}
\usepackage{url}
\usepackage{float}
\usepackage[full]{textcomp}
\makeatletter
\newcommand{\citecomment}[2][]{\citen{#2}#1\citevar}
\newcommand{\citeone}[1]{\citecomment{#1}}
\newcommand{\citetwo}[2][]{\citecomment[,~#1]{#2}}
\newcommand{\citevar}{\@ifnextchar\bgroup{;~\citeone}{\@ifnextchar[{;~\citetwo}{]}}}
\newcommand{\citefirst}{\@ifnextchar\bgroup{\citeone}{\@ifnextchar[{\citetwo}{]}}}

\makeatother

\usepackage{tikz,tikz-cd,tikz-3dplot}
\usepackage{pgfplots}
\pgfplotsset{compat=1.18}
\usetikzlibrary{calc}
\usetikzlibrary{fadings}
\usetikzlibrary{decorations.pathmorphing}
\usetikzlibrary{decorations.pathreplacing}
\usetikzlibrary{patterns}
\usetikzlibrary{arrows,shadows,positioning, calc, decorations.markings, 
	hobby,quotes,angles,decorations.pathreplacing,intersections,shapes}
\usepgflibrary{shapes.geometric}
\usetikzlibrary{fillbetween,backgrounds}

\usepackage{xcolor}

\definecolor{highlight}{HTML}{3465a4}

\usepackage{anyfontsize}
\usepackage[lining]{libertine}%
\usepackage{courier}
\usepackage{csquotes}
\makeatletter
\renewcommand*\libertine@figurestyle{LF}
\makeatother
\usepackage[libertine,libaltvw,liby]{newtxmath}
\usepackage{microtype}

\usepackage{array}
\usepackage{slashbox}
\newcolumntype{H}{>{\setbox0=\hbox\bgroup}c<{\egroup}@{}}

\BibSpec{article}{%
+{}{\sc \PrintAuthors} {author}
+{,}{ } {title}
+{,}{ \textit} {journal}
+{,}{ } {volume}
+{}{ \IfEmptyBibField{journal}{}{\PrintDate{date}}} {transition}
+{,}{ } {pages}
+{,}{ } {status}
+{.}{} {transition}
+{}{ \, Preprint available at \tt} {eprint}
+{.}{} {transition}
}

\BibSpec{book}{%
+{}{\sc \PrintAuthors} {author}
+{,}{ } {title}
+{.}{ \textit} {series}
+{,}{ } {volume}
+{.}{ } {publisher}
+{,}{ } {place}
+{,}{ } {date}
+{.}{} {transition}
}

\BibSpec{incollection}{
+{} {\sc \PrintAuthors} {author}
+{,} { } {title}
+{,} { in \textit} {booktitle}
+{,}{ } {series}
+{,}{ } {volume}
+{}{ \PrintDate } {date}
+{,} { pp.~} {pages}
+{,}{ } {publisher}
+{ }{ Preprint available at \tt} {eprint}
+{.}{} {transition}
}

\BibSpec{inproceedings}{
+{} {\sc \PrintAuthors} {author}
+{,} { } {title}
+{,} { in \textit} {booktitle}
+{,} { pp.~} {pages}
+{,} { \PrintDateB} {date}
+{,} { available at \eprint} {eprint}
+{.} {} {transition}
}

\BibSpec{collection.article}{
+{} {\sc \PrintAuthors} {author}
+{,} { } {title}
+{,} { in \it \PrintConference} {conference}
+{,} { pp.~} {pages}
+{.} {\PrintBook} {book}
+{,} { \PrintDateB} {date}
+{.} {} {transition}
}

\BibSpec{innerbook}{
+{.} { \emph} {title}
+{.} { } {part}
+{:} { \emph} {subtitle}
+{.} { } {series}
+{,} { } {volume}
+{.} { Edited by \PrintNameList} {editor}
+{.} { Translated by \PrintNameList}{translator}
+{.} { \PrintContributions} {contribution}
+{.} { } {publisher}
+{.} { } {organization}
+{,} { } {address}
+{,} { \PrintEdition} {edition}
+{,} { \PrintDateB} {date}
+{.} { } {note}
}

\tikzset{
	commutative diagrams/.cd, 
	arrow style=tikz, 
	diagrams={>=stealth}
}
\tikzset{
	arrow/.pic={\path[tips,every arrow/.try,->,>=#1] (0,0) -- +(0,4pt);},
	pics/arrow/.default={triangle 90}
}
\tikzset{->-/.style={decoration={
			markings,
			mark=at position .6 with {\arrow{latex}}},postaction={decorate}}
}
\tikzset{
	c/.style={every coordinate/.try}
}
\usepackage{enumerate}
\usepackage{enumitem}
\setlist[enumerate,1]{%
	label=(\roman*)	,itemsep=0.3em
	}
\setlist[itemize]{itemsep=0.3em}
\setlist[description]{itemsep=0.3em,labelindent=0.5cm}

\newcommand{\mc}{\mathcal}
\newcommand{\mr}{\mathrm}
\newcommand{\mf}{\mathfrak}

\DeclareMathOperator{\Aut}{Aut}

\DeclareMathOperator*{\Tot}{Tot}
\DeclareMathOperator{\Rep}{Rep}

\DeclareMathOperator{\SL}{SL}
\DeclareMathOperator{\GL}{GL}

\newcommand{\define}[1]{\textbf{#1}}

\newcommand{\mgp}[1]{\mathsf{#1}} 
\newcommand{\mchar}[1]{\mathsf{#1}} 
\newcommand{\invariantsfont}[1]{\mathsf{#1}}
\newcommand{\modulifont}[1]{#1}
\newcommand{\Mbar}{\overline{\modulifont{M}}}

\newcommand{\Nbar}{\overline{\modulifont{N}}}

\newcommand{\vir}{\mathrm{virt}}

\newcommand{\GW}{\invariantsfont{GW}}
\newcommand{\PT}{\invariantsfont{PT}}
\newcommand{\BPS}{\Omega}
\newcommand{\MT}{\invariantsfont{MT}}
\newcommand{\refGWgb}[4]{\GW_{#1, #4}(#2, #3 )}
\newcommand{\refGWb}[3]{\GW_{#3}(#1, #2)}
\newcommand{\refGWdiscb}[3]{\GW^{\bullet}_{#3}(#1, #2)}
\newcommand{\refGWg}[3]{\GW_{#1}(#2, #3)}
\newcommand{\refBPSb}[3]{\BPS_{#3}(#1, #2)}
\newcommand{\refBPSbarb}[3]{\overline{\BPS}_{#3}(#1, #2)}
\newcommand{\refMTb}[2]{\MT_{#2}(#1)}
\newcommand{\refPTb}[3]{\PT_{#3}(#1 , #2)}
\newcommand{\refPTnb}[4]{\PT_{#1,#4}(#2 , #3)}

\newcommand{\kfield}{\mathbb{C}}
\newcommand{\Speck}{\operatorname{Spec}\kfield}

\newcommand{\Aone}{\mathbb{C}}
\newcommand{\Aaff}[1]{\mathbb{C}^{#1}}
\usepackage{xifthen}
\newcommand{\Gm}[1][]{%
	\ifthenelse{\isempty{#1}}%
	{\mathbb{C}^\times}
	{(\mathbb{C}^\times)^{#1}}
}

\newcommand{\intEquiv}[1]{{}^{#1}\!\!\!\int}
\newcommand{\ch}[1][]{\mathrm{ch}_{#1}}
\newcommand{\Chow}{\mathsf{CH}}
\newcommand{\hhh}{\mathsf{H}}
\newcommand{\HodgeLambda}[2][]{\Lambda (#2)}

\newcommand{\Tmax}{\mgp{T}_{\!\mr{M}}}
\newcommand{\TCY}{\mgp{T}_{\!\mr{CY}}}

\newcommand{\Kappa}{K}

\newcommand{\dd}{\mathrm{d}}
\newcommand{\ri}{\mathrm{i}}
\newcommand{\re}{\mathrm{e}}

\newcommand{\de}{{\partial}}

\newcommand{\bbM}{\mathbb{M}}
\newcommand{\bbN}{\mathbb{N}}

\newcommand{\bbL}{\mathbb{L}}
\newcommand{\bbZ}{\mathbb{Z}}
\newcommand{\bbR}{\mathbb{R}}
\newcommand{\bbC}{\mathbb{C}}
\newcommand{\bbP}{\mathbb{P}}
\newcommand{\bbE}{\mathbb{E}}

\newcommand{\bbQ}{\mathbb{Q}}

\newcommand{\cY}{\mathcal{Y}}
\newcommand{\cZ}{\mathcal{Z}}
\newcommand{\cO}{\mathcal{O}}

\newcommand{\cE}{\mathcal{E}}

\newcommand{\cC}{\mathcal{C}}
\newcommand{\DD}{\mathcal{D}}
\newcommand{\LL}{\mathcal{L}}

\newcommand{\cN}{\mathcal{N}}
\newcommand{\cW}{\mathcal{W}}
\newcommand{\cG}{\mathcal{G}}

\newcommand{\HH}{\mathcal{H}}

\newcommand{\cI}{\mathcal{I}}

\newcommand{\cV}{\mathcal{V}}

\newcommand{\MM}{\mathcal{M}}
\newcommand{\cM}{\mathcal M}
\newcommand{\cQ}{\mathcal Q}

\renewcommand{\l}{\left}
\renewcommand{\r}{\right}

\def\beq{\begin{equation}}                     %
	\def\eeq{\end{equation}}                       %
\def\bea{\begin{eqnarray}}                     
	\def\eea{\end{eqnarray}}
\def\bary{\begin{array}} 
	\def\eary{\end{array}} 
\def\ben{\begin{enumerate}} 
	\def\een{\end{enumerate}}
\def\bit{\begin{itemize}} 
	\def\eit{\end{itemize}}
\def\nn{\nonumber}

\theoremstyle{plain}
\newtheorem{thm}{Theorem}[section]
\newtheorem{lem}[thm]{Lemma}

\newtheorem{prop}[thm]{Proposition}
\newtheorem*{prop*}{Proposition}
\newtheorem{conj}[thm]{Conjecture}

\newtheorem*{conj*}{Conjecture}

\newtheorem{cor}[thm]{Corollary}
\newtheorem*{cor*}{Corollary}

\theoremstyle{definition}
\newtheorem{defn}[thm]{Definition}
\newtheorem{rmk}[thm]{Remark}
\newtheorem{notation}[thm]{Notation}
\newtheorem{example}[thm]{Example}

\theoremstyle{plain}
\newtheorem{innercustomthm}{Theorem}
\newenvironment{customthm}[1]
{\renewcommand\theinnercustomthm{#1}\innercustomthm}
{\endinnercustomthm}

\theoremstyle{plain}

\theoremstyle{plain}
\newtheorem{innercustomconj}{Conjecture}
\newenvironment{customconj}[1]
{\renewcommand\theinnercustomconj{#1}\innercustomconj}
{\endinnercustomconj}

\theoremstyle{definition}

\theoremstyle{plain}

\crefname{equation}{Eq.}{Eqs.}
\crefname{eqnarray}{Eq.}{Eqs.}
\crefname{algo}{algorithm}{algorithms}
\crefname{conj}{conjecture}{conjectures}
\crefname{lem}{lemma}{lemmas}
\crefname{thm}{theorem}{theorems}
\crefname{claim}{claim}{claims}
\crefname{rmk}{remark}{remarks}
\crefname{prop}{proposition}{propositions}
\crefname{section}{section}{sections}
\crefname{appendix}{appendix}{appendices}
\crefname{cor}{corollary}{corollaries}
\crefname{figure}{figure}{figures}
\crefname{table}{table}{tables}
\crefname{example}{example}{examples}
\crefname{prob}{problem}{problems}
\crefname{assm}{assumption}{assumptions}
\crefname{defn}{definition}{definitions}
\crefname{notation}{notation}{notations}
\crefname{speculation}{speculation}{speculations}
\crefname{construction}{construction}{constructions}
\crefname{observation}{observation}{observations}
\crefname{innercustomthm}{Theorem}{Theorems}
\crefname{innercustomconj}{Conjecture}{Conjectures}
\crefname{innercustomassumption}{assumption}{Assumption}
\crefname{innerpracticalresult}{practical result}{practical results}

\newcommand{\bra}{\left\langle}
\newcommand{\ket}{\right\rangle}

\crefname{equation}{Eq.}{Eqs.}
\crefname{eqnarray}{Eq.}{Eqs.}
\crefname{algo}{Algorithm}{Algorithms}
\crefname{conj}{Conjecture}{Conjectures}
\crefname{lem}{Lemma}{Lemmas}
\crefname{thm}{Theorem}{Theorems}
\crefname{customthm}{Theorem}{Theorems}
\crefname{claim}{Claim}{Claims}
\crefname{rmk}{Remark}{Remarks}
\crefname{prop}{Proposition}{Propositions}
\crefname{section}{Section}{Sections}
\crefname{appendix}{Appendix}{Appendices}
\crefname{cor}{Corollary}{Corollaries}
\crefname{figure}{Figure}{Figures}
\crefname{table}{Table}{Tables}
\crefname{example}{Example}{Examples}
\crefname{prob}{Problem}{Problems}
\crefname{assm}{Assumption}{Assumptions}
\crefname{defn}{Definition}{Definitions}
\crefname{customconj}{Conjecture}{Conjectures}

\numberwithin{equation}{section}

\title{Refined Gromov--Witten invariants}
\author{Andrea Brini$^{1,2}$ and Yannik Schuler$^{1,3,4}$}
\address{$^1$ University of Sheffield, School of Mathematical and Physical Sciences, Hounsfield Road,  Sheffield S3 7RH, United Kingdom.}
\address{$^2$ On leave from CNRS, DR~13, Montpellier, France}
\address{$^3$ University of Cambridge, Department of Pure Mathematics and Mathematical Statistics, Wilberforce Road, Cambridge CB3 0WB, United Kingdom}
\address{$^4$ ETH Z\"{u}rich, Department of Mathematics, R\"{a}mistrasse 101, 8092 Z\"{u}rich, Switzerland}
\email{a.brini@sheffield.ac.uk, yannik.schuler@math.ethz.ch}
\begin{document}

\begin{abstract}
We study the enumerative geometry of stable maps to Calabi--Yau 5-folds $Z$ with a group action preserving the Calabi--Yau form. In the central case $Z=X \times \bbC^2$, where $X$ is a Calabi--Yau 3-fold with a group action scaling the holomorphic volume form non-trivially, we conjecture that the disconnected equivariant Gromov--Witten generating series of $Z$ returns the Nekrasov--Okounkov equivariant K-theoretic PT partition function of $X$ and, under suitable rigidity conditions, its refined BPS index. We show that in the unrefined limit the conjecture reproduces known statements about the higher genus Gromov--Witten theory of $X$;  we prove it for $X$ the resolved conifold; and we establish a refined cycle-level local/relative correspondence for local del Pezzo surfaces, implying the Nekrasov--Shatashvili limit of the conjecture when $X$ is the local projective plane. We further establish B-model physics predictions of Huang--Klemm for refined higher genus mirror symmetry for local $\bbP^2$. In particular, we prove that our refined Gromov--Witten generating series obey extended holomorphic anomaly equations, are quasi-modular functions of $\Gamma_1(3)$, have leading asymptotics at the conifold point given by the logarithm of the Barnes double-Gamma function, and satisfy a version of the higher genus Crepant Resolution Correspondence with the refined orbifold Gromov--Witten theory of $[\bbC^3/\mu_3]$. This refines results, and partially proves conjectures, of Lho--Pandharipande, Coates--Iritani, and Bousseau--Fan--Guo--Wu.
\end{abstract}

\maketitle

\setcounter{tocdepth}{1}
\tableofcontents

\hfill
\pagebreak

\section{Introduction}

\subsection{Motivation}
Let $X$ be a smooth projective Calabi--Yau threefold. For a fixed genus $g\in \bbZ_{\geq 0}$ and class $\beta \in H_2(X, \bbZ)$, the Gromov--Witten (GW) invariant $\GW_{g,\beta}(X)$ is the virtual count
of 
parametrised curves
%
\[
\GW_{g,\beta}(X) \coloneqq \int_{[\overline{M}_{g}(X,\beta)]^{\vir}}1 \in \bbQ\,,
\] 
where $\overline{M}_{g}(X, \beta)$ is the moduli space of genus $g$, degree $\beta$ unpointed stable maps to $X$. Based on physical considerations, 
Gopakumar--Vafa \cite{GV98:MthTopStrI,GV98:MthTopStrII} were led to express the generating function of GW invariants as
\beq
\sum_{g \geq 0} \sum_{\beta \in H_2(X, \bbZ)} (-\epsilon^2)^{g-1} \mathcal{Q}^\beta \, \GW_{g,\beta}(X)  = \sum_{k \geq 1} \sum_{\beta \in H_2(X, \bbZ)} \frac{\mathcal{Q}^{k\beta}}{k [\mchar{q}^k]^2} \, \Omega^{\rm un}_\beta(X)(\mchar{q}^k)\,,
\label{eq:GWGV}
\eeq
where we denoted 
\[ \mchar{q} \coloneqq \re^{\epsilon}\,, \quad [\mchar{q}] \coloneqq \mchar{q}^{1/2}-\mchar{q}^{-1/2}\,,\] and $\Omega_{\mathrlap{\beta}}{\mathstrut}^{\rm un}(X)$ is the integral Laurent polynomial \cite{IP18:GV,DIW21:GVfinite}
\beq
\Omega^{\rm un}_\beta(X)(\mchar{q}) \coloneqq \sum_{g \geq 0} (-1)^{g-1} [\mchar{q}]^{2g} \, \invariantsfont{GV}_{g,\beta}(X)
\in \bbZ[\mchar{q}^{\pm 1}] \,.
\label{eq:Omegabeta}
\eeq
The integers $\invariantsfont{GV}_{g,\beta}(X)$ are the \emph{Gopakumar--Vafa} (GV) invariants of $X$. These are computed numerically via the Gromov--Witten invariants through \eqref{eq:GWGV}, and independent definitions in algebraic geometry have since been proposed using perverse sheaves of vanishing cycles \cite{MT18:GV} or, more recently, unramified stable maps in the case of primitive curve classes \cite{Nesterov:2024bdv}.

In the string theory heuristics, the GV invariants are determined by a finer physical count of the number $N^{\mathrlap{\beta}}{\mathstrut}_{j_-,j_+} (X)$ of stable BPS particles with charge $\beta$ and Lorentz weight $(j_-, j_+)$,
\[
\Omega^{\rm un}_\beta(X)(\mchar{q}) = 
\sum_{j_-, j_+ \in \frac{1}{2}\bbZ_{\geq 0}} N_{j_-,j_+}^\beta(X)\, (-1)^{2 j_-+2j_+} (2j_++1) \sum_{d_-=-j_-}^{j_-} \mchar{q}^{2 d_-} \in \bbZ[\mchar{q}^{\pm 1}]\,.
\]
Unlike the GV invariants, the BPS numbers 
are not expected to be deformation-invariant, and might generally depend on the complex structure of $X$.
However, when $X$ is a rigid quasi-projective Calabi--Yau threefold, 
it would be reasonable to expect that the individual BPS numbers 
could receive a geometric definition in terms of {\it bona fide} curve-counting invariants of $X$, simultaneously refining its ordinary Gromov--Witten and Gopakumar--Vafa invariants: for example, for local del Pezzo surfaces,
the numbers $N^{\mathrlap{\beta}}{\mathstrut}_{j_-,j_+} (X)$ are identified with the graded dimensions of  the cohomology of moduli spaces of pure one-dimensional sheaves on the surface, viewed as $\mathrm{Spin}(4)=\mathrm{SU}(2)\times \mathrm{SU}(2)$-representations for the hard Lefschetz actions
induced by 
the Hilbert--Chow morphism \cite{CKK14:refBPS,KPS23:refBPSP2,KLMP24:CohOneDimSheavesP2}. 
The corresponding generating function is the ``refined BPS index''
\beq
\sum_{k \geq 1} \sum_{\beta \in H_2(X, \bbZ)} \frac{\mc{Q}^{k\beta}}{k [\mchar{q}_1^k] [\mchar{q}_2^k]} \, \Omega_{\beta}(X)(\mchar{q}_+^k, \mchar{q}_-^k)\,,
\label{eq:GVref}
\eeq
where $\mchar{q}_1 \mchar{q}_2 =\mchar{q}_+^2$, $\mchar{q}_1/\mchar{q}_2=\mchar{q}_-^2$, and now
\beq
\Omega_\beta(X)(\mchar{q}_+,\mchar{q}_-) \coloneqq \sum_{j_-, j_+ \in \frac{1}{2}\bbZ_{\geq 0}} N_{j_-,j_+}^\beta(X)\, (-1)^{2 j_-+2j_+} \sum_{d_-=-j_-}^{j_-} \!\!\mchar{q}_-^{2d_-}\sum_{d_+=-j_+}^{j_+} \!\!\mchar{q}_+^{2d_+} \in \bbZ[\mchar{q}_+^{\pm 1}, \mchar{q}_-^{\pm 1}]\,,
\label{eq:Omegabetaref}
\eeq
from which the unrefined GV generating function $\Omega_{\mathrlap{\beta}}{\mathstrut}^{\rm un}(X)$ in \eqref{eq:Omegabeta} is recovered when $\mchar{q}_+ = 1$, $\mchar{q}_-=\mchar{q}$.\\
More generally, when $X$ is a Calabi--Yau threefold supporting a non-trivial $\Gm$-action on its canonical bundle,  
an influential suggestion of Nekrasov--Okounkov\footnote{With earlier incarnations in work of Nekrasov \cite{Nek03:SWprePInstCounting,Nek05:Ztheory,MR2491283}; see also \cite{Ok15:lectures}.} interprets the refined BPS index as a putative K-theoretic index of a moduli space of M2-branes of the extended target geometry $X \times \Aaff{2}$ 
\cite{NO14:membranes}, which effectively resums the K-theoretic $\bbC^\times$-equivariant Pandharipande--Thomas (PT) generating function of $X$. From this 5-dimensional vantage point, the variables $\mchar{q}_1$ and $\mchar{q}_2$ acquire a geometric interpretation as the coordinates of the torus $\Gm[2]$ scaling $\Aaff{2}$. 
\\

A natural question is 
how to circle back the BPS/sheaf-theoretic refinement \eqref{eq:Omegabetaref} to an equally meaningful refinement of the original \emph{stable maps} counting problem -- i.e.~to give a geometric definition of \emph{refined Gromov--Witten invariants}. On this point,  both the physics (see \cite{AFHNZ13:WorldsheetRefTopStr, NO12:WorldsheetOmegaBg} for attempts in this direction) and, rather more so, the mathematics of the refinement have however so far been elusive.

In this paper we address this question by  considering the enumerative geometry of stable maps to Calabi--Yau fivefolds with a reductive group action preserving their holomorphic volume form. We loosely state the core idea here, which we will qualify more precisely in \cref{conj: refined GV integrality intro,conj: GWPT intro}. 

\begin{customconj}{M}
    \label{speculation: link to Mtheory}
    For $Z$  a smooth quasi-projective complex Calabi--Yau 5-fold with a Calabi-Yau reductive group action,
    the equivariant disconnected Gromov--Witten generating function of $Z$ equates the Chern character of its equivariant  K-theoretic Nekrasov--Okounkov membrane index.
\end{customconj}

When $Z=X \times \Aaff{2}$, 
\cref{speculation: link to Mtheory} 
will posit that the equivariant enumeration of stable maps in $Z$ provides a definition of refined GW invariants of $X$;
under suitable properness conditions, we claim that the resulting equivariant invariants indeed return the expected refinement for both GV and PT invariants.
The setup, main statements, and broader context of the paper are described below.

\subsection{Setup}
\label{sec:introsetup}
Throughout this section, we let $X$ be a smooth quasi-projective complex Calabi--Yau threefold, endowed with the action of an algebraic torus $\mgp{T}'$ having compact fixed locus and scaling non-trivially the canonical bundle $\omega_X$. We will take $Z$ to be the Calabi--Yau five-fold $Z \coloneqq X \times \bbC^2$ endowed with a $\mgp{T} \simeq \Gm[2]$ action which simultaneously lifts the $\mgp{T}'$ action on $X$ and the canonical linear 2-torus action on $\bbC^2$, and which furthermore acts trivially on $\omega_Z$. We will write $\mchar{q}_1$, $\mchar{q}_2$ for the coordinates of
the torus action on $\Gm[2]$. By possibly passing to a cover, we will assume that there exist $\mgp{T}$-characters $\mchar{q}_\pm \in \mathrm{Rep}(\mgp{T}$) with $\mchar{q}_+^2=\mchar{q}_1 \mchar{q}_2$, $\mchar{q}_-^2=\mchar{q}_1/\mchar{q}_2$.

The main quantity we will consider is the generating series of $\mgp{T}$-equivariant, unpointed higher genus Gromov--Witten invariants of $Z$.
\begin{defn}[=\Cref{defn: refined GW invariant} \& \labelcref{defn: refined GW series}]
We define
\beq
\refGWb{Z}{\mgp{T}}{\beta} \coloneqq \sum_{g \geq 0} ~ \intEquiv{\mgp{T}}_{[\Mbar_{g}(Z,\beta)]^{\vir}_{\mgp{T}}} 1 \,.
\label{eq:GWref}
\eeq     
\label{def:refGWintro}
\end{defn}
In \eqref{eq:GWref}, the right-hand side 
is defined as a $\mgp{T}$-equivariant residue \cite{GP97:virtloc} and under suitable properness assumptions on $X$ the sum can be regarded as an element of a completion $(\epsilon_1 \epsilon_2)^{-1} \widehat{H(B \mgp{T}, \bbQ)} \simeq (\epsilon_1 \epsilon_2)^{-1}\bbQ \llbracket \epsilon_1, \epsilon_2 \rrbracket$ of the $\mgp{T}$-equivariant cohomology of the point, where $\epsilon_{i}=c_1^{\mgp{T}}(\mchar{q}_i)$ are 
the equivariant parameters of the action on $\bbC^2$. 
By Mumford's formula for the Chern character of the Hodge bundle \cite{Mu83}, 
the equivariant parameter $\epsilon_-=c_1^{\mgp{T}}(\mchar{q}_-)$ of the anti-diagonal torus scaling $\bbC^2$ is naturally identified with the genus expansion variable (\Cref{prop: unref limit}):
\begin{equation*}
    \refGWb{Z}{\mgp{T}}{\beta} ~ \Big\vert_{\epsilon_1 = -\epsilon_2 = \epsilon_-} = \sum_{g \geq 0} (-\epsilon_-^2)^{g-1} \GW_{g,\beta}(X)\,.
\end{equation*}

The definition in \eqref{eq:GWref}  depends \emph{a priori} on the specific lift $\mgp{T}$ of the 2-torus action on $\bbC^2$: for example, if $X$ is a toric Calabi--Yau threefold, there is a $(\Gm)^2$-worth of such choices, on which the invariants will in principle depend. An assumption we will often make is to restrict our attention to target varieties $X$ for which the moduli space of stable maps is proper in all genera and (non-zero) degrees: we call such geometries \textit{equirigid}.
This is for example the case for $X$ the resolved conifold, or a local del Pezzo surface, whenever $\beta$ is a non-zero effective curve class; examples where this is not the case will be considered in \cref{sec: local P1}. For equirigid $X$, we claim in \cref{conj: rigidity} 
that the resulting invariants are independent of the particular lift to $Z$ of the 2-torus action on $\bbC^2$, and provide ample evidence for this phenomenon in examples. Lacking a rigidity assumption, we will generally observe a dependence on the choice of torus action\footnote{This is expected from the point of view of GV/PT refinement: for non-equirigid targets there are known discrepancies between the K-theoretic and motivic refinement of PT theory (see for instance \cite{Ob24:EnriquesKth,Ob24:EnriquesMotivic} or \cite[Rmk.~0.2]{Th24:RefinedK3}). Relatedly, when $X$ is a toric CY3 which is not equirigid, the outcome of the refined topological vertex is contingent on choices of ``preferred directions'' \cite{IKV09:RefVert}.}: this observation suggests to work equivariantly with respect to the largest symmetry group $\mgp{G}$ acting on $X\times \Aaff{2}$, which indeed is the perspective we take in the main body of the paper.

\subsection{Main statements}
\label{sec:introstatements}

We will articulate the general expectation of \Cref{speculation: link to Mtheory} in the specific setup where $Z=X\times \bbC^2$ as two main correspondences. We formulate them somewhat informally here, but include pointers to the main body of the paper for the precise statements.

The first correspondence asserts that, for equirigid $X$, the fixed-degree equivariant Gromov--Witten generating series of $X \times \Aaff{2}$ in \cref{def:refGWintro} are governed by a finite number of integer invariants. Furthermore, when $X=K_S$ is a local del Pezzo surface, these integers encode the cohomological refinement \cite{IKV09:RefVert,CKK14:refBPS,KPS23:refBPSP2,KLMP24:CohOneDimSheavesP2} of the Gopakumar--Vafa invariants defined by Maulik--Toda \cite{MT18:GV} using the perverse filtration associated to the Hilbert--Chow map.

\begin{conj}[Refined BPS integrality and GW/GV correspondence, =\Cref{conj: refined GV integrality,conj: MT type BPS local surfaces}; see \cref{sec: rationality conjecture,sec: integrality conjecture,sec: geometric interpretation conjecture} for details]
    \label{conj: refined GV integrality intro}
    If $X$ is equirigid and $\beta \neq 0$, there exist integral Laurent polynomials $\Omega_{\beta}(X)(\mchar{q}_+, \mchar{q}_-) \in \bbZ[\mchar{q}_+^{\pm 1}, \mchar{q}_-^{\pm 1}]$ satisfying
    \[
    \refGWb{Z}{\mgp{T}}{\beta} = \ch[\mgp{T}] \sum_{k | \beta} \frac{1}{k [\mchar{q}_1^k][\mchar{q}_2^k]} \,\Omega_{\beta/k}(X)(\mchar{q}_+^k, \mchar{q}_-^k)\,,
    \]
    with the equivariant Chern character acting via $\ch[\mgp{T}] \, \mchar{q}_i = \re^{\epsilon_i}$. Moreover, if $X=K_S$ is a local del Pezzo surface, we have
    \[
    \Omega_{\beta}(K_S) = \MT_\beta(K_S)
    \,,
    \]
    where
    $\MT_\beta(K_S)(\mchar{q}_+,\mchar{q}_-)$ are the Poincar\'e polynomials \eqref{eq:refMTb} for the perverse filtration on the cohomology of the moduli space of 1-dimensional semi-stable sheaves on $S$ with support $\beta$ and Euler characteristic one (see \cref{sec: geometric interpretation conjecture}).
\end{conj}

The second correspondence combines \Cref{speculation: link to Mtheory} with the membranes/sheaves correspondence of \cite{NO14:membranes} to provide 
a natural refinement of the  numerical GW/PT correspondence conjecture \cite{PT09:StabPairs}
 upon identifying the coordinate  of the anti-diagonal torus $\Gm_{\mchar{q}_-}$ acting on $\bbC^2$ with the box-counting parameter in PT theory.

\begin{conj}[Refined GW/PT correspondence, =\Cref{conj: refined PT GW}; see \cref{sec:KthPT,sec:KthPTconj} for a detailed treatment]
    \label{conj: GWPT intro}
Let $\beta \neq 0$ and
\begin{align*}
 \refGWdiscb{Z}{\mgp{T}}{\beta} & \coloneqq \sum_{g\in \bbZ} \intEquiv{\mgp{T}}_{[\Mbar^{\bullet}_{g}(Z,\beta)]^{\vir}_{\mgp{T}}} 1\,, \\
\refPTb{X}{\mgp{T}^{\prime}}{\beta} & \coloneqq \sum_{n \in \bbZ} (-\mchar{q}_-)^n \, \chi_{ \mgp{T}^{\prime}} \!\left( P_n(X,\beta), \, \widehat{\mc{O}}^{\vir}_{P_n(X,\beta)} \right)
\end{align*}
denote, respectively, the fixed degree-$\beta$  disconnected $\mgp{T}=\mgp{T}' \times \Gm_{\mchar{q}_-}$-equivariant Gromov--Witten generating series of $X \times \Aaff{2}$ and the $\mgp{T}'$-equivariant K-theoretic PT generating series of $X$, defined in both cases by virtual localisation. Then,
\beq
 \refGWdiscb{Z}{\mgp{T}}{\beta} = 
 \ch[\mgp{T}] \, \refPTb{X}{\mgp{T}^{\prime}}{\beta}\,.
\label{eq:GWPTintro}
\eeq
\end{conj}
We stress that, unlike for \Cref{conj: refined GV integrality intro}, we do not impose any rigidity assumption in \Cref{conj: GWPT intro}. As K-theoretic and motivic PT invariants generally disagree when the stable pairs moduli space is non-compact, this in particular means that we can only expect an agreement with motivically refined sheaf counting when the target is equirigid. In that case, the Rigidity \cref{conj: rigidity} implies that \eqref{eq:GWPTintro} would hold with $\mgp{T}$ replaced by any Calabi--Yau 2-torus on $Z$ covering the $\mgp{T}'$ action on $X$ and the $\Gm[2]$-action on $\Aaff{2}$.
\\


We will show in \cref{sec:unrefGW,sec: PT GW unref,sec:refGVevidence} that \cref{conj: refined GV integrality intro,conj: GWPT intro} reduce in the unrefined limit ($\epsilon_1=-\epsilon_2$) to the standard notions of GV integrality and GW/PT correspondence for the higher genus GW theory of $X$. Away from this limit, the validity of \cref{conj: refined GV integrality intro,conj: GWPT intro}, and therefore of \cref{def:refGWintro} as a potential definition of refined GW invariants of $X$, can be broadly subjected to two types of litmus tests. 

\bit
\item In the A-model, one would like to directly contrast the equivariant deformation in \eqref{eq:GWref} 
with, respectively, the Poincar\'e refinement of GV invariants \cite{MT18:GV} arising from the double Lefschetz action on the moduli of pure 1-dimensional sheaves on $X$, and with the
K-theoretic equivariant  refinement \cite{NO14:membranes} (equivalently, when $X$ is equirigid, the motivic refinement  \cite{MR2851153,CKK14:refBPS}) in PT theory. 
\item In the B-model, \cref{def:refGWintro} should satisfy the ``refined BCOV axioms'' of \cite{HK10:OmegaBG}, including the refined holomorphic (or modular) anomaly equations, a finite generation/qua\-si-modularity property, as well as a set of boundary conditions at suitable divisors in the stringy K\"ahler moduli space of $X$.  
\eit

Our results will be, accordingly, structured along these two main strands. 

\subsubsection{A-model}

Our first results here concern genus zero local curves. We start by considering the basic equirigid example of the resolved conifold.


\begin{customthm}{A1}[=\cref{prop:ResConifold}]
Let $X=\mathrm{Tot}(\cO_{\bbP^1}^{\oplus 2}(-1))$ and $\mgp{T}$ be the fibrewise 2-torus acting with opposite weights on each $\cO_{\bbP^1}(-1) \oplus \cO_{\bbP^1}$ fibre of $Z=X\times \Aaff{2} =\mathrm{Tot}(\cO_{\bbP^1}^{\oplus 2}(-1)\oplus \cO_{\bbP^1}^{\oplus 2})$.  
    Then, the refined BPS integrality in \cref{conj: refined GV integrality intro} holds as 
    \beq
    \refGWb{Z}{\mgp{T}}{d [\bbP^1]} = \ch[\mgp{T}] ~ \frac{1}{d [\mchar{q}_1^d] [\mchar{q}_2^d]}
    \,.
    \label{eq:GVrescon}
    \eeq
In particular, $\Omega_{d [\bbP^1]}(X) = \delta_{d,1}$.
\label{thm:A1}
\end{customthm}
The result follows from a direct calculation of \eqref{eq:GWref} to all genera and degrees via the degeneration formula in Gromov--Witten theory, which reduces the calculation to known results about the local GW theory of curves \cite{BP08:GWLocCurves,BP05:CurvesCY3TQFT}. In accordance with the discussion below \cref{conj: GWPT intro}, the equivariant GW partition function of $Z$ should return the K-theoretic $\mgp{T}'$-equivariant PT partition function of $X$, where $\mgp{T}'$ acts on $\omega_X^{-1}$ with character $\mgp{q}_+^{2}$ \cite{IKV09:RefVert,Arb21:KthDT,NO14:membranes}. Indeed, by 
\eqref{eq:GVrescon}, we have
\[
1+\sum_{d =1}^\infty \refGWdiscb{Z}{\mgp{T}}{d [\bbP^1]} \,\cQ^d =\ch[{\mgp{T}}]
 \prod_{m \geq 1} \prod_{j=0}^{m-1} \l(1-\mchar{q^{2j-m-1}_+} \mchar{q_-^m} \cQ\r)
 =1+\ch[{\mgp{T}}] \sum_{d =1}^\infty \refPTb{X}{\mgp{T}'}{d [\bbP^1]}\, \cQ^d 
 \,.
\]
Aspects of non-equirigid local $\bbP^1$ geometries, along with calculational evidence for \cref{conj: GWPT intro} and a comparison with the Nekrasov--Okounkov membrane index in that setup, are further discussed in \cref{sec: local P1 shifted}.\\

Our second series of results considers the substantially more involved case of a local Calabi--Yau threefold with compact divisors. We take here $X=K_S$ to be the total space of the canonical bundle over a smooth del Pezzo surface $S$, and denote $\mgp{T}_{\rm F}$ the 2-torus that scales the fibres of $Z=K_S \times \bbC^2$ while acting trivially on $\omega_Z=\omega_{K_S \times \bbC^2}$. We do not, in particular, assume that $S$ itself is toric.

We will first relate the equivariantly refined generating GW function of $Z$ to the generating series of 
relative GW invariants of $S$ with maximal contact along a smooth anti-canonical curve $D$, together with an insertion of the top Chern class of the Hodge bundle:
\begin{equation*}
    \GW_{\beta}(S / D) \coloneqq \frac{(-1)^{D\cdot \beta +1}}{D\cdot \beta} \sum_{g\geq 0} \epsilon_1^{2g-1} \int_{[\Mbar_{g,0,(D\cdot \beta)}(S / D,\beta)]^{\vir}} \lambda_g\,.
\end{equation*}

\begin{customthm}{A2}[=\Cref{thm: NS vs log cycle,thm: NS vs log num}]
    \label{thm: NS vs log num intro}
    Let $S$ be a smooth del Pezzo surface and $D$ a smooth connected anticanonical curve. Then
    \begin{equation}
        \epsilon_2 \, \refGWb{Z}{\mgp{T}_{\rm F}}{\beta} ~\Big|_{\epsilon_2=0} = \GW_{\beta}(S / D)\,.
    \label{eq:NSlogGW}
    \end{equation}
\label{thm:A2}
\end{customthm}

The simple equality in \eqref{eq:NSlogGW} follows from a stronger cycle-valued identity (\cref{thm: NS vs log cycle}) relating the $\mgp{T}_{\mr{F}}$-equivariant virtual class on $\overline{M}_{g,n}(\mathrm{Tot}(\cO_S(-D) \oplus \cO^{\oplus 2}_S), \beta)$ to the virtual classes of the moduli spaces of stable maps 
to the base surface $S$ relative to $D$: the key technical tool employed in its proof is the Kim--Lho--Ruddat degeneration formula \cite{KLR21:degen}. 
\cref{thm: NS vs log cycle} provides in particular a refined version of the log-local correspondence of \cite{vGGR19} for surfaces in the so-called Nekrasov--Shatashvili (NS) limit ($\epsilon_2=0$); this is to be contrasted with the higher genus log-local correspondence of \cite[Thm.~1.2]{BFGW21:HAE}, which corresponds instead to the unrefined limit, $\epsilon_1=-\epsilon_2=\epsilon_-$. 

When $S=\bbP^2$, we deduce from \cref{thm:A2} a proof of \cref{conj: refined GV integrality intro} in the NS  limit. For this, we will use a result of Bousseau
\cite[Thm.~0.4.4 and 0.4.6]{Bou19:Takahashi}
which equates the right-hand side of \eqref{eq:NSlogGW} to a plethystic generating series of Poincar\'e polynomials of the moduli spaces of pure one-dimensional sheaves on $\bbP^2$.


\begin{customthm}{A3}[=\cref{cor: MT BPS local P2 NS}]
\label{thm:A3}
    The refined GW/GV correspondence holds  in the Nekrasov--Shatashvili (NS) limit for $X=K_{\bbP^2}$ with fibrewise $\mgp{T}_{\rm F}$-action,
\beq
\epsilon_2 \, \refGWb{Z}{\mgp{T}_{\rm F}}{\beta} ~\Big|_{\epsilon_2=0} = \sum_{ k \vert \beta} \frac{1}{k^2 [\mchar{q}_1^k]} \, \MT^{\rm NS}_{\beta/k}(K_{\bbP^2})(\mchar{q}_1^k)~\bigg|_{\mchar{q}_1=\re^{\epsilon_1}} \,.
\label{eq:GWGVrefNS}
\eeq
\end{customthm}

In \eqref{eq:GWGVrefNS}, we denoted by $\MT^{\rm NS}_\beta(K_{\bbP^2})$ the refined Maulik--Toda generating function in the NS limit,
\[ \MT^{\rm NS}_\beta(K_{\bbP^2})(\mchar{q}) \coloneqq (-\mchar{q}^{-1/2})^{\dim(M_{\beta})} \sum_{j=0}^{\dim(M_{\beta})} b_{2j}(M_{\beta}) \, (-\mchar{q})^j\,,
\]
given by the symmetrised Poincar\'{e} polynomial of the moduli space $M_{\beta}$ of pure one-dimensional Gieseker semi-stable sheaves on $\bbP^2$ with second Chern character $\beta$ and Euler characteristic one. The right-hand side of \eqref{eq:GWGVrefNS} should be then viewed as the limit
\[
(\mchar{q}_2-1) \sum_{ k \vert \beta} \frac{1}{k [\mchar{q}_1^k][\mchar{q}_2^k]} \, \MT_{\beta/k}(K_{\bbP^2})(\mchar{q}_+^k, \mchar{q}_-^k) ~\Bigg|_{\mchar{q}_1=\re^{\epsilon_1}\,,\mchar{q}_2=1}\,.
\] 
\cref{thm:A2,thm:A3} further elucidate the surprising observation made in \cite{Bou20:QuTropVert, BFGW21:HAE,Bou19:Takahashi} that the relative Gromov--Witten theory of a del Pezzo surface $S$ with one $\lambda_g$-insertion encodes the NS limit ($\mchar{q}_2 \to 1$) of the BPS generating function the local surface geometry $K_S$. From \eqref{eq:NSlogGW}, this is conceptually explained by the former arising as the geometric $\epsilon_2=0$ ($\mchar{q}_2 \to 1$) restriction of the equivariant five-fold GW refinement in \cref{def:refGWintro}.
Away from the specialisation $\mchar{q}_2=1$, and in order to be compatible with the Maulik--Toda proposal in the unrefined limit \cite{MT18:GV}, the refined BPS index should additionally keep track of the perverse filtration induced by the Hilbert--Chow morphism: we formulate this expectation as a conjecture for general del Pezzo surfaces in \Cref{conj: MT type BPS local surfaces}.


\subsubsection{B-model}
\label{sec:Bmodel}
There is, at the time of writing, no conceptual definition of refined higher genus mirror symmetry\footnote{More is known in the case of refined \emph{genus zero} mirror symmetry, i.e.~the Nekrasov--Shatashvili limit, in which case the refinement should be related to some non-commutative deformation of the local mirror geometry \cite{Aganagic:2011mi}.} for a local Calabi--Yau threefold, akin e.g.~to the construction of the unrefined higher genus B-model through  the Kodaira--Spencer theory of gravity \cite{Bershadsky:1993cx,Costello:2012cy}. There is, however, a {\it calculational} definition of refined B-model free energies in terms of the ``refined direct integration'' formalism of \cite{Grimm:2007tm,HK10:OmegaBG,KW11:ExtHAE}, which is conjectured to reproduce the motivic/equivariant refinement in PT theory under the mirror map \cite{HK10:OmegaBG,HKK13:OmegaBmodelRigidN2,CKK14:refBPS}. Viewed solely from an A-model vantage point, one would regard the refined direct integration method as  providing a list of (remarkable) structural properties of the refined GW generating series, which furthermore determine them uniquely. We recall what these are below, specialising the entire discussion to the case $X=K_{\bbP^2}$.

The mirror family of local $\bbP^2$ is described by the family of elliptic curves over the weighted projective line $\bbP(3,1)$ \cite{Hori:2000ck, Coates:2018hms},

\[
\cZ = \l\{\l([z_0:z_1:z_2], y\r) \in \bbP^2 \times \bbP(3,1)\, \big|\, z_0^3+ z_1^3+ z_2^3 + y^{-1/3} z_0 z_1 z_2=0 \r\} \longrightarrow \bbP(3,1)\,.
\]

The family has singular fibres at the two values $y=-1/27$ (the {\it conifold point}) and $y=0$ (the {\it large radius point}), away from which the fibres are smooth elliptic curves. On a regular fibre $\cZ_{y}$ over $y \neq 0$, $y \neq -1/27$, write
\[
\Pi_\gamma(y) = \int_\gamma \log z_1 \dd \log z_2
\]
for the (single-valued) period of the (multi-valued) 1-form $\log z_1 \dd \log z_2$ over a homology cycle $\gamma \in H_1(\cZ_y, \bbZ)$. The periods $\Pi_\gamma$ have unipotent (resp.~maximal unipotent) monodromy around the singular fibre at the conifold (resp.~the large radius point), and cyclic $\mu_3$-monodromy around $y=\infty$. Around the large radius point $y=0$, there is a unique $\gamma \in H_1(\cZ_y, \bbZ)$ such that $\Pi_{\gamma} \sim \log y + \cO(y)$ diverges logarithmically \cite{CKYZ99:locMS}. Around the conifold point $y=-1/27$, there is a unique 
$\gamma_{\rm cf} \in H_1(\cZ_y, \bbZ)$ such that $\Pi_{\gamma_{\rm cf}} \sim (y-1/27) + \cO(y-1/27)^2$ vanishes linearly \cite{HKR:HAE}; and around the orbifold point, there is a unique $\gamma_{\rm orb} \in H_1(\cZ_y, \bbC)$ such that $\Pi_{\gamma_{\rm orb}} \sim y^{-1/3} + \cO(y^{-4/3})$ \cite{Aganagic:2006wq}. 

The complement of the conifold and large radius points is the coarse moduli space of the modular curve $[\mathbb{H}/\Gamma_1(3)]$; the discriminant $X(y)$, Jacobian $A(y)$, and propagator $F(y)$, defined as
\[
X(y) \coloneqq \frac{1}{1+27 y}\,, \qquad A(y)= y \de_y \log \Pi_{\gamma}\,, \qquad F(y) \coloneqq y\de_y \log A(y)
\]
are quasi-modular functions of $\Gamma_1(3)$ of weight and depth equal to  $(0,0)$, $(1,0)$, and $(0,1)$  respectively. In the following, we will similarly define
\[
 A_{\rm orb}(y)= y \de_y \log \Pi_{\gamma_{\rm orb}}\,, \quad F_{\rm orb}(y) \coloneqq y\de_y \log A_{\rm orb}(y)\,, \quad A_{\rm cf}(y)= y \de_y \log \Pi_{\gamma_{\rm cf}}\,, \quad F_{\rm cf}(y) \coloneqq y\de_y \log A_{\rm cf}(y)\,,
\]
and denote by \[\big[x^n\big] f\] the $n^{\rm th}$ formal Laurent  coefficient of $f \in \bbQ(\!(x)\!)$. 

\begin{defn}[Refined direct integration for local $\bbP^2$, \cite{HK10:OmegaBG}] 
\label{def:dirint}
For all $k+g \geq 2$, the refined generating series are uniquely and recursively determined from $\mathscr{F}_{0,0}$, $\mathscr{F}_{0,1}$ and $\mathscr{F}_{1,0}$ by imposing:    

\begin{description}
 \item[(Finite generation)] $\mathscr{F}_{k,g}$ is a weight-zero, depth $3(g+k-1)$ quasi-modular function of $\Gamma_1(3)$ with
  \[ \mathscr{F}_{k,g} \in \bbQ[X^{\pm 1}][F]\,, \quad \deg_F \mathscr{F}_{k,g} = 3(g+k-1)\,.\] 
 \item[(Modular anomaly)] the refined modular anomaly equations hold:
 \[
\frac{X}{A^2}
\frac{\de \mathscr{F}_{k,g}}{\de F}=\frac{1}{2} (\cQ \partial_\cQ)^2 \mathscr{F}_{k,g-1} + \frac{\cQ^2}{2} \sum_{\substack{0\leq k' \leq k\,, 0 \leq g' \leq g \\(0,0)\neq (k',g')\neq (k,g) }} ~ \de_\cQ \mathscr{F}_{k',g'} \,  \partial_\cQ \mathscr{F}_{k-k',g-g'}\,.\]
 \item[(Orbifold regularity)]  near $y=\infty$, the following asymptotic behaviour holds
 \[ \mathscr{F}_{k,g}\Big|_{F \to F_{\rm orb}} = \cO(y^{-1})\,. \]
  \item[(Conifold asymptotics)] near $y=-1/27$, the following asymptotic behaviour holds
 \[ \mathscr{F}_{k,g}\Big|_{F \to F_{\rm cf}} =  \mathfrak{b}_{k,g} \Pi^{2-2g+2k}_{\gamma_{\rm cf}} + \cO(1)\,,
  \]
  where
  \beq \mathfrak{b}_{k,g} =   (-3)^{g-1} \big[t^{2g} (\epsilon_1+\epsilon_2)^k(\epsilon_1\epsilon_2)^{g-1}\big] \l(\frac{ (2g-1)!}{[\mchar{q}_1^t] [\mchar{q}_2^t]}\r)\,.
  \label{eq:bkgcf}
  \eeq 
\end{description}
\end{defn}

With a suitable mirror symmetry prescription for the initial data $\mathscr{F}_{k,g}$ with $k+g=0$ (the special geometry prepotential) and $k+g=1$ (the refined Ray--Singer torsion), the authors of \cite{HK10:OmegaBG,CKK14:refBPS} conjecture and verify in low genus and degree that
\[
\mathscr{F}_{k,g} = \big[(\epsilon_1+\epsilon_2)^k (\epsilon_1 \epsilon_2)^{g-1}\big]  \sum_{k \geq 1} \l[\sum_{0 \neq \beta} \frac{\cQ^{k\beta}}{k [\mchar{q}_1^k][\mchar{q}_2^k]} \, \MT_\beta(K_{\bbP^2})( \mchar{q}_+^k, \mchar{q}_-^k)\r]\Bigg|_{\mchar{q}_i=\re^{\epsilon_i}}\,.
\]

According to \cref{conj: refined GV integrality intro}, it would then be expected that the $(\epsilon_1+\epsilon_2, \epsilon_1 \epsilon_2)$ expansion of the refined GW generating series of \cref{def:refGWintro} satisfies the refined direct integration properties in \cref{def:dirint}, with the same initial conditions.


\begin{customthm}{B}[Refined higher genus mirror symmetry for local $\bbP^2$, =\cref{thm:haegen,thm:haeref,thm:crc,thm:cfasym}]
\label{thm:B}
Let $X=K_{\bbP^2}$, $Z=K_{\bbP^2}\times \Aaff{2}$ and $\mgp{T}_{\rm F}$ be the fibrewise action as in \cref{thm:A3}. Writing $H$ for the class of a line in $\bbP^2$, define
\[
\mathscr{F}_{k,g} \coloneqq \sum_{d=1}^\infty  \l(\big[(\epsilon_1+\epsilon_2)^k (\epsilon_1 \epsilon_2)^{g-1}\big]\, \refGWb{Z}{\mgp{T}_{\rm F}}{d [H]}\r) \cQ^d \in \bbQ\llbracket \cQ \rrbracket\,.
\]
Then, under the mirror map $\cQ \to \re^{\Pi_\gamma(y)}$:
\bit
\item 
for $k+g<2$, $\mathscr{F}_{k,g}$ coincides with the Huang--Klemm initial data  \cite{HK10:OmegaBG,CKK14:refBPS};
\item for $k+g\geq 2$, $\mathscr{F}_{k,g}$ satisfies {\bf (Finite generation)},  {\bf (Modular anomaly)}, {\bf (Orbifold regularity)}, and the following weak version of {\bf (Conifold asymptotics)}:
 \[ \mathscr{F}_{k,g}\Big|_{F \to F_{\rm cf}} = \mathfrak{b}_{k,g} \Pi^{2-2g+2k}_{\gamma_{\rm cf}}\l(1 + \cO(\Pi_{\gamma_{\rm cf}})\r)\,.
  \]
\eit
\end{customthm}

The proof of {\bf (Finite generation)} and ({\bf Modular anomaly}) is identical in spirit to, although in practice more challenging than, the proof of the two analogous statements in the unrefined limit by \cite{LP18:HAE,Coates:2018hms} via the Givental--Teleman reconstruction theorem. The main complication here  arises from the need to consider a non-Calabi--Yau torus action on $X$, leading to a much enlarged set of generators in {\bf (Finite generation)}: equivariance with respect to the $\mathrm{PGL(3,\bbC)}$ automorphism group of the target is imposed to
dramatically reduce this to the claimed $\Gamma_1(3)$-quasimodularity. The unrefined limit $\epsilon_1=-\epsilon_2$ yields a slightly stronger finite generation property than the one proved in \cite[Thm.~1(i)]{LP18:HAE}, with the NS limit $\epsilon_2 \to 0$ recovering the quasi-modularity and modular anomaly results of \cite[Thm.~1.7-1.8]{BFGW21:HAE}. \\
The proof of ({\bf Orbifold regularity)} follows as a consequence of a stronger statement, establishing a higher genus refined Crepant Resolution Correspondence (CRC) for the resolution $\pi: K_{\bbP^2} \times \bbC^2 \longrightarrow \bbC^5/\mu_3$: this is a relatively straightforward generalisation 
of 
the proof of the threefold higher genus CRC for $\pi: K_{\bbP^2} \longrightarrow \bbC^3/\mu_3$ in \cite{MR3960668}. \\
 Lastly, {\bf (Conifold asymptotics)} follows from showing that, at the leading order in an expansion in the conifold flat coordinate, the conifold-frame GW refined generating series reduce to the Givental--Teleman potential of the special cubic Hodge CohFT  \cite{Okounkov:2003aok}, albeit with a subtly deformed translation action. The closed-form summation of the resulting stable graph expansion to the expression in \eqref{eq:bkgcf} is implied by the comparison with an analogous R-matrix calculation for the resolved conifold (\cref{thm:A1}). When $\epsilon_1=-\epsilon_2$ (resp., $\epsilon_2\to 0$), this in particular establishes, to all genera, the conjectured behaviour for the leading order of the conifold asymptotics considered in the unrefined case in \cite[Sec.~10.8]{Coates:2018hms} (resp., in the NS limit in \cite[Conj.~1.9]{BFGW21:HAE}). 

The full version of {(\bf Conifold asymptotics)}, including the vanishing of the sub-leading polar part (the ``conifold gap''), appears to require a significant increase in sophistication in the methods. We could, however, verify it in low genus (see \cref{exmpl:F2,rmk:cfgap}), where we also verify that the initial conditions of \cite{HK10:OmegaBG} for the recursion in \cref{def:dirint} are satisfied 
for $k+g< 2$ (\cref{exmpl:F1}).

\subsection{Context and implications}

\subsubsection{Relation to the work of Nekrasov--Okounkov} 
As has already been mentioned, although the actual mathematics of this paper is in practice much different, our main conceptual thread bears a very significant intellectual debt to the work of Nekrasov--Okounkov \cite{NO14:membranes}, in which the K-theoretic PT theory of complex threefolds $X$ is geometrised through a putative equivariant index of a  compactification of the moduli of immersed curves on an associated local Calabi--Yau fivefold $Z$ (the \emph{stable membranes} moduli space).  \cref{speculation: link to Mtheory} is indeed the natural Gromov--Witten counterpart of the membranes/sheaves correspondence of \cite[Conj.~1]{NO14:membranes} where one insists on computing the membrane contributions using stable maps as advised (against) in \cite[Sec.~4.1.1]{NO14:membranes}: as boundedness of the moduli spaces of fixed-degree immersed curves is lost in a stable map compactification, the membrane index should be reproduced by a manual summation of GW contributions to all genera, in a specular manner to the box-counting summation of \cite{NO14:membranes} on the PT side. This is exactly the prescription given in \cref{def:refGWintro}: the absence, or rather redundance (see \cref{sec:unrefGW}), of a dedicated genus counting parameter in \eqref{eq:GWref} closely reflects the M-theory perspective of \cite[Sec.~2.2.10]{NO14:membranes}, whereby no membrane degrees of freedom can couple to the genus of the immersed curves. \\ An upside of the use of stable maps is, of course, enhanced computability. A consequence of \cref{speculation: link to Mtheory} is that
structural properties of GW theory (such as virtual localisation, degeneration, and $R$-matrix quantisation) would consequently apply to the K-theoretic membrane index, and in particular 
to the refined BPS index of equirigid local CY threefolds: this is at the heart of the proof of \cref{thm:A1,thm:A2,thm:A3} and \cref{thm:B}.

\subsubsection{Rationality, integrality, rigidity, and the vertex}

\Cref{conj: GWPT intro} has a number of consequences. Firstly, for the change-of-variables under application of the Chern character to be well-defined, both the PT and GW generating series need to admit lifts to a rational function on $\mgp{T}'\times \Gm_{\mchar{q}_-}$. This can be viewed as a weak form of the integrality part of \Cref{conj: refined GV integrality intro} which we expect to hold for arbitrary (possibly non-equirigid) $X$. This common lift is expected to admit a modular interpretation as an index in localised K-theory with integral coefficients (\Cref{speculation: link to Mtheory}).\\
Secondly, for equirigid $X$, Nekrasov--Okounkov prove that the K-theoretic PT generating series only depends on the choice of $\mgp{T}'$-action on $X$ through the character of the induced action on $\omega_X$ \cite[Thm.~1]{NO14:membranes}. By \cref{conj: GWPT intro}, the same should apply to the equivariant Gromov--Witten generating series: this is the GW Rigidity Conjecture (\cref{conj: rigidity}), predicting the independence of $\refGWb{Z}{\mgp{T}}{\beta}$ on the particular choice of two-torus action.\\
Finally, in the equirigid case the PT generating series can efficiently be computed using the refined topological vertex \cite{NO14:membranes,Arb21:KthDT,IKV09:RefVert}, and \Cref{conj: GWPT intro} entails that the same method could be applied to obtain refined GW invariants. For non-equirigid geometries, the last statement is more subtle, and we will illustrate it with a simple example in \Cref{sec: refined vertex}. A refined GW vertex formalism for toric CY threefolds would require a better handle on quintuple Hodge integrals, and establishing \cref{conj: rigidity} appears to be key in this respect, as the rigidity principle strongly constrains these integrals. 


\subsubsection{Higher rank curve counts on threefolds}

In \cite[Prop.~5.1]{NO14:membranes}, the authors link rank-$r$ DT theory on $X$ to the membrane index on the minimal resolution of $Z=X \times (\bbC^2/\mu_r)$. Through \Cref{speculation: link to Mtheory}, one should accordingly expect a connection to the equivariant GW theory of $Z$, which in particular should be governed by the 7d master formula \cite{dZNPZ22:PlayingMtheory}. Some subtleties in the correspondence arise from the fact that  GW generating series  detect perturbative terms of the membrane index which the DT series do not see (cf.~\Cref{rmk: rank 2 DT C3} for a discussion of the case $X=\Aaff{3}$). It would be interesting to compare the correspondence with the recent reconstruction of rank-$r$ invariants from rank one, which also establishes a link to the Gromov--Witten theory of $X$ \cite{FT23:rankrrankoneDT}.

\subsubsection{The refined ``remodelled B-model''}

A long-standing problem in local mirror symmetry has been to find a place for the GV/PT refinement within the topological recursion formalism of \cite{Eynard:2007kz,Bouchard:2007ys,Fang:2016svw},
as early proposals to model the refinement around the theory of solutions of loop equations for the $\beta$-deformed matrix model ensembles \cite{Dijkgraaf:2009pc} were found to be invalid away from conifold loci where filling fractions vanish \cite{Brini:2010fc}. 
On the other hand, the identification in \cref{def:refGWintro} of refined GW invariants with the potential of a semi-simple CohFT automatically comes with a refinement of the ``remodelled B-model'' formalism of 
\cite{Bouchard:2007ys}, at least for the closed free energies: by \cite{Eynard:2011ga,Dunin-Barkowski:2012vii}, this takes the form of the Eynard--Orantin local topological recursion for superpotential germs defined, near each critical point, via the Laplace transform of the identity component of the equivariant 
$R$-matrix of the CY 5-fold $Z$.
Although the local superpotentials may not in general glue to a global spectral curve as in \cite{Bouchard:2007ys}, in the special cases where they do they can be seen to correctly reproduce the expressions for ``refined spectral curves'' heuristically found in  examples  by \cite{Eynard:2011vs} 
in the analysis of matrix integrals arising in the theory of the refined topological vertex \cite{IKV09:RefVert}: see \cref{rmk:refmm}. We defer to future work to find a  description for the open sector of the topological recursion in our setting, with a view towards 
generalising \eqref{eq:GWref} to refined open GW theory \cite{IKV09:RefVert,Aganagic:2011sg,Kozcaz:2018ndf}. 
\label{sec:refTR}

\subsubsection{The worldsheet definition of refined topological strings}

From a physics point of view, the identification of refined invariants of threefolds with equivariant invariants of fivefolds yields a first-principles definition of refined topological strings as a twist of an $\cN=(2,2)$ worldsheet theory in the original sense of \cite{Wit88:TopSigMod,Witten:1991zz}, without appealing to target-space string dualities \cite{AFHNZ13:WorldsheetRefTopStr}. The tree-level refined A-model on $X$ is the A-twist of the non-linear $\sigma$-model on $Z=X\times\Aaff{2}$, deformed by  potential terms generated by gauging the $\mgp{T}$-flavour symmetry and freezing the associated vector multiplets to complex twisted masses in the weak gauge coupling limit \cite{Alvarez-Gaume:1983uye,Labastida:1991yq,Hori:2000kt}. The corresponding refined topological string is then defined from the usual coupling to worldsheet topological gravity \cite{Dijkgraaf:1990qw}. For toric $X$,  by the Hori--Vafa construction \cite{Hori:2000ck,MR4047552}, the tree-level B-model mirror of the refined A-model on $X$ (=$\mgp{T}$-equivariant A-model on $Z$) would be a Landau--Ginzburg\footnote{As mentioned in \cref{sec:refTR}, a geometric CY B-model description (akin to the global spectral curve mirrors of \cite{Hori:2000kt}) may not be generally available in the refined setting.}  model on $\Aaff{5}$ deformed by logarithmic superpotential terms associated to the twisted masses of the torus $\mgp{T}$. It would also  be interesting to make direct contact
to the proposals for worldsheet refinement appearing in some of the physics literature \cite{AFHNZ13:WorldsheetRefTopStr,NO12:WorldsheetOmegaBg}, and also to shed light on the Rigidity \cref{conj: rigidity} and the refined direct integration (\cref{def:dirint}) from a physics standpoint.

\subsection{Organisation of the paper}
We describe the setup and largely fix notation for our treatment of refined Gromov--Witten invariants in \cref{sec: refined GW}, concluding with a discussion of the unrefined limit and the Rigidity Conjecture. In \cref{sec: local P1} we consider bundles over $\bbP^1$ as a first test case: we establish 
\cref{thm:A1} in \cref{sec: resolved conifold}, and consider aspects of non-equirigid examples in \cref{sec: local P1 shifted}. In \cref{sec:localSurf} we consider the specialisation to the NS limit of local del Pezzo surfaces, and prove \cref{thm:A2}. In \cref{sec: ref mirror sym KP2,sec: ref KP2 beyond LV} we study refined higher genus mirror symmetry for local $\bbP^2$. We show {\bf (Quasi-modularity}) and {\bf (Modular anomaly)} in \cref{thm:haegen,thm:haeref}, 
 prove in \cref{thm:crc} the Refined Crepant Resolution Correspondence  implying {\bf (Orbifold regularity)}, 
and establish the leading order of the conjectured universal relation to the Barnes double Gamma function of {\bf (Conifold asymptotics)} in \cref{thm:cfasym}, which collectively show \cref{thm:B}.
The refined GW/GV and GW/PT correspondences and general evidence in their favour are finally explored in \cref{sec:further implications}, including  
the proof of \cref{thm:A3}.
\cref{sec:resconRmat} contains a proof of the refined Harer--Zagier formula \eqref{eq:betagiv},
which is in turn key to the weak {(\bf Conifold Asymptotics)} statement of \cref{thm:B}.
\subsection*{Acknowledgements} 
A.~B.~ thanks the Galileo Galilei Institute for Theoretical Physics and the Institut Mittag-Leffler for hospitality, and INFN for partial support during the completion of this work. 
We acknowledge with thanks discussions and correspondence with P.~Bousseau, T.~Bridgeland, D.~Genlik, H.~Iritani, C.~Manolache, M.~Moreira, G.~Oberdieck, and W.~Pi. A.~B.~ was supported by the EPSRC Fellowship EP/S003657/2. Y.~S.~ was supported by the SNF grant SNF-200020-21936.
\section{Equivariant GW theory of Calabi--Yau fivefolds}
\label{sec: refined GW}

\subsection{Primer on equivariant Chow groups}

We start off by setting some notation. 
Let $\mgp{G}$ be a connected reductive group with a maximal torus $\mgp{T}$ and Weyl group $\mgp{W}$ acting on a Deligne--Mumford stack $M$. The $\mgp{G}$-equivariant Chow groups of $M$, as introduced in \cite{EG00:equivRR}, satisfy
\begin{equation}
	\label{eq: Chow G versus T W}
	\Chow_\bullet^{\mgp{G}}(M) \cong \Chow_\bullet^{\mgp{T}}(M)^{\mgp{W}}\,.
\end{equation}
We take Chow groups with coefficients in $\bbQ$. We will say that the $\mgp{G}$-action factors through the action of a group $\mgp{G}'$ on $M$ if there exists a group homomorphism $\mgp{G} \rightarrow \mgp{G}'$ compatible with the respective actions on $M$, i.e.~the following diagram commutes:
\begin{equation*}
	\begin{tikzcd}[column sep=1.2em]
		\mgp{G} \times M \ar[rr] \ar[dr] & & M\,.	\\
		& \mgp{G}' \times M \ar[ur]  & 
	\end{tikzcd}
\end{equation*}
In this situation we get an induced morphism of Chow groups
\begin{equation*}
	\Chow_\bullet^{\mgp{G}'}(M) \longrightarrow \Chow_\bullet^{\mgp{G}}(M)\,,
\end{equation*}
which is compatible with push-forwards, pull-backs, Chern classes and  (virtual) fundamental classes.\footnote{Following Edidin and Graham \cite{EG98:equivIntTh}, the morphism is given by flat pull-back along the approximation of $M\times_{\mgp{G}} \bbE\mgp{G} \times \bbE\mgp{G}' \rightarrow M\times_{\mgp{G}'} \bbE\mgp{G}'$. Compatibility can be argued as in \cite[Sec.~2.3 \& 2.4]{EG98:equivIntTh}.}
In the following, we will use the shorthand \[R_{\mgp{G}} \coloneqq \Chow^\bullet_{\mgp{G}}(\mr{pt})\] to denote the equivariant Chow cohomology of a point. Fixing an isomorphism $\mgp{T}\cong \Gm[m]$ we may identify
\begin{equation}
	\label{eq: Chow point trivialisation}
	R_{\mgp{T}} \cong \bbQ [\alpha_1,\ldots,\alpha_m]\,,
\end{equation}
where $\alpha_i$ is the first Chern class $c_1^{\mgp{T}}(\mchar{t}_i)$ of the natural one-dimensional $\Gm[m]$-representation
\begin{equation*}
	\mchar{t}_i \in K_{\mgp{T}}(\mr{pt}) = \Rep \mgp{T}
\end{equation*}
given by multiplication by the $i^{\rm th}$ coordinate of $\Gm[m]$. 
We will further write $R^{\mr{loc}}_{\mgp{G}}$ for the localisation of $R_{\mgp{G}}$ to the multiplicative subset of homogeneous elements of strictly positive degree: $R^{\mr{loc}}_{\mgp{G}}$ inherits the grading of $R_{\mgp{G}}$, and we shall denote $R_{\mgp{G},k}$ and $R^{\mr{loc}}_{\mgp{G},k}$ the respective $k^{\rm th}$ graded pieces, and
\begin{equation*}
    \widehat{R}_{\mgp{G}} \qquad \text{and} \qquad \widehat{R}^{\,\mr{loc}}_{\mgp{G}} 
\end{equation*}
the respective completions w.r.t.~the filtration induced by the grading. Under the isomorphism \eqref{eq: Chow point trivialisation}, we may identify
\begin{equation*}
	\widehat{R}_{\mgp{G}} \cong \big(\widehat{R}_{\mgp{T}}\big)^{\mgp{W}} \cong \bbQ \llbracket \alpha_1,\ldots,\alpha_m \rrbracket^{\mgp{W}}\,.
\end{equation*}

Suppose now that $M$ is equipped with a $\mgp{G}$-equivariant perfect obstruction theory, providing us with a virtual fundamental class $[M]^{\vir}_{\mgp{G}} \in \Chow_{\mr{vd}\,M}^{\mgp{G}}(M)$. If
the $\mgp{T}$-fixed locus of $M$ is proper, we may define the degree of this cycle via virtual torus localisation.
\begin{notation}
    \label{notation: defn via localisation}
    If $M^{\mgp{T}}$ is proper, we set\footnote{Of course, when $M$ is proper this relation is a theorem rather than a definition \cite{GP97:virtloc}, as the l.h.s.~is well-defined as an element in $R_{\mgp{G}} \subseteq R^{\mr{loc}}_{\mgp{G}}$.} 
    \begin{equation}
    	\label{eq: defn via localisation}
    	\intEquiv{\mgp{G}}_{[M]^{\vir}_{\mgp{G}}} \gamma \coloneqq \intEquiv{\mgp{T}}_{[M^{\mgp{T}}]^{\vir}_{\mgp{T}}} \frac{\gamma\vert_{M^{\mgp{T}}}}{e^{\mgp{T}}(N^{\vir})} \qquad \in \big(R^{\mr{loc}}_{\mgp{T}, \mr{vd}\,M}\big)^{\mgp{W}} \cong R^{\mr{loc}}_{\mgp{G}, \mr{vd}\,M}
    \end{equation}
    where $N^{\vir}$ is the virtual normal bundle of $M^{\mgp{T}}\hookrightarrow M$.
\end{notation} 
Note that $M^{\mgp{T}}$ is preserved by the action of the Weyl group, and the moving and fixed part of the perfect obstruction theory are $\mgp{W}$-invariant. Hence, via \eqref{eq: Chow G versus T W} the right-hand side integral  yields an element in $R^{\mr{loc}}_{\mgp{G}}$.

\subsection{Calabi--Yau fivefolds}
Let now $Z$ be a smooth quasi-projective Calabi--Yau fivefold, and 
suppose that $\mgp{G}$ is a connected reductive group acting effectively on $Z$. This action naturally lifts to $\Mbar_g(Z,\beta)$, and the usual perfect obstruction theory is readily checked to be $\mgp{G}$-equivariant. Let $\mgp{T}$ be the maximal torus of $\mgp{G}$, and suppose that the $\mgp{T}$-fixed locus $\Mbar_g (Z,\beta)^\mgp{T}$ is proper.

\begin{defn}
	\label{defn: refined GW invariant}
	We define the \textbf{genus-$g$ $\mgp{G}$-equivariant Gromov--Witten invariants of $Z$ in class} $\beta$ as
	\begin{equation}
		\label{eq: term refined GW}
		\refGWgb{g}{Z}{\mgp{G}}{\beta} \coloneqq \intEquiv{\mgp{G}}_{[\Mbar_g (Z,\beta)]^{\vir}_{\mgp{G}}} 1 \coloneqq \intEquiv{\mgp{T}}_{[\Mbar_g (Z,\beta)^{\mgp{T}}]^{\vir}_{\mgp{T}}} \frac{1}{e^{\mgp{T}}(N^{\vir})} \qquad \in R^{\mr{loc}}_{\mgp{G},2g-2} \,.
	\end{equation}
\end{defn}


Since each equivariant Gromov--Witten invariant is an element of the $(2g-2)$-graded piece of $R^{\mr{loc}}_{\mgp{G}}$, summing these invariants over all genera we obtain a well-defined element in the completed ring. 

\begin{defn}
	\label{defn: refined GW series}
	We define the fixed degree $\beta$ (resp.~fixed genus-$g$) \define{$\mgp{G}$-equivariant Gromov--Witten generating series} as
	\begin{alignat*}{2}
		\refGWb{Z}{\mgp{G}}{\beta} \coloneqq & \sum_{g\geq 0} ~ \intEquiv{\mgp{G}}_{[\Mbar_g (Z,\beta)]^{\vir}_{\mgp{G}}} 1 && \qquad \in \widehat{R}^{\,\mr{loc}}_{\mgp{G}}\,, \nn \\
		\refGWg{g}{Z}{\mgp{G}} \coloneqq & \sum_{\beta \neq 0} ~\intEquiv{\mgp{G}}_{[\Mbar_g (Z,\beta)]^{\vir}_{\mgp{G}}} \cQ^\beta && \qquad \in R^{\mr{loc}}_{\mgp{G},2g-2}\llbracket \cQ \rrbracket\,, 
	\end{alignat*}
	where the sum in the last equality is taken over all non-zero effective curve classes in $Z$.
\end{defn}

By design, the invariants enjoy the following functorality property.

\begin{lem}[Contravariance] \label{lem: GW contravariance}
	Let $\mgp{G}$ and $\mgp{G}'$ act on $Z$ and $\mgp{G} \rightarrow \mgp{G}'$  a group homomorphism compatible with the actions on $Z$. Then
	\begin{equation*}
		\refGWgb{g}{Z}{\mgp{G}'}{\beta} \longmapsto \refGWgb{g}{Z}{\mgp{G}}{\beta}
	\end{equation*}
	under the induced morphism $R^{\mr{loc}}_{\mgp{G}'} \rightarrow R^{\mr{loc}}_{\mgp{G}}$ of Chow rings.\qed
\end{lem}

\subsection{Threefolds in local Calabi--Yau fivefolds}
\label{sec: refined threefold setup}
Let us now specialise our discussion to our 
main case of interest. We will take the Calabi--Yau fivefold $Z$ to be a product
\begin{equation*}
	Z = X \times \Aaff{2}
\end{equation*}
with $X$ a smooth complex quasi-projective Calabi--Yau threefold. Throughout the rest of this section, we will assume that we are given an action of a connected reductive group $\mgp{G}$ with maximal torus $\mgp{T}$ on $X$ and $\Aaff{2}$ providing us with a $\mgp{G}$-action on $Z$. 
We will write $\kappa$ and $\Kappa$ for the $\mgp{G}$-character (equivalently, its associated one-dimensional representation) associated to the induced action on $\omega_{X}^{-1}\cong \mc{O}_X$, respectively $\omega_{Z}^{-1}\cong \mc{O}_Z$. We say that the $\mgp{G}$-action on $X$ (resp.~$Z$) is \textbf{Calabi--Yau} if $\kappa$ (resp.~$\Kappa$) are trivial $\mgp{G}$-representations. For instance, the natural $\SL(2,\bbC)$-action on $\Aaff{2}$ lifts to a Calabi--Yau action on $Z$. We will usually refer to its maximal torus
\begin{alignat*}{2}
	\Gm_{\mchar{q}_-} & \lhook\joinrel\longrightarrow {}&& \SL(2,\bbC) \, , \\ \mchar{q}_- &\longmapsto {}&& \begin{psmallmatrix}
		\mchar{q}_-^{-1} & \, \\
		\, & \mchar{q}_-
	\end{psmallmatrix}
\end{alignat*}
as the \define{anti-diagonal torus} acting on $\Aaff{2}$ and $Z$.

Moreover, we are always going to assume that $\mgp{T}$ acts on the two affine factors in $X\times \Aaff{2}$ via characters 
\begin{equation*}
\mchar{q}_1^{-1} \qquad \text{and}\qquad \mchar{q}_2^{-1}
\end{equation*}
and we will write
\begin{equation*}
	\epsilon_1 \coloneqq c^{\mgp{T}}_1(\mchar{q}_1) \qquad \text{and}\qquad \epsilon_2 \coloneqq c^{\mgp{T}}_1(\mchar{q}_2)
\end{equation*}
for the associated elements in $R_{\mgp{T}}$. For convenience, we also introduce
\begin{equation*}
	\epsilon_\pm \coloneqq \frac{\epsilon_1 \pm \epsilon_2}{2}\,.
\end{equation*}
This notation is chosen so that, whenever the action of the anti-diagonal torus factors through the action of $\mgp{G}$ on $\Aaff{2}$, then under the induced morphism $R_{\mgp{G}} \rightarrow R_{\Gm_{\mchar{q}_-}}$ we have
\begin{equation*}
	\epsilon_- \mapsto c_1^{\Gm} \! (\mchar{q}_-)\,.
\end{equation*}

\subsubsection{Maps to local threefolds}
Since the composition of a stable map to $Z=X\times \Aaff{2}$ with the projection to the affine plane is constant, we have
\begin{equation*}
	\Mbar_{g}(Z,\beta) = \Mbar_{g}(X,\beta) \times \Aaff{2}\,,
\end{equation*}
and for the equivariant Gromov--Witten invariants of this target to be well-defined we shall need to assume that the $\mgp{T}$-action on the affine directions is non-trivial, or in other words $\epsilon_1$ and $\epsilon_2$ are non-zero.

Writing $N_{-}^{\vir}$ for the virtual normal bundle of $\Mbar_{g}(X,\beta)\times \{(0,0)\}\hookrightarrow \Mbar_{g}(Z,\beta)$ we can rewrite 
\eqref{eq: term refined GW} as
\begin{equation*}
	\refGWgb{g}{Z}{\mgp{G}}{\beta} = \intEquiv{\mgp{T}}_{[\Mbar_{g}(X,\beta)]^{\vir}_{\mgp{T}}} \frac{1}{e^{\mgp{T}}(N^{\vir}_{-})}\,.
\end{equation*}
As an element of $\mgp{T}$-equivariant K-theory, the virtual normal bundle is 
\begin{equation*}
	N^{\vir}_{-} = \big(\mchar{q}_1 - \bbE^{\vee}_g \mchar{q}_1 \big) + \big(\mchar{q}_2 - \bbE^{\vee}_g \mchar{q}_2 \big)\,,
\end{equation*}
where $\bbE_g$ is the Hodge bundle on $\Mbar_{g}(X,\beta)$. Denoting
by $\lambda_k\coloneqq c_k (\bbE_g)$ its Chern classes, we have
\begin{equation*}
	\frac{1}{e^{\mgp{T}}\big(\mchar{q}_i - \bbE_g^\vee \mchar{q}_i\big)} = \sum_{k\geq 0} \lambda_k (-1)^k \epsilon_i^{g-k-1} \eqqcolon \HodgeLambda[g]{\epsilon_i}\,,
\end{equation*}
and therefore
\begin{equation}
	\label{eq: ref GW simplified}
	\refGWgb{g}{Z}{\mgp{G}}{\beta} = \intEquiv{\mgp{T}}_{[\Mbar_{g}(X,\beta)]^{\vir}_{\mgp{T}}} \HodgeLambda[g]{\epsilon_1}\,  \HodgeLambda[g]{\epsilon_2} \,.
\end{equation}
The above presentation has two immediate consequences which we are going to discuss below.

\subsection{The unrefined limit}
\label{sec:unrefGW}
The first observation will be that the equivariant Gromov--Witten invariants of $X\times \Aaff{2}$ naturally refine the ones of $X$. Since the Hodge bundle vanishes in genus zero, we have the equality
\begin{equation*}
	\refGWgb{0}{Z}{\mgp{G}}{\beta} = (\epsilon_1 \epsilon_2)^{-1} ~ \intEquiv{\mgp{G}}_{[\Mbar_{0}(X,\beta)]^{\vir}_{\mgp{G}}}1\,.
\end{equation*}
Put another way, there is \textit{no} refinement in genus zero. In higher genus, we make a simple observation concerning the special cases when the action is either trivial on $\omega_{\Aaff{2}}$, or pointwise trivial on $X$.

\begin{prop}
	\label{prop: unref limit}
	If $\epsilon_+=0$, or if $\mgp{T}$ acts trivially on $X$, we have
	\begin{equation*}
		\refGWgb{g}{Z}{\mgp{G}}{\beta} = (\epsilon_1 \epsilon_2)^{g-1} ~~ \intEquiv{\mgp{G}}_{[\Mbar_{g}(X,\beta)]^{\vir}_{\mgp{G}}} 1\,.
	\end{equation*}
\end{prop}

The first part of the \namecref{prop: unref limit} is an immediate consequence of the following classical statement.

\begin{lem}[Mumford's relation]
	\label{lem: Mumfords relation}
	 We have \cite[Sec.~5]{Mu83}
	\begin{equation*}
		\HodgeLambda[g]{\epsilon} \, \HodgeLambda[g]{-\epsilon} = (-\epsilon^2)^{\mathrm{rk}\, \bbE_{\cdot} -1}\,.
	\end{equation*}
\end{lem}

\begin{proof}[Proof of \Cref{prop: unref limit}]
	The first part of the \namecref{prop: unref limit} ($\epsilon_+=0$) follows from applying \cref{lem: Mumfords relation} to \eqref{eq: ref GW simplified}. As for the second part, note that if $\mgp{T}$ acts trivially on $X$, then $[\Mbar_{g}(X,\beta)]^{\vir}_{\mgp{T}}$ coincides with its non-equivariant limit under the isomorphism
	\begin{equation*}
		\Chow_\bullet^{\mgp{T}}\big(\Mbar_{g}(X,\beta)\big) \cong \Chow_\bullet\big(\Mbar_{g}(X,\beta)\big) \otimes R_{\mgp{T}}\,.
	\end{equation*}
	As $\mr{vd}\, \Mbar_{g}(X,\beta) = 0$, only the leading order $\HodgeLambda[g]{\epsilon_1}\,  \HodgeLambda[g]{\epsilon_2} = (\epsilon_1 \epsilon_2)^{g-1} + \ldots$ may contribute non-trivially to the integral, from which the claim follows.
\end{proof}

By contravariance, one corollary of \Cref{prop: unref limit} is that the all-genus equivariant Gromov--Witten generating series of the fivefold $X\times \Aaff{2}$ is a one-parameter deformation of the higher genus GW generating series of the threefold $X$. To see this, suppose that $\epsilon_1,\epsilon_2$ are linearly independent. Fixing an isomorphism $R_{\mgp{G}}\cong \bbQ[\epsilon_1,\epsilon_2,\alpha_1,\ldots, \alpha_m]^{\mgp{W}}$ we may view the all-genus generating series as a formal power series in the equivariant parameters.

\begin{cor}
	\label{cor: GW unref limit}
	Suppose $\epsilon_1$, $\epsilon_2$ are linearly independent and set $\mgp{G}' \coloneqq \ker \mchar{q}_1 \mchar{q}_2$. Then, as formal power series,
	\begin{flalign*}
		&&\refGWb{Z}{\mgp{G}}{\beta} ~ \Big\vert_{\epsilon_1 = -\epsilon_2 = \epsilon_-} =  \sum_{g \geq 0} (-\epsilon_-^2)^{g-1} ~~  \intEquiv{\mgp{G}^\prime\!}_{[\Mbar_{g}(X,\beta)]^{\vir}_{\mgp{G}^\prime}} 1\,.&&\qed
	\end{flalign*}
\end{cor}

Motivated by this observation, we will refer to the restriction induced by the inclusion $\mgp{G}^\prime \hookrightarrow \mgp{G}$ as the \textbf{unrefined limit}.

\begin{rmk}
	\label{rmk: rigid unref limit}
	In the statement of the \namecref{cor: GW unref limit} we implicitly assumed that $\Mbar_{g}(X,\beta)^{\mgp{T}^\prime}$ is proper for all $g\geq 0$ in order for the integral on the r.h.s.~to be well-defined. This is of course trivially the case if $\Mbar_{g}(X,\beta)$ is proper itself, in which case
	\begin{equation*}
		\intEquiv{\mgp{G}^\prime\!}_{[\Mbar_{g}(X,\beta)]^{\vir}_{\mgp{G}^\prime}} 1 = \int_{[\Mbar_{g}(X,\beta)]^{\vir}} 1\,.
	\end{equation*}
\end{rmk}

\begin{rmk}
	Another lesson we learn from \Cref{prop: unref limit} is that the equivariant Gromov--Witten invariants of $X\times \Aaff{2}$ only refine the ones of the 
 threefold $X$ if the latter admits a non-trivial torus action. Therefore, a GW-theoretic refinement for compact Calabi--Yau threefolds would require new ideas. 
\end{rmk}

\subsection{Rigidity}
Another consequence of the presentation \eqref{eq: ref GW simplified} concerns the general form of the equivariant Gromov--Witten invariants of local Calabi--Yau threefolds. The insertion $\HodgeLambda[g]{\epsilon_1} \HodgeLambda[g]{\epsilon_2}$ on the right-hand side of equation \eqref{eq: ref GW simplified} is polynomial in $\epsilon_1$ and $\epsilon_2$, up to an overall factor $(\epsilon_1 \epsilon_2)^{-1}$. If the moduli stack of stable maps to $X$ is proper, this implies a polynomial dependence on all equivariant parameters.
\begin{lem}
	\label{lem: shape of refGW series}
	If $\Mbar_{g}(X,\beta)$ is proper we have
	\begin{flalign*}
		&&\epsilon_1 \epsilon_2 \, \refGWgb{g}{Z}{\mgp{G}}{\beta} \in R_{\mgp{G},2g}\,.&&\qed
	\end{flalign*}
\end{lem}

The \namecref{lem: shape of refGW series} motivates the following terminology.

\begin{defn}
    \label{defn: rigid}
	We call an effective curve class $\beta$ in $X$ \textbf{rigid} if $\Mbar_g(X,\beta)$ is proper for all $g\geq 0$. We say $X$ is \textbf{equirigid} if all non-zero effective curve classes in $X$ are rigid.
\end{defn}


The class of Calabi--Yau threefolds which are equirigid, but non-proper, includes the resolved conifold and all local del Pezzo surfaces.

Now, suppose 
that $\epsilon_1$ and $\epsilon_2$ are linearly independent in $R_{\mgp{T}}$. Then $\mgp{T}\cong \Gm[2+m]$ for some $m\geq 0$ and we may choose additional generators such that $R_{\mgp{G}}\cong \bbQ[\epsilon_1,\epsilon_2,\alpha_1,\ldots,\alpha_m]^{\mgp{W}}$. We may reformulate the statement of \Cref{lem: shape of refGW series} as
\begin{equation*}
	\epsilon_1 \epsilon_2 \, \refGWgb{g}{Z}{\mgp{G}}{\beta} \in \bbQ [\epsilon_1,\epsilon_2,\alpha_1,\ldots,\alpha_{m} ]_{2g}^{\mgp{W}}\,.
\end{equation*}
For $m=0$, the dependence on equivariant parameters takes a particularly simple form.

\begin{lem}
	\label{lem: shape of refGW series 2}
	Suppose $\Mbar_{g}(X,\beta)$ is proper and $\mgp{G}$ acts such that
	\begin{enumerate}
		\item $\mgp{T}\cong \Gm[2]$;
		\item the $\mgp{G}$-action on $X\times \Aaff{2}$ is Calabi--Yau;
		\item the $\mgp{G}$-action on $X$ is \textit{not} Calabi--Yau;
		\item \label{item: contains anti diag} the action of the anti-diagonal torus on $X\times \Aaff{2}$ factors through the $\mgp{G}$-action.
	\end{enumerate}
	Then
	\begin{equation*}
		\epsilon_1 \epsilon_2 \, \refGWgb{g}{Z}{\mgp{T}}{\beta} \in \bbQ [\epsilon_-^2, \epsilon_+^2 ]_{2g} \cong \bbQ [ \epsilon_1\cdot  \epsilon_2 , \epsilon_+^2 ]_{2g} \,.
	\end{equation*}
\end{lem}
\begin{proof}
	Assumptions (ii) and (iii) guarantee that $\epsilon_-$ and $\epsilon_+$ are linearly independent. Hence, by (i) these elements generate $R_{\mgp{T}}$. Moreover, by (iv) the only dependence on $\epsilon_-$ comes from the insertion $\HodgeLambda[g]{\epsilon_1}\,  \HodgeLambda[g]{\epsilon_2}$. The \namecref{lem: shape of refGW series 2} then follows as the latter is invariant under the exchange $(\epsilon_1,\epsilon_2) \mapsto (\epsilon_2,\epsilon_1)$ (which is the Weyl group action related to the full $\GL(2,\bbC)$ symmetry) and therefore can only depend on $\epsilon_-^2$.
\end{proof}

If $m>0$ one should in principle expect a dependence on the additional equivariant parameters. 
We claim that, remarkably, this is not the case if $\beta$ is rigid and $\mgp{G}$ preserves $\omega_Z$.

\begin{conj}[Rigidity]
	\label{conj: rigidity}
	Suppose $\beta$ is a rigid curve class in $X$ and the $\mgp{G}$-action on $Z$ is Calabi--Yau. Then, as a formal Laurent series, $\refGWb{Z}{\mgp{G}}{\beta}$ only depends on $\epsilon_1$ and $\epsilon_2$. In other words,
	\begin{equation*}
		\begin{split}
			\epsilon_1 \epsilon_2 \, \refGWb{Z}{\mgp{G}}{\beta} \quad \in \quad \bbQ \llbracket \epsilon_-^2 , \epsilon_+^2 \rrbracket_{2g}  \quad \subsetneq \quad \widehat{R}_{\mgp{G},2g}\,.
		\end{split}
	\end{equation*}
\end{conj}

Evidence for this Rigidity Conjecture comes from the unrefined limit.

\begin{prop}
    \label{prop: rigidity unref limit}
	\Cref{conj: rigidity} holds if furthermore 
 the $\mgp{G}$-action on $X$ is Calabi--Yau.
\end{prop}
\begin{proof}
The $\mgp{G}$-action being separately Calabi--Yau on $X$ and $Z$ implies in particular that it is Calabi--Yau on $\bbC^2$, hence $\epsilon_+=0$. The claim follows from combining \Cref{prop: unref limit} with \Cref{rmk: rigid unref limit}.
\end{proof}
A few comments are in order regarding the two conditions (the rigidity of the curve class, and the triviality of the action on $\omega_Z$) in  the statement of \cref{conj: rigidity}. 

\begin{rmk}
To start with, it is worth emphasising that \emph{both} of the stated conditions in \cref{conj: rigidity} are 
necessary for the statement to hold: in \Cref{sec: resolved conifold} we will discuss an example where, for non-Calabi--Yau actions, we indeed observe a dependence on additional equivariant parameters, and in \Cref{sec: local P1 shifted} we will show that the statement of the Conjecture is generally false if one drops the condition that the curve class is rigid.
\end{rmk}

\begin{rmk}
From a Gromov--Witten perspective, it does not appear {\it a priori} immediate to conceptually justify  the importance of the Calabi--Yau condition for the action on $Z$ in \Cref{conj: rigidity}. On the other hand, this requirement will become much more natural in relation to the refined GW/PT correspondence presented in \Cref{sec: refined PT-GW}, and indeed the reader familiar with \cite{NO14:membranes} may have been reminded of a similar rigidity phenomenon occurring in K-theoretic Pandharipande--Thomas theory. In particular, in \Cref{sec: rigidity implication} we will see that \cite[Prop.~7.5]{NO14:membranes} implies \Cref{conj: rigidity}, assuming the refined GW/PT correspondence stated in \cref{conj: refined PT GW}. 
%
 Our discussion of refined PT theory in \Cref{sec: refined PT-GW} will also make the rigidity condition for the curve class conceptually more transparent.
\end{rmk} 

\begin{rmk}
	A proof of \Cref{conj: rigidity} also seems desirable from a practical point of view, as it implies that the restriction to \textit{any} two-dimensional subtorus $\Gm[2] \hookrightarrow \mgp{G}$ satisfying the properties\footnote{The last condition \labelcref{item: contains anti diag} can actually be relaxed. It is sufficient to demand that the action of the anti-diagonal torus on $\Aaff{2}$ factors through $\mgp{G}$. In other words $\Gm_{\mchar{q}_-}$ may act non-trivially on $X$.} listed in \Cref{lem: shape of refGW series 2} recovers the fully equivariant generating series
	\begin{equation*}
		\refGWb{Z}{\mgp{G}}{\beta} = \refGWb{Z}{\Gm[2]}{\beta}
	\end{equation*}
	where both sides are viewed as formal power series in $\epsilon_1$ and $\epsilon_2$. The flexibility in the choice of the restriction is often key to some dramatic simplifications in the calculations: we will witness two instances of such a phenomenon in this paper in \cref{prop:ResConifold,thm:haeref}.
\end{rmk}

\section{Local \texorpdfstring{$\bbP^1$}{P1} geometries}
\label{sec: local P1}
We now move on to a discussion of 
some 
concrete case studies. In this section we will explicitly compute the equivariant Gromov--Witten generating series of two local $\bbP^1$ geometries.

\subsection{The resolved conifold}
\label{sec: resolved conifold}

Let us consider the product of the resolved conifold and $\Aaff{2}$
\begin{equation*}
	\begin{tikzcd}
		Z &[-2.5em] = &[-2.5em] \Tot \big( \mc{O}_{\bbP^1}(-1) \, \oplus \, \mc{O}_{\bbP^1}(-1) \, \oplus \, \mc{O}_{\bbP^1} \, \oplus \,  \mc{O}_{\bbP^1} \big) \arrow[d] \\
		& & \bbP^1
	\end{tikzcd}
\end{equation*}
together with the action of the dense open torus $\Tmax \cong \Gm[5]$ scaling the fibres over the fixed point $0\in \bbP^1$ with weights $\alpha_1,\alpha_2,-\epsilon_1,-\epsilon_2 \in R_{\Tmax}$ respectively. Together with the tangent weight of $0\in \bbP^1$, which we are going to denote by $\alpha_0$, these weights freely generate $R_{\Tmax}$.

Since our target is toric, for low degree and low genus one can efficiently compute its equivariant Gromov--Witten invariants via virtual localisation \cite{GP97:virtloc}. For instance in primitive degree $\beta = [\bbP^1]$ a computer calculation reveals that
\begin{equation*}
	\begin{split}
		\epsilon_1 \epsilon_2 \, & \refGWb{Z}{\Tmax}{[\bbP^1]}\\
		= 1 & - \frac{1}{12} \big( \epsilon_+^2 +\epsilon_-^2  \big) + \frac{1}{720} \big( 3 \epsilon_+^4 + 8 \epsilon_+^2\epsilon_-^2 + 3 \epsilon_-^4\big) \\
		& + (\epsilon_1 + \epsilon_2 - \textstyle{\sum_i} \displaystyle \alpha_i) \, \epsilon_+\, \left(\frac{1}{12} + \frac{1}{1440} \big(2 \alpha_1 \alpha_2 + (\alpha_1+\alpha_2) \alpha_0 +\alpha_0^2 - 7 (\alpha_1 + \alpha_2 + \alpha_0) \epsilon_+ -14 \epsilon_-^2\big)\right)  + \ldots
	\end{split}
\end{equation*}
up to contributions from genus $g>2$, and where we collected terms proportional to
\begin{equation*}
	c_1^{\Tmax}(\Kappa) = \epsilon_1 + \epsilon_2 - \textstyle{\sum_i} \displaystyle \alpha_i
\end{equation*}
in the second row. Passing to the four dimensional Calabi--Yau subtorus $\TCY \subseteq \ker \Kappa$ this weight is set to zero, and in agreement with \Cref{conj: rigidity} the coefficients of the expansion of $\refGWb{Z}{\TCY}{[\bbP^1]}$ only depend on $\epsilon_1$ and $\epsilon_2$. The fact that the term in the second line equally vanishes when passing to a subtorus on which $\epsilon_+=0$ is an instance of \Cref{prop: rigidity unref limit}.

Upon restriction to an appropriate three-dimensional subtorus of $\TCY$, we can give a closed-form solution for the equivariant Gromov--Witten generating series of $X$.

\begin{prop}
	\label{prop:ResConifold}
	Let $\mgp{T} \hookrightarrow \TCY$ be a subtorus for which $\alpha_i = \epsilon_j$ for some $i,j\in\{1,2\}$ in $R_{\mgp{T}}$. Then, for all $d>0$, we have
	\begin{equation}
		\refGWb{Z}{\mgp{T}}{d[\bbP^1]}  = \frac{1}{d} \left(2 \sinh  \frac{d \epsilon_1}{2}\right)^{-1} \left(2 \sinh  \frac{d \epsilon_2}{2}\right)^{-1} \,.
    \label{eq:resconGWref}
	\end{equation}
\end{prop}

\begin{proof}
	We will only discuss the case $\alpha_1=\epsilon_1$ here since all other cases can be treated similarly. We observe that under this assumption, the CY condition $c_1^{\mgp{T}}(\Kappa)=0$ fixes the $\mgp{T}$-weight on the fibre over the fixed point $\infty\in\bbP^1$ of the second $\mc{O}_{\bbP^1}(-1)$ bundle to be $-\epsilon_2$. We degenerate the base curve into a union of two $\bbP^1$s glued along a point $p$ in a way that is compatible with the $\mgp{T}$-action; in particular, $p$ is a $\mgp{T}$-fixed point. We extend the two negative line bundles to the special fibre of the degeneration as displayed below:
	\begin{equation*}
		\begin{tikzpicture}
			\draw[thick] (3,0.5) -- (-0.6,-0.1);
			\draw[thick] (-3,0.5) -- (0.6,-0.1);
			\draw[thick] (-8.6,0) -- (-5,0);
			\node[above] at (-6.8,1) {\small $\mc{O}_{\bbP^1}(-1)\oplus\mc{O}_{\bbP^1}(-1)$};
			\draw[->] (-6.8,1) -- (-6.8,0.2);
			\node[above] at (-1.4,1) {\small $\mc{O}_{\bbP^1}\oplus\mc{O}_{\bbP^1}(-1)$};
			\draw[->] (-1.4,1) -- (-1.4,0.4);
			\node[above] at (1.4,1) {\small $\mc{O}_{\bbP^1}(-1)\oplus\mc{O}_{\bbP^1}$};
			\draw[->] (1.4,1) -- (1.4,0.4);
			\node at (-4,0.5) {$\rightsquigarrow$};
			%
			\node[above] at (0,0) {\small$p$};
		\end{tikzpicture}
	\end{equation*}
	The degeneration formula \cite{Li02:Degen} expresses $\refGWb{Z}{\mgp{T}}{d[\bbP^1]}$ as a convolution of the relative Gromov--Witten invariants of the irreducible components of the special fibre:
	\begin{equation*}
		\begin{split}
			\refGWb{Z}{\mgp{T}}{d[\bbP^1]} = \sum_{\Gamma_1,\Gamma_2}  \frac{\prod_i d_i}{|\Aut(\Gamma_1,\Gamma_2)|} & \, (\epsilon_1 \epsilon_2)^{2 |\mathbf{d}|} \\
			& \times \intEquiv{\text{$\mgp{T}$}}_{[\Mbar_{\Gamma_1}(\bbP^1/ p)]^{\vir}_{\text{$\mgp{T}$}}} \HodgeLambda{-\epsilon_1} ~ e^{\text{$\mgp{T}$}}\!\big(\mathbf{R}^{1} \pi_*f^*\mc{O}_{\bbP^1}(-1)\big)\,  \HodgeLambda{\epsilon_1} \, \HodgeLambda{\epsilon_2}\\
			& \times\intEquiv{\text{$\mgp{T}$}}_{[\Mbar_{\Gamma_2}(\bbP^1/ p)]^{\vir}_{\text{$\mgp{T}$}}} e^{\text{$\mgp{T}$}}\!\big(\mathbf{R}^{1} \pi_*f^*\mc{O}_{\bbP^1}(-1)\big) \, \HodgeLambda{-\epsilon_2}\, \HodgeLambda{\epsilon_1} \, \HodgeLambda{\epsilon_2}\,.
		\end{split}
	\end{equation*}
	Here, $\Mbar_{\Gamma_j}(\bbP^1/ p)$ is Li's moduli stack of type $\Gamma_j$ stable maps to expanded degenerations of $(\bbP^1/ p)$ \cite{Li01:RelStabMaps}. The sum is over all possible types $(\Gamma_1,\Gamma_2)$ a stable map can split over the degenerated target and $\mathbf{d}$ is a partition of $d$ encoding the ramification profile over $p$. The factor $(\epsilon_1 \epsilon_2)^{2|\mathbf{d}|}$ is due to the normalisation exact sequence at the contact markings. Now, by Mumford's relation (\Cref{lem: Mumfords relation}), we have
	\begin{equation*}
		\HodgeLambda{-\epsilon_i} \, \HodgeLambda{\epsilon_i} = (-\epsilon_i^2)^{-\chi_j}\,,
	\end{equation*}
	where $\chi_j$ is the holomorphic Euler characteristic of the domain curves of relative stable maps of type $\Gamma_j$: this is where we crucially use the assumption that $\alpha_1=\epsilon_1$. The expression hence simplifies to
	\begin{equation*}
		\begin{split}
			\refGWb{Z}{\mgp{T}}{d[\bbP^1]} = \sum_{\Gamma_1,\Gamma_2} \frac{\prod_i d_i}{|\Aut(\Gamma_1,\Gamma_2)|} (-1)^{-\chi_1}\,\epsilon_1^{2|\mathbf{d}|-2\chi_1} &~  \intEquiv{\text{$\mgp{T}$}}_{[\Mbar_{\Gamma_1}(\bbP^1/ p)]^{\vir}_{\text{$\mgp{T}$}}} e^{\text{$\mgp{T}$}}\!\big(\mathbf{R}^{1} \pi_*f^*\mc{O}_{\bbP^1}(-1)\big)\, \HodgeLambda{\epsilon_2}\\
			\times (-1)^{-\chi_2}\,\epsilon_2^{2|\mathbf{d}|-2\chi_2} &~ \intEquiv{\text{$\mgp{T}$}}_{[\Mbar_{\Gamma_2}(\bbP^1/ p)]^{\vir}_{\text{$\mgp{T}$}}} e^{\text{$\mgp{T}$}}\!\big(\mathbf{R}^{1} \pi_*f^*\mc{O}_{\bbP^1}(-1)\big)\, \HodgeLambda{\epsilon_1} \,.
		\end{split}
	\end{equation*}
	A quick dimension count as in the proof of \cite[Lem.~6.3]{BP08:GWLocCurves} shows that all terms in the above sum vanish except those with a single contact marking, i.e.~$\mathbf{d}=(d)$. Therefore, the expression further simplifies to
	\begin{equation*}
		\label{eq: res conf intermed}
		\begin{split}
			\refGWb{Z}{\mgp{T}}{d[\bbP^1]} = d \sum_{g_1 \geq 0} (-1)^{g_1-1}\,\epsilon_1^{2g_1} &~  \intEquiv{\text{$\mgp{T}$}}_{[\Mbar_{g_1,(d)}(\bbP^1/ p,d[\bbP^1])]^{\vir}_{\text{$\mgp{T}$}}} e^{\text{$\mgp{T}$}}\!\big(\mathbf{R}^{1} \pi_*f^*\mc{O}_{\bbP^1}(-1)\big)\,  \HodgeLambda[g]{\epsilon_2}\\
			\times \sum_{g_2 \geq 0} (-1)^{g_2-1}\,\epsilon_2^{2g_2}  &~ \intEquiv{\text{$\mgp{T}$}}_{[\Mbar_{g_2,(d)}(\bbP^1/ p,d[\bbP^1])]^{\vir}_{\text{$\mgp{T}$}}} e^{\text{$\mgp{T}$}}\!\big(\mathbf{R}^{1} \pi_*f^*\mc{O}_{\bbP^1}(-1)\big)\, \HodgeLambda[g]{\epsilon_1} \,.
		\end{split}
	\end{equation*}
	The Gromov--Witten invariants appearing on the r.h.s.~have been computed in \cite[Thm.~5.1]{BP05:CurvesCY3TQFT} and read
	\[
		\intEquiv{\text{$\mgp{T}$}}_{[\Mbar_{g,(d)}(\bbP^1/ p,d[\bbP^1])]^{\vir}_{\text{$\mgp{T}$}}} e^{\text{$\mgp{T}$}}\!\big(\mathbf{R}^{1} \pi_*f^*\mc{O}_{\bbP^1}(-1)\big) \, \HodgeLambda[g]{\epsilon_i} = \frac{(-1)^{d-1} }{d \epsilon_i} ~ \big[u^{2g-1}\big] \left(2 \sin \frac{d u}{2} \right)^{-1}\,.
	\]
	Inserting this into our last expression for $\refGWb{Z}{\mgp{T}}{d[\bbP^1]}$ we arrive at formula \eqref{eq:resconGWref}.
\end{proof}

\begin{rmk}
	The proof heavily relied on the assumption that $\alpha_i = \epsilon_j$ for some $i,j\in\{1,2\}$. However, \Cref{conj: rigidity}, applied to the geometry at hand, predicts that the formula on the right-hand side of \eqref{eq:resconGWref} should also more generally hold for $\refGWb{Z}{\TCY}{d[\bbP^1]}$, and we verified this claim in a computer calculation for $d\leq 3$ up to contributions from genus five stable maps. Also, note that in the unrefined limit $\epsilon_1 = -\epsilon_2 = \epsilon_-$ we recover the well-known formula for the Gromov--Witten invariants of the resolved conifold \cite[Thm.~3]{MR1728879}:
	\begin{equation*}
		\label{eq: res conifold unref}
		\refGWb{Z}{\Gm_{\mchar{q}_-}}{d [\bbP^1]} = -\frac{1}{d} \left(2 \sinh \frac{d \epsilon_-}{2}\right)^{-2}\,.
	\end{equation*}
\end{rmk}

\subsection{A shifted example}
\label{sec: local P1 shifted}
The ideas of the proof of \Cref{prop:ResConifold} apply equally well to the equivariant Gromov--Witten theory of the threefold
\begin{equation*}
	X = \Tot \mc{O}_{\bbP_1}\! (-2) \oplus \mc{O}_{\bbP_1}\!(0)\,.
\end{equation*}
We again denote by $\TCY$ the Calabi--Yau subtorus of the dense open torus $\Tmax$ acting on $X\times \Aaff{2}$. For conformity, let us write $\epsilon_0$ for the tangent $\TCY$-weight at the origin of the $\Aone$-factor of $X$. Under restriction to appropriate three-dimensional subtori $\mgp{T} \hookrightarrow \TCY$ we are again able to compute the equivariant generating series explicitly.

\begin{prop}
	\label{prop: shifted local P1}
	Let $d>0$ and $\mgp{T} \hookrightarrow \TCY$ a subtorus of the Calabi--Yau torus acting on $Z$.
    \begin{itemize}
        \item If there exist 
        $i, j\in\{0,1,2\}$ pairwise distinct such that 
        $\epsilon_i =-\epsilon_j$ in $R_{\mgp{T}}$, then
	    \begin{equation*}
		      \refGWb{Z}{\mgp{T}}{d[\bbP^1]} = \frac{1}{d}\left(2 \sinh  \frac{d \epsilon_i}{2}\right)^{-2}\,.
	    \end{equation*}

        \item If $\epsilon_0 + \epsilon_1 + \epsilon_2= 0$ in $R_{\mgp{T}}$, then
        \begin{equation*}
            \refGWb{Z}{\mgp{T}}{d[\bbP^1]} = 0\,.
        \end{equation*}
    \end{itemize}
\end{prop}

\begin{proof}
	The first part of the \namecref{prop: shifted local P1} can be proven with arguments similar to the ones we used in the proof of \Cref{prop:ResConifold}, and we leave it
 as an exercise to the reader. The second part follows from the observation that the weight of the induced $\TCY$-action on the holomorphic two-form of $\Tot \mc{O}_{\bbP^1}(-2)$ is $-\epsilon_0 - \epsilon_1 - \epsilon_2$. If this weight is set to zero, we
 get that
    \begin{equation*}
        [\Mbar_{g}(\Tot \mc{O}_{\bbP^1}(-2), d[\bbP^1])]^{\vir}_{\mgp{T}}=0\,,
    \end{equation*}
    as the obstruction bundle features an equivariantly trivial factor.
\end{proof}

The case of the full $\TCY$-equivariant Gromov--Witten generating series of $\Tot \mc{O}_{\bbP^1} (-2) \times \Aaff{3}$ is harder, but it appears to be amenable to a general closed-form solution. 

\begin{conj}
    \label{conj: ref GW A1xC3}
	For all $d>0$ we have
	\begin{equation}
    \label{eq: ref GW A1xC3}
		\refGWb{Z}{\TCY}{d[\bbP^1]} = - \frac{1}{d} \frac{2 \sinh  \frac{d (\epsilon_0 + \epsilon_1 + \epsilon_2)}{2} }{ \prod_{i=0}^2 2 \sinh  \frac{d \epsilon_i}{2} }\,.
	\end{equation}
\end{conj}

The formula recovers the special cases computed in \Cref{prop: shifted local P1}, and we verified it  in a localisation calculation 
for $d\leq 3$, up to contributions from genus five stable maps.
	Notice that 
 in this case we do find a dependence on equivariant parameters other than $\epsilon_1$ and $\epsilon_2$: this is consistent with \Cref{conj: rigidity}, as $X$ does not support any rigid curve classes.

\begin{rmk}[The membrane perspective]
	\Cref{conj: ref GW A1xC3} is also motivated by \Cref{speculation: link to Mtheory}, as the formula in \eqref{eq: ref GW A1xC3} can be alternatively deduced from the heuristic proposal for the moduli space of M2-branes in \cite[Sec.~4]{NO14:membranes} for this particular geometry. As membranes are expected to only stack up in normal directions with strictly positive degree, for this target 
 the membrane index is only non-zero in primitive curve degree, and 
 the $\TCY$-fixed locus consists of a single point --- the embedding of the zero section $\bbP^1 \hookrightarrow Z$. Its contribution to the membrane index is
	\begin{equation}
		\label{eq:M2 index C3A1}
		\frac{-(\mchar{q}_0 \mchar{q}_1 \mchar{q}_2)^{1/2} + (\mchar{q}_0 \mchar{q}_1 \mchar{q}_2)^{-1/2}}{\prod_{i=0}^2 (\mchar{q}_i^{1/2} - \mchar{q}_i^{-1/2})} ~ \mc{Q}^{[\bbP^1]}\,,
	\end{equation}
	where the denominator comes from deformations in the trivial directions and the numerator captures the obstructions towards deforming into the $\mc{O}_{\bbP^1}(-2)$-direction \cite[Sec.~4]{NO14:membranes}. 
 Taking the logarithm 
 of the plethystic exponential of \eqref{eq:M2 index C3A1} and applying the Chern character $\mchar{q}_i \mapsto \re^{\epsilon_i}$ we recover the right-hand side of \eqref{eq: ref GW A1xC3}. The same logic could be applied to deduce formula \eqref{eq:resconGWref} for the resolved conifold.
\end{rmk}

\begin{rmk}
    \label{rmk: rank 2 DT C3}
		In \cite[Prop.~5.1]{NO14:membranes}, the rank-$r$ Donaldson--Thomas generating function of a Calabi--Yau threefold $X$ is conjectured to equate to the M2-brane index of $X \times \mc{A}_{r-1}$, modulo perturbative contributions, where the equivariant framing parameters get identified with the variables recording the exceptional curve class degrees. Combined with \Cref{speculation: link to Mtheory}, and using that $\mc{A}_{1} = \Tot \mc{O}_{\bbP^1} (-2)$, we should therefore expect a relationship between the generating series of rank two DT invariants of $\Aaff{3}$ and \Cref{conj: ref GW A1xC3}. There is an apparent tension between these, as it was shown in \cite{FMR21:HigherRankDT} that the former has no dependence on any of the framing parameters, while \eqref{eq: ref GW A1xC3} clearly depends non-trivially on the curve class variable. The contradiction is resolved by observing that the left-hand side of \eqref{eq: ref GW A1xC3} is \textit{purely perturbative} (cf e.g.\ \eqref{eq:M2 index C3A1} with \cite[Eq.~(18)]{NO14:membranes}), and therefore contributes trivially to the rank two DT theory of $\Aaff{3}$ as expected. The only non-trivial contribution from the equivariant Gromov--Witten theory of $\Aaff{3}\times \mc{A}_{1}$ comes from constant maps, which  
  will be discussed elsewhere.
\end{rmk}


\section{Local surfaces and the Nekrasov--Shatashvili limit}
\label{sec:localSurf}
We now shift the complexity up one notch and consider local Calabi--Yau threefolds with compact divisors. The simplest example is given by $X=K_S$ being the total space of the canonical bundle over a projective surface $S$. For these targets, there is a preferred torus which acts fibrewise preserving the zero section: this action scales non-trivially the holomorphic threeform, so that ~$\kappa\neq 1$. 
The main result of this section concerns local del Pezzo surfaces, for which we identify the so called NS limit ($\epsilon_2=0$) with the GW theory of the surface relative to a smooth anticanonical curve. 

\subsection{Statement of the correspondence}
\subsubsection{Fivefold invariants}
\label{sec: local surface prelim}
Let $S$ be a smooth projective surface and $D$ a smooth curve in $S$. We consider the action of $\mgp{T}_{\rm F}= \Gm[2]$ on
\begin{equation}
	\label{eq: local surface fibre action}
	\begin{tikzcd}
		Z &[-2.5em] = &[-2.5em] \Tot \big( \mc{O}_S(-D) \, \oplus \, \mc{O}_{S} \, \oplus \,  \mc{O}_{S} \big) \arrow[d] \\
		& & S
	\end{tikzcd}
\end{equation}
leaving the base invariant and scaling the fibres with weight $(1,1)$, $(-1,0)$ and $(0,-1)$ respectively. If $Z$ is Calabi--Yau (i.e.~if $D$ is  anticanonical), this torus action is 
Calabi--Yau. In a slight deviation from our setup in \Cref{sec: refined GW}, we will not necessarily assume $Z$ to be Calabi--Yau in what follows.

Under the assumption that $D^2\geq 0$, all stable maps to $\Tot \mc{O}_S(-D)$ of class $\beta$ intersecting $D$ positively factor through the zero section $S \hookrightarrow \Tot \mc{O}_S(-D)$. Hence
\begin{equation*}
	\label{eq:TotSpModuliStMaps}
	\Mbar_{g,n}(\mc{O}_S(-D),\beta) = \Mbar_{g,n}(S,\beta)\,,
\end{equation*}
which in particular implies that all curve classes $\beta$ with $D\cdot \beta >0$ are rigid. The respective virtual fundamental classes, however, differ by an insertion
\begin{equation*}
	[\Mbar_{g,n}(\text{$\mc{O}$}_S(-D),\beta)]^{\vir}_{\text{$\mgp{T}$}_{\rm F}} = 
	e^{\text{$\mgp{T}$}_{\rm F}} \big(\mathbf{R}^{1}\pi_* f^* \mc{O}_S(-D)\big) \cap [\Mbar_{g,n}(S,\beta)]^{\vir}\,.
\end{equation*}
We introduce the notation
\begin{equation*}
	[\Mbar_{g,n}(Z,\beta)]^{\vir_{\sigma}}_{\text{$\mgp{T}$}_{\rm F}} \coloneqq \HodgeLambda[g]{\epsilon_1} \, \HodgeLambda[g]{\epsilon_2} \, 
	e^{\text{$\mgp{T}$}_{\rm F}} \big(\mathbf{R}^{1}\pi_* f^* \mc{O}_S(-D)\big) \cap [\Mbar_{g,n}(S,\beta)]^{\vir}
\end{equation*}
for the virtual fundamental class of stable maps to the extended target $Z$ localised at the cocharacter $\sigma: \Gm_{\mchar{q}_-} \hookrightarrow \mgp{T}_{\rm F} , \, \mchar{q}_- \mapsto (\mchar{q}_-,\mchar{q}_-^{-1})$. We will be especially concerned with the \emph{Nekrasov--Shatahvili limit} of this cycle, i.e.~
the restriction to the inclusion of the first factor
\begin{equation*}
	\Gm \hookrightarrow \Gm[2]=\mgp{T}_{\rm F}\,.
\end{equation*}
As the $\mgp{T}_{\rm F}$-action on $\Mbar_{g,n}(S,\beta)$ is trivial, the restriction morphism essentially sets $\epsilon_2$ equal to zero:
\begin{equation*}
	\epsilon_2 \, [\Mbar_{g,n}(Z,\beta)]^{\vir_{\sigma}}_{\text{$\mgp{T}$}_{\rm F}}~ \Big\lvert_{\epsilon_2 = 0} = (-1)^g \lambda_g \, \HodgeLambda[g]{\epsilon_1} \, 
	e^{\Gm}\!\big(\mathbf{R}^{1}\pi_* f^* \mc{O}_S(-D)\big) \cap [\Mbar_{g,n}(S,\beta)]^{\vir} \,.
\end{equation*}

\subsubsection{Relative invariants}
\label{sec: notation relative GW}
Let $\mathbf{d} = (d_1,\ldots,d_m)$ be an ordered partition of $D\cdot \beta >0$, where we assume that $d_i>0$ for all $i\in \{1,\ldots,m\}$. We write
\begin{equation*}
	\Mbar_{g,n,\mathbf{d}}(S / D,\beta)
\end{equation*}
for Kim's moduli stack of genus $g$, class $\beta$ stable logarithmic maps to expanded degenerations of $(S/ D)$ with $n$ interior markings and $m$ markings with tangency profile $\mathbf{d}$ along $D$ \cite{Kim10:LogStabMaps}. This moduli stack comes with a forgetful morphism
\begin{equation*}
	\phi:\Mbar_{g,n,\mathbf{d}}(S/D,\beta) \longrightarrow \Mbar_{g,n}(S,\beta)\,.
\end{equation*}
If $D$ is anticanonical, the virtual dimension of the moduli problem is $g+n+m-1$. When $n=0$ and $m=1$, we may cap with the top Chern class of the Hodge bundle and form the following generating series in $\epsilon_1$:
\begin{equation*}
	\GW_{\beta}(S/ D) \coloneqq \frac{(-1)^{D\cdot \beta +1}}{D\cdot \beta} \sum_{g\geq 0} \epsilon_1^{2g-1} \int_{[\Mbar_{g,0,(D\cdot \beta)}(S / D,\beta)]^{\vir}} \lambda_g\,.
\end{equation*}

\subsubsection{The correspondence}
The main result of this section is the following cycle-valued identity relating the NS limit of the equivariant Gromov--Witten theory of $\Tot\mc{O}_S(-D) \times \Aaff{2}$ with the relative Gromov--Witten theory of $(S/ D)$.

\begin{thm}
	\label{thm: NS vs log cycle}
	Suppose $D$ is a smooth connected genus one curve with $D^2\geq 0$ and $\beta$ is an effective curve class satisfying $D\cdot \beta >0$. Then, for all $g,n\geq 0$, we have
	\begin{equation}
		\label{eq:NS vs log cycle}
		\epsilon_2 [\Mbar_{g,n}(Z,\beta)]^{\vir_{\sigma}}_{\text{$\mgp{T}$}_{\rm F}} ~  \Big\lvert_{\epsilon_2 = 0} = \sum_{m>0} \epsilon_1^{2g+m-2} \!\!\! \sum_{\substack{d_1,\ldots,d_m > 0 \\ \sum_i d_i = D\cdot\beta}} \frac{(-1)^{D\cdot \beta  - 1}}{m! \, \prod_{i=1}^m d_i} ~ \phi_{*} \big( \lambda_g \cap [\Mbar_{g,n,\mathbf{d}}(S/D,\beta)]^{\vir} \big)\,.
	\end{equation}
\end{thm}

Specialising to the case where $D$ is an anticanonical curve in $S$, and pushing forward to the point the identity \eqref{eq:NS vs log cycle}, we immediately obtain the following

\begin{cor}
	\label{thm: NS vs log num}
	Suppose $S$ supports a smooth connected anticanonical curve $D$ with $D^2 \geq 0$. Then, for all curve classes $\beta$ with $D\cdot \beta >0$, we have
	\begin{equation*}
		\epsilon_2 \, \refGWb{Z}{\mgp{T}_{\rm F}}{\beta} ~ \Big\lvert_{\epsilon_2 = 0} = \GW_{\beta}(S/ D)\,.
	\end{equation*}
\end{cor}


The rest of this section is dedicated to the proof of \Cref{thm: NS vs log cycle}.

\subsection{Degeneration}
Our proof follows the same approach as \cite{vGGR19} and uses a degeneration to the normal cone argument. Since this technique has become a well-established approach (see for instance \cite{BNTY23:LocOrbSNC,BFGW21:HAE,TY20:HigherGenusRelOrb,CMS24:QuasiMapDegen}) we economise on details.

Given a surface $S$ together with a smooth connected curve $D \subset S$, we consider the degeneration to the normal cone of $D$ in $S$,
\begin{equation*}
	\mc{S} = \mr{Bl}_{D\times\{0\}}(S\times \Aone) \rightarrow \Aone\,.
\end{equation*}
This is a family over $\Aone$ with general fibre $S$ and special fibre $\mc{S}_0$ obtained by gluing $S$ along $D$ with $P=\bbP_D(\mc{O}_D \oplus N_D S)$ along the zero section, which we denote by $D_0$.

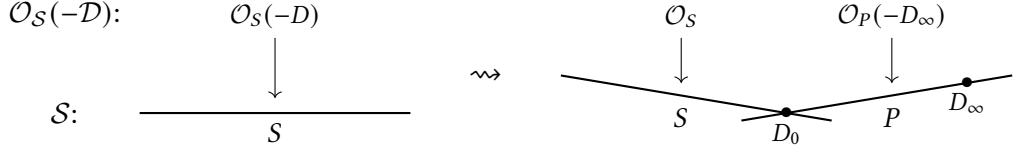
\begin{figure}[t]
    \centering
    \begin{tikzpicture}
			\draw[thick] (3,0.5) -- (-0.6,-0.1);
			\draw[thick] (-3,0.5) -- (0.6,-0.1);
			\draw[thick] (-8.6,0) -- (-5,0);
			\node[above] at (-6.8,1) {\small $\mc{O}_{S}(-D)$};
			\draw[->] (-6.8,1) -- (-6.8,0.2);
            \node[below] at (-6.8,0) {$S$};
			\node[above] at (-1.4,1) {\small $\mc{O}_{S}$};
			\draw[->] (-1.4,1) -- (-1.4,0.4);
            \node[below] at (-1.4,0.2) {$S$};
			\node[above] at (1.4,1) {\small $\mc{O}_{P}(-D_{\infty})$};
			\draw[->] (1.4,1) -- (1.4,0.4);
            \node[below] at (1.4,0.2) {$P$};
			\node at (-4,0.5) {$\rightsquigarrow$};
			\node at (0,0) {\small$\bullet$};
			\node[below] at (0,0) {\small$D_0$};
			\node at (2.4,0.4) {\small$\bullet$};
			\node[below] at (2.4,0.4) {\small$D_\infty$};
            \node at (-9.6,0) {$\mc{S}$:};
            \node[above] at (-9.6,1) {$\mc{O}_{\mc{S}}(-\mc{D})$:};
		\end{tikzpicture}
    \caption{Degeneration to the normal cone of $D$ in $S$.}
    \label{fig: degeneration to normal cone}
\end{figure}

Denote by $\mc{D}$ the proper transform of $D\times \Aone$. Its intersection with a general fibre is $D\hookrightarrow S$, while its intersection with the special fibre is the infinity section $D_\infty$ of $P$. We will write $\mc{Z}$ for the total space of 
\begin{equation*}
	\mc{O}_{\mc{S}}(-\mc{D})\oplus\mc{O}_{\mc{S}}\oplus \mc{O}_{\mc{S}}
\end{equation*}
and consider the same fibrewise $\mgp{T}_{\rm F}$-action on $\mc{Z}$ as in \eqref{eq: local surface fibre action}. The whole setup is illustrated in \Cref{fig: degeneration to normal cone}.

\subsection{Decomposition}
Composition with the blowup morphism $\mc{S}_0 \rightarrow S$ contracting the component $P$ gives a pushforward morphism
\begin{equation*}
	\rho : \Mbar_{g,n}(\mc{S}_0,\beta) \longrightarrow \Mbar_{g,n}(S,\beta)
\end{equation*}
under which
\begin{equation}
	\label{eq: general vs special fibre}
	[\Mbar_{g,n}(Z,\beta)]^{\vir_{\sigma}} = \rho_{*} [\Mbar_{g,n}(\mc{Z}_0,\beta)]^{\vir_{\sigma}}\,.
\end{equation}
We will use the degeneration formula for stable logarithmic maps of Kim, Lho and Ruddat \cite{KLR21:degen} to decompose the cycle on the right-hand side. For purely expository reasons, we will be working with their formula rather than Li's for relative stable maps \cite{Li02:Degen}; both essentially yield the same result \cite{ACFW13:ExpDegenPairs}.

The right-hand side of \eqref{eq: general vs special fibre} decomposes into contributions labelled by certain decorated graphs $\Gamma$, which we call splitting types. This is the datum of a bipartite graph with edges $E(\Gamma)$, legs $L(\Gamma)$ and vertices $V(\Gamma)=V_S(\Gamma) \sqcup V_P(\Gamma)$ split into $S$- and $P$-vertices. Moreover, each edge $e$ is decorated with a contact order $d_e\in \bbZ_{>0}$ and to each vertex $v$ we assign a genus $g_v \in \bbZ_{\geq 0}$ and a curve class $\beta_v$ in either $S$ or $P$. Lastly, the datum of a splitting type also includes a bijection $L(\Gamma)\cong \{1,\ldots,n\}$ and, for convenience, also a labelling $E(\Gamma)\cong\{1,\ldots,m\}$. This datum is subject to the usual compatibility conditions with the discrete datum $(g,\beta)$. 

Following \cite{KLR21:degen}, for each such $\Gamma$ there is an associated moduli stack $\Mbar_\Gamma$ parametrising stable logarithmic maps to $\mc{S}_0$ of type $\Gamma$. It admits a finite morphism $\iota: \Mbar_\Gamma \longrightarrow \Mbar_{g,n}(\mc{S}_0,\beta)$ and an \'{e}tale morphism $\nu$ to the fibre product
\begin{equation*}
	\begin{tikzcd}
		\Mbar_\Gamma \arrow[r,"\nu"] & \bigodot_{v} \Mbar_v \arrow[r] \arrow[d]  \ar[rd,phantom,"\square",start anchor=center,end anchor=center] & \prod_{v} \Mbar_v \arrow[d] \\
		& \prod_e D_0 \arrow[r,"\Delta"] & \prod_v \prod_{e \ni v}  D_0\,.
	\end{tikzcd}
\end{equation*}
Given a vertex $v$ of $\Gamma$, in the above diagram we write $\Mbar_v$ to denote the moduli stack of stable logarithmic maps to either $S(\mr{log}\, D_0)$ or $P(\mr{log}\, D_0)$ depending on whether $v$ is an $S$- or $P$-vertex and whose topological type is given the star of $\Gamma$ at $v$. This moduli stack admits evaluation morphisms to $D_0$ associated to edges adjacent to $v$ and we take the right vertical arrow in the diagram to be the collection of all these evaluation morphisms. Then $\bigodot_{v} \Mbar_v$ is defined to be the fibre product along the inclusion of the diagonal $\Delta$.

By the degeneration formula for stable logarithmic maps \cite[Thm~1.5]{KLR21:degen} we have
\begin{equation}
	\label{eq: degen formula}
	[\Mbar_{g,n}(\mc{S}_0,\beta)]^{\vir} = \sum_{\Gamma} \frac{\mr{lcm}(d_e)}{m!} \iota_{*} \nu^{*} \Delta^! \textstyle{\prod_v} [\Mbar_v]^{\vir}\,,
\end{equation}
where the sum is over all splitting types. Now observe that there is a morphism $\theta: \bigodot_{v} \Mbar_v \rightarrow \Mbar_{g,n}(S,\beta)$ which composes a stable map with the blowup morphism and glues the domain curve along edge markings. This morphism makes the diagram
\begin{equation*}
	\begin{tikzcd}
		\Mbar_\Gamma \ar[d,"\iota"] \ar[r,"\nu"] & \bigodot_{v} \Mbar_v \ar[d,"\theta"] \\ 
		\Mbar_{g,n}(\mc{S}_0,\beta) \ar[r,"\rho"] &\Mbar_{g,n}(S,\beta)
	\end{tikzcd}
\end{equation*}
commute. In combination with the degeneration formula \eqref{eq: degen formula} we obtain
\begin{equation*}
	\rho_{*} [\Mbar_{g,n}(\mc{S}_0,\beta)]^{\vir} = \sum_{\Gamma} \frac{\prod_e d_e}{m!} \theta_{*} \Delta^! \textstyle{\prod_v} [\Mbar_v]^{\vir}\,,
\end{equation*}
using that the degree of $\nu$ is $(\prod_e d_e)/\mr{lcm}(d_e)$ by  \cite[Eq.~(1.4)]{KLR21:degen}.

Now, the cycles $[\Mbar_{g,n}(\mc{S}_0,\beta)]^{\vir}$ and $[\Mbar_{g,n}(\mc{Z}_0,\beta)]^{\vir_{\sigma}}$ differ by an insertion
\begin{equation*}
	\HodgeLambda[g]{\epsilon_1} \, \HodgeLambda[g]{\epsilon_2} \, e^{\mgp{T}_{\rm F}} \big(\mathbf{R}^{1}\pi_{*} f^{*} \mc{O}_{\mc{S}_0}(-D_{\infty})\big)\,.
\end{equation*}
Pulling this insertion back to the product $\prod_v \Mbar_v$ yields
\begin{equation}
	\label{eq: degeneration final}
	[\Mbar_{g,n}(Z,\beta)]^{\vir_{\sigma}} = \sum_{\Gamma} \big(\epsilon_1\epsilon_2(-\epsilon_1-\epsilon_2)\big)^{|E(\Gamma)|} ~ \frac{\prod_e d_e}{m!} \theta_{*} \Delta^! \big(\textstyle{\prod_v} C_v \cap [\Mbar_v]^{\vir}\big)\,,
\end{equation}
where we introduce the notation
\begin{equation*}
	C_v \coloneqq \begin{cases}
		\HodgeLambda[g_v]{\epsilon_1}\,  \HodgeLambda[g_v]{\epsilon_2} \, \HodgeLambda[g_v]{-\epsilon_1-\epsilon_2} & \text{if $v$ is an $S$-vertex,}\\[0.5em]
		\HodgeLambda[g_v]{\epsilon_1} \, \HodgeLambda[g_v]{\epsilon_2} \, e^{\mgp{T}_{\rm F}} \big(\mathbf{R}^{1}\pi_{*} f^{*} \mc{O}_{\mc{S}_0}(-D_{\infty})\big) & \text{if $v$ is a $P$-vertex.}
	\end{cases}
\end{equation*}
Note that the factor $\big(\epsilon_1\epsilon_2(-\epsilon_1-\epsilon_2)\big)^{|E(\Gamma)|}$ is due to the normalisation exact sequence coming from gluing along the edges of $\Gamma$.

From the description of the insertions $C_v$ we observe that the contribution of an $S$-vertex is polynomial in $\epsilon_1$, $\epsilon_2$ up to a factor $(\epsilon_1\epsilon_2(\epsilon_1+\epsilon_2))^{-1}$. Similarly, the contribution of a $P$-vertex is polynomial in the equivariant parameters up to a factor $(\epsilon_1\epsilon_2)^{-1}$. This means that each term in the sum \eqref{eq: degeneration final} is a cycle-valued polynomial in $\epsilon_1$, $\epsilon_2$, times an overall factor
\begin{equation*}
	\big(\epsilon_1 \epsilon_2 (\epsilon_1 + \epsilon_2)\big)^{|E(\Gamma)|} ~ \big(\epsilon_1\epsilon_2(\epsilon_1+\epsilon_2)\big)^{-|V_S(\Gamma)|} ~ \big(\epsilon_1\epsilon_2\big)^{-|V_P(\Gamma)|} = (\epsilon_1 \epsilon_2)^{h_1(\Gamma)-1} (\epsilon_1 + \epsilon_2)^{|E(\Gamma)|-|V_S(\Gamma)|} \,.
\end{equation*}
As a consequence, the only splitting types $\Gamma$ that can contribute non-trivially to
\begin{equation*}
	\epsilon_2\, [\Mbar_{g,n}(Z,\beta)]^{\vir_{\sigma}}~\Big\lvert_{\epsilon_2 = 0}
\end{equation*}
must have $h_1(\Gamma)=0$, or in other words, the graph underlying $\Gamma$ must be a tree. Furthermore, in the limit $\epsilon_2=0$, the insertions simplify to
\begin{equation*}
	C^{\mr{NS}}_v \coloneqq \begin{cases}
		- \epsilon_1^{2g_v-2}\, \lambda_{g_v} & \text{if $v$ is an $S$-vertex,}\\[0.5em]
		(-1)^{g_v} \lambda_{g_v} \,\HodgeLambda[g_v]{\epsilon_1} \, e^{\mgp{\Gm}} \! \big(\mathbf{R}^{1}\pi_{*} f^{*} \mc{O}_{\mc{S}_0}(-D_{\infty})\big) & \text{if $v$ is a $P$-vertex}\,,
	\end{cases}
\end{equation*}
where in the top line we used \Cref{lem: Mumfords relation}. We arrive at
\begin{equation}
	\label{eq: NS degeneration final}
	\epsilon_2[\Mbar_{g,n}(Z,\beta)]^{\vir_{\sigma}}~ \Big\lvert_{\epsilon_2 = 0} = \sum_{\text{tree type }\Gamma} \big(-\epsilon_1^2\big)^{|E(\Gamma)|} ~ \frac{\prod_e d_e}{m!} \theta_{*} \Delta^! \big( \textstyle{\prod_v} C^{\mr{NS}}_v \cap [\Mbar_v]^{\vir} \big)\,.
\end{equation}

\subsection{Vanishing}
In this section we will prove that most splitting types contribute trivially to \eqref{eq: NS degeneration final}.

\begin{prop}
	\label{prop: vanishing}
	The cycle $\theta_{*} \Delta^! \big(\textstyle{\prod_v} C^{\mr{NS}}_v \cap [\Mbar_v]^{\vir} \big)$ associated to a splitting type $\Gamma$ is non-zero only if
	\begin{itemize}  
		\item the graph underlying $\Gamma$ is of star shape with a unique $S$-vertex $\tilde{v}$;
		
		\item $\tilde{v}$ carries all marking legs and is decorated with $g_{\tilde{v}}=g$ and $\beta_{\tilde{v}}=\beta$;
        
		\item all $P$-vertices $v$ are decorated with $g_v=0$ and $\beta_v$ is the multiple of fibre class in $P$:
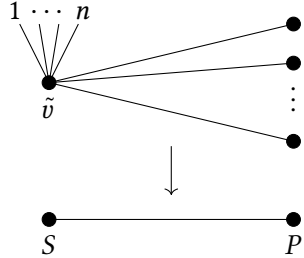
        \begin{equation*}
        \begin{tikzpicture}[baseline=(current  bounding  box.center),scale=1.3]
			\draw[fill=black] (-2.5,1.4) circle[radius=2pt];
			\draw (-2.5,1.35) node[below]{$\tilde{v}$};
			
			\draw (-2.5,1.4) -- (0,0.8);
			\draw[fill=black] (0,0.8) circle[radius=2pt];
			
			\draw (0,1.31) node{$\vdots$};
			
			\draw (-2.5,1.4) -- (0,1.6);
			\draw[fill=black] (0,1.6) circle[radius=2pt];		
			
			\draw (-2.5,1.4) -- (0,2);
			\draw[fill=black] (0,2) circle[radius=2pt];		

			\draw (-2.5,1.4) -- (-2.8,2);
			\draw (-2.85,1.95) node[above]{$1$};
			\draw (-2.5,1.95) node[above]{$\cdots$};
			\draw (-2.5,1.4) -- (-2.6,2);
			\draw (-2.5,1.4) -- (-2.4,2);
			\draw (-2.5,1.4) -- (-2.2,2);
			\draw (-2.15,1.95) node[above]{$n$};
   
			\draw[->] (-1.25,0.75) -- (-1.25,0.25);
			
			\draw (-2.5,0) -- (0,0);
			
			\draw[fill=black] (-2.5,0) circle[radius=2pt];
			\draw (-2.5,-0.05) node[below]{$S$};
			
			\draw[fill=black] (0,0) circle[radius=2pt];
			\draw (0,-0.05) node[below]{$P$};
		\end{tikzpicture}
        \end{equation*}
	\end{itemize}
\end{prop}

\subsubsection{First reduction: confined curve class and genus label}
Let $\Gamma$ be a splitting type whose underlying graph is a tree. We will constrain its shape in several steps. First we will establish the last condition in \Cref{prop: vanishing}.

Let $v\in V(\Gamma)$ be a $P$-vertex. The projection $p:P\rightarrow D_0$ induces a pushforward morphism
\begin{equation*}
	\Mbar_v \longrightarrow \Mbar_{g_v,n_v\cup m_v}(D_0,p_*\beta_v)
\end{equation*}
where $n_v\subseteq L(\Gamma)$ denotes the set of marking legs and $m_v \subseteq E(\Gamma)$ the set of edges adjacent to $v$. Here of course we need to assume that the moduli stack on the right-hand side is non-empty. Following \cite{vGGR19}, we factor the morphism through the fibre product
\begin{equation}
	\label{eq: factor projection morphism}
	\begin{tikzcd}
		\Mbar_v \ar[r,"u"] \ar[dr] & \widetilde{\modulifont{M}}_v \ar[r] \ar[d] \ar[rd,phantom,"\square",start anchor=center,end anchor=center]  & \Mbar_{g_v,n_v\cup m_v}(D_0,p_*\beta_v) \ar[d] \\
		& \mf{M}^{\mr{log}}_{g_v,n_v\cup m_v,H_2(P,\bbZ)^+} \ar[r,"\nu"] & \mf{M}_{g_v,n_v\cup m_v,H_2(D_0,\bbZ)^+}\,,
	\end{tikzcd}
\end{equation}
where the moduli stacks in the bottom row are the Artin stacks parametrising prestable (log) curves with irreducible components labelled by effective curve classes. Pulling the obstruction theory of $\Mbar_{g_v,n_v\cup m_v}(D_0,p_*\beta_v)$ back to $\widetilde{\modulifont{M}}_v$, we find
\begin{equation*}
	[\widetilde{\modulifont{M}}_v]^{\vir}= \nu^! [\Mbar_{g_v,n_v\cup m_v}(D_0,p_*\beta_v)]^{\vir}
\end{equation*}
as proven in \cite[Sec.~4]{vGGR19}. Moreover, the short exact sequence
\begin{equation*}
	0 \longrightarrow \mc{T}^{\mr{log}}_{P(\mr{log}\, D_0) / D_0} \longrightarrow  \mc{T}^{\mr{log}}_{P(\mr{log}\, D_0)} \longrightarrow  \mc{T}_{ D_0} \longrightarrow 0
\end{equation*}
induces a compatible triple for the  commuting triangle on the left-hand side of \eqref{eq: factor projection morphism}. Hence, by \cite[Cor.~4.9]{Ma11:VirtPull}, we have
\begin{equation}
	\label{eq: pullback virt class}
	[\Mbar_v]^{\vir} = u^![\widetilde{\modulifont{M}}_v]^{\vir}= u^! \nu^! [\Mbar_{g_v,n_v\cup m_v}(D_0,p_*\beta_v)]^{\vir}\,.
\end{equation}
We are now equipped to prove the following vanishing statement.

\begin{lem}
	\label{lem: first vanishing}
	Suppose $D$ is a genus one curve. Then, for every $P$-vertex $v$ with $p_* \beta_v\neq 0$ or $g_v>0$, we have
	\begin{equation*}
		\lambda_{g_v} \cap [\Mbar_v]^{\vir}=0\,.
	\end{equation*}
\end{lem}

\begin{proof}
	Suppose we have $p_* \beta_v\neq 0$. In this case $\Mbar_{g_v,n_v\cup m_v}(D_0,p_*\beta_v)$ is non-empty and therefore \eqref{eq: pullback virt class} holds. Moreover, since lambda classes are preserved under stabilisation, we have
	\begin{equation}
		\label{eq: pullback virt class lambda}
		\lambda_{g_v} \cap[\Mbar_v]^{\vir} =  u^! \nu^! \big( \lambda_{g_v} \cap [\Mbar_{g_v,n_v\cup m_v}(D_0,p_*\beta_v)]^{\vir} \big)\,.
	\end{equation}
	It therefore suffices to demonstrate the vanishing of the cycle on the right-hand side. Observe that, since we assume $D_0$ to be a genus one curve, $D_0\times \Aone$ is a surface with trivial canonical bundle: consequently, the virtual fundamental class $[\Mbar_{g_v,n_v\cup m_v}(D_0\times \Aone,p_*\beta_v)]^{\vir}$ vanishes. Writing 
	\begin{equation*}
		\iota_0 : \Mbar_{g_v,n_v\cup m_v}(D_0,p_*\beta_v) \hookrightarrow \Mbar_{g_v,n_v\cup m_v}(D_0\times \Aone,p_*\beta_v)
	\end{equation*}
	for the inclusion of the zero section, a comparison of obstruction theories yields
	\begin{equation*}
		\lambda_{g_v} \cap [\Mbar_{g_v,n_v\cup m_v}(D_0,p_*\beta_v)]^{\vir} = (-1)^{g_v} \, \iota_0^! [\Mbar_{g_v,n_v\cup m_v}(D_0\times \Aone,p_*\beta_v)]^{\vir} = 0\,,
	\end{equation*}
	which proves the first part of the \namecref{lem: first vanishing}.
	
	Suppose now that $p_* \beta_v = 0$ and $g_v>0$. We again use
 \eqref{eq: pullback virt class lambda}, where now
	\begin{equation*}
		[\Mbar_{g_v,n_v\cup m_v}(D_0,0)]^{\vir} = (-1)^{g_v} \lambda_{g_v} \cap [D_0 \times \Mbar_{g_v,n_v\cup m_v}]
	\end{equation*}
	since the Chern class of the tangent bundle of $D_0$ is trivial. The statement of the \namecref{lem: first vanishing} then follows from the fact that $\lambda_{g_v}^2 =0$ whenever $g_v>0$, which is the leading order vanishing in \Cref{lem: Mumfords relation}.
\end{proof}

\subsubsection{Second reduction: confined markings}
We continue by assuming that all $P$-vertices of $\Gamma$ satisfy the third property stated in \Cref{prop: vanishing}. In this Section we will show that all $P$-vertices carry at most one marking.

Let $v$ be a $P$-vertex of $\Gamma$. The image of a stable map parametrised by $\Mbar_v$ gets contracted to a point under composition with the projection $P \rightarrow D_0$ since we assume that $\beta_v$ is a multiple of a fibre class. The morphism $\theta$, which composes with the blowup morphism and glues the domain curves, factors as follows:
\begin{equation*}
	\label{eq:BigDiag}
	\begin{tikzcd}
		\Mbar_{g,n}(S,\beta) & \ar[l,swap,"\theta"] \bigodot_{v} \Mbar_{v} \arrow[r] \arrow[d] \ar[rd,phantom,"\square",start anchor=center,end anchor=center] & \prod_{v\in V_P(\Gamma)} \Mbar_{v} \times \prod_{\tilde{v}\in V_S(\Gamma)} \Mbar_{\tilde{v}} \arrow[d] \arrow[d,"\Pi_v\psi_v \times \mr{Id}"]\\
		& \ar[ul,"\phi"] \Nbar_\Gamma \arrow[r] \arrow[d] \ar[rd,phantom,"\square",start anchor=center,end anchor=center] &  \prod_{v\in V_P(\Gamma)} \big(D_0 \times \Mbar_{0,m_v \cup n_v} \big) \times \prod_{\tilde{v}\in V_S(\Gamma)} \Mbar_{\tilde{v}} \arrow[d]\\
		& \prod_e D_0 \arrow[r,"\Delta"] & \prod_v \prod_{e\ni v} D_0\,.
	\end{tikzcd}
\end{equation*}
Here, we denote by $\psi_v$ the forgetful morphism $\Mbar_{v}\rightarrow D_0 \times \Mbar_{0,m_v \cup n_v}$ which only remembers the evaluation of the edge marking and the stabilised domain curve, and we adopt the convention that $\Mbar_{0,n}=\Speck$ if $|n|\leq 2$. 
Since the diagram commutes and the top square is Cartesian we find that
\begin{equation}
	\label{eq: term simplified 1}
	\theta_{*} \Delta^! \big(\textstyle{\prod_v} C^{\mr{NS}}_v \cap [\Mbar_v]^{\vir} \big) = \phi_{*} \Delta^! \Big(\prod_{v\in V_P(\Gamma)} \psi_{v*} \big(C^{\mr{NS}}_v \cap [\Mbar_{v}]^{\vir}\big) \times \prod_{\tilde{v}\in V_S(\Gamma)} C^{\mr{NS}}_{\tilde{v}} \cap [\Mbar_{\tilde{v}}]^{\vir}\Big)\,.
\end{equation}
A quick dimension count reveals that for all $v\in V_P(\Gamma)$ the cycle $C^{\mr{NS}}_v \cap [\Mbar_{v}]^{\vir}$ is an element of
\begin{equation*}
	\Chow^{\mgp{\Gm}}_{|m_v|+|n_v|}(\Mbar_{v})
\end{equation*}
times a factor $\epsilon_1^{-1}$. Therefore, since $\Mbar_{0,m_v \cup n_v}\times D_0$ is of dimension $\min(|m_v|+|n_v| - 2,1)$, the pushforward
\begin{equation*}
	\psi_{v*} \big(C^{\mr{NS}}_v \cap [\Mbar_{v}]^{\vir}\big)
\end{equation*}
vanishes whenever $|m_v|+|n_v| \geq 2$. Since necessarily $|m_v| \geq 1$, we deduce the following.
\begin{lem}
	\label{lem: confined markings}
	$\theta_{*} \Delta^! \big(\textstyle{\prod_v} C^{\mr{NS}}_v \cap [\Mbar_v]^{\vir} \big)=0$ unless every $P$-vertex has exactly one edge and carries no legs.\qed
\end{lem}

In combination with \Cref{lem: first vanishing} the above indeed proves all the vanishings claimed in \Cref{prop: vanishing}.

\subsection{Remaining terms}
Let us analyse the remaining terms of the sum \eqref{eq: NS degeneration final}. We fix a star shaped splitting type $\Gamma$ as specified in \Cref{prop: vanishing} and observe that the formula \eqref{eq: term simplified 1} for its associated term simplifies to
\begin{equation*}
	\label{eq: term simplified 2}
	\theta_{*} \Delta^! \big(\textstyle{\prod_v} C^{\mr{NS}}_v \cap [\Mbar_v]^{\vir} \big) = -\epsilon_1^{2g-2} ~ \phi_{*} \Delta^! \Big(\left(\prod_{v\in V_P(\Gamma)} \psi_{v*} \big(C^{\mr{NS}}_v \cap [\Mbar_{v}]^{\vir}\big)\right) \times \big(\lambda_g \cap [\Mbar_{\tilde{v}}]^{\vir}\big)\Big)\,,
\end{equation*}
where we write $\tilde{v}$ for the unique $S$-vertex of $\Gamma$. Now, for every $P$-vertex $v$, a dimension count analogous to the one we carried out above reveals that
\begin{equation*}
	\psi_{v*} \big(C^{\mr{NS}}_v \cap [\Mbar_{v}]^{\vir}\big) = \epsilon_1^{-1} a_v [D_0]
\end{equation*}
for some $a_v\in\bbQ$. We compute $a_v$ as in \cite[Prop.~2.4]{vGGR19} by considering the fibre square
\begin{equation*}
	\begin{tikzcd}
		\Mbar_{0,0,(d_e), d_e} (\bbP^1 / 0) \ar[r,hook] \ar[d] \ar[rd,phantom,"\square",start anchor=center,end anchor=center] & \Mbar_{v} \ar[d,"\psi_v"] \\
		\Speck \ar[r,hook,"\xi"] & D_0\,.
	\end{tikzcd}
\end{equation*}
By \cite[Thm.~5.1]{BP05:CurvesCY3TQFT} we find that
\begin{equation*}
	\epsilon_1^{-1} a_v = \xi^!\psi_{v*} \big(C^{\mr{NS}}_v \cap [\Mbar_{v}]^{\vir}\big) = \epsilon_1^{-1} \int_{[\Mbar_{0,0,(d_e)} (\mc{O}_{\bbP^1}(-1) / 0 , d_e[\bbP^1])]^{\vir}_{\mgp{\Gm}}} 1 =  \frac{(-1)^{d_e+1}}{\epsilon_1\, d_e^2}\,,
\end{equation*}
where $e$ is the unique edge adjacent to $v$. Inserting this relation into \eqref{eq: term simplified 2} gives
\begin{equation*}
	\theta_{*} \Delta^! \big({\textstyle\prod_v} \displaystyle C^{\mr{NS}}_v \cap [\Mbar_v]^{\vir} \big) = \epsilon_1^{2g-|E(\Gamma)|-2} (-1)^{D\cdot \beta + |E(\Gamma)| +1} \frac{1}{\prod_e d_e^2} ~ \phi_{*} \big( \lambda_g \cap [\Mbar_{\tilde{v}}]^{\vir}\big)\,.
\end{equation*}
Inserting this into \eqref{eq: NS degeneration final} yields \eqref{eq:NS vs log cycle} upon identifying $\Mbar_{\tilde{v}}$ with $\Mbar_{g,n,\mathbf{d}}(S/ D,\beta)$, concluding the proof of \Cref{thm: NS vs log cycle}.

\section{Refined mirror symmetry for local \texorpdfstring{$\bbP^2$}{P2}}
\label{sec: ref mirror sym KP2}

We will now specialise to the case when the surface $S$ is the projective plane
and study the relation between the equivariantly refined GW invariants of \cref{sec: refined GW} and the refined mirror symmetry formalism of \cite{HK10:OmegaBG}. The main result of this Section will be the proof of the first two parts, {\bf (Finite generation)} and {\bf (Modular anomaly)}, of \cref{thm:B}. In order to do that, we shall need to set up quite a bit of notation first.

\subsection{Setup and notation}

As a toric variety, $Z=K_{\bbP^2} \times \bbC^2$ is a GIT quotient of $\bbC^6$ with coordinates $(x_0, x_1, x_2; y_0, y_1, y_2)$,
\beq
Z \coloneqq \bbC^6/\!\!/_{\xi>0}  \Gm = (\bbC^6\setminus\Delta)/\Gm\,,
\label{eq:GIT}
\eeq
where
\[
\Delta \coloneqq V\l(\bra x_0,x_1,x_2\ket\r)\,,
\]
and the quotient torus acts as
\bea
\Gm \times \bbC^6 & \longrightarrow & \bbC^6\,,  \nn \\
(\lambda; x_0, x_1, x_2, y_0,y_1,y_2) & \to & (\lambda x_0, \lambda x_1, \lambda x_2, \lambda^{-3} y_0, y_1, y_2)\,.
\label{eq:GITact}
\eea
%
From this, we see that
$[x_0:x_1:x_2]$ serve as homogeneous coordinates for the base $\bbP^2$, and $(y_0, y_1,y_2)$ are fibre coordinates. 
The maximal torus $\Tmax \simeq (\Gm)^5 \circlearrowright Z$ descends from an action on $\bbC^6$ with geometric weights \[(\alpha_0, \alpha_1, \alpha_2, 0, -\epsilon_1, -\epsilon_2)\,.\]  This lifts to $\pi:Z \to \bbP^2$ the $(\Gm)^3$ action on $\bbP^2$ having $P_0=[1:0:0]$, $P_1 = [0:1:0]$ and $P_2=[0:0:1]$ as fixed points, and whose natural lift to $\cO_{\bbP^2}(1)$ has weight $-\alpha_i$ on $P_i$, $i=0,1,2$. Denoting $H \coloneqq \pi^*c_1^{\Tmax}(\cO_{\bbP^2}(1))$ the equivariant line class, 
the ${\Tmax}$-equivariant 
Chow cohomology of $Z$ is
%
\[
\Chow_{\Tmax}(Z) \simeq
\frac{\bbC[H,\alpha_0,\alpha_1,\alpha_2, \epsilon_1,\epsilon_2]}{\bra \prod_{i=0}^2(H+\alpha_i) \ket}
\]
%
%
%
which we are going to take with complex coefficients in this and the following section. In the following we will write
\[
Q_{\Tmax} \coloneqq \mathrm{Frac}\,\Chow^{}_{\Tmax}(\mathrm{pt}) \simeq \bbC(\alpha_0, \alpha_1, \alpha_2, \epsilon_1, \epsilon_2)\,, \quad 
%
\cV_{\Tmax} \coloneqq \Chow^{}_{\Tmax}(Z) \otimes_{\Chow^{}_{\Tmax}\!(\mathrm{pt})} Q_{\Tmax}  \simeq  \frac{\bbC(\alpha_0, \alpha_1, \alpha_2, \epsilon_1, \epsilon_2)[H]}{\bra \prod_{i=0}^2(H+\alpha_i)\ket} \,.\] 

The ${\Tmax}$-characters $w^{(j)}_i$ ($j=0, \dots, 4$, $i=0,1,2$) on the tangent spaces $T_{P_i}Z$ are 
\[
w^{(j)}_i = \left\{\bary{ll} \alpha_j-\alpha_i\,, & i \neq j,~ i,j \in \{0,1,2\}\,, \\ 3 \alpha_i\,, & i=j\,, \\ -\epsilon_{j-2}\,, & j=3,4\,.
\eary\right.
\]
Writing $\mathsf{w}_i\coloneqq \prod_{j=0}^4 w^{(j)}_i$ for the ${\Tmax}$-character on $\wedge^5 T_{P_i}Z$, the ${\Tmax}$-equivariant intersection pairing on $Z$ is
\[
\eta(\phi, \psi) \coloneqq \intEquiv{\Tmax}_{[Z]_{\Tmax}} \phi \cup \psi = \sum_{k=0}^2 \frac{\phi|_{P_k} \psi|_{P_k}}{\mathsf{w}_k}\,.
\]
%
%
Let
\[e_i(x) \coloneqq \big[y^{3-i}\big] \prod_{j=0}^2 \big(y+x_j\big)\] be the $i^{\rm th}$~elementary symmetric polynomial in $3$ variables
$x=(x_0, x_1, x_{2})$. In components, the Gram matrix $(\eta_{ab})_{a,b \in \{0,1,2\}}$
in the basis $\{H^a\}_{a=0,1,2}$ for $\Chow_{{\Tmax}}(Z)$ is
%
\beq
\eta_{ab} \coloneqq \eta(H^a, H^b) = \frac{1}{3\epsilon_1 \epsilon_2}
\left(
\begin{array}{ccc}
 \frac{1}{e_3(\alpha)} & 0 & 0 \\
0 & 0 & -1  \\
0 & -1 &  e_1(\alpha)
\end{array}
\right)_{ab}\,, 
\quad  \eta^{ab}= (\eta^{-1})_{ab}\,.
\label{eq:etaflat}
\eeq
As in \cref{sec:localSurf}, of particular interest to us will be the 2-dimensional subtorus $\mgp{T}_{\rm F} \subset {\Tmax}$ fixing the base and rotating the fibres with weights $(1,1)$, $(-1,0)$ and $(0,-1)$. From \eqref{eq:GITact}, this is given by the diagonal restriction 
%
\beq 
\alpha_i = \frac{\epsilon_1+\epsilon_2}{3}\,,
\label{eq:diagweights}
\eeq 
which covers the trivial action on $\bbP^2$. 

\subsubsection{Gromov--Witten invariants}
%
We will employ the usual correlator notation for $\mgp{T}$-e\-qui\-va\-riant Gro\-mov--Witten invariants: for non-negative integers ${a_1}, \dots, {a_n}$, $r_1, \dots, r_n \in \bbZ_{\geq 0}$, and $t \in \cV_{\Tmax}$ we write
\bea
\bra H^{a_1} \psi_1 ^{r_1}, \dots, H^{a_n}\psi_n ^{r_n} \ket_{g,n,d}^Z
& \coloneqq & \intEquiv{\Tmax}_{[\Mbar_{g,n}(Z,
    d [H])]^{\vir}_{\Tmax}} \prod_{i=1}^n \mathrm{ev}^*_i (H^{a_i})\,
\psi_i^{r_i}, \nn
\\
\bra \! \bra H^{a_1} \psi_1 ^{r_1}, \dots, H^{a_n} \psi_n^{r_n}\ket \! \ket_{g,n}^Z & 
\coloneqq &
\sum_{d>0, m\geq 0} \big\langle 
H^{a_1} \psi_1 ^{r_1}, \dots, H^{a_n}\psi_n ^{r_n},  \overbrace{t,t,\ldots,t}^{\text{$m$
        times}} \big\rangle_{g,n+m,d}^Z\,.
\label{eq:gwprim}
\eea
Expanding $t=\sum t_i H^i$, 
the big quantum cohomology product of $Z$ is defined as
\beq
\eta\Big(H^{a} \circ_t H^b, H^c \Big) \coloneqq  \Big\langle \! \Big\langle H^a, H^b, H^c \Big\rangle \! \Big\rangle_{0,3}^Z = \eta\Big(H^{a+b}, H^c \Big)+ \sum_{d>0, m\geq 0} \Big\langle H^{a}, H^b, H^c, 
  \overbrace{H^2,H^2,\ldots,H^2}^{\text{$m$
        times}} \Big\rangle_{0,m+3,d}^Z\frac{t_2^m} {m!}\cQ^d\,,
\label{eq:qprod}
\eeq
where we 
used the Divisor Axiom, and set 
\[
\cQ \coloneqq \re^{t_1}\,.\]
The structure constants in \eqref{eq:qprod} have a non-zero radius of convergence $\cQ_{\rm cf} \in \bbR^+$  \cite{Klemm:1999gm} and define, in
the region
$B_{\cQ_{\rm cf}} \coloneqq \{0<|\cQ|<\cQ_{\rm cf}\}$
with coordinates $(t_0, \cQ, t_2)$ centred at $\cQ=0$, 
an analytic family of semi-simple associative rings 
over $B_{\cQ_{\rm cf}}$. 
We will write $\mathfrak{e}_i$ for the idempotents of the algebra \eqref{eq:qprod}, and $\mathfrak{f}_i$ for their normalised counterparts:
\[
\mathfrak{e}_i \circ_t \mathfrak{e}_j = \delta_{ij} \mathfrak{e}_i, \quad \mathfrak{f}_i \coloneqq \frac{\mathfrak{e}_i}{\sqrt{\Delta_i}}\,,
\]
where \[\Delta_i \coloneqq \eta(\mathfrak{e}_i,\mathfrak{e}_i)\,.\] 
%
Formal flat sections $\chi \in \bbC(\!(z)\!) \otimes \Omega^1(\cV_{\Tmax})$ of the Dubrovin connection
\[
\nabla^{(z)}_X \chi \coloneqq  \dd -z^{-1} X \circ_t \chi =0\,,
\label{eq:defconn1}
\]
form a 3-dimensional vector space 
over $\bbC(\!(z)\!)\otimes Q_{\Tmax}$.
A canonical set of flat coordinates for the Dubrovin connection  is given by the components of the big $J$-function of $Z$,
\beq
\label{eq:resj}
J(t,z) \coloneqq (z+t_0)\mathbf{1} +t_1 H + t_2 H^2 + \sum_{a,b=0}^{2}H^a \eta^{ab}\bra \! \bra
\frac{H^b}{z-\psi}\ket \! \ket_{0,1}^Z
\,.
\eeq
%
The fundamental solution of $\nabla^{(z)}\chi=0$ is the corresponding flat frame $S(t,z)\coloneqq \mathrm{d} J(t, z)$.
%
%
Explicitly,
\beq
S(t,z)(H^a)
= \partial_{t_a}  J(t, z) = H^a +\sum_{c,b=0}^{2}H^c \eta^{cb} \bra \! \bra H^a,\frac{H^b}{z-\psi}\ket \! \ket_{0,2}^Z
\,.
\label{eq:fundsol}
\eeq

\subsubsection{Small quantum cohomology}
The small quantum cohomology locus is defined by 
\[
\mathsf{B}_{\cQ_{\rm cf}}= \big\{t \in B_{\cQ_{\rm cf}} \, \big| \, t_0=t_2=0 \big\}\,.
\]
Throughout this Section and the next, for a function $f : B_{\cQ_{\rm cf}} \to Q_{\Tmax}$, we will denote in 
sans-serif
typeface $\mathsf{f}: \mathsf{B}_{\cQ_{\rm cf}} \to Q_{\Tmax}$ its restriction to small quantum cohomology
\[ 
\mathsf{f}(\cQ) \coloneqq f(t)\big|_{t_0=t_2=0; t_1=\log \cQ}\,.
\]
Restricting \eqref{eq:resj} to this locus, $t_0=t_2=0$,  and applying the divisor axiom gives the small $J$-function
\beq
\mathsf{J}(\cQ,z)  
=z \cQ^{H/z}
\l(\mathbf{1}+ \sum_{d,k} \cQ^d H^a \eta^{ab} \bra
\frac{H^b}{z(z-\psi_{1})}\ket_{0,1,d}^Z\r).
\label{eq:Jred}
\eeq
The small $J$-function determines entirely the restriction $\mathsf{S}(\cQ,z)$ of the fundamental solution to small quantum cohomology. From \eqref{eq:fundsol} and the string and divisor equation, we have that
\beq
\mathsf{S} (\mathbf{1}) =  \mathsf{J}/z\,, \quad
\mathsf{S} (H) = \cQ \partial_{\cQ} \mathsf{J}\,.
\label{eq:SJ1}
\eeq
To compute the restriction to the small quantum cohomology locus of the remaining component, $\mathsf{S}(H^2) $,
write
\[
H^a \circ_t H^b = \sum_{c=0}^2 C_{ab}^c (t) H^c.
\]
%
The $\nabla^{(z)}$-flatness condition for the big $J$-function implies that
\[
S(H^2) = \partial_{t_2} J = 
\frac{1}{C_{11}^2}\l( z \partial_{t_1}^2 J -C_{11}^1 \partial_{t_1} J-C_{11}^0 J/z \r)\,, \]
so that, restricting to small quantum cohomology,
%
\bea
\mathsf{S}(H^2) &=& \frac{z \cQ \partial_{\cQ}\mathsf{S}(H) -\mathsf{C}^1_{11}\mathsf{S}(H) -\mathsf{C}^0_{11}\mathsf{S}(\mathbf{1})}{\mathsf{C}^2_{11}} \,.
\label{eq:SJ2}
\eea
%

\subsection{Finite generation and quasi-modularity: genus zero}

Let $A$, $F$, $X$, $(Y_i,W_i^{1/2})_{i=0}^2$ be formal symbols, and consider the ring
\[ 
G_{\Tmax} \coloneqq \frac{Q_{\Tmax}[\{Y_i^{\pm 1},W_i^{\pm 1/2}\}_{i=0}^2]}{\bra \big\{\widehat{\LL}(Y_i), \widehat{\MM}(Y_i,W_i)\big\}_{i=0}^2 \ket}\,,
\] 
where
\beq
\widehat{\LL}(Y) \coloneqq \prod_{j=0}^2 (Y+\alpha_j) + 27 y Y^3 \,, \quad \widehat{\MM}(Y,W) \coloneqq W- \sum_{j=1}^3 j Y^{3-j} e_{j}(\alpha) \,.
\label{eq:LLMM}
\eeq
We will be interested in the  polynomial ring $G_{\Tmax}[A,F,X],$ when the formal variables $A$, $F$, $X$ (respectively, $Y_i$ and $W_i$) are specialised to certain terms in the expansion of a system of flat coordinates for the Dubrovin connection around $z=\infty$ (respectively, $z=0$).
This ring is equipped with an $S_3$-action permuting the indices of the ${\Tmax}$-weights $\alpha_i$ and the generators $Y_i$ and $W_i$. We will show that,  under the restriction \eqref{eq:diagweights} to the diagonal subtorus $\mgp{T}_{\rm F} \subset \Tmax$, 
the invariant ring 
$G^{S_3}_{\mgp{T}_{\rm F}}[F]$ for this action  is closely related to the ring of quasi-modular forms of $\Gamma_1(3)$. This prepares the ground for the next Section, where we shall show that the higher genus $\mgp{T}_{\rm F}$-equivariant Gromov--Witten potentials of $X$ will be elements of one of its graded components.
\\

\subsubsection{Generators at \texorpdfstring{$z=\infty$}{z = infinity}}
\label{sec:genzinf}
To start off, define the $I$-function of $Z$ as
\beq
\label{eq:Ifun}
I(y,z) \coloneqq  z y^{H/z} \sum_{d\geq 0} y^d \frac{\prod_{m=0}^{3d-1} (-3 H-m z)}{\prod_{j=0}^2 \prod_{m=1}^d (H+\alpha_j+mz)}\,.
\eeq
The sum in \eqref{eq:Ifun} is a convergent power series when $0\leq |y|<1/27$. It satisfies the equivariant Picard--Fuchs equation
\beq
\l[\prod_{j=0}^2 (z \theta_y+\alpha_j) + z^3 y \prod_{j=0}^2 (3\theta_y+ j)  \r] I(y,z) =0\,,
\label{eq:pf}
\eeq
where we write $\theta_y \coloneqq 
y \partial_y$ for the logarithmic $y$-derivative. Equivariant mirror symmetry for $\bbP^2$ and the quantum Riemann--Roch theorem \cite{CG07:QuantumRR} recover the small $J$-function in terms of \eqref{eq:Ifun}. Writing
\[
I^{[a]}_{n}(y) \coloneqq [H^a z^{-n}] I(y,z)\,,
\] 
we have
\[
I^{[a]}_{0} = \delta_{a,1}\l(\log y + 3 \sum_{d=1}^\infty \frac{(3d-1)!}{(d!)^3} (-y)^d\r)\,, \quad
I^{[0]}_{1} = 0\,, \quad
I^{[1]}_{1} = -3 e_1(\alpha)\sum_{d=1}^{\infty} \frac{ (3 d-1)! H_d (-y)^d}{(d!)^3}\,, \]
\beq
I^{[2]}_{1} = -\frac{\log^2 y}{2}+ \log y~ I_{1}^{[0]} +9\sum_{d=1}^{\infty}  \frac{(3 d-1)! \left(\psi ^{(0)}(3 d)-\psi^{(0)}(d+1) \right)}{(d!)^3}(-y)^d\,,
\label{eq:Ian}
\eeq
where \[\psi^{(0)}(x)\coloneqq \partial_x \log \Gamma(x)\,, \quad H_d \coloneqq \sum_{k=1}^d \frac{1}{k}\,.\]
In particular, the mirror map is 
\[
\cQ \longrightarrow \cQ(y) = \re^{I^{[1]}_{0}(y)}\,,
\]
and genus zero local mirror symmetry for $Z$ is the statement that \cite{MR2486673}
%
\beq
\mathsf{J}(\cQ(y),z) = I(y,z)\,.
\label{eq:J=I}
\eeq
%
We will define the following analytic functions in $|y|<1/27$ from the first few Taylor coefficients of the $I$-function around $z=\infty$:
\[ 
X(y)  \coloneqq  \frac{1}{1+27 y}\,,\quad
A(y)  \coloneqq  \theta_y \log \cQ = \theta_y I^{[1]}_{0} = \, _2F_1\left(\frac{1}{3},\frac{2}{3};1;-27 y\right) \,, 
\]
\beq
D(y)  \coloneqq  \theta_y(A^{-1} \theta_y I^{[2]}_{1})\,, \quad
E(y) \coloneqq  \frac{\theta_y(A^{-1} \theta_y I^{[1]}_{1})}{e_1(\alpha)}\,,\quad
F(y)  \coloneqq   \theta_y \log A \,.
\label{eq:ADEFXdef}
\eeq
The ring $\bbQ[A,D,E,F,X]$ is not freely generated, as the following Lemma shows.
\begin{lem}
The following relations hold:
\beq
A^2 D=X\,, \quad 
E = \frac{1}{2} \l(D- X\r)\,.
\label{eq:DErel}
\eeq 
In particular, $\bbQ[A,D,E,F,X] \simeq \bbQ[A^{\pm 1}, F,X] \simeq \bbQ[A,D,F]$.
\label{lem:modrel}
\end{lem}
\begin{proof}[Proof (sketch)]
The first relation in \eqref{eq:DErel} is the Zagier--Zinger relation, proved\footnote{With previous appearances in the physics literature (see e.g.~\cite[Eq.~6.17]{Aganagic:2006wq}): it is equivalent to the mirror symmetry expectation that, in B-model coordinates, Yukawa couplings of 1-parameter Calabi--Yau threefolds are given by the inverse of the discriminant of the mirror family.} in \cite{MR2454324}. The second relation follows from
\beq
\frac{2}{e_1(\alpha)}\theta_y I_1^{[1]} - \theta_y I_1^{[2]} +A \log\l(\frac{y}{1+27 y}\r)  = 0\,,
\label{eq:DEXrel}
\eeq 
upon dividing by $A(y)$ and acting 
with $\theta_y$. To show \eqref{eq:DEXrel}, note that, by expanding around $y=0$ and using \eqref{eq:Ian}, it is equivalent to the  convolution identity
\beq 
\sum_{k=0}^{d-1}\frac{1}{d-k}\frac{\left(\frac{1}{3}\right)_k \left(\frac{2}{3}\right)_k}{(k!)^2}=\frac{(3 d)! \left(2 H_d+\frac{1}{d}+3 \psi ^{(0)}(3 d)-3 \psi ^{(0)}(d+1)\right)}{3^{3d} (d!)^3}
\,,
\label{eq:limitconv}
\eeq 
where $(x)_k = \Gamma(x+k)/\Gamma(x)$ is the Pochhammer symbol. 
This identity can be proved by regularising the convolution kernel as $\frac{1}{d-k} \to \frac{(\epsilon)_{d-k}}{(1+\epsilon)_{d-k}}$, 
and considering instead  the hypergeometric identity
\bea
\sum_{k=0}^{d-1} \frac{\left(\frac{1}{3}\right)_k \left(\frac{2}{3}\right)_k (\epsilon )_{d-k}}{(k!)^2 (\epsilon +1)_{d-k}}
&=&
\frac{\epsilon  \, _3F_2\left(\frac{1}{3},\frac{2}{3},-d-\epsilon ;1,1-d-\epsilon;1\right)}{d+\epsilon }  \nn \\ & &
-\frac{(3d)! 
   \, _4F_3\left(1,d+\frac{1}{3},d+\frac{2}{3},-\epsilon ;d+1,d+1,1-\epsilon ;1\right)}{3^{3d}(d!)^3}\,, \nn
\eea 
which follows from writing the sum on the l.h.s.~as a difference of two infinite sums $\sum_{k=0}^\infty - \sum_{k=d}^\infty$. 
%
%
%
%
The left-hand side manifestly reduces to the l.h.s.~of \eqref{eq:limitconv} in the limit $\epsilon\to 0$. The summands in the r.h.s.~are individually divergent\footnote{This is a consequence of the Pfaff--Saalsch\"utz formula \cite[Sec.~2.3.1]{MR0201688}.}  in the limit $\epsilon \to 0$ with simple poles as $\epsilon\to 0$; however, writing the second term as a $d^{\rm th}$ primitive of a $_3 F_2$ hypergeometric function evaluated at the unity, and repeated use of the Thomae relation, allows to show that the that the divergent terms cancel out between the summands, and the $\cO(\epsilon^0)$ coefficient returns the r.h.s.~of  \eqref{eq:limitconv}. We omit the details here; these are available to the curious reader upon request. 
\end{proof}

\subsubsection{Quasi-modular forms of \texorpdfstring{$\Gamma_1(3)$}{Gamma1(3)}}
Consider the leading order term in $1/z$ of the equivariant Picard--Fuchs equation \eqref{eq:pf},
\beq
\l[\theta_y^3 + y \prod_{j=0}^2 (3\theta_y+ j)  \r] \Pi =0\,.
\label{eq:pfinf}
\eeq
Applying the Frobenius method around the maximal unipotent monodromy point $y=0$ shows that the space of solutions of \eqref{eq:pfinf} is spanned by $\{1, I_0^{[1]}, I_1^{[2]}\}$. The two non-constant solutions are furthermore known to be periods of the universal family of elliptic curves 
\[ 
\big\{([x_0 : x_1 : x_2],y) \in \bbP^2 \times \bbP(3,1): x_0^3 + x_1^3 + x_2^3 + y^{-1/3} x_0 x_1 x_2 = 0\big\}\,,
\] 
parametrised by a point on the modular curve $y \in \cC_1(3) \coloneqq \bbP(3,1) \setminus \{[0:1], [1:0],[-1/27:1]\}$. The base $\cC_1(3)$ of the family is the coarse moduli space of the modular curve $[\mathbb{H}/\Gamma_1(3)]$ for the congruence subgroup
$\Gamma_1(3)<\mathrm{SL}_2(\bbZ)$, where
\[ 
\Gamma_1(3) \coloneqq
\l\{
\l(\begin{matrix}
a & b \\
c & d
\end{matrix}\r)
\in \mathrm{SL}(2, \bbZ)\, \bigg| 
\l(\begin{matrix}
a & b \\
c & d
\end{matrix}\r)
\equiv
\l(\begin{matrix}
1 & \star \\
0 & 1
\end{matrix}\r)\,\mathrm{mod}~3
\r\}
\] 
and $\mathbb{H} = \{\tau \in \bbC \, | \, \mathrm{Im}(\tau)>0\}$ is the upper half-plane. The modular parameter is related to the periods on the family and the Hauptmodul as 
\beq
\tau = \frac{1}{2}
 + \frac{1}{2\pi\ri}
\frac{\theta_y I_1^{[2]}}{\theta_y I_0^{[1]}}\,, \quad y= -\frac{\eta (3 \tau )^{12}}{\eta (\tau )^{12}+27 \eta (3 \tau )^{12}}\,,
\label{eq:tau}
\eeq
where $\eta(\tau) \coloneqq 
\re^{\pi \ri \tau/12} \prod_{n=1}^\infty (1-\re^{2\pi \ri n \tau})$ is the Dedekind eta function. The boundary points $\tau=\ri \infty$ and $\tau = 0$ correspond to the large radius ($y=0$) and conifold ($y=-1/27$) points in the extended K\"ahler moduli space of $Z$. We will write \[
\mathrm{Mod}^{[k]}(\Gamma_1(3))\,, \qquad 
\mathrm{QMod}^{[k,l]}(\Gamma_1(3))\,,
\]
for the complex vector spaces of modular forms  for $\Gamma_1(3)$ of weight $k$ (respectively, quasi-modular forms  for $\Gamma_1(3)$ of weight $k$ and depth $\leq l$), and
\[
\mathrm{Mod}(\Gamma_1(3))
\,, \qquad \mathrm{QMod}(\Gamma_1(3)) 
\,,
\]
for the associated rings graded by weight.
We will use Gothic typeface
\[
\mathfrak{Mod}^{[k]}(\Gamma_1(3))\,, \qquad 
\mathfrak{QMod}^{[k,l]}(\Gamma_1(3))\,, \qquad 
\mathfrak{Mod}(\Gamma_1(3))\,, \qquad 
\mathfrak{QMod}(\Gamma_1(3))\,,\]
to denote the corresponding spaces of modular functions, where we drop the requirement of holomorphicity 
throughout $\mathbb{H}$.

Let 
\[
B(y) \coloneqq  \frac{6 F(y)-X(y)+2}{D(y)}\,, \quad C(y) \coloneqq \frac{A(y)}{D(y)}\,.
\]
Expressed as functions of the elliptic modulus \eqref{eq:tau}, it is well-known (see e.g.~\cite[Sec.\,6.3]{Aganagic:2006wq}, \cite[Sec.\,4.1]{BFGW21:HAE}) that 
\beq A = \l(\frac{\eta(\tau)^9}{
\eta(3\tau)^3}+27\frac{\eta(3\tau)^9}{
\eta(\tau)^3}\r)^{1/3} \,, \qquad B = \frac{E_2(\tau)+3 E_2(3\tau)}{4}\,, \qquad
 C= \l(\frac{\eta(\tau)^3}{\eta(3\tau)}\r)^3 \,,
\label{eq:ABCmod}
\eeq
where $E_2(\tau)$ is the second Eisenstein series \[E_2(\tau)=1 - 24 \sum_{n=1}^\infty
\frac{n\re^{2\pi\ri n \tau}}{1 - \re^{2\pi\ri n \tau}}\,.\] In particular, we have 
\[
A\in \mathrm{Mod}^{[1]}(\Gamma_1(3))\,, \quad B\in \mathrm{QMod}^{[2,1]}(\Gamma_1(3))\,, \quad C \in \mathrm{Mod}^{[3]}(\Gamma_1(3))\,,
\]
\[D^{-1} \in \mathfrak{Mod}^{[2]}(\Gamma_1(3))\,, \quad X \in \mathfrak{Mod}^{[0]}(\Gamma_1(3))\,, \quad F \in \mathfrak{QMod}^{[0,1]}(\Gamma_1(3))\,.\] The following statement was proved in \cite[Lem.~4.2]{BFGW21:HAE}:
\begin{prop}[\cite{BFGW21:HAE}]
The quasi-modular forms $A$, $B$ and $C$ are algebraically independent over $\bbC$, and they freely generate $\mathrm{QMod}(\Gamma_1(3))$,
\[
\mathrm{Mod}(\Gamma_1(3)) \simeq \bbC[A,C]\,, \quad \mathrm{QMod}(\Gamma_1(3)) \simeq \bbC[A,B,C]\,.
\]
\label{prop:algindABC}
\end{prop}
%
%
%
%
%

\subsubsection{Generators at \texorpdfstring{$z=0$}{z=0}}
Expand now the small J-function \[\mathsf{J}(\cQ,z) \coloneqq \sum_{i=0}^{2}\mathsf{J}_i(\cQ,z) \mathfrak{f}_i\] along the basis of normalised idempotents of the quantum product. Its coefficients
%
%
 admit
\cite
{Givental:1997hm,
MR1653024,
MR1901075}
a formal asymptotic expansion around $z=0$,
\beq
\mathsf{J}_i \simeq z\re^{\mathsf{u}_i/z} \sum_{n \geq 0} \mathsf{S}^{[n]}_{0i} z^n\,,
\label{eq:asymp}
\eeq
where $\mathsf{u}_i +\alpha_i \log y$ 
are the restriction of the canonical coordinates to the small quantum cohomology locus, and $\mathsf{S}^{[n]}_{0i}$ 
is analytic in $\mathsf{B}_{\cQ_{\rm cf}}$. 
Substituting the expansion \eqref{eq:asymp} into \eqref{eq:pf}, and solving order by order in $z$, fixes $\mathsf{u}_i(y)$ up to a constant and gives rise to a set of recursive differential relations for the formal Taylor coefficients $\mathsf{S}^{[n]}_{0i}$.
For example, at the leading order in $z$ we find that
%
\beq
\theta_y \mathsf{u}_i=Y_i\,,
\label{eq:uili}
\eeq 
where $Y_i(y)$ is the root of the polynomial $\widehat{\LL}$ in \eqref{eq:LLMM} satisfying $Y_i(0)=-\alpha_i$. At $\cO(z)$, we find
\beq
\mathsf{S}^{[0]}_{i,0} = c W_i^{-1/2}
\label{eq:Ji0}
\eeq 
up to a multiplicative constant $c \in \Gm$, where
\beq 
W_i(y) \coloneqq \sum_{j=1}^3 j Y_i(y)^{3-j} e_{j}(\alpha)
\label{eq:delta1}
\eeq 
%
%
and we periodically identify $Y_{i+3}=Y_i$. 
%
%
Let then
\beq
G_{\Tmax} \coloneqq \frac{Q_{\Tmax}[(Y_i^{\pm}, W_i^{\pm 1/2})_{i=0}^2]}{\bra \{ \widehat{\LL}(Y_i), \widehat{\MM}(Y_i, W_i)\}_{i=0}^2\ket}\,.
\label{eq:GTquot}
\eeq
%
%
%

The mirror theorem \eqref{eq:J=I} allows to 
express the structure constants $\mathsf{C}_{ab}^c$ and the Jacobian matrix of the change-of-variables from canonical to flat coordinates in terms of the generators $A$, $X$, $Y_i$ and $W_i^{-1/2}$.

\begin{prop}
Let
\[\Phi_{ai} \coloneqq (\de_{t_a} u_i)\big|_{\mathsf{B}_{\cQ_{\rm cf}}}\,, \quad \Psi_{ai}\coloneqq \Phi_{ai} \Delta_i^{-1/2}\big|_{\mathsf{B}_{\cQ_{\rm cf}}}\]
be, respectively, the change-of-basis matrices from the canonical frame (resp., normalised canonical frame) to the flat frame, restricted to the small quantum cohomology locus. Then,
\[
\Phi_{ia} \in G_{\Tmax}(A)\,, \quad \Psi_{ia} \in G_{\Tmax}(A)\,, \quad  \mathsf{C}_{ab}^c \in Q_{\Tmax}[A^{\pm 1}, X^{\pm 1}]\,.
\]
Explicitly, the change-of-basis matrices are given by
\beq
\Phi_{0i} = 1\,,  \quad
\Phi_{1i}
= \frac{Y_i}{A}, \quad 
\Phi_{2i}
= \frac{\mathsf{C}_{12}^1Y_i+A \mathsf{C}_{12}^0 }{Y_i - A \mathsf{C}_{12}^2}\,, \quad
\Delta_i^{-1} = \frac{1}{3\epsilon_1 \epsilon_2 W_i}\,, 
\label{eq:phi}
\eeq
and the structure constants of the small quantum product are:
\beq
\begin{gathered}
\mathsf{C}_{0a}^b = \delta_a^b, \quad \mathsf{C}_{ab}^0 = \eta^{00}\eta_{ab}\,, \\
\mathsf{C}^1_{11} =
\frac{e_1(\alpha)X(1-A^2)}{2A^3}
\,,\quad
\mathsf{C}^2_{11} = 
\frac{X}{A^3}
\,, \quad
\mathsf{C}^2_{12} = 
-\frac{e_1(\alpha )(A^2+1)X}{2A^3}
\,, \\
\mathsf{C}^1_{12} = \frac{\left(A^4-1\right) X e_1(\alpha ){}^2}{4 A^3}-A e_2(\alpha )  \,, \quad \mathsf{C}^2_{22} = \frac{\left(A^2+1\right)^2 X e_1(\alpha ){}^2-4 A^4 e_2(\alpha )}{4 A^3}\,, \\
 \mathsf{C}^1_{22} = -\frac{A^3 e_3(\alpha )}{X}+\frac{1}{2} \left(A^2+1\right) A e_1(\alpha ) e_2(\alpha )-\frac{\left(A^2-1\right) \left(A^2+1\right)^2 X e_1(\alpha ){}^3}{8 A^3}\,. 
 \end{gathered}
\label{eq:allc}
\eeq  
\end{prop}

\begin{proof}
The two equalities in the first row of \eqref{eq:allc} follow from the String Axiom, together with \eqref{eq:etaflat}.
%
%
Write now
\[
\mathsf{S}(H^b) =  H^b  + \sum_{a=0}^2 \frac{\mathfrak{s}^a_{b}}{z} + \cO\l(\frac{1}{z^2}\r)
\] 
for the leading terms of the expansion of \eqref{eq:fundsol} at $z=\infty$. In particular, by \eqref{eq:SJ1}, \eqref{eq:SJ2} and \eqref{eq:J=I}, we get
\beq 
\mathfrak{s}^a_0 = \log \cQ~\delta^{a}_1\,, \quad 
\mathfrak{s}^a_1 = A^{-1}\theta_y I_1^{[a]}\,.
\label{eq:ssH1}
\eeq 
Recall that the fundamental solution satisfies the symplectic condition
\beq 
\eta\big(S(t, z) H^a, S(t, -z) H^b\big) = \eta_{ab}\,,
\label{eq:sympl}
\eeq 
as can be shown by differentiating the l.h.s.~with respect to $t$. Imposing  \eqref{eq:sympl} then further sets
\beq 
\mathfrak{s}^0_2 = \log \cQ~e_3(\alpha)\,, \quad
\mathfrak{s}^2_2 = \mathfrak{s}^1_1 - e_1(\alpha)\mathfrak{s}^2_1\,,
\label{eq:ssH2}
\eeq 
thus determining $\mathfrak{s}^a_b$ in terms of $\{I_1^{[c]}\}_{c=0,1,2}$ for all $(a,b)$ except $(a,b)=(1,2)$. From the definition \eqref{eq:fundsol} and the Divisor Axiom, we have that
\[
\mathsf{C}^a_{1b} = \cQ \partial_\cQ \mathfrak{s}^a_b = A^{-1} \theta_y \mathfrak{s}^a_b\,,
\]
and the 
second  row of \eqref{eq:allc}
follows from \eqref{eq:ssH1}--\eqref{eq:ssH2}.
%
In order to compute the remaining structure constants 
note that, by definition of idempotency,
\beq
\Phi_{ai} \Phi_{bi}
= \sum_c \mathsf{C}_{ab}^c  \Phi_{ci}\,.
\label{eq:duidta}
\eeq
%
The claimed expressions for $\Phi_{ai}$ in \eqref{eq:phi} are implied by \eqref{eq:uili} together with the $(a,b)=(0,b)$ and $(a,b)=(1,2)$ cases of \eqref{eq:duidta}. The inverse change-of-basis matrix, $(\Phi)^{-1}_{ia}$, is directly computed using \eqref{eq:delta1} to be
\[
 (\Phi^{-1})_{i0}  =  \frac{e_3(\alpha)}{W_i}\,, \quad
 (\Phi^{-1})_{i1}  = -\l(e_1(\alpha )\frac{\left(A^2+1\right) }{2 A }+\frac{A Y_i}{X}\r) \frac{Y_i}{W_i} \,, \quad
 (\Phi^{-1})_{i2}  = -\frac{Y_i  }{A W_i}\,,
\]
%
from which the square norm of the canonical idempotents can be read off as
\[ 
\Delta_i^{-1} =
\eta(\mathfrak{e}_i,\mathfrak{e}_i) = 
\eta(\mathfrak{e}_i,\mathbf{1}) =  (\Phi^{-1})_{0i} \eta(\mathbf{1},\mathbf{1}) = \frac{1}{3\epsilon_1 \epsilon_2 W_i} \,.
\] 
It remains to compute the last three structure constants in \eqref{eq:allc}.
%
%
The $(a,b)=(1,1)$ and $(2,2)$ cases of \eqref{eq:duidta}, together with the single non-trivial WDVV relation,
\[ 
0= \sum_{a,b=0}^2 \eta_{ab} \l(\mathsf{C}^a_{12}\mathsf{C}^b_{12}-\mathsf{C}^a_{11}\mathsf{C}^b_{22}\,
\r)\,,
\]
give a rank-three linear system over $G_{\Tmax}(A,X)$ for the three  unknowns $\mathsf{C}^1_{12}$, $\mathsf{C}^1_{22}$ and $\mathsf{C}^2_{22}$. Imposing  
the Picard--Fuchs relations $\{\widehat{\LL}(Y_i)=0\}_{i=0,1,2}$ on the solution eliminates all dependence on $Y_i$ and
returns the last two rows of \eqref{eq:allc},  completing the proof.
%
%
%
\end{proof}
%
%

\subsubsection{Finite generation of the fundamental solution}

There are sufficient conditions \cite[Sec.\,3.2]{lho2018note} to ensure that all coefficients in the formal expansion \eqref{eq:asymp} belong to $G_{\Tmax}$.

\begin{defn}
  We say that $\mgp{T}_{\rm R} \subset \Tmax$ is a \define{rational} sub-torus of $\Tmax$ if the weights $\alpha=(\alpha_0, \alpha_1, \alpha_2)$ of its action on $X$ satisfy
\ben
\item $e_2(\alpha)^2 = 3 e_1 (\alpha) e_3(\alpha)$,
\item $\alpha_i\neq \alpha_j$ for $i\neq j$.
\een
\label{def:adm}
\end{defn}
\noindent The first condition  in \cref{def:adm} implies that the expression of $W_i$ in \eqref{eq:delta1} factorises as a square
\beq
W_i = \frac{\big(e_2(\alpha) Y_i + 3e_3(\alpha)\big)^2 }{3e_3(\alpha)}\,,
\label{eq:Wadm}
\eeq 
hence $W_i^{-1/2}$ is a rational function of $Y_i$. The second  condition is the (generic) statement that the torus action on $Z$ has only 0-dimensional point-wise fixed strata. A one-parameter family of weights satisfying the constraints in \cref{def:adm} is given by
\beq
\alpha_i= \frac{\epsilon _1+\epsilon _2}{3}\l(\frac{1}{1- \delta  \re^{\frac{2 \ri \pi  i}{3}}}\r) 
\label{eq:adelta}
\eeq
for $\delta \in \Gm$. Two important boundary points are obtained when $\delta \to \infty$, where this reduces to the special Calabi--Yau torus chosen by \cite{LP18:HAE} with weights formally specialised to cubic roots of unity; and when $\delta \to 0$, when it corresponds to the fibrewise torus $\mgp{T}_{\rm F}$.

The following Lemma was proved in \cite[Sec.~3.2]{lho2018note} for the restriction to rational sub-tori of $\Tmax$. 

\begin{lem}[\cite{lho2018note}]
Let  $\mgp{T}_{\rm R} \subset \Tmax$ be a rational sub-torus of $\Tmax$. Then~
%
$\mathsf{S}^{[n]}_{0i} \in G_{\mgp{T}_{\rm R}}\,.$
%
\label{lem:lho}
\end{lem}
By combining \cref{lem:lho} with \eqref{eq:SJ2}, we can deduce a statement akin to \cref{lem:lho} for the $z=0$ expansion of the components $\mathsf{S}(H^a)$ of the fundamental solution with $a>0$. %
Note that, from \eqref{eq:SJ1}, \eqref{eq:SJ2}, \eqref{eq:J=I}, \eqref{eq:ADEFXdef} and \eqref{eq:allc}, the components of the fundamental solution in B-model coordinates $\cQ \to \cQ(y)$ read
\beq
\mathsf{S}(\mathbf{1})  = I(y)/z\,, \quad
\mathsf{S}(H)  = \frac{\theta_y I(y)}{A}\,,  \quad
\mathsf{S}(H^2) =  \frac{z \theta_y\mathsf{S}(H)- e_1(\alpha )E(y) \mathsf{S}(H) }{D(y)}\,.  
\label{eq:SHb}
\eeq
In particular, the expression of $\mathsf{S}(H^k)(y)$ is entirely reconstructed from the knowledge of the $I$-function, and likewise for their formal expansion around $z=0$. It will be convenient, following \cite{LP18:HAE}, to express the expansion of 
\[\mathsf{S}_{aj}\coloneqq \re^{-\mathsf{u}_j/z}\mathfrak{f}_j^*\mathsf{S}(H^a)\] in the form
\beq
\mathsf{S}_{0j} = \sum_{n=0}^\infty\,
\mathsf{S}^{[n]}_{0j} z^n\,, \quad
\mathsf{S}_{1j} = \frac{Y_j}{A} \sum_{n=0}^\infty 
\mathsf{S}^{[n]}_{1j} z^n\,, \quad
\mathsf{S}_{2j} = \frac{Y_j^2 }{A D} \sum_{n=0}^\infty 
\mathsf{S}^{[n]}_{2j} z^n\,.
\label{eq:asympS}
\eeq
Plugging \eqref{eq:asympS} into the second and third equalities of \eqref{eq:SHb} and solving order by order in $z$ gives the following recursive relations for $\mathsf{S}^{[n]}_{aj}$, $a>0$. 
\begin{prop}
The following equalities hold:
\beq 
\begin{split}
\mathsf{S}^{[n+1]}_{1i} &=
\mathsf{S}^{[n+1]}_{0i} + \frac{\theta_y\mathsf{S}^{[n]}_{0i}}{Y_i}\,, \\
\mathsf{S}^{[n+1]}_{2i} &= \mathsf{S}^{[n+1]}_{1i}+ \frac{\theta_y\mathsf{S}^{[n]}_{1i}}{Y_i}+ \l[\frac{\theta_y Y_i}{Y_i^2}-\frac{F}{Y_i}\r] \mathsf{S}^{[n]}_{1i} + \frac{e_1(\alpha) X (A^2-1)}{2Y_i A^2 }\mathsf{S}^{[n+1]}_{1i}\,.
\end{split}
\label{eq:Sain}
\eeq
In particular, for $\mgp{T}_{\rm R} \subset \Tmax$ a rational sub-torus, we have that
\[
\mathsf{S}^{[n]}_{0i}\in G_{\mgp{T}_{\rm R}}\,, \quad 
\mathsf{S}^{[n]}_{1i} \in G_{\mgp{T}_{\rm R}}\,, \quad 
\mathsf{S}^{[n]}_{2i} \in G_{\mgp{T}_{\rm R}} + \frac{1}{A^2}G_{\mgp{T}_{\rm R}}+F\, G_{\mgp{T}_{\rm R}} \,.
\] 
%
\label{prop:fingenS}
\end{prop}
\begin{proof}
The first part of the Proposition is a straightforward substitution of \eqref{eq:asympR} into \eqref{eq:SHb}, using \eqref{eq:ADEFXdef}. The second part follows from \eqref{eq:Sain} together with \cref{lem:modrel,lem:lho}.
\end{proof}

\subsubsection{The fibrewise action}
\label{sec:fwaction}
The Weyl group $\cW(\mathrm{PGL}_3(\bbC)) \simeq S_3$ of the automorphism group of $\bbP^2$ induces an action of the symmetric group in three elements on $G_{\Tmax}$, 
\[
\bary{rcl}
S_3 \times G_{\Tmax} & \longrightarrow & G_{\Tmax} \\
\Big(\sigma, f\big((\alpha_i, Y_i, W_i^{1/2})_{i=0}^2\big)\Big) & \longmapsto & f((\alpha_{\sigma(i)}, Y_{\sigma(i)}, W_{\sigma(i)}^{1/2})_{i=0}^2)
\eary 
\]
permuting the fixed points $\{P_0, P_1, P_2\}$ and the associated localisation data.
%
%
%
%
We will be interested in the quotient ring under this action when we restrict to the diagonal (fibrewise) torus $\mgp{T}_{\rm F}$.
%
By \eqref{eq:Wadm}, under the specialisation  \eqref{eq:diagweights} we have 
\[
W_i = -\frac{(\epsilon_1+\epsilon_2)(3Y_i+\epsilon_1+\epsilon_2)^2}{9}\,,
\]
so that
\[
G_{\mgp{T}_{\rm F}} = \frac{\bbC(\epsilon_1, \epsilon_2)[(Y_i^{\pm 1}, W_i^{\pm1/2})_{i=0}^2]}{\bra \{ \LL(Y_i),\MM(Y_i,W_i) \}_{i=0}^2  \ket} \simeq 
\frac{\bbM[(Y_i^{\pm 1}, (3Y_i+\epsilon_1+\epsilon_2)^{\pm 1})_{i=0}^2]}{\bra \{ \LL(Y_i) \}_{i=0}^2  \ket}\,,
\]
where 
\[
\LL(x) \coloneqq \widehat{\LL}(x)|_{\alpha_i=(\epsilon_1+\epsilon_2)/3}\,, \quad \MM(x,y) \coloneqq \widehat{\MM}(x,y)|_{\alpha_i=(\epsilon_1+\epsilon_2)/3}\,, \quad
\bbM \coloneqq \frac{\bbC(\epsilon_1, \epsilon_2, \epsilon_3)}{\bra \epsilon_3^2 - (\epsilon_1+ \epsilon_2) \ket}\,.
\]
Since $\alpha_0=\alpha_1=\alpha_2$, 
an $S_3$-invariant element $f\in G^{S_3}_{\mgp{T}_{\rm F}}$ is just a symmetric Laurent polynomial in $Y_i$ and $3Y_i+\epsilon_1+\epsilon_2$:
\[ 
f=\frac{g(Y_0,Y_1,Y_2)}{ e_3(Y)^m e_3(3Y+\epsilon_1+\epsilon_2)^n}
\] 
for some $m$, $n \in \bbZ_{\geq 0}$ and  
$g \in \bbM \otimes_\bbZ \Lambda_3$ an $\bbM$-valued symmetric polynomial in $Y=(Y_0, Y_1, Y_2)$.
From \eqref{eq:LLMM}, we find 
%
\[
G_{\mgp{T}_{\rm F}}^{S_3} \simeq \frac{G_{\mgp{T}_{\rm F}}}{\bra \cN(Y)\ket}
\,,
\]
with 
\[ 
\cN(Y)=\l\{e_k(Y) + \frac{(-1)^{k} (\epsilon_1+\epsilon_2)}{3} \binom{3}{k}   X\r\}_{k=1}^3 \bigcup \l\{e_3\big(3Y+\epsilon_1+\epsilon_2\big) +
(X-1) \left(\epsilon _1+\epsilon _2\right)^3
\r\}\,,
\] 
and therefore
\beq
G_{\mgp{T}_{\rm F}}^{S_3} \simeq \bbM[X^{\pm 1}, (1-X)^{-1}]\,.
\label{eq:GTFS3}
\eeq
\begin{defn}
Let $\mgp{T} \subset \Tmax$ be a sub-torus of the maximal torus of $Z$. 
An element of $G_{\mgp{T}}$, viewed as a polynomial in $Y_i^{\pm 1}$ and $W_i^{-1/2}$, will be called {\it regular}, if 
\ben
\item it is even in $W_i^{-1/2}$;
\item its image under the specialisation map 
\[ 
G_{\mgp{T}} \to Q_{\mgp{T}}(\!(y)\!)\,, \quad (Y_i, W_i) \to (Y_i(y),W_i(y)) 
\]
is contained in $(\epsilon_1 \epsilon_2)^{-1} \bbQ[\epsilon_1, \epsilon_2,\alpha_0,\alpha_1,\alpha_2]\llbracket y\rrbracket$. 
\een
We will denote by $G_{{\mgp{T}}}^{\rm reg} \subset G_{{\mgp{T}}}$ the subring of regular elements of $G_{\mgp{T}}$.
\label{def:regel}
 \end{defn}

Let now $f \in G_{\mgp{T}_{\rm F}}^{\rm reg}$
be an $S_3$-invariant regular element of $G_{\mgp{T}_{\rm F}}$. By \eqref{eq:GTFS3} and the first condition in \cref{def:regel}, $f$ is a polynomial in $X^{\pm 1}$ and $(1-X)^{-1}$ with coefficients in $\bbQ(\epsilon_1, \epsilon_2)$, as only even powers of $(\epsilon_1+\epsilon_2)^{1/2}$ appear in the expression of $W_i^{-1}$. By the second condition in \cref{def:regel}, it has to be regular at $y=0$,
meaning that it must be pole-free at $X=1$.
We therefore have
\beq
\l(G_{\mgp{T}_{\rm F}}^{\rm reg}\r)^{S_3} \simeq \frac{1}{\epsilon_1 \epsilon_2}\bbQ[\epsilon_1,\epsilon_2][X^{\pm 1}]\,.
\label{eq:Greg}
\eeq 

\subsection{Finite generation and quasi-modularity: higher genus}

We will use the Givental--Teleman quantised symplectic loop group action on the space of semi-simple cohomological field theories  \cite{MR1866444, MR1901075, MR2917177}  
to leverage the finite generation results of the previous Section from genus zero to higher genus. \\ We briefly review the relevant background.
%
%
%
%
By semi-simplicity, there exists a unique formal one-parameter family of sections of $\mathrm{End}(\cV_{\Tmax}) \simeq T^*\cV_{\Tmax}$,  $\{ R_j \in \Omega^1(\cV_{\Tmax}) \otimes \re^{u^j/z}Q_{\Tmax}\llbracket z\rrbracket \}_{i=0}^2$,
satisfying:

\medskip

\begin{description}[itemsep = 0.5em]
\item[(Asymptotics)] \smash{$R_j(t,z)(H^a) =\sum_{a,i} R_{aj}(t, z)\re^{u_j/z}$}, with \smash{$R_{aj} \in \bbC\llbracket z\rrbracket$} for fixed $t$\,;
\item[(Flatness)] \smash{$\nabla_Z^{(\eta,z)} R_j = 0$}\,;
\item[(Symplecticity)] \smash{$\sum_j R_{ai}(t,z) \eta^{ab} R_{bj}(t, -z) = \delta_{ij}$}\,,
\item[(Normalisation)] the small quantum cohomology restriction $R(t) \to \mathsf{R}(\cQ)$ satisfies
\[ 
\mathsf{R}_{aj}\big|_{\cQ=0} 
= \Psi_{aj}\big|_{\cQ=0} \re^{-c_j(z)}\,,
\] 
where
\[ 
c_j(z) \coloneqq \sum_{n=1}^\infty 
p_{2n-1}\l(w_j^{(0)},\dots, w_j^{(4)}\r) \frac{B_{2n} z^{2n-1}}{2n(2n-1)}\,,
\] 
and $p_n(x_1, \dots, x_k) = \sum_{i=1}^k x_i^n$ is the $n^{\rm th}$ Newton polynomial in $k$ variables.
\end{description}



%
%
%
We will refer to $\mathsf{R}_{aj}$ as the {\it canonical R-matrix of $Z$}, which we will now compute. 
Let $\{\mathfrak{f}_i^*\}_{i=0,1,2}$ be the dual basis at $t\in B_{\rm cf}$ of $\{\mathfrak{f}_i\}_{i=0,1,2}$, and define
\[
R_{j}(H^a) \coloneqq 
\mathfrak{f}_j^* S(H^a) N_j(z)
\]
for some $N_j \in Q_{\Tmax}\llbracket z\rrbracket$ to be specified. Then $R_j$ satisfies ({\bf Flatness}) by \eqref{eq:fundsol}, and ({\bf Symplecticity}) by \eqref{eq:sympl}. Its small quantum cohomology restriction, 
\[
\mathsf{R}_j(H^a) \coloneqq \re^{\mathsf{u}_j/z} \mathsf{R}_{aj} = \re^{\mathsf{u}_j/z}\mathsf{S}_{aj} N_j\,,
\]
has, by \eqref{eq:asymp} and \eqref{eq:SHb}, an asymptotic expansion at $z=0$ of the form
\beq 
\mathsf{R}_j(H^a) \simeq \re^{\mathsf{u}_j/z} \mathsf{R}_{aj} = \re^{\mathsf{u}_j/z} \sum_{n \geq 0} \rho^{[a]}_{j,n} z^n \,.
\label{eq:asympRrho}
\eeq 
%
Therefore, $\mathsf{R}_j$ also satisfies ({\bf Asymptotics}). Noting that, by \eqref{eq:Jred} and \eqref{eq:SHb},
\[ 
\mathsf{R}_{aj} \big|_{\cQ=0} = N_j \Psi_{aj}\big|_{\cQ=0}
\] 
and imposing ({\bf Normalisation}) we conclude by uniqueness that $N_j= \re^{-c_j}$. The canonical R-matrix therefore reads
\beq 
\mathsf{R}_{aj}(\cQ,z) = \mathsf{S}_{aj} \re^{-c_j(z)} \,.
\label{eq:canR}
\eeq 
In particular, from \eqref{eq:SHb}, we have
\beq
\mathsf{R}_{1j}  = \frac{z \theta_y R_{0j}+Y_j R_{0j}}{A}\,, \quad
\mathsf{R}_{2j} =  \frac{z \theta_y\mathsf{R}_{1j}+ (Y_j - e_1(\alpha )E(y) )\mathsf{R}_{1j} }{D(y)}\,.  
\label{eq:RHb}
\eeq
As we did in \eqref{eq:asympS}, it will be convenient to rewrite the coefficients of the asymptotic expansion \eqref{eq:asympRrho} as
\beq
\rho^{[0]}_{j,n} = 
\mathsf{R}^{[0]}_{j,n}\,, \quad
\rho^{[1]}_{j,n} = \frac{Y_j}{A}  
\mathsf{R}^{[1]}_{j,n} \,, \quad
\rho^{[2]}_{j,n} = \frac{Y_j^2}{A D}
\mathsf{R}^{[2]}_{j,n}\,.
\label{eq:asympR}
\eeq
%
A verbatim repetition of the argument used for \cref{prop:fingenS} then shows the following:
\begin{prop}
The following equalities hold:
\bea
\mathsf{R}^{[1]}_{i,n+1} &=& 
\mathsf{R}^{[0]}_{i,n+1} + \frac{\theta_y\mathsf{R}^{[0]}_{i,n}}{Y_i}\,, \nn \\
\mathsf{R}^{[2]}_{i,n+1} &=& \mathsf{R}^{[1]}_{i,n+1}+ \frac{\theta_y\mathsf{R}^{[1]}_{i,n}}{Y_i}+ \l[\frac{\theta_y Y_i}{Y_i^2}-\frac{F}{Y_i}\r] \mathsf{R}^{[1]}_{i,n} + \frac{e_1(\alpha) X (A^2-1)}{2Y_i A^2 }\mathsf{R}^{[1]}_{i,n+1}\,. \nn
\eea
In particular, for $\mgp{T}_{\rm R} \subset \Tmax$ a rational sub-torus, we have that
\begin{flalign*}
&&\mathsf{R}^{[0]}_{i,n}\in G_{\mgp{T}_{\rm R}}\,, \quad 
\mathsf{R}^{[1]}_{i,n} \in G_{\mgp{T}_{\rm R}}\,, \quad 
\mathsf{R}^{[2]}_{i,n} \in G_{\mgp{T}_{\rm R}} + \frac{1}{A^2}G_{\mgp{T}_{\rm R}}+F\, G_{\mgp{T}_{\rm R}} \,. &&\qed
\end{flalign*}
\label{prop:fingenR}
\end{prop}
\label{sec:Rmatrix}
\vspace{-.5em}
\subsubsection{Higher genus reconstruction}

%
%
%
%
The main output of the Givental--Teleman reconstruction theorem \cite{MR1866444, MR1901075, MR2917177}
is an explicit formula for 
$\refGWg{g}{K_{\bbP^2}}{\Tmax}$
as a Feynman graph sum determined by the canonical R-matrix. We fix some notation first.


\begin{defn}
    
An $m$-decorated stable graph $\Gamma=\left(\cV_\Gamma, \HH_\Gamma, \cE_\Gamma, \LL_\Gamma, \DD_\Gamma, \mathsf{g},\mathsf{i}, \mathsf{l}, \mathsf{v}, \mathsf{q}\right)$ is the datum of
\bit
\item $\cV_\Gamma$: a set of vertices, endowed with
\bit 
\item 
a {\it genus assignment} $\mathsf{g}: \cV_\Gamma \rightarrow \mathbb{Z}_{\geq 0}$;
\item a {\it marking assignment}  $\mathsf{q}: \cV_\Gamma \rightarrow \{0,1,\dots, m\}$;
\eit 
\item $\HH_\Gamma$: a set of half-edges, equipped with 
\bit
\item an involution $\iota: \HH_\Gamma \rightarrow \HH_\Gamma$; 
\item $\mathsf{v}: \HH_\Gamma \rightarrow \cV_\Gamma$: the assignment of the attaching vertex to the half-edge;
\eit
\item $\cE_\Gamma$: a set of compact edges, with self-edges being allowed, defined by the 2-orbits $e=\{h_1, h_2\}$ of $\iota: \HH_\Gamma \rightarrow \HH_\Gamma$;
\item $\LL_\Gamma \sqcup \DD_\Gamma$: a splitting of the set of fixed points of $\iota$ (the set of leaves $\HH_\Gamma \setminus \cE_\Gamma$) into a set of {\it marked leaves} $\LL_\Gamma$ and {\it dilaton leaves} $\DD_\Gamma$;
\item  $\mathsf{q}: \LL_\Gamma \rightarrow \{0,1,\dots, m\}$: a {\it marking assignment} for the marked leaves;
\item  $\mathsf{a}: \HH_\Gamma \rightarrow \bbZ_{\geq 0}$: an {\it ancestor assignment} for the half-edges;
\eit

and further satisfying the following two conditions:

\ben 
\item the pair $\left(\cV_\Gamma, \cE_\Gamma\right)$ defines a connected graph;
\item for each vertex ${v} \in \cV_\Gamma$, denoting by $\HH(v)$ the set of half-edges incident to $v$, 
the stability condition
\begin{equation*}
2\mathsf{g}({v}) - 2 + |\HH(v)| > 0\,
\end{equation*}
is satisfied.
\end{enumerate}

\end{defn}

We define the genus of $\Gamma$ to be
$\mathsf{g}(\Gamma)=h^1(\Gamma)+\sum_{{v}\in\cV_\Gamma}\mathsf{g}({v})$. We will denote by $\cG^{[m]}_{g,n}$ the set of all $m$-decorated stable graphs of genus $g$ and with $n$ marked legs, up to isomorphism. 
%
%

\begin{figure}[t]
\centering
\begin{tikzpicture}[baseline=(current  bounding  box.center),scale=1.2]
	\draw[thick] (-0.5,0) circle[radius=0.5];
	\draw[thick] (0,0) -- (2,0) -- (4,0);
	\draw[thick] (2,2) -- (2,0);
	\draw[thick] (1,0) ellipse (1 and 0.5);
	\draw[fill=white,thick] (0,0) circle[radius=4pt];
	\draw[color=red] (0,0) node{\footnotesize 2};
	\draw[fill=white,thick] (2,0) circle[radius=4pt];
	\draw[color=red] (2,0) node{\footnotesize 3};
	\draw[fill=white,thick] (2,2) circle[radius=4pt];
	\draw (2,2) node{\Large $\times$};
	\draw[color=blue] (-0.25,0.25) node{\tiny\textbf{\texttt{2}}};
	\draw[color=blue] (-0.25,-0.25) node{\tiny\textbf{\texttt{1}}};
	\draw[color=blue] (0.35,0.15) node{\tiny\textbf{\texttt{3}}};
	\draw[color=blue] (0.35,0.52) node{\tiny\textbf{\texttt{0}}};
	\draw[color=blue] (0.35,-0.25) node{\tiny\textbf{\texttt{2}}};
	\draw[color=blue] (1.65,0.15) node{\tiny\textbf{\texttt{2}}};
	\draw[color=blue] (1.65,0.52) node{\tiny\textbf{\texttt{3}}};
	\draw[color=blue] (1.65,-0.25) node{\tiny\textbf{\texttt{4}}};
	\draw[color=blue] (2.35,0.15) node{\tiny\textbf{\texttt{0}}};
	\draw[color=blue] (2.15,0.4) node{\tiny\textbf{\texttt{2}}};
	\draw (3.9,0) node[above]{$\mr{pt}$};
\end{tikzpicture}
\caption{An example of a stable graph in $\cG^{[n]}_{8,1}$, with cohomological decorations at the vertices omitted. For each vertex the genus is marked in red; for each half-edge the power of the corresponding $\psi$-class is marked in blue, so that their sum at a vertex $v$ is equal to $\dim \Mbar_{\mathrm{g}(v), |\HH(v)|}$. On the leaves, crossed vertices indicate dilaton markings; the remaining leaves carry markings by  cohomology classes (the point class for the only marked leaf in the picture).}
\label{fig:stgrgen}
\end{figure}
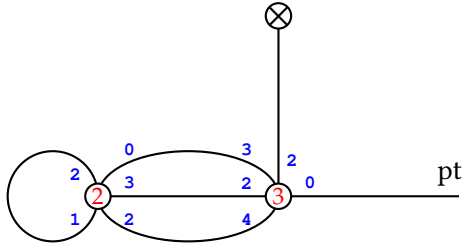

\subsubsection{Givental--Teleman reconstruction}
\label{sec:givtel}

Let $\Gamma\in \cG^{[2]}_{g,n}$ be a $2$-decorated stable graph. 
%
We will assign to each leaf, edge, and vertex in $\Gamma$ a set of $Q_{\Tmax}$-valued functions on $\mathsf{B}_{\cQ_{\rm cf}}$, as follows.
Given the canonical $R$-matrix \eqref{eq:canR}, we define for $i,j\in\{0,1,2\}$, $n,m\in\bbZ_{\geq 0}$ the dilaton and marking factors
%
\beq
\mathsf{D}_{i,n}(z) \coloneqq [z^{n-1}] \l(\Psi_{0i}-\sum_j \mathsf{R}(-z)_{0i}\r)\,, \quad \mathsf{L}_{i,n}^{[a]} \coloneqq [z^n] \mathsf{R}_{ia} (-1)^n\,,
\label{eq:leafterms}
\eeq 
and the edge factor
\beq 
\mathsf{E}^{ij}_{n,m} \coloneqq [z^{n} w^m]\l( \frac{\delta_{ij}-\sum_{ab}\mathsf{R}_{ai}(-z) \eta^{ab} \mathsf{R}_{bj}(-w)}{z+w}\r)\,.
\label{eq:edgeterm}
\eeq
For each leaf $h \in \DD_\Gamma \sqcup \LL_\Gamma $, we define
\[
\mathrm{Cont}(h) \coloneqq 
\begin{cases}
    \mathsf{D}_{\mathsf{q}(\mathsf{v}(h)),\mathsf{a}(h)}\,, & h \in \DD_\Gamma\,,\\[0.3em]
    \mathsf{L}_{\mathsf{q}(\mathsf{v}(h)),\mathsf{a}(h)}^{[\mathsf{q}(h)]}\,, & h \in \LL_\Gamma\,.
\end{cases}
\] 
For each edge $e=\{h_1,h_2\} \in \cE_\Gamma$, we define
\[
\mathrm{Cont}(e) \coloneqq
\mathsf{E}^{\mathsf{q}(\mathsf{v}(h_1)), \mathsf{q}(\mathsf{v}(h_2)}_{\mathsf{a}(h_1),\mathsf{a}(h_2)}\,.
\] 
For each vertex $v \in \cV_\Gamma$, fixing a labelling $\{h_1, \dots, h_{|\HH(v)|}\}$ of the half-edges incident to $v$, 
we define 
\[
\mathrm{Cont}(v) = (3 W_i)^{\frac{2\mathsf{g}(v)-2+\mathsf{|H(v)|}}{2}} 
\int_{[\Mbar_{\mathsf{g}(v), |\HH(v)|}]}\prod_{k=1}^{|\HH(v)|} 
\psi_k^{\mathsf{a}(h_k)}\,.
\]
%
%
Recalling that $\Delta_i^{-1}=\epsilon_1 \epsilon_2 (3 W_i)$ from \eqref{eq:delta1}, the Givental--Teleman reconstruction formula then states that, for $2g+n-2>0$,
\beq 
\bra \! \bra H^{a_1}, \dots, H^{a_n} \ket \! \ket_{g,n}^Z =
\sum_{\Gamma \in \cG^{[2]}_{g,n}} \frac{(\epsilon_1 \epsilon_2)^{g-1}}{|\mathrm{Aut}(\Gamma)|} 
\prod_{e\in \cE_\Gamma}\mathrm{Cont}(e)
\prod_{h\in \LL_\Gamma 
\sqcup \DD_\Gamma
}\mathrm{Cont}(h) \prod_{v\in \cV_\Gamma} 
\mathrm{Cont}(v)\,.
\label{eq:gtH}
\eeq 
Note that, by the degree condition
\[ 
\sum_{k} \mathsf{a}(h_k)  = 3 \mathsf{g}(v)-3+|\HH(v)|\,,
\] 
the sum in \eqref{eq:gtH} is finite.

\begin{prop}
Let $\mgp{T}_{\rm R} \subset \Tmax$ be a rational sub-torus. Then,
\[
\mathsf{D}_{i,n} \in G_{\mgp{T}_{\rm R}}\,, \quad \mathsf{L}_{i,n}^{[1]} \in \frac{1}{A}G_{\mgp{T}_{\rm R}}\,, \quad \mathsf{E}^{i,j}_{n,m} \in G_{\mgp{T}_{\rm R}}  + F\,G_{\mgp{T}_{\rm R}}\,.
\]
\label{prop:fingenTLE}
\end{prop}
\begin{proof}
The first two inclusions are immediate as, from  \eqref{eq:asympR} and \eqref{eq:leafterms}, \[
\mathsf{D}_{i,n} = (-1)^n \mathsf{R}_{i,n}^{[0]}\,, \quad \mathsf{L}_{i,n}^{[1]} =  \frac{(-1)^n Y_i \mathsf{R}_{i,n}^{[1]}}{A}\,,\]
and therefore $\mathsf{D}_{i,n}\in G_{\mgp{T}_{\rm R}}$ and $\mathsf{L}_{i,n}^{[1]} \in A^{-1}  G_{\mgp{T}_{\rm R}}$ by \cref{prop:fingenR}.
%
%
Expanding out the summation in the numerator of \eqref{eq:edgeterm}, we find
\begin{equation}
\begin{split}
\sum_{a,b} \mathsf{R}_{ai}(-z) \eta^{ab} \mathsf{R}_{bj}(-w) & =
\frac{\epsilon _1 \epsilon _2 e_3(\alpha)}{9} \mathsf{R}_{0i}(-z) \mathsf{R}_{0j}(-w) -3 \epsilon _1 \epsilon _2 \big( \mathsf{R}_{1i}(-z) \mathsf{R}_{2j}(-w)+ \mathsf{R}_{2i}(-z) \mathsf{R}_{1j}(-w)\big) 
\\ 
&- 3\epsilon _1 \epsilon _2 e_1(\alpha) \mathsf{R}_{1i}(-z) \mathsf{R}_{1j}(-w) 
\,.
\end{split}
\label{eq:preedge}
\end{equation}
%
Recall that, by \cref{prop:fingenR}, we have
\begin{align*}
[z^n w^m] \mathsf{R}_{0i}(-z) \mathsf{R}_{0j}(-w) \in\, & G_{\mgp{T}_{\rm R}}\,, \nn \\
[z^n w^m] \mathsf{R}_{1i}(-z) \mathsf{R}_{1j}(-w) \in\, & A^{-2} G_{\mgp{T}_{\rm R}}\,,\nn \\
[z^n w^m] \mathsf{R}_{2i}(-z) \mathsf{R}_{1j}(-w) \in\,  &  G_{\mgp{T}_{\rm R}}+A^{-2} G_{\mgp{T}_{\rm R}}+F\, G_{\mgp{T}_{\rm R}}\,. 
\end{align*}
Let us analyse the polar part in $A$ of \eqref{eq:preedge}. From   \eqref{eq:DErel}, \eqref{eq:asympS}, and \eqref{eq:RHb}, $\mathsf{R}_{2i}(z)$ has a linear term in $A^{-1}$ equal to $e_1(\alpha)\mathsf{R}_{1i}(z)/2$, which arises from the second summand in the numerator of \eqref{eq:RHb},
\[
\frac{e_1(\alpha) E(y) \mathsf{R}_{1i}(z)}{D(y)} =
e_1(\alpha)\mathsf{R}_{1i}(z) \l( \frac{1}{2}- \frac{A^2}{2}\r)\,.
\]
We then have that
\[
[z^n w^m A^{-2}] \mathsf{R}_{2i}(-z) \mathsf{R}_{1j}(-w) = [z^n w^m A^{-2}]\frac{e_1(\alpha)}{2}\mathsf{R}_{1i}(-z) \mathsf{R}_{1j}(-w)\,,
\]
which implies that the terms linear in $A^{-2}$ cancel between the first and the second line of \eqref{eq:preedge}. Therefore,
\[
[z^n w^m]\sum_{ab}\mathsf{R}_{ai}(-z) \eta^{ab} \mathsf{R}_{bj}(-w) \in G_{\mgp{T}_{\rm R}} + F\,G_{\mgp{T}_{\rm R}}\,,
\]
from which the claim follows by \eqref{eq:edgeterm}.
\end{proof}

\subsection{The modular anomaly equation}

\cref{prop:fingenR,prop:fingenTLE} give all the ingredients for the proof of the holomorphic (or modular) anomaly equation. In the following, for any torus $\mgp{T}\subseteq \Tmax$, we will write $\GW_g^{\mgp{T}}$ as a short-hand for \[\GW_g^{\mgp{T}} \coloneqq\refGWg{g}{K_{\bbP^2}\times \Aaff{2}}{\mgp{T}}\,.\]
\begin{thm}
Let $\mgp{T}_{\rm R} \subset \Tmax$ be a rational sub-torus of $\Tmax$. Then, for $g \geq 2$:
\ben
\item  the fixed-genus generating series $\GW_g^{\mgp{T}_{\rm R}}$ satisfy the finite generation property
\[ 
(\cQ \de_\cQ)^n \GW_g^{\mgp{T}_{\rm R}} \in A^{-n} G_{\mgp{T}_{\rm R}}[F]\,;
\] 
\item the modular anomaly equation holds,
\beq
\frac{X}{A^2}
\frac{\de \GW^{\mgp{T}_{\rm R}}_g}{\de F}=\frac{\epsilon_1 \epsilon_2}{2} \l[(\cQ \partial_\cQ)^2 \mathsf{GW}^{\mgp{T}_{\rm R}}_{g-1} + \cQ^2 \sum_{g'=1}^{g-1}   \de_\cQ \mathsf{GW}^{\mgp{T}_{\rm R}}_{g'} \, \partial_\cQ \mathsf{GW}^{\mgp{T}_{\rm R}}_{g-g'}\r]\,,
\label{eq:HAEF}
\eeq
as an identity in $A^{-2}G_{\mgp{T}_{\rm R}}[F]$;
\item $\deg_F \GW^{\mgp{T}_{\rm R}}_g \leq 3g-3$.
\een
\label{thm:haegen}
\end{thm}
\begin{proof}

We will break down the proof according to the three parts of the statement.\\

\ben
\item From the Givental--Teleman formula, \eqref{eq:gtH}, combined with \cref{prop:fingenTLE}, we have that 
\[
\bra \! \bra H, \dots, H \ket \! \ket_{g,n}^Z \in A^{-n} G_{\mgp{T}_{\rm R}}[F]\,.
\]
The Divisor Axiom then imposes
\[ 
(\cQ \de_\cQ)^n \GW_g^{\mgp{T}_{\rm R}} = \bra \! \bra H, \dots, H \ket \! \ket_{g,n}^Z\big|_{t_a=I_0^{[a]}(y)}
\]
as an identity in $Q_{\mgp{T}_{\rm R}}\llbracket y\rrbracket$, where the right hand side is now viewed as the specialisation of an element in $A^{-n}G_{\mgp{T}_{\rm R}}[F] \to Q_{\mgp{T}_{\rm R}}\llbracket\cQ\rrbracket$ via $A \to A(y)$, $F \to F(y)$. By the algebraic independence of $A$ and $F$ \cite[Lem.~4.2]{BFGW21:HAE}, the left hand side can be lifted to a unique representative in $ A^{-n}\,G_{\mgp{T}_{\rm R}}[F]$.\\

\item Let $\Gamma \in \cG^{[2]}_{g,n}$, $e\in \cE_\Gamma$. Following \cite{LP18:HAE}, we adopt the following conventions:
\bit
\item
if $\Gamma \setminus e$ is connected, we write $\Gamma^0_e \in \cG^{[2]}_{g-1,2}$ for the resulting decorated stable graph;
\item if it is disconnected, we write $\Gamma^i_e \in \cG^{[2]}_{g_i-1,1}$, $i=1,2$, for the two resulting decorated stable graphs, where $g_1+g_2=g$.
\eit
\medskip
\noindent By the Leibnitz rule applied to \eqref{eq:gtH}, we have 
\[
\frac{X}{A^2}
\frac{\de \GW_g^{\mgp{T}_{\rm R}}}{\de F} = \sum_{\Gamma \in \cG^{[2]}_{g,0}}
\frac{(\epsilon_1 \epsilon_2)^{g-1}}{|\mathrm{Aut}(\Gamma)|}  \sum_{e\in \cE_\Gamma} \frac{X}{A^2} \frac{\de \mathrm{Cont}(e)}{\partial F}  \prod_{f \neq e\in \cE_\Gamma}\mathrm{Cont}(f) \prod_{h\in \LL_\Gamma \sqcup \DD_\Gamma}\mathrm{Cont}(h) \prod_{v\in \cV_\Gamma}  \mathrm{Cont}(v)\,. 
\]
Now, from \cref{prop:fingenR}, we have that 
\beq 
\de_F \mathsf{R}_{i,n}^{[2]} = -\frac{\mathsf{R}_{i,n}^{[1]}}{6 Y_i}\,,
\label{eq:dFR}
\eeq 
and therefore, from \eqref{eq:edgeterm} and denoting $i_k=\mathsf{q}(\mathsf{v}(h_k))$, $n_k=\mathsf{a}(h_k)$,
\beq 
\de_F \mathrm{Cont}(e) =
\de_F \mathsf{E}^{i_1, i_2}_{n_1+1, n_2+1}=\frac{(-1)^{n_1+n_2} Y_i Y_j \mathsf{R}^{[1]}_{i_1,n_1} \mathsf{R}^{[1]}_{i_2,n_2}
}{X}\,. 
\label{eq:dFedge}
\eeq 
From \eqref{eq:asympR} and \eqref{eq:dFedge}, acting with $\de_F$ converts edge terms corresponding to $e \in \Gamma$ to insertions of the hyperplane class at the corresponding new markings in $\Gamma^0_e$ (if $\Gamma \setminus e$ is connected) or $\Gamma^i_e$, $i=1,2$ (if $\Gamma \setminus e$ is disconnected). As these are not ordered, when we sum over labellings we need to compensate for that by a symmetry factor of 2. Explicitly, substituting \eqref{eq:dFedge} into \eqref{eq:gtH}, decomposing the resulting sum over graphs according to the connectedness of $\Gamma \setminus e$, and using \eqref{eq:gtH} for each of these separately, we have
\bea 
& &\frac{X}{A^2}
 \frac{\de \GW_g^{\mgp{T}_{\rm R}}}{\de F} = \frac{X}{2A^2} \Bigg[\sum_{\Gamma \in \cG^{[2]}_{g-1,2}}
\frac{(\epsilon_1 \epsilon_2)^{g-1}}{|\mathrm{Aut}(\Gamma)|}  \sum_{e\in \cE_\Gamma} \frac{\de \mathrm{Cont}(e)}{\partial F} \prod_{f \neq e\in \cE_\Gamma}\mathrm{Cont}(f) \prod_{h\in \LL_\Gamma \sqcup \DD_\Gamma}\mathrm{Cont}(h) \prod_{v\in \cV_\Gamma} 
\mathrm{Cont}(v) \nn \\
& & +   \sum_{g'=1}^{g-1} \sum_{\Gamma \in \cG^{[2]}_{g',1}} \sum_{\Gamma \in \cG^{[2]}_{g-g',1}}
\frac{(\epsilon_1 \epsilon_2)^{g-1}}{|\mathrm{Aut}(\Gamma)|}  \sum_{e\in \cE_\Gamma} \frac{\de \mathrm{Cont}(e)}{\partial F} \prod_{f \neq e\in \cE_\Gamma}\mathrm{Cont}(f) \prod_{h\in \LL_\Gamma \sqcup \DD_\Gamma}\mathrm{Cont}(h) \prod_{v\in \cV_\Gamma} 
\mathrm{Cont}(v)\Bigg]\nn \\
& & =
\frac{\epsilon_1 \epsilon_2}{2} \l(\bra \! \bra H,H \ket \! \ket^Z_{g-1,2} + \sum_{g'=1}^{g-1} \bra \! \bra H \ket \! \ket^Z_{g',1}\bra \! \bra H \ket \! \ket^Z_{g-g',1}\r)\bigg|_{t_a=I_0^{[a]}} \nn \\
& & =
\frac{\epsilon_1 \epsilon_2}{2} \l[(\cQ \partial_\cQ)^2 \mathsf{GW}^{\mgp{T}_{\rm R}}_{g-1} + \cQ^2 \sum_{g'=1}^{g-1}   \de_\cQ \mathsf{GW}^{\mgp{T}_{\rm R}}_{g'}  \, \partial_\cQ \mathsf{GW}^{\mgp{T}_{\rm R}}_{g-g'}\r]\,,
\eea 
where in the last line we have used the Divisor Axiom. By the first part of the theorem, the identity holds as an equality in $A^{-2}G_{\mgp{T}_{\rm R}}[F]$. 
\item From \eqref{eq:gtH} and \eqref{eq:dFR}, since $\GW_g^{\mgp{T}_{\rm R}}$ is polynomial in the edge terms, and these are linear in $F$, the degree of $\GW_g^{\mgp{T}_{\rm R}}$ in $F$ is bounded by the number of edges in the graph. We have $\max_{\Gamma \in \cG^{[2]}_{g,0}} |\cE_\Gamma|=3g-3$, corresponding to the degenerate stratum of $\Mbar_{g,0}$ of rational curves with $3g-3$ nodes.\qedhere
\een
\end{proof}

\subsubsection{The fibrewise torus}
Consider now the specialisation to the fibrewise torus $\mgp{T}_{\rm F}$. The diagonal restriction in \eqref{eq:diagweights} obviously violates the (open) Condition~(ii) in \cref{def:adm}, so $\mgp{T}_{\rm F}$ is not a rational sub-torus of $\Tmax$. It is immediate to check, however, that it satisfies the (closed) Condition~(i), hence the fibrewise action can be regarded as a limiting version a rational torus $\mgp{T}_{\rm R}$: one such example is given by taking the $\delta\to 0$ limit of \eqref{eq:adelta}.

Restricting the equivariant generating series $\refGWg{g}{K_{\bbP^2}\times \Aaff{2}}{\Tmax}$ to $\mgp{T}_{\rm F}$,
\[
\GW_g^{\mgp{T}_{\rm F}} \coloneqq \GW_g^{\Tmax}\Big|_{\alpha_i=\frac{\epsilon_1+\epsilon_2}{3}}\,,
\]
we define the refined generating series
\beq 
\mathscr{F} \coloneqq \sum_{k \geq 0,g\geq 0} (\epsilon_1+\epsilon_2)^{2k} (\epsilon_1 \epsilon_2)^{g-1} \mathscr{F}_{k,g} \coloneqq \sum_{g \geq 0}   \GW_g^{\mgp{T}_{\rm F}} \in \widehat{R}^{\,\mr{loc}}_{\mgp{T}_{\rm F}}\llbracket\cQ\rrbracket  \,. 
\label{eq:Fref}
\eeq 

By \cref{lem: shape of refGW series}, an expansion of the form \eqref{eq:Fref}, involving only even powers of $2\epsilon_+ = (\epsilon_1+\epsilon_2)$, is guaranteed to exist.

\begin{thm}
For $k+g \geq 2$, the refined generating series are weight $0$, depth $3g+3k-3$ quasi-modular functions of $\Gamma_1(3)$,
\[ 
\mathscr{F}_{k,g} \in \bbQ[X^{\pm 1}][ F]\subset \mathfrak{QMod}^{[0,3g+3k-3]}(\Gamma_1(3))\,.
\] 
%
Moreover, the refined modular anomaly equation holds,
\[
\frac{X}{A^2}
\frac{\de \mathscr{F}_{k,g}}{\de F}=\frac{1}{2} (\cQ \partial_\cQ)^2 \mathscr{F}_{k,g-1} + \frac{\cQ^2}{2} \sum_{\substack{0\leq k' \leq k\,, 0 \leq g' \leq g \\(0,0)\neq (k',g')\neq (k,g) }}   \de_\cQ \mathscr{F}_{k',g'} \, \partial_\cQ \mathscr{F}_{k-k',g-g'}\,,
\]
as an identity in $A^{-2}\,\bbQ[X^{\pm 1},F]$.
\label{thm:haeref}
\end{thm}

\begin{proof}

%

By \cref{lem: shape of refGW series}, the Gromov--Witten invariant
\beq 
\epsilon_1\epsilon_2 \, \refGWgb{g}{K_{\bbP^2} \times \Aaff{2}}{\mgp{T}_{\rm F}}{d[H]}
\label{eq:gwref}
\eeq 
is a homogenous polynomial in $(\epsilon_1,\epsilon_2)$ of degree $2g$, and therefore
\beq 
\GW^{\mgp{T}_{\rm F}}_g = (\epsilon_1 \epsilon_2)^{g-1}
\sum_{n=1-g}^{g-1}
\GW^{\mgp{T}_{\rm F}}_{g,n} (\epsilon_1/\epsilon_2)^n\,, \quad \GW^{\mgp{T}_{\rm F}}_{g,n}\in \bbQ\llbracket\cQ\rrbracket\,.
\label{eq:GWLpol}
\eeq 
%
By properness of the $\mgp{T}_{\rm F}$-fixed locus $\Mbar^{\mgp{T}_{\rm F}}_{g,0}(Z,
    d [H])$, the equivariant GW generating function  for the maximal torus, $\GW^{\Tmax}_g$,
is a regular function of the weights along the diagonal \eqref{eq:diagweights}; and as argued below \eqref{eq:diagweights}, the fibrewise action is in the closure of the set of rational weights in \cref{def:adm}. Therefore, by continuity, \cref{thm:haegen}(i) extends to the fibrewise torus:
\beq 
\GW^{\mgp{T}_{\rm F}}_{g} \in \,G_{\mgp{T}_{\rm F}}[F]\,.
\label{eq:fingenTf}
\eeq 
%
The $\mathrm{PGL}_3(\bbC)$ automorphism group of the target further induces an $S_3$-e\-qui\-va\-riance of the Gromov--Witten invariants \eqref{eq:gwprim}, where $\sigma \in S_3$ acts by simultaneous permutation of the weights at the fixed points and the canonical data associated to the respective idempotents, \[\{\alpha_i \to \alpha_{\sigma(i)}, 
Y_i \to Y_{\sigma(i)}, W_i \to W_{\sigma(i)}\}\,.\] When $\alpha_0=\alpha_1=\alpha_2$, by \eqref{eq:Wadm}, this in particular constrains $\GW^{\mgp{T}_{\rm F}}_g$ to be a  symmetric rational function of $(Y_0, Y_1, Y_2)$
\[
\GW^{\mgp{T}_{\rm F}}_g \in G_{\mgp{T}_{\rm F}}^{S_3}[F]\,.
\]
Furthermore, as a power series in $\cQ$,
$\GW_g^{\mgp{T}_{\rm F}}$ is by definition regular at $\cQ=y=0$; and by \eqref{eq:GWLpol}, the coefficients of the $\cQ$-expansion are homogeneous $\bbQ$-Laurent polynomials in $(\epsilon_1, \epsilon_2)$.
We then conclude, from \eqref{eq:Greg} and \eqref{eq:GWLpol}, that
\bea
\GW_g^{\mgp{T}_{\rm F}} & \in & \l(G_{\mgp{T}_{\rm F}}^{\rm reg}[F]\r)^{S_3} = (\epsilon_1 \epsilon_2)^{g-1} \bbQ\l[(\epsilon_1 / \epsilon_2)^{\pm 1}\r][X^{\pm 1}, F]\,, \nn \\ 
\GW^{\mgp{T}_{\rm F}}_{g,n} & \in & \bbQ[X^{\pm 1},F]\,,\nn
\eea 
with $\deg_F \GW_{g,n}^{\mgp{T}_{\rm F}}\leq \deg_F \GW_g^{\mgp{T}_{\rm F}}\leq 3g-3$. In terms of the refined generating series \eqref{eq:Fref}, we have
\[ 
\GW_{g,n}^{\mgp{T}_{\rm F}}= \sum_{g'=0}^g \binom{2(g-g')}{n+g-g'} \mathscr{F}_{g-g',g'}\,,
\]
and therefore
\[ 
\mathscr{F}_{k,g} \in \bbQ[X^{\pm 1},F]
\] 
with $\deg_F \mathscr{F}_{k,g}\leq 3(g+k-1)$. Finally, note that for a rational torus~$\mgp{T}_{\rm R}$, the modular anomaly equation \eqref{eq:HAEF} holds as a closed condition in the weights $(\alpha_0, \alpha_1, \alpha_2)$, which is furthermore regular on the diagonal \eqref{eq:diagweights}. Arguing as for \eqref{eq:fingenTf}, it will therefore extend to the specialisation to the fibrewise torus.
Expanding  \eqref{eq:HAEF} in $\epsilon_1 \epsilon_2$ and $\epsilon_1/\epsilon_2$ (respectively, $\epsilon_1 \epsilon_2$ and $\epsilon_1+\epsilon_2$), we obtain
\[
\frac{X}{A^2}
\frac{\de \GW^{\mgp{T}_{\rm F}}_{g,n}}{\de F}=\frac{1}{2} (\cQ \partial_\cQ)^2 \GW^{\mgp{T}_{\rm F}}_{g-1,n} + \frac{\cQ^2}{2} \sum_{\substack{0\leq n' \leq n\,, 0 \leq g' \leq g \\(0,0)\neq (g',n')\neq (g,n) }}   \de_\cQ \GW^{\mgp{T}_{\rm F}}_{g',n'} \, \partial_\cQ \GW^{\mgp{T}_{\rm F}}_{g-g',n-n'}\,,
\label{eq:HAEFdiag}
\]
and
\[
\frac{X}{A^2}
\frac{\de \mathscr{F}_{k,g}}{\de F}=\frac{1}{2} (\cQ \partial_\cQ)^2 \mathscr{F}_{k,g-1} + \frac{\cQ^2}{2} \sum_{(0,0)\neq (k',g')\neq (k,g) }   \de_\cQ \mathscr{F}_{k',g'} \, \partial_\cQ \mathscr{F}_{k-k',g-g'}\,,
\label{eq:HAEFref}
\]
as identities in $A^{-2}\bbQ[X^{\pm 1},F]$.
\end{proof}

\begin{rmk}
By taking $\epsilon_1=-\epsilon_2$, 
\cref{prop: unref limit,thm:haegen} give a quasi-modularity statement for the unrefined generating series,
\[
\mathscr{F}_{0,g} \in \bbQ[X^{\pm 1}][F] \subset \mathfrak{QMod}^{[0,3g-3]}(\Gamma_1(3))\,.
\]
This slightly strengthens \cite[Thm.~1(i)]{LP18:HAE}, where it was proved that $\mathscr{F}_{0,g} \in \bbQ[X^{\pm 1/3}][F]$. That this stronger form of the result of \cite{LP18:HAE} should hold was anticipated in \cite{HKR:HAE, BFGW21:HAE}; it also follows from the quasi-modularity theorem for local $\bbP^2$ as presented in \cite{Coates:2018hms}.
\end{rmk}

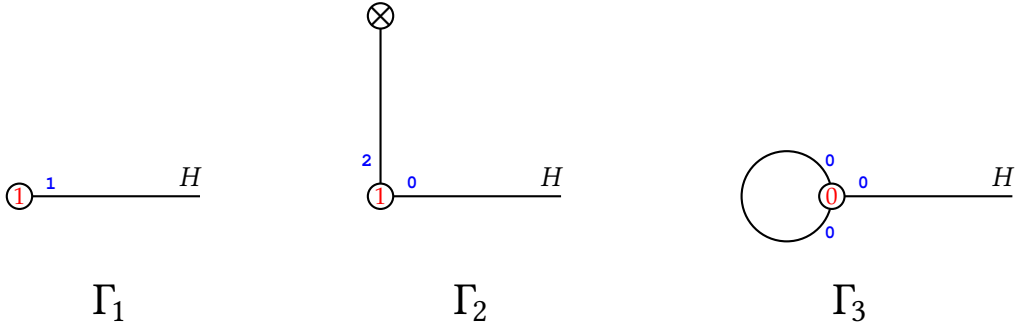
\begin{figure}[t]
\centering
\begin{tikzpicture}[baseline=(current  bounding  box.center),scale=1.2]
    \draw[thick] (0,0) -- (2,0);
    \draw[fill=white,thick] (0,0) circle[radius=4pt];
    \draw[color=red] (0,0) node{\footnotesize 1};
    \draw[color=blue] (0.35,0.15) node{\tiny\textbf{\texttt{1}}};
    \draw (1.9,0) node[above]{$H$};
    \draw (1,-1.5) node[above]{\LARGE $\Gamma_1$};
    
    \draw[thick] (4,2) -- (4,0) -- (6,0);
    \draw[fill=white,thick] (4,2) circle[radius=4pt];
    \draw (4,2) node{\Large $\times$};
    \draw[color=blue] (3.85,0.4) node{\tiny\textbf{\texttt{2}}};
    \draw[fill=white,thick] (4,0) circle[radius=4pt];
    \draw[color=red] (4,0) node{\footnotesize 1};
    \draw[color=blue] (4.35,0.15) node{\tiny\textbf{\texttt{0}}};
    \draw (5.9,0) node[above]{$H$};
    \draw (5,-1.5) node[above]{\LARGE $\Gamma_2$};
    
    \draw[thick] (8.5,0) circle[radius=0.5];
    \draw[thick] (9,0) -- (11,0);
    \draw[fill=white,thick] (9,0) circle[radius=4pt];
    \draw[color=red] (9,0) node{\footnotesize 0};
    \draw[color=blue] (9.35,0.15) node{\tiny\textbf{\texttt{0}}};
    \draw[color=blue] (8.98,0.4) node{\tiny\textbf{\texttt{0}}};
    \draw[color=blue] (8.98,-0.4) node{\tiny\textbf{\texttt{0}}};
    \draw (10.9,0) node[above]{$H$};
    \draw (9.2,-1.5) node[above]{\LARGE $\Gamma_3$};
\end{tikzpicture}
\caption{The three stable graphs in $\cG^{[2]}_{1,1}$, with cohomological decorations omitted and with powers of $\psi$ markings indicated in blue as in \cref{fig:stgrgen}, that give a non-zero contribution to $\cQ \de_\cQ \GW^{\mgp{T}_{\rm F}}_1$.}
\label{fig:stgr11}
\end{figure}

\begin{example}[Genus one]
\label{exmpl:F1}
The genus one Gromov--Witten potential with one line insertion is computed by \eqref{eq:gtH} as a Givental sum over the three weighted stable graphs in \cref{fig:stgr11}. We have
\[
\cQ \partial_\cQ \GW^{\mgp{T}_{\rm F}}_1 = \mathrm{Cont}_{\Gamma_1}+\mathrm{Cont}_{\Gamma_2}+\mathrm{Cont}_{\Gamma_3}\,,
\]
with 
\begin{alignat*}{2}
\mathrm{Cont}_{\Gamma_1} &= \frac{1}{24} \sum_{i=0}^2 \mathsf{L}_{i,1}^{[1]} \Delta_i^{1/2} &&= \frac{\epsilon _1^2+\epsilon _1 \epsilon _2+\epsilon _2^2}{288 \epsilon _1 \epsilon _2 } \frac{X}{A}\,,\\
\mathrm{Cont}_{\Gamma_2} &= \frac{1}{24} \sum_{i=0}^2 \mathsf{D}_{i,2}\mathsf{L}_{i,0}^{[1]} \Delta_i &&= -\frac{\epsilon _1^2-11 \epsilon _1 \epsilon _2+\epsilon _2^2}{288  \epsilon _1 \epsilon _2} \frac{X}{A}\,,\\
\mathrm{Cont}_{\Gamma_3} &= \frac{1}{2} \sum_{i=0}^2 \mathsf{L}_{i,0}^{[1]} \mathsf{E}^{i,i}_{0,0} \Delta_i^{1/2} &&= 
-\frac{F}{2 A}-\frac{X \left(\epsilon _1^2+\epsilon _2 \epsilon _1+\epsilon _2^2\right)+4 \epsilon _1 \epsilon _2}{24 A \epsilon _1 \epsilon _2}\,,
\end{alignat*}
so that
\[
\cQ \partial_\cQ \mathscr{F}_{0,1} = -\frac{F}{2 A}+ \frac{X-2}{12 A}\,, \quad
\cQ \partial_\cQ \mathscr{F}_{1,0} =  -\frac{X}{24 A}\,.
\]
\end{example}

    Up to a constant (unstable) term, the unrefined genus-one generating series, $\mathscr{F}_{0,1}$, coincides with the expression calculated in \cite{LP18:HAE}. Likewise, the Nekrasov--Shatashvili 
    generating series in genus one, $\mathscr{F}_{1,0}$, in turn recovers \cite[Thm.~1.4]{BFGW21:HAE}.

\begin{example}[Genus two]

The $g=2$ refined Gromov--Witten potential is computed by \eqref{eq:gtH} as a Givental sum over the seventeen decorated weighted stable graphs in $\cG^{[2]}_{2,0}$. To compute their individual contributions, we restrict first to a rational sub-torus $\mgp{T}_{\rm R}$  whose fixed locus is zero-dimensional; one example is the choice in \eqref{eq:adelta} with $\delta \neq 0$. We then look for a representative in the quotient ring \eqref{eq:GTquot} of an asymptotic solution of \eqref{eq:pf} in the form \eqref{eq:asympRrho},
\[ \mathsf{R}_i(\mathbf{1})(y,z)   \simeq  a(z)\re^{\mathsf{u}_i/z} \sum_{n \geq 0} \rho_{i,n}^{[0]}(y) z^n\,,\] with the ansatz
\[
\rho_{i,n}^{[0]}(y)=\frac{\sum_{l=0}^2\sum_{m=0}^{2n} \sigma_{nlm} Y_i^l y^m}{ (1+27 y)W_i^{\frac{n+2}{2}}  (27(1+27y)e_3(\alpha)^2-e_2(\alpha)^3)^n}\,,
\]
where $\sigma_{nlm} \in Q_{\mgp{T}_{\rm R}}$,  $a(z) \in Q_{\mgp{T}_{\rm R}}\llbracket z\rrbracket$. A direct inspection for the first few values of $n$ reveals that \eqref{eq:pf} becomes a full-rank linear system for $\sigma_{nlm}$, uniquely determining $\rho_{i,n}^{[0]}$. Imposing the symplecticity condition and the quantum Riemann--Roch normalisation \eqref{eq:normNLR} fixes $a(z)$ uniquely, and the restriction  of the $R$-matrix to the fibrewise torus ${\mgp{T}_{\rm R}} \to \mgp{T}_{\rm F}$ is obtained by taking the limit $\delta \to 0$ in \eqref{eq:adelta}, which exists for $y \neq \{0,-1/27\}$. 

With a closed-form expression for the $R$-matrix at hand, the resulting graph sum can be evaluated first as an element of
$Q_{\mgp{T}_{\rm F}}[\{Y_i^{\pm}, W_i^{\pm 1/2}\}_{i=0}^2][F]$ and then quotiented down to the invariant ring $G_{\mgp{T}_{\rm F}}^{S_3}[F]$
by polynomial division w.r.t.~a Groebner basis for the ideal $\bra \{\{\LL(Y_i), \cM(Y_i,W_i)\}_{i=0}^2, \cN(Y)\} \ket$. The full calculation is available upon request to the interested reader, and we will only reproduce the final result here: 
\begin{equation}
\begin{split}
\mathscr{F}_{0,2} &= 
\frac{5 F^3}{8 X}+ \left(\frac{5}{8 X}-\frac{1}{2}\right)F^2+\frac{1}{96}  \left(13 X+\frac{20}{X}-32\right)F-\frac{27 X^3-98 X^2+121 X-50}{2160 X}
\,, \\
\mathscr{F}_{1,1} &= \frac{F^2}{16}+\frac{1}{96} (8-7 X)F+\frac{84 X^2-196 X+121}{5760}
\,,  \\
\mathscr{F}_{2,0} &= \frac{X}{384}F-\frac{21 X^2-29 X+8}{5760}
\,.
\end{split}
\label{eq:F2ref}
\end{equation}
\label{exmpl:F2}
\end{example}

The unrefined generating series $\mathscr{F}_{0,2}$ coincides with the expression found in \cite{LP18:HAE}, up to the coefficient $-121/2160$ of the term of $\cO(F^0 X^0)$. This discrepancy is expected, as the constant term in $\mathscr{F}_{0,2}$ corresponds to the degree zero term contribution for which the corresponding moduli space is non-compact, and different torus actions (as is the case here compared to the choice of \cite{LP18:HAE}) will typically lead to different numerical invariants. The Nekrasov--Shatashvili generating series $\mathscr{F}_{2,0}$ matches the expression found in \cite[Thm.~1.4]{BFGW21:HAE} from the higher genus relative/local correspondence applied to the Gromov--Witten theory of $\bbP^2$ relative to a smooth cubic.

\section{Refined local \texorpdfstring{$\bbP^2$}{P2} beyond large volume}
\label{sec: ref KP2 beyond LV}

We will study here the behaviour of the refined Gromov--Witten generating series 
at the boundary points $y \in \bbP(3,1) \setminus \cC_1(3)$ of the $\Gamma_1(3)$-modular curve other than $y=0$ (the large volume limit): these are $y=\infty$ (the orbifold point) and $y=-1/27$ (the conifold point). For that, we shall need analytical control on the behaviour of the coefficients of the $R$-matrix away from
$y=0$: 
this will be provided by an integral representation of the flat coordinates of the Dubrovin connection 
as steepest-descent integrals with phase function 
given by Givental's equivariant mirror \cite{MR1653024}. We recall some generalities, and fix notation, about steepest-descent asymptotics of $n$-dimensional oscillating integrals below.\\

Let $\phi(q)$, $\gamma(q):\bbC^n \to \bbC$ be holomorphic functions on $\bbC^n$ and
let $q^{\rm cr}$ be a critical point of $\phi$, $\de_a \phi(q^{\rm cr})=0$, with non-degenerate Hessian $\HH_{ab} = \de^2_{ab} \phi(q^{\rm cr})$. Consider the oscillatory integral
\[
\cI(z) \coloneqq \int_{\Gamma} \re^{-\frac{\phi(q)}{z}} \gamma(q) \dd^n q
\]
where $z>0$ and $\Gamma\subset \bbC^n$ be the union of the negative gradient flow lines of $\mathrm{Re}(\phi)$ emerging from $q^{\rm cr}$: by the Cauchy--Riemann equations, this is an $n$-dimensional contour. In the sense of Poincar\'e expansions (i.e.~neglecting exponentially sub-leading terms), the asymptotic behaviour of $\omega(z)$ as $z\to 0^+$ can be computed by the classical Laplace method: we have
\[
\cI(z) \sim \frac{(2 \pi z)^{n/2}\re^{-\phi(q^{\rm cr})/z}}{\sqrt{\det_{ab}(H_{ab})}} \omega(z)\,,
\]
where 
\[
\omega(z) = \exp\l[\frac{z}{2} \sum_{ab} (\HH^{-1})_{ab} \de^2_{ab}\r] \l(\re^{-\phi^{\rm ng}/z} \gamma\r)\,\Bigg|_{q=q^{\rm cr}}\,.
\]
Here, \[\phi^{\rm ng} \coloneqq \phi(q)-\phi(q^{\rm cr}) - \frac{1}{2}\sum_{a,b=1}^n \HH_{ab} (q_a-q_a^{\rm cr}) (q_b-q_b^{\rm cr})\] is the non-Gaussian tail of the Taylor expansion of $\phi(q)$ at $q=q^{\rm cr}$. Following \cite{MR4047552}, we will denote
\[
\mathsf{Asym}_{q=q^{\rm cr}}\big[ \re^{-\phi/z} \gamma \dd^n q \big] \coloneqq \frac{\omega(z)}{\sqrt{\det_{ab}(H_{ab})}} \in \bbC\llbracket z\rrbracket\,.
\]
\begin{example}
Take $n=1$, $\phi(q)=q$, $\gamma(q)=q^{\alpha/z+\beta-1}$ with $\alpha \in \bbR^+$, $\beta \in \bbC$, and $\cC$ the steepest descent path emerging from the unique critical point at $q^{\rm cr} = \alpha$. Then 
\[
\int_{\cC} \re^{-\frac{q}{z}}
q^{\frac{\alpha}{z} + \beta-1} \dd q= \Gamma\l(\frac{\alpha}{z} + \beta\r) z^{\frac{\alpha}{z} + \beta}\,,
\] 
and the Stirling approximation for the Gamma function $\Gamma(x)$ for $\mathrm{Re}(x) \gg 1$  gives
\[
\mathsf{Asym}_{q=q^{\rm cr}}\big[ \re^{-\frac{q}{z}}
q^{\frac{\alpha}{z} + \beta-1} \dd q \big] = \alpha^{\beta-1/2} \exp\l[-\sum_{k > 0} \frac{B_{k+1}(\beta)}{k(k+1)} \l(-\frac{z}{\alpha}\r)^k\r]\,,
\] 
with $B_n(x) = [t^n] \frac{t n! e^{t x}}{e^t-1}$ denoting the $n^{\rm th}$ Bernoulli polynomial. For future use, we will introduce the short-hand notations
%
\[ \Theta(z,x)  \coloneqq  ~
\exp\l[-\sum_{k > 0} \frac{B_{k+1}(x) (-z)^k}{k(k+1)}\r] \in 1+ z \bbQ[x]\llbracket z\rrbracket  \,, \quad 
\Theta(z) \coloneqq \Theta(z,0)  \in 1+ z \bbQ\llbracket z\rrbracket\,, \]
\[
\Upsilon(z_1,z_2) \coloneqq \frac{\Theta(z_1) \Theta(z_2)}{\Theta(z_1+z_2)} \in \bbQ\llbracket z_1,z_2\rrbracket\,.\]
\end{example}

\subsection{Equivariant mirror symmetry for local \texorpdfstring{$\bbP^2$}{P2}}
\label{sec:oscintp2}

Let 
\beq 
\cW(y, \alpha,q) \coloneqq \sum_{i=0}^2 (q_i-\alpha_i \log q_i) + \l(\frac{\prod_{i=0}^2 q_i}{y}\r)^{1/3}
\,.
\label{eq:spotP2}
\eeq 
By \cite[Sec.~0]{MR1653024}, $\cW$ is a mirror superpotential for the $(\Gm)^3$-equivariant small quantum cohomology of $X$, with primitive form given by \[\dd \Omega \coloneqq \frac{1}{3} \prod_{i=0}^2\dd \log q_i\,.\] In particular, there exist Morse-theoretic cycles $\{\Gamma_i\}_{i=0}^2$ for $\mathrm{Re}(\cW)$ such that the periods of $\re^{\cW/z}\dd\Omega$ along $\Gamma_i$ are the restriction to small quantum cohomology of a flat coordinate chart for the Dubrovin connection. Then, taking as representative for $\Gamma_i$ the stable manifold of $\mathrm{Re}(\cW)$ associated to a critical point $q^{\rm cr}_i$ with $\cW(y, \alpha, q^{\rm cr}_i) = \mathsf{u}_i(y)$, the periods equate, up to a $y$-independent factor, the localised component of the $I$-function at $P_i$.
%
%
The steepest descent asymptotics of this equality is
\beq
\mathsf{Asym}_{q_i^{\rm cr}}[\re^{-\cW/z} \dd \Omega] = N_i \mathsf{R}_{0i} \,,
\label{eq:oscintP}
\eeq 
for a formal Taylor series $N_i \in Q_{\Tmax}[(w_i^{(j)})^{1/2}]\llbracket z\rrbracket$. This relative normalisation term  was computed in \cite[Prop.~6.9]{MR4047552}:  specialised to our case, it reads
\beq
N_i(z)= \sqrt{\epsilon_1 \epsilon_2}\,
\Theta\l(\frac{z}{\epsilon_1}\r)\Theta\l(\frac{z}{\epsilon_2}\r)
\,.
\label{eq:normNLR}
\eeq
%
%
%
In particular, when $\alpha_i \neq \alpha_j $ for $i \neq j$, the large radius limit of \eqref{eq:oscintP} is given by
\[
\lim_{y \to 0} \mathsf{Asym}_{q_i^{\rm cr}}[\re^{-\cW/z} \dd \Omega] = \prod_{j=0}^2
\frac{1}{\sqrt{w_i^{(j)}}}
\Theta\l(\frac{z}{w^{(j)}_i}\r)\,.
\] 
When $\alpha_i=\alpha_j$ for some $i\neq j$, meaning that $w_i^{j}=0$, the r.h.s.~of the last equality is not defined.
On the other hand, for $y \neq 0$ the coefficients in the saddle-point asymptotics of \eqref{eq:oscintP} are regular, order by order in $z$, at $\alpha_i=\frac{\epsilon_1+\epsilon_2}{3}$, and \eqref{eq:oscintP} then holds by continuity as an equality in $G_{\mgp{T}_{\rm F}}$. We shall restrict henceforth our analysis of \eqref{eq:oscintP} to the case of the fibrewise torus.\\

When $\alpha_i=\frac{\epsilon_1+\epsilon_2}{3}$, write $s_n \coloneqq e_n(q^{1/3})$ for the $n^{\rm th}$ elementary symmetric polynomial in $(q_0^{1/3}, q_1^{1/3},q_2^{1/3})$. The morphism
\[ 
s : (\Gm)^3 \to \bbC
\] 
is finite 
and unramified away from the diagonals $q_i=q_j$, $i \neq j$. Under the change-of-variables $q\to q(s)$, we have
\bea 
\cW &=& s_1^3 -3 s_1 s_2+ \l(3+\frac{1}{y^{1/3}} \r)s_3 - \l(\epsilon_1+\epsilon_2\r) \log s_3\,, \nn \\ 
\dd\Omega &=&  \sqrt{s_1^2 s_2^2-4 s_1^3  s_3+18 s_1 s_2 s_3-4 s_2^3-27 s_3^2}~\frac{\dd s_1 \dd s_2 \dd s_3}{s_3^3}\,. \nn
\eea
There are two possible sources of critical points of $\cW : (\bbC^{\times})^3 \to \bbC$. The first is given, obviously, by the vanishing of its gradient in $s$-coordinates, 
$\de_{s_n} \cW =0$, which occurs for 
$ (s_1,s_2,s_3)= \l(0,0,\frac{y^{1/3} (\epsilon_1+\epsilon_2)}{1+3 y^{1/3}}\r) 
$. Its pre-images under $q \to s(q)$ are given by the permutations of the components of 
\[
q^{\rm cr}_{0} = \frac{y^{1/3} (\epsilon_1+\epsilon_2)}{1+3 y^{1/3}} (1, \omega, \omega^2)\,,
\]
with $\omega=\re^{2\pi \ri/3}$: by $S_3$-symmetry of $\cW$ in $q$, these all share the same critical value and Laplace expansion. 
The second source of critical points is given by the values of $q$ such that $\de_s \cW$ is in the kernel of the differential of $s(q)$. This is non-zero on the discriminant $q_i=q_j$, $i \neq j$, where it is spanned by $(q_i^{2/3},-q_i, 1)$. Imposing $\de_s \cW(q) = \lambda (q_i^{2/3},-q_i^{1/3}, 1)$ for $\lambda \in \Gm$ yields, on an appropriate branch of $\cW$, two solutions given by:
\beq 
q_{i}^{\rm cr}= \frac{{y}^{1/3} \left(\epsilon _1+\epsilon _2\right)}{\omega^{2i}+3   {y}^{1/3}} \l(1,1, 1\r)\,, \quad i=1,2\,.
\label{eq:qi}
\eeq
We further verify that, for $i=0,1,2$,
\[
y \de_y \cW(q_{i}^{\rm cr}) =
\theta_y u_i = Y_i\,, 
\]
\beq
\det_{jk} \l(q^{1/3}
_j q^{1/3}_k \de^2\cW_{q^{1/3}_j q^{1/3}_k}\r)\bigg|_{q=q^{\rm cr}_{i}} = W_i\,.
\label{eq:crituiWi}
\eeq 
All other critical points of \eqref{eq:spotP2} are each related to one of the three critical points $q_0^{\rm cr}$, $q_1^{\rm cr}$ or $q_2^{\rm cr}$ under a subgroup of the monodromy action on $\cW$, and share the same asymptotic expansion data \eqref{eq:crituiWi} with it. This is a consequence of the invariance of the superpotential \eqref{eq:spotP2} under 
\[
(q_0^{1/3},q_1^{1/3},q_2^{1/3}) \longrightarrow (\omega^{j_0} q_0^{1/3}, \omega^{j_1}q_1^{1/3},\omega^{j_2}q_2^{1/3})
\] 
with $\sum_i {j_i} \equiv 0$ mod $3$. In particular, we have that
\[
\mathsf{Asym}_{q_i^{\rm cr}} \l[\re^{-\cW/z} \dd\Omega \r]=
\mathsf{Asym}_{\widetilde{q_i^{\rm cr}}} \l[\re^{-\cW/z} \dd\Omega \r]\,,
\] 
with 
\beq
\widetilde{q_i^{\rm cr}} =\frac{y^{1/3} (\epsilon_1+\epsilon_2)}{\omega^{2i}+3  y^{1/3}} (1, 1, 1)\,, \quad i=0,1,2\,.
\label{eq:qtilde}
\eeq 
%

\subsection{The orbifold point and the refined Crepant Resolution Conjecture}

Consider the variation of GIT obtained by moving the parameter $\xi$ across the wall $\xi=0$ in \eqref{eq:GIT}. The resulting quotient is
\[
X \coloneqq \bbC^6/\!\!/_{\xi<0}  \Gm = (\bbC^6\setminus\Delta_{\rm orb})/\Gm\,,
\]
where
\[
\Delta_{\rm orb} \coloneqq V\l(\bra x_4\ket\r).
\]
with the same quotient torus as in \eqref{eq:GITact}. Quotienting out $x_4$ yields a residual $\mu_{3}$ action on $\bbC^5$
with weights $(1,1,1,0,0)$. The resulting orbifold $\cY=[\bbC^5/\mu_{3}]$ is the ${\Tmax}$-equivariant vector bundle over the classifying stack of $\mu_3$
\beq
\cY = \textstyle{\bigoplus}_{i=0}^2 \cO_{1} \oplus \textstyle{\bigoplus}_{k=1,2} \cO_{0}  \rightarrow B\mu_{3}\,.
\label{eq:YbunBmu3}
\eeq
The ${\Tmax}$-equivariant Chen--Ruan cohomology is spanned by the fundamental classes $\{\mathbf{1}_0, \mathbf{1}_{1/3}, \mathbf{1}_{2/3}\}$ of each connected component of the inertia stack of $\cY$. 
We will henceforth restrict to the diagonal torus $\mgp{T}_{\rm F}$ by setting $\alpha_i= (\epsilon_1+\epsilon_2)/3$. \\

For $\tau^{\rm orb} \in \Chow_{\Tmax}^{\rm orb}(\cY)$ an orbifold cohomology class, write 
$\tau^{\rm orb}=\sum_{\nu=0}^{2} \tau^{\rm orb}_{\nu} \mathbf{1}_{\frac{\nu}{3}}$ for its expansion in the basis of fundamental classes of the components of the inertia stack $I\cY$. Orbifold Gromov--Witten invariants of $\cY$ are defined as the equivariant Gromov--Witten invariants of $B\mu_3$, twisted by the inverse $\mgp{T}_{\rm F}$-equivariant Euler class of \eqref{eq:YbunBmu3}. 
We will be specifically interested in the 
generating series of 
primary invariants with  
insertions 
in age-shifted degree two:
\begin{equation} 
\begin{split}
\refGWg{g}{\cY}{\mgp{T}_{\rm F}} \coloneqq & \sum_{n \geq 0}
\frac{\sigma^n}{n!} \Big\langle 
\overbrace{\mathbf{1}_{\frac{1}{3}}, \dots, \mathbf{1}_{\frac{1}{3}}}^{n~\mathrm{times}}  \Big\rangle^\cY_{g,n}\,, \\
\mathscr{F}^{\rm orb}(\epsilon_1, \epsilon_2, \sigma) \coloneqq & \sum_{k \geq 0,g\geq -1} (\epsilon_1+\epsilon_2)^{2k} (\epsilon_1 \epsilon_2)^{g-1} \mathscr{F}^{\rm orb}_{k,g} \coloneqq \sum_{g \geq 0} \refGWg{g}{\cY}{\mgp{T}_{\rm F}}
\,, 
\end{split}
\label{eq:Forbref}
\end{equation}
where
\[
\Big\langle 
\overbrace{\mathbf{1}_{\frac{\nu_1}{3}}, \dots, \mathbf{1}_{\frac{\nu_n}{3}}}^{n~\mathrm{times}}  \Big\rangle^\cY_{g,n}  \coloneqq  
\intEquiv{\mgp{T}_{\rm F}}_{[\Mbar_{g,n}(\cY)]^{\vir}_{\mgp{T}_{\rm F}}}
\prod_{k=1}^n \mathrm{ev}_k^* \mathbf{1}_{\frac{\nu_k}{3}} \in \frac{1}{\epsilon_1 \epsilon_2}\bbQ[\epsilon_1, \epsilon_2]\,,  
\]
\[
[\Mbar_{g,n}(\cY)]^{\vir}  \coloneqq e^{\mgp{T}_{\rm F}}(\mathbf{R}^1\pi_* f^* (\cO_{1}))^3 \cap e^{\mgp{T}_{\rm F}}\big( \mathbf{R}\pi_* f^* \cO_{0}\big)^2 \cap [\Mbar_{g,n}(B\mu_3)]^{\vir} \,. 
\]
When $g=0$ and $n=2$, the correlators give the $\mgp{T}_{\rm F}$-equivariant Chen--Ruan intersection pairing on $\cY$:
\[
\eta_{\rm orb}\big(\mathbf{1}_{\frac{\mu}{3}},\mathbf{1}_{\frac{\nu}{3}}\big) \coloneqq 
(\eta_{\rm orb})_{\mu\nu} \coloneqq
\intEquiv{\mgp{T}_{\rm F}}_{[I\cY]_{\mgp{T}_{\rm F}}} \mathbf{1}_{\frac{\mu}{3}} \cup \mathbf{1}_{\bra\frac{3-\nu}{3}\ket} = \frac{1}{3 \epsilon_1 \epsilon_2} 
\l(
\bary{ccc}
\frac{27}{(\epsilon_1+\epsilon_2)^3} & 0 & 0 \\
0 & 0 & 1 \\
0 & 1 & 0 \\
\eary
\r)_{\mu\nu}\,,
\]
with inverse $\eta_{\rm orb}^{\mu\nu} \coloneqq (\eta_{\rm orb}^{-1})_{\mu\nu}$.

\subsubsection{Finite generation in genus zero}
The Gromov--Witten invariants of $\cY$ can be computed
using the orbifold quantum Riemann--Roch theorem \cite{MR2578300}.
In genus zero, the $I$-function $I^{\rm orb}(x,z)$ of $\cY$ reads, for $|x|<3$,
\[
I^{\rm orb}(x,z) \coloneqq   \sum_{k\geq 0} \frac{
\prod_{\stackrel{\bra b \ket = \bra k/3 \ket}{0\leq
      b<\frac{k}{3}}}(\epsilon_1+\epsilon_2-3b z)^3}{27 z^k}\frac{x^k}{k!}\mathbf{1}_{\bra k/{3}\ket}\,.
\]
Denoting $I^{\rm orb}_{\nu,n} \coloneqq [{\mathbf{1}_{\frac{\nu}{3}} z^{-n}}] I^{\rm orb}(x,z) $ as in \cref{sec:genzinf}, we find for example that
\[
I^{\rm orb}_{1,0} = \sum_{m =0}^{\infty}\frac{(-1)^m x^{3 m+1} \l(\frac{1}{3}\r)_m^3}{(3 m+1)!}\,, \qquad I^{\rm orb}_{2,1} = \sum_{m =0}^{\infty}\frac{(-1)^m x^{3 m+2} \l(\frac{2}{3}\r)_m^3}{(3 m+2)!}\,,
\]
\[
I^{\rm orb}_{1,1} =(\epsilon_1+\epsilon_2) \sum_{m=0}^\infty \frac{(-1)^m  x^{3 m+1}  \left(\frac{1}{3}\right)_m^3 \left(\psi ^{(0)}\left(\frac{4}{3}\right)-\psi ^{(0)}\left(m+\frac{1}{3}\right)-3\right)}{27
 (3 m+1)!}\,.
\]
 The $I$-function is equal \cite{MR2510741,MR2486673} to the restriction of the big $J$-function $J^{\rm orb}(\tau,z)$ of $\cY$ to the locus $\tau^{\rm orb}_0=\tau^{\rm orb}_2=0$:
\[
J^{\rm orb}(\tau,z)\big|_{\tau^{\rm orb}_{\nu}=  \delta_{\nu,1} \sigma}
= I^{\rm orb}(x(\sigma),z)\,,
\]
where the orbifold mirror map is 
\[
\sigma \rightarrow \sigma(x) = I^{\rm orb}_{1,0}\,.
\]
In complete analogy with \eqref{eq:ADEFXdef}, and writing  $\theta_x \coloneqq x \de_x$, we define the following  analytic functions in the disk $|x|<3$ from the first few Taylor coefficients of $I^{\rm orb}(x,z)$ around $z=\infty$:
\[
A_{\rm orb}(x)  \coloneqq  
\theta_x I^{\rm orb}_{1,0} = x \, _2F_1\left(\frac{1}{3},\frac{1}{3};\frac{2}{3};-\frac{x^3}{27}\right) \,, \quad F_{\rm orb}(x)  \coloneqq   \theta_x \log A_{\rm orb} \,,
%
%
%
\]
\[
D_{\rm orb}(x)  \coloneqq  \theta_x\l(\frac{1}{A_{\rm orb}} \theta_x I^{\rm orb}_{2,1}\r)\,, \quad
E_{\rm orb}(x)  \coloneqq \quad 
\frac{1}{\epsilon_1+\epsilon_2}\theta_x\l(\frac{1}{A_{\rm orb}} \theta_x I^{\rm orb}_{1,1}\r)\,. \quad
\]
Recall that
\[
X(y)\Big|_{y=x^{-3}} = \frac{x^3}{x^3+27}\,. 
\]
As in \cref{lem:modrel}, the ring $\bbQ[A_{\rm orb},D_{\rm orb},E_{\rm orb},F_{\rm orb},X]$ is not freely generated.
\begin{lem}
The following relations hold in $\bbQ \llbracket x \rrbracket$:
\beq
A_{\rm orb}^2 D_{\rm orb}=27 X\,, \quad 
E_{\rm orb} = \frac{3}{2}X\,.
\label{eq:DEorbrel}
\eeq 
In particular, $\bbQ[A_{\rm orb},D_{\rm orb},E_{\rm orb},F_{\rm orb},X] \simeq \bbQ[A_{\rm orb}^{\pm 1}, F_{\rm orb},X] \simeq \bbQ[A_{\rm orb},D_{\rm orb},F_{\rm orb}]$.
\label{lem:modrelorb}
\end{lem}
\begin{proof}[Proof (sketch)]
The proof is largely identical to that of \cref{lem:modrel}. The first relation in \eqref{eq:DEorbrel} is proved using the methods of \cite{MR2454324} (see also \cite{MR3960668}). The second relation is implied by
\[
 A_{\rm orb}  \log \left(1+\frac{x^3}{27}\right)=\frac{2}{\epsilon_1+\epsilon_2}\theta_x I^{\rm orb}_{1,1}\,,
\]
upon dividing by $A_{\rm orb}(x)$ and acting 
which $\theta_x$. This is in turn equivalent to the  convolution identity
\[ 
\frac{1}{2}\sum_{k=0}^{d-1}\frac{1}{d-k}
\frac{\left(\frac{1}{3}\right)_k^2}{k!\left(\frac{2}{3}\right)_k}
=
\frac{3^{3d} \left(\frac{1}{3}\right)_d^3 \left(\psi ^{(0)}\left(d+\frac{1}{3}\right)-\psi ^{(0)} \left(\frac{4}{3}\right)+3\right)}{(3 d)!}\,,
\]
whose proof is, {\it mutatis mutandis}, the same as the one for \eqref{eq:limitconv} in \cref{lem:modrel}.
\end{proof}

\noindent The structure constants of the orbifold quantum product,
\[
(\mathsf{C}_{\rm orb})^{\lambda}_{\mu\nu}
\coloneqq \sum_{\rho} (\eta_{\rm orb})^{\lambda \rho} \bra \! \bra \mathbf{1}_{\frac{\mu}{3}}, \mathbf{1}_{\frac{\nu}{3}},\mathbf{1}_{\frac{\rho}{3}}\ket \! \ket_{0,3}^\cY\Big|_{\tau^{\rm orb}_\nu = \delta_{\nu 1} \sigma(x)}
\]
are easily expressed as Laurent polynomials in $A_{\rm orb}$ and $X$. The String Equation and the form of the Chen--Ruan pairing imply that
\[
(\mathsf{C}_{\rm orb})^{0}_{\mu\nu} = \eta_{\rm orb}^{00} (\eta_{\rm orb})_{\mu\nu}\,, \quad 
(\mathsf{C}_{\rm orb})^{\mu}_{0\nu} = (\mathsf{C}_{\rm orb})^{\mu}_{\nu 0} = \delta^{\mu}_\nu\,, \quad (\mathsf{C}_{\rm orb})^{1}_{12} = (\mathsf{C}_{\rm orb})^{1}_{21} =
(\mathsf{C}_{\rm orb})^{2}_{11}\,,
\]
\[
(\mathsf{C}_{\rm orb})^{1}_{22} = (\mathsf{C}_{\rm orb})^{2}_{21} =
(\mathsf{C}_{\rm orb})^{2}_{12} =
(\mathsf{C}_{\rm orb})^{1}_{11}\,,
\]
while the WDVV equation gives
\[
(\mathsf{C}_{\rm orb})^{2}_{22}=
\frac{\left(\epsilon _1+\epsilon _2\right){}^3}{27
   (\mathsf{C}_{\rm orb})^1_{11}}+(\mathsf{C}_{\rm orb})^2_{11}\,,
\]
meaning that all structure constants are computed from $(\mathsf{C}_{\rm orb})^\mu_{11} = \de^2_{\sigma} I^{\rm orb}_{\mu,1}$, $\mu=1,2$. Using \cref{lem:modrelorb}, we get
\[
(\mathsf{C}_{\rm orb})^1_{1,1} = \l(\frac{1}{A_{\rm orb}}\theta_x\r)^2 I^{\rm orb}_{1,1} = \frac{3(\epsilon_1+\epsilon_2)}{2} \frac{X}{A_{\rm orb}}\,, \quad 
(\mathsf{C}_{\rm orb})^2_{1,1} = \l(\frac{1}{A_{\rm orb}}\theta_x\r)^2 I^{\rm orb}_{2,1} = \frac{27 X}{A_{\rm orb}^3}\,. 
\]

\subsubsection{Finite generation in higher genus}

Armed with this, the $\mgp{T}_{\rm F}$-equivariant Gromov--Witten R-matrix of $\cY$, restricted to $\tau^{\rm orb}_{\nu}=\delta_{\nu,1}\sigma(x)$, satisfies 
\begin{equation}
\begin{split}
\mathsf{R}_{\rm orb}(\mathbf{1}_{\frac{1}{3}}) &= \, z \de_\sigma \mathsf{R}_{\rm orb}(\mathbf{1}_0)  = \frac{z}{A_{\rm orb}}\theta_x \mathsf{R}_{\rm orb}(\mathbf{1}_0)\,,  \\
\mathsf{R}_{\rm orb}(\mathbf{1}_{\frac{2}{3}}) &= 
\frac{
z \de_\sigma \mathsf{R}_{\rm orb}(\mathbf{1}_{\frac{1}{3}}) - \sum_{\mu \neq 2} (\mathsf{C}_{\rm orb})^\mu_{1,1} \mathsf{R}_{\rm orb}(\mathbf{1}_{\frac{\mu}{3}})}
{(\mathsf{C}_{\rm orb})^2_{1,1}} 
= 
\frac{A_{\rm orb}^2}{27 X}\l[
z\theta_x - \frac{3(\epsilon_1+\epsilon_2)}{2} A_{\rm orb} X\r]
\mathsf{R}_{\rm orb}(\mathbf{1}_{\frac{1}{3}})\,,
\end{split}
\label{eq:R12orb}
\end{equation}
where, in the same vein as for \eqref{eq:RHb}, we have used the fact that $\sigma = \tau^{\rm orb}_{1}$ in the first equality, and the WDVV equation in the second. As in \cref{sec:Rmatrix}, write \[\mathsf{R}_{\nu i}^{\rm orb}(x,z) \coloneqq \re^{-\mathsf{u}_i/z} \mathsf{R}^{\rm orb}(\mathbf{1}_{\frac{\nu}{3}})\,,\]
where $\Psi^{\rm orb}_{\nu i}$ is the Jacobian matrix of the change-of-basis \[\{\mathbf{1}_{\frac{\nu}{3}}\}_{\nu=0,1,2} \to \{ \mathsf{f}_i\}_{i=0,1,2}\] restricted to $\tau_\mu^{\rm orb}=\delta_{\mu, 1} \sigma(x)$. By \eqref{eq:R12orb}, $\mathsf{R}_{\nu i}^{\rm orb}(x,z)$ is determined by $\mathsf{R}_{0i}^{\rm orb}(x,z)$.
\begin{prop}
The following identity holds in $\sqrt{\epsilon_1 \epsilon_2}\, Q_{\mgp{T}_{\rm F}}\llbracket x,z\rrbracket$:
\[
\mathsf{Asym}_{\widetilde{q}_i^{\rm cr}}[\re^{-\cW/z} \dd\Omega]\Big|_{y=x^{-3}}
={\sqrt{\epsilon_1 \epsilon_2}}\,\Theta\l(\frac{z}{\epsilon_1} \r) \Theta\l(\frac{z}{\epsilon_2} \r)\mathsf{R}^{\rm orb}_{0i} \,. 
\]
\label{prop:Rnorm}
\end{prop}

\begin{proof}
Upon identifying $y=x^{-3}$, the kernel of the equivariant Picard--Fuchs operator \eqref{eq:pf} is spanned by the components of $I^{\rm orb}(x,z)$  \cite{MR2486673,Brini:2014fea},
\[
\l[x^3(\epsilon_1+\epsilon_2-z \theta_x) - 27 z^3\prod_{j=0}^2 (\theta_x- j)  \r] I^{\rm orb}(x,z) =0\,.
\]
This implies that the components of the I-function $I(x^{-3},z)$ of $Z$ are a flat chart for the Dubrovin connection of $\cY$ restricted to the locus $\tau_\nu^{\rm orb}=\delta_{\nu, 1} \sigma(x)$. By \eqref{eq:normNLR},
the same is true for the steepest descent asymptotics of the periods of $\re^{-\cW/z} \dd \Omega$, with $y=x^{-3}$, to all orders in a formal expansion at $z=0$. As such, the latter equate the identity component of the canonical $R$-matrix of $\cY$ up to right-multiplication by a diagonal factor $N_i^{\rm orb} \in Q_{\mgp{T}_{\rm F}}[\sqrt{\epsilon_1},\sqrt{\epsilon_2}]\llbracket z\rrbracket$ which is constant in $x$:
\[ 
\mathsf{Asym}_{\widetilde{q}_i^{\rm cr}}[\re^{-\cW/z} \dd\Omega]\bigg|_{y=x^{-3}}
=\mathsf{R}^{\rm orb}_{0i} N_i^{\rm orb}\,. 
\] 
As explained in  \cite{Brini:2013zsa}, Tseng's orbifold quantum Riemann--Roch theorem \cite{MR2578300}, applied to the  Gromov--Witten theory of $B\mu_3$ \cite{MR1950944} twisted by the $\mgp{T}_{\rm F}$-equivariant inverse Euler class of \eqref{eq:YbunBmu3}, fixes uniquely the $R$-matrix in terms of its boundary behaviour at $\sigma=0$. We have \cite[Lem.~6.5]{Brini:2013zsa}
\beq
\mathsf{R}^{\rm orb}_{\nu i}\big|_{\sigma=0} =  
\frac{\omega^{\nu i}}{3\sqrt{\epsilon_1 \epsilon_2}} \l(\frac{\epsilon_1+\epsilon_2}{3}\r)^{\nu-3/2}\exp\left[
\sum_{k>0} \l(-\frac{3B_{k+1}\left(\frac{\nu}{3}\right) (-3z)^k}{k(k+1)(\epsilon_1+\epsilon_2)^k}+ 
\frac{B_{2k} z^{2k-1}(\epsilon_1^{1-2k}+\epsilon_2^{1-2k})}{2k(2k-1)}\r) \r]\,. 
\label{eq:Rorbnorm}
\eeq
On the other hand, at $x=0$ we have that $\cW=\sum_{j=0}^2 (q_i-(\epsilon_1+\epsilon_2)/3 \log q_i)$, $\dd \Omega = 1/3 \prod_{j=0}^2 \dd \log q_i$, and $\widetilde{q^{\rm cr}_i} = (\epsilon_1+\epsilon_2)/3$. Then, the oscillating integral \eqref{eq:oscintP} factorises at $x=0$ into a product of three Gamma integrals:
\[
\int_{\Gamma_i} \re^{-\cW/z}\dd \Omega\bigg|_{x=0} = \frac{1}{3}\prod_{j=0}^2 \int_0^\infty \re^{-q_j/z} q_j^{\frac{\epsilon_1+\epsilon_2}{3 z}-1} \dd q_i = \frac{z^{\frac{\epsilon_1+\epsilon_2}{z}}}{3} 
\Gamma\l( \frac{\epsilon_1+\epsilon_2}{3z} \r)^3
\,,
\label{eq:Gammaorb}
\]
so that
\[ 
\mathsf{Asym}_{\widetilde{q}_i^{\rm cr}}[\re^{-\cW/z} \dd\Omega]\bigg|_{x=0} = \sqrt{\frac{3}{(\epsilon_1+\epsilon_2)^3 }}\exp\left[
-\sum_{k>0}\frac{3B_{k+1}\left(0\right) (-3z)^k}{k(k+1)(\epsilon_1+\epsilon_2)^k}
 \r]\,.
\] 
%
%
%
Since 
$X = 27x^3 + \cO(x^6)$ and $A_{\rm orb} = x + \cO(x^2)$, 
by \eqref{eq:R12orb} the components $R^{\rm orb}_{\nu i}$ of the R-matrix satisfy 
\[
\frac{\de^\nu}{\de x^\nu}\mathsf{Asym}_{\widetilde{q}_i^{\rm cr}}[\re^{-\cW/z} \dd\Omega]\bigg|_{x=0}
=\mathsf{R}^{\rm orb}_{\nu i} N_i^{\rm orb}\bigg|_{x=0}\,. 
\]
Using that, from \eqref{eq:qtilde},
\begin{align}
\frac{\de^\nu}{\de x^\nu}
\int_{\Gamma_i}  \re^{-\cW/z}\dd \Omega\bigg|_{x=0} & =
\omega^{\nu i} \frac{\de^\nu}{\de x^\nu}
\int_{\Gamma_0}  \re^{-\cW/z}\dd \Omega\bigg|_{x=0}=
\frac{1}{3}\prod_{j=0}^2 \int_0^\infty \re^{-q_j/z} q_j^{\frac{\epsilon_1+\epsilon_2}{3 z}+\frac{\nu}{3}-1} \dd q_i \nn \\ & =  \omega^{\nu i}\frac{z^{\frac{\epsilon_1+\epsilon_2}{z} + \nu}}{3} 
\Gamma\l( \frac{\epsilon_1+\epsilon_2}{3z} +\frac{\nu}{3} \r)^3
\,, \nn
\end{align}
we find, for $\nu \in \{0,1,2\}$, that
\[
\frac{\de^\nu}{\de x^\nu}\mathsf{Asym}_{\widetilde{q}_i^{\rm cr}}[\re^{-\cW/z} \dd\Omega]\bigg|_{x=0}=
\frac{\omega^{\nu i}}{3}
\l(\frac{\epsilon_1+\epsilon_2}{3}\r)^{\nu-3/2}\exp\left[
\sum_{k>0} \l(-\frac{3B_{k+1}\left(\frac{\nu}{3}\right) (-3z)^k}{k(k+1)(\epsilon_1+\epsilon_2)^k}\r) \r]\,. 
\]
Comparing with \eqref{eq:Rorbnorm}, we conclude that 
\[
N_i^{\rm orb} = {\sqrt{\epsilon_1 \epsilon_2}}\,\Theta\l(\frac{z}{\epsilon_1} \r) \Theta\l(\frac{z}{\epsilon_2} \r)\,.
\]
\end{proof}

\medskip 
Define a homomorphism of rings
\[
\mathsf{Orb} :  G_{\mgp{T}_{\rm F}}[A,F]  \longrightarrow   G_{\mgp{T}_{\rm F}}[A_{\rm orb},F_{\rm orb}]
\]
induced by the map on generators
\[
\mathsf{Orb}(A)\coloneqq \frac{A_{\rm orb}}{3}\,, \quad  
\mathsf{Orb}(F) \coloneqq -\frac{F_{\rm orb}}{3}
\,.
\]
%
As we did in \eqref{eq:asympR}, we will express the asymptotic expansion of $\mathsf{R}^{\rm orb}_{aj}$ in the form
\[
\mathsf{R}^{\rm orb}_{0j} = 
\sum_{n=0}^\infty\,
\mathsf{R}^{\rm orb, [0]}_{j,n} z^n\,, \quad
\mathsf{R}^{\rm orb}_{1j} = \frac{Y_j}{A_{\rm orb}} \sum_{n=0}^\infty 
\mathsf{R}^{{\rm orb},[1]}_{j,n} z^n\,, \quad
\mathsf{R}^{\rm orb}_{2j} = \frac{Y_j^2}{A_{\rm orb} D_{\rm orb}} \sum_{n=0}^\infty 
\mathsf{R}^{{\rm orb},[2]}_{j,n} z^n\,.
\]
%
%
\begin{prop}
    Upon identifying $y=x^{-3}$, the R-matrices of $Z$ and $\cY$ are related as
\[
\mathsf{R}^{\rm orb}_{0i} =  
\mathsf{Orb}(\mathsf{R}_{0i})\,, 
\quad \mathsf{R}^{\rm orb}_{1i} =  
\mathsf{Orb}(-\mathsf{R}_{1i})\,, 
\quad \mathsf{R}^{\rm orb}_{2i} =  
\mathsf{Orb}\l(\mathsf{R}_{2i}
+\frac{\epsilon_1+\epsilon_2}{2} \mathsf{R}_{1i}
\r)\,.
\]  
Moreover, the Taylor coefficients of $\mathsf{R}^{\rm orb}_{\nu i}$ satisfy
\[
\mathsf{R}^{{\rm orb}, [0]}_{i,n}\in G_{\mgp{T}_{\rm F}}\,, \quad 
\mathsf{R}^{{\rm orb}, [1]}_{i,n} \in G_{\mgp{T}_{\rm F}}\,, \quad 
\mathsf{R}^{{\rm orb}, [2]}_{i,n} \in G_{\mgp{T}_{\rm F}} +F_{\rm orb} G_{\mgp{T}_{\rm F}} \,.
\] 
\label{prop:fingenRorb}
\end{prop}
\begin{proof}
The first equality follows from 
\[
\mathsf{Asym}_{\widetilde{q}_i^{\rm cr}}[\re^{-\cW/z} \dd\Omega] = \mathsf{R}_{0i} N_i = \mathsf{R}^{\rm orb}_{0i} N^{\rm orb}_{i}\,, \quad
N_i=N_i^{\rm orb} = \sqrt{\epsilon_1 \epsilon_2}\, \Theta\l(\frac{z}{\epsilon_1}\r)\Theta\l(\frac{z}{\epsilon_2}\r).
\]
The other two follow from comparing \eqref{eq:RHb} and \eqref{eq:R12orb}, using \cref{lem:modrel,lem:modrelorb}  and bearing in mind that for $f\in G_{\mgp{T}_{\rm F}}[A_{\rm orb}, F_{\rm orb}]$, 
\beq
\de_\sigma f = \frac{\theta_x f}{A_{\rm orb}} = \frac{-3 \theta_y f}{3 \mathsf{Orb}(A)} = -\mathsf{Orb}\l(\frac{ \theta_y f}{A}\r) = -\mathsf{Orb}\l(\cQ \de_\cQ f\r)\,.
\label{eq:desigmacrc}
\eeq
From this, the second part of the statement follows from \eqref{eq:R12orb}, \cref{lem:modrelorb}, and the finite-generation statement of \cref{prop:fingenR}.
\end{proof}

\subsubsection{The Refined Crepant Resolution Correspondence}

We are now in a position to compute $\refGWg{g}{\cY}{\mgp{T}_{\rm F}}$ 
from the 
orbifold $R$-matrix. Dilaton and edge factors are
\[ 
\mathsf{D}^{\rm orb}_{i,n} \coloneqq [z^{n-1}]\l(\Psi^{\rm orb}_{0i}-\sum_j \mathsf{R}_{\rm orb}(-z)_{0i}\r)\,, \quad 
(\mathsf{E}^{\rm orb})^{ij}_{n,m} \coloneqq [z^{n} w^m]\l( \frac{\delta_{ij}-\sum_{\mu\nu}\mathsf{R}^{\rm orb}_{\mu i}(-z) (\eta_{\rm orb})^{\mu\nu} \mathsf{R}^{\rm orb}_{\nu j}(-w)}{z+w}\r)\,,
\]
with dilaton leaf, edge, and vertex contributions given as
\[
\mathrm{Cont}_{\rm orb}(d) \coloneqq
\mathsf{D}^{\rm orb}_{\mathsf{q}(\mathsf{v}(d)), \mathsf{a}(d)}\,, \quad
\mathrm{Cont}_{\rm orb}(e) \coloneqq
(\mathsf{E}^{\rm orb})^{\mathsf{q}(\mathsf{v}(h_1)), \mathsf{q}(\mathsf{v}(h_2)}_{\mathsf{a}(h_1),\mathsf{a}(h_2)}\,,
\]
\[
\mathrm{Cont}_{\rm orb}(v) = (3 W_i)^{\frac{2\mathsf{g}(v)-2+\mathsf{|H(v)|}}{2}} 
\int_{[\Mbar_{\mathsf{g}(v), |\HH(v)|}]}\prod_{k=1}^{|\HH(v)|} 
\psi_k^{\mathsf{a}(h_k)}\,.
\]
%
The Givental--Teleman reconstruction formula then reads
\beq
\refGWg{g}{\cY}{\mgp{T}_{\rm F}}
=
\sum_{\Gamma \in \cG^{[2]}_{g,0}} \frac{(\epsilon_1 \epsilon_2)^{g-1}}{|\mathrm{Aut}(\Gamma)|} 
\prod_{e\in \cE_\Gamma}\mathrm{Cont}_{\rm orb}(e)
\prod_{d\in \DD_\Gamma 
}\mathrm{Cont}_{\rm orb}(d) \prod_{v\in \cV_\Gamma} 
\mathrm{Cont}_{\rm orb}(v)\,.
\label{eq:gtorb}
\eeq 
Define now \[
\mathscr{F}^{\rm orb}_{k,g} = [(\epsilon_1+\epsilon_2)^{2k} (\epsilon_1 \epsilon_2)^{g-1}] \refGWg{k+g}{\cY}{\mgp{T}_{\rm F}}
\,.\]

\begin{thm}[Refined Crepant Resolution Correspondence]
For all $k,g\geq 0$ with $k+g \geq 2$, we have 
%
\begin{align*}
\mathscr{F}^{\rm orb}_{k,g} = 
\mathsf{Orb}\l(\mathscr{F}_{k,g}\r)\,.
\end{align*}
\label{thm:crc}
\end{thm}
\begin{proof}
From \cref{prop:fingenRorb} we obtain immediately that 
\[
\mathsf{D}_{i,n}^{\rm orb}(d)=\mathsf{Orb}(\mathsf{D}_{i,n}(d))\,.
\]
For the edge contributions we find, using \cref{prop:fingenRorb}, that
\begin{align}
& \frac{1}{3 \epsilon_1 \epsilon_2}\sum_{\mu\nu}\mathsf{R}^{\rm orb}_{\mu i}(-z) (\eta_{\rm orb})^{\mu\nu} \mathsf{R}^{\rm orb}_{\nu j}(-w) =  
\frac{(\epsilon_1+\epsilon_2)^3}{27}\mathsf{R}^{\rm orb}_{0i}(-z)  \mathsf{R}^{\rm orb}_{0j}(-w)+ \mathsf{R}^{\rm orb}_{1i}(-z)  \mathsf{R}^{\rm orb}_{2j}(-w)+\mathsf{R}^{\rm orb}_{2i}(-z)  \mathsf{R}^{\rm orb}_{1j}(-w) \nn \\
& =  
\mathsf{Orb}\l[\frac{(\epsilon_1+\epsilon_2)^3}{27}
\mathsf{R}_{0i}(-z)  \mathsf{R}_{0j}(-w)\r]  + \mathsf{Orb}\l[
- \mathsf{R}_{1i}(-z)  \l(\mathsf{R}_{2j}(-w)+ \frac{\epsilon_1+\epsilon_2}{2} \mathsf{R}_{1j}(-w)\r) + \big( i \leftrightarrow j \big) \r] \nn \\
& =  \mathsf{Orb}\l[\frac{1}{3 \epsilon_1 \epsilon_2}\sum_{ab}\mathsf{R}_{a i}(-z) \eta^{ab} \mathsf{R}_{b j}(-w)\r]\,, \nn
\end{align}
where in the last line we have used \eqref{eq:preedge}. This in turn implies that 
\[
\mathrm{Cont}_{\rm orb}(v)=\mathsf{Orb}(\mathrm{Cont}(v))\,, \quad \mathrm{Cont}_{\rm orb}(d)=\mathsf{Orb}(\mathrm{Cont}(d))\,, \quad \mathrm{Cont}_{\rm orb}(e)=\mathsf{Orb}(\mathrm{Cont}(e))\,.
\]
where the homomorphism $\mathsf{Orb}$ acts as the identity in the first two equalities. The claim then follows from \eqref{eq:gtH} and \eqref{eq:gtorb}. 
\end{proof}

\begin{cor}
For $k+g \geq 2$, the refined generating series of $\cY$ satisfy the finite generation property 
\[ 
\mathscr{F}^{\rm orb}_{k,g} \in \bbQ[X^{\pm 1}][ F_{\rm orb}]\,,
\] 
with $\deg_F\mathscr{F}_{k,g}\leq 3(g+k-1)$. Moreover,  the refined modular anomaly equation holds,
\[
-\frac{27 X}{A_{\rm orb}^2}
\frac{\de \mathscr{F}^{\rm orb}_{k,g}}{\de F_{\rm orb}}=\frac{1}{2}  \partial_\sigma^2 \mathscr{F}^{\rm orb}_{k,g-1} + \frac{1}{2} \sum_{\substack{0\leq k' \leq k\,, 0 \leq g' \leq g \\(0,0)\neq (k',g')\neq (k,g) }}   \de_\sigma\mathscr{F}^{\rm orb}_{k',g'} \, \partial_\sigma \mathscr{F}^{\rm orb}_{k-k',g-g'}\,,
\]
as an identity in $A_{\rm orb}^{-2}\,\bbQ[X^{\pm 1}][F_{\rm orb}]$.
\label{cor:orbHAE}
\end{cor}

\begin{proof}
    This follows from \cref{thm:crc,thm:haeref} using \eqref{eq:desigmacrc}.
\end{proof}

\subsubsection{Orbifold regularity}
An important corollary of \cref{thm:crc} is the ``orbifold regularity'' property \cite{HK10:OmegaBG}. Let $\tilde F \coloneqq F+1/3 $, and define a grading on the ring $\bbQ[X^{\pm 1}, \tilde F]$ by $\deg X = \deg \tilde F=1$. We will denote by $\mathfrak{g}_k$ the subspaces in the corresponding filtration, generated by monomials of homogeneous degree less than or equal to $k \in \bbZ_{\geq 0}$, and by $\mathfrak{g}_k^{\rm reg}$ the subspace of $\mathfrak{g}_k$ satisfying $3 \deg X + \deg \tilde F \geq 0$.

\begin{prop}[Orbifold regularity]
For $k+g \geq 2$, we have
    $
     \mathscr{F}_{k,g} \in X^{1-k-g} ~\mathfrak{g}^{\rm reg}_{3g+3k-3}\,.$
\label{cor:orbreg}
\end{prop}

\begin{proof}
Under the specialisation 
\beq
F_{\rm orb} \longrightarrow F_{\rm orb}(x) \in \bbQ\llbracket x\rrbracket\,, \quad X \longrightarrow X(x^{-3}) \in \bbQ\llbracket x\rrbracket\,,
\label{eq:specFX}
\eeq
the orbifold refined generating series
 are by definition regular (in fact analytic) at $\sigma=x=0$, from \eqref{eq:Forbref}. By  \cref{cor:orbHAE},  
$\mathscr{F}^{\rm orb}_{k,g}$ is a polynomial in $F_{\rm orb}$ and a Laurent polynomial in $X$: the fact that $F_{\rm orb}(x)=1 + \cO(x)$, $X(x^{-3})=\cO(x^3)$, and $\deg_{F_{\rm orb}}\mathscr{F}^{\rm orb}_{k,g}=3(g+k-1)$ collectively imply that the degree of $\mathscr{F}^{\rm orb}_{k,g}$ in $X^{-1}$ cannot exceed $g+k-1$. More in detail, 
\[
\mathscr{F}^{\rm orb}_{k,g} = X^{1-k-g} 
\sum_{n=0}^{3(g+k-1)} \sum_{m=0}^{d_n} a_{mn} X^m F_{\rm orb}^n
\,,
\]
for some $d_n \in \bbZ_{\geq 0}$ and coefficients $a_{mn} \in \bbQ$ such that the polar part of $\mathscr{F}^{\rm orb}_{k,g} $ in $X$, under the specialisation \eqref{eq:specFX}, vanishes at $x=0$.
Shifting $F_{\rm orb} \to \tilde{F}_{\rm orb} = F_{\rm orb} -1=-3 \mathsf{Orb}(\tilde F)$, so that $\tilde{F}_{\rm orb}(x)=x+\cO(x^2)$, entails that
\[
\mathscr{F}^{\rm orb}_{k,g} = X^{1-k-g} 
\sum_{n=0}^{3(g+k-1)} \sum_{m=0}^{d_n} b_{mn} X^m \tilde F_{\rm orb}^n
\,
\]
for some $b_{mn} \in \bbQ$ such that $b_{mn}=0$ if $3m+n<0$, again by requiring that the r.h.s.~is $\cO(1)$ at $x=0$ under \eqref{eq:specFX}. Therefore $\mathscr{F}^{\rm orb}_{k,g} \in X^{1-k-g} ~ \mathsf{Orb}(\mathfrak{g}^{\rm reg}_{3g+3k-3})$, and by \cref{thm:crc},
\[
\mathscr{F}_{k,g} =  X^{1-k-g} 
\sum_{n=0}^{3(g+k-1)} \sum_{m=0}^{d_n} (-3)^{n} b_{mn} X^m \tilde F^n \in X^{1-k-g}\,\mathfrak{g}_{3g+3k-3}^{\rm reg}\,. \qedhere
\]
\end{proof}

\begin{table}[t]
    \centering
\def\arraystretch{1.5}
\begin{tabular}{|c|c|c|c|c|c|c|}
\hline
 \backslashbox{$(k,g)$}{$m$}  & 1 & 2 & 3 & 4 & 5 & 6 \\
 \hline
$(0,1)$  & 0 & $\frac{1}{243}$ & $-\frac{14}{243}$ & $\frac{13007}{6561}$ & $-\frac{8354164}{59049}$ & $\frac{9730293415}{531441}$ \\  \hline
$(1,0)$  & $\frac{1}{108}$ & $-\frac{5}{324}$ & $\frac{1319}{8748}$ & $-\frac{114983}{26244}$ & $\frac{66572075}{236196}$ & $-\frac{24043411967}{708588}$ \\  \hline
$(0,2)$  & $\frac{1}{19440}$ & $-\frac{13}{11664}$ & $\frac{20693}{524880}$ & $-\frac{12803923}{4723920}$ & $\frac{314291111}{944784}$ & $-\frac{8557024202467}{127545840}$ \\  \hline
$(1,1)$  & $-\frac{1}{1620}$ & $\frac{7}{972}$ & $-\frac{8933}{43740}$ & $\frac{1628851}{131220}$ & $-\frac{331638709}{236196}$ & $\frac{2836254577177}{10628820}$ \\  \hline
$(2,0)$  & $\frac{1}{1080}$ & $-\frac{79}{11664}$ & $\frac{29}{180}$ & $-\frac{4656751}{524880}$ & $\frac{1332163447}{1417176}$ & $-\frac{804769240343}{4723920}$ \\  \hline
    \end{tabular}
    \vspace{.3cm}
    \caption{Refined orbifold Gromov--Witten invariants 
    of $[\bbC^3/\mu_3]$.}
    \label{tab:orbrefGW}
\end{table}

A uniform bound $d_n \leq 3(g+k-1)$ holds from the conifold asymptotics, as we will shall see in \cref{sec:cfpt}, meaning that as a Laurent polynomial in $X$ we have $\deg_{X}\mathscr{F}_{k,g} \leq 2g+2k-2$. 
\cref{cor:orbreg} then implies the following enhanced quasi-modularity result about $\mathscr{F}_{k,g}$.

\begin{cor}
For $k+g \geq 2$, the refined generating series satisfy
    \[
     \mathscr{F}_{k,g} \in C^{2(1-g-k)} \mathrm{QMod}^{[6(g+k-1),3(g+k-1)]}\l(\Gamma_1(3)\r)\,.
     \]    
\end{cor}

\begin{proof}
Let $\bbC[A,C]_n$ (respectively, $\bbC[A,B,C]_n$) denote the weight-$n$ graded component of $\bbC[A,C]$ (respectively, $\bbC[A,B,C]$). By \cite[Lem.~4.1]{BFGW21:HAE},
we have that $\bbC[A,C]_{n} \simeq \mathrm{Mod}^{[n]}(\Gamma_1(3))$, hence 
\[
\bbC[A,B,C]_{2n} \simeq 
\mathrm{QMod}^{[2n,n]}(\Gamma_1(3))\,.
\]
Expressing $X$ and $F$ in terms of the quasi-modular forms $A$, $B$, $C$, there is an isomorphism of $\bbC$-vector spaces \cite[Prop.~4.3]{BFGW21:HAE}
\[
\l\{\phi \in X^{1-k-g} \mathfrak{g}^{\rm reg}_{3g+3k-3} \Big| \deg_X \phi \leq 2(g+k-1)\,,  \deg_F \phi \leq 3g+3k-3\r\} \simeq C^{2(1-g-k)} \bbC[A,B,C]_{6(g+k-1)}\,,
\]
from which the claim follows.
\end{proof}

In \cref{tab:orbrefGW}, we tabulate the refined orbifold Gromov--Witten invariants 
\[\de^{3m}_\sigma \mathscr{F}^{\rm orb}_{k,g}\Big|_{\sigma=0}\] for $k+g \leq 2$ and $m \leq 6$. These confirm\footnote{The row $(k,g)=(0,2)$ corrects typographical mistakes in the last row of \cite[Eqs.~7.16]{HKK13:OmegaBmodelRigidN2}.} the numerical predictions  of \cite{HKK13:OmegaBmodelRigidN2} based on B-model physics.

\pagebreak
\subsection{The conifold point and the Barnes double Gamma function}
\label{sec:cfpt}


\subsubsection{Asymptotics of the double Gamma function}

For $n \geq -2$, define Laurent polynomials $\beta_n(\epsilon_1,\epsilon_2) \in (\epsilon_1 \epsilon_2)^{-1} \bbQ[\epsilon_1, \epsilon_2]$ by
\beq
\beta_n(\epsilon_1,\epsilon_2) \coloneqq 
\begin{cases}
\frac{1}{\epsilon_1 \epsilon_2}\,, & n=-2\,,\\
0 & n=-1\,,\\
- \frac{\epsilon_1^2+\epsilon_2^2}{24\epsilon_1 \epsilon_2} & n=0\,,\\
(n-1)![t^n]
\l(4\sinh\l(\frac{\epsilon_1 t}{2}\r) \sinh\l(\frac{\epsilon_2 t}{2}\r)\r)^{-1} & n>0\,,
\end{cases}
\label{eq:betapol}
\eeq
so that e.g. 
$$\beta_2 = \frac{7 \epsilon _1^4+10 \epsilon _2^2 \epsilon _1^2+7 \epsilon _2^4}{5760 \epsilon _1 \epsilon _2}\,, \quad \beta_4 = -\frac{31 \epsilon _1^6+49 \epsilon _1^2 \epsilon _2^4+49 \epsilon _1^4
   \epsilon _2^2+31 \epsilon _2^6}{161280 \epsilon _1 \epsilon _2}\,, \quad \dots $$ and $\beta_{2n+1}=0$. The polynomials \eqref{eq:betapol} arise as the coefficients of the asymptotic expansion of the logarithm of the (shifted) Barnes double Gamma function 
   \cite{Nekrasov:2003rj,Brini:2010fc,MR2528052},
\begin{align*}
\log \Gamma_2\l(x+\frac{\epsilon_1+ \epsilon_2}{2},\epsilon_1,\epsilon_2\r) & \coloneqq  \frac{\dd}{\dd s} \l[\frac{1}{\Gamma(s)} \int_0^\infty \frac{\dd t}{t}\, \frac{ t^{s}\re^{-tx}}{(\re^{t\epsilon_1}-1)(\re^{t\epsilon_2}-1)}\r]\Bigg|_{s=0} \\
& =
\frac{3-2  \log (x)}{4 \epsilon _1 \epsilon _2}x^2
+\frac{\left(\left(\epsilon _1-1\right) \epsilon _1+\left(\epsilon _2-1\right) \epsilon
   _2\right) \log (x)+ \left(\epsilon _1+\epsilon _2\right)}{2 \epsilon _1 \epsilon _2}x \\
   &   +\frac{\left(\epsilon _1+\epsilon _2\right) \left(\epsilon _1 \left(6 \epsilon _1-5\right)+\epsilon _2 \left(6 \epsilon _2-5\right)\right) \log x}{24
   \epsilon _1 \epsilon _2}
+\kappa(\epsilon_1,\epsilon_2)+\sum_{n=1}^{\infty}  \beta_{n}(\epsilon_1,\epsilon_2) x^{-n}\,,
\end{align*}
where the constant term $\kappa(\epsilon_1,\epsilon_2)$, whose exact expression we will not need, is related to the derivative of the Riemann--Barnes double zeta function at the origin \cite{MR2528052}. Two notable cases are obtained when $\epsilon_2=-\epsilon_1$ and $\epsilon_1 \to 0$: for $n\geq 1$, we get respectively
\beq
\begin{split}
\beta_n(\epsilon,-\epsilon) &= (n-1)![t^n]
\frac{-1}{4\sinh^2(\epsilon_1 t/2)} =  \frac{B_{n+2} \epsilon^n}{n(n+2)}\,, \\
\lim_{\epsilon_1 \to 0} \epsilon_1 \beta_n(\epsilon_1,\epsilon_2) &= (n-1)![t^n]
\frac{1}{2\sinh(\epsilon_2 t/2)} = \frac{B_{n+2}(1/2)\epsilon_2^{n+1}}{n(n+1)(n+2)}\,.
\end{split}
\label{eq:betalimits}
\eeq
\subsubsection{Double-Gamma asymptotics from stable graph sums}
For $n>0$, the Laurent polynomials $\beta_{2n}(\epsilon_1, \epsilon_2)$ can be expressed in terms of a certain sum over stable graphs, which we will prove in \cref{sec:resconRmat}. Define
\beq
\begin{split}
\mathrm{Cont}_{\rm cf}(d) & \coloneqq  [z^{\mathsf{a}(d)-1}]\l(1- \Upsilon\l(-\frac{z}{\epsilon_1},-\frac{z}{\epsilon_2}\r)\l(1-\frac{z}{\epsilon_1+\epsilon_2}\r)^{-1}\r)\,, \\
\mathrm{Cont}_{\rm cf}(e) & \coloneqq  [z^{\mathsf{a}(h_1)} w^{\mathsf{a}(h_2)}]\l( \frac{1- \Upsilon\l(-\frac{z}{\epsilon_1},-\frac{z}{\epsilon_2}\r)\Upsilon\l(-\frac{w}{\epsilon_1},-\frac{w}{\epsilon_2}\r)}{z+w}\r)\,, \\
\mathrm{Cont}_{\rm cf}(v) & \coloneqq (-(\epsilon_1+\epsilon_2)^3\epsilon_1 \epsilon_2)^{\frac{2\mathsf{g}(v)-2+\mathsf{|H(v)|}}{2}} 
\int_{[\Mbar_{\mathsf{g}(v), |\HH(v)|}]}
\prod_{k=1}^{|\HH(v)|} 
\psi_k^{\mathsf{a}(h_k)}\,.
\end{split}
\label{eq:contbetagiv}
\eeq
The edge, vertices, and dilaton factors are identical to the corresponding terms appearing in the Givental--Teleman higher genus potentials of the special cubic Hodge-CohFT \cite{Okounkov:2003aok}, but for the appearance of a linear factor $(1-2z/\epsilon_+)$ rescaling the generating function for the dilaton leaf contributions. Deferring the proof to \cref{sec:resconRmat}, we claim that the following equality holds:
\begin{prop}[= \cref{prop:cfgraphsum}]
\beq 
\beta_{2g}(\epsilon_1, \epsilon_2)=
    \sum_{\Gamma \in \cG^{[0]}_{g,0}} \frac{1}{|\mathrm{Aut(\Gamma)}|} \prod_{e\in \cE_\Gamma} \mathrm{Cont}_{\rm cf}(e)\prod_{d\in \DD_\Gamma} 
\mathrm{Cont}_{\rm cf}(d) \prod_{v\in \cV_\Gamma}  \mathrm{Cont}_{\rm cf}(v) 
\,.
\label{eq:betagiv}
    \eeq
    \label{prop:cfgraphsummainbody}
\end{prop}
\cref{prop:cfgraphsummainbody} is a consequence of a closed-form formula for the all-order asymptotics of the R-matrix of the resolved conifold near its conifold point, which we derive in \cref{sec:resconRmat}, together with \cref{prop:ResConifold}.
\begin{rmk}
   The unrefined limit $\epsilon_1=-\epsilon_2$ of $\beta_{2g}(\epsilon_1, \epsilon_2)$ is well--known  to be related to the orbifold Euler characteristic of the moduli space of genus $g$ smooth curves by the Harer--Zagier formula \cite{MR848681},
    \[
    \chi(M_g) = \lambda^{2-2g} \beta_{2g-2}(\lambda, -\lambda) = \frac{B_{2g}}{2g(2g-2)}\,.
\]
The graph sum \eqref{eq:betagiv} then gives a new expression for $\chi(M_g)$ in terms of intersection numbers of tautological classes on $\Mbar_{g,n}$, for which formulas of similar flavour have recently appeared in \cite{MR4562996}. It would be interesting to trace the exact connection between those and \eqref{eq:betagiv}, and to see  how the refined case ($\epsilon_1\neq -\epsilon_2$) considered here could be related geometrically to cohomological information about $M_g$.

\end{rmk}

\subsubsection{Generators at the conifold point}
Let now
\[
t_{\rm cf} \coloneqq 
\frac{\sqrt{3}}{2\pi} I_{2,1}(y)\,, \quad
A_{\rm cf}(y) \coloneqq  \theta_y t_{\rm cf}\,, \quad F_{\rm cf} \coloneqq \theta_y \log A_{\rm cf}(y)\,.
\]
\noindent There is a homomorphism of rings,
\[
\mathsf{Cf} :  G_{\mgp{T}_{\rm F}}[A,F]  \longrightarrow   G_{\mgp{T}_{\rm F}}[A_{\rm cf},F_{\rm cf}]\,,
\]
induced by the map on generators
\[
\mathsf{Cf}(A)\coloneqq A_{\rm cf}\,, \quad  
\mathsf{Cf}(F)\coloneqq F_{\rm cf}
\,.
\]
In terms of quasi-modular forms \eqref{eq:ABCmod}, this map is induced by the change-of-variables
\[
\mathsf{Cf}(f)(\tau) \coloneqq f\l(-\frac{1}{3 \tau}\r)\,.
\]
We will also define
\[
D_{\rm cf} \coloneqq \mathsf{Cf}(D) = A_{\rm cf}^2 X\,, \quad E_{\rm cf}\coloneqq\mathsf{Cf}(E) = \frac{1}{2}(D_{\rm cf}-X)\,, 
\]
Finally, the conifold R-matrix and conifold refined generating series are defined as 
\[
\mathsf{R}^{\rm cf}_{ai} \coloneqq \mathsf{Cf} (\mathsf{R}_{ai})\,, \quad \GW_{g}^{\rm cf}=\mathsf{Cf}(\GW_{g}^{\mgp{T}_{\rm F}})\,, \quad \mathscr{F}^{\rm cf}_{k,g}=\mathsf{Cf}(\mathscr{F}_{k,g})\,.
\]
By \eqref{eq:gtH}, we have
\beq 
\GW^{\rm cf}_{g}
= \sum_{\Gamma \in \cG^{[2]}_{g,n}} \frac{(\epsilon_1 \epsilon_2)^{g-1}}{|\mathrm{Aut}(\Gamma)|} 
\mathsf{Cf}\l[\prod_{e\in \cE_\Gamma}\mathrm{Cont}(e)
\prod_{h\in \LL_\Gamma 
}\mathrm{Cont}(h) \prod_{v\in \cV_\Gamma} 
\mathrm{Cont}(v)\r]\,.
\label{eq:gtcf}
\eeq 
The refined conifold generating series have a singular behaviour as Laurent series in $t_{\rm cf}$ near $t_{\rm cf}=0$ \cite{Coates:2012vs}, with a pole of order $2g-2$:
\[
\GW^{\rm cf}_{g}\sim \frac{c_g(\epsilon_1,\epsilon_2)}{t_{\rm cf}^{2g-2}}+ \cO\l(\frac{1}{t_{\rm cf}^{2g-3}}\r)\,.
\]
We will determine the precise form of the leading term $c_g(\epsilon_1, \epsilon_2)$ in their Laurent expansion in the next section. 

\subsubsection{Oscillating integrals at the conifold point}
Let $\xi=1/X = 1+27 y$, and consider the asymptotic behaviour of the oscillating integrals \eqref{eq:oscintP} near the conifold point $\xi = 0$. When $i \neq 0$, the critical points \eqref{eq:qi}, the corresponding critical values, and the associated coefficients of the saddle-point expansion are analytic and non-zero at $\xi= 0$:
\[
\mathsf{Asym}_{q_i^{\rm cr}}[\re^{-\cW/z}\dd\Omega] = \cO(\xi^0)\,.
\]
%
%
On the other hand, the critical point $q^{\rm cr}_0 \to \infty$ has a simple pole at $\xi=0$, and the corresponding critical value $\mathsf{u}_0$ diverges logarithmically at $\xi=0$. This divergence similarly affects the coefficients of the Laplace expansion of \eqref{eq:oscintP} around $q^{\rm cr}_0$, and therefore all divergence in the Gromov--Witten potential at the conifold point stems from quantities arising as coefficients from the steepest descent expansion around the critical point $q_{\rm cr}^{[0]}$. In $s$-coordinates, this is the steepest descent integral
\beq 
\mathsf{Asym}_{s^{\rm cr}}[\re^{-\cW/z} \dd \Omega] = \frac{1}{\sqrt{\det(-h)}} \l[ \re^{\frac{z}{2} \sum_{m,n} h^{mn} \partial^2_{s_m s_n}} \re^{-(\cW^{(1)}_{\rm ng}+\cW^{(2)}_{\rm ng})/z} g(s) \r]\Bigg|_{s=s_{\rm cr}}\,,
\label{eq:asympRp2}
\eeq 
where
\[
s^{\rm cr} = \l(0,0, \frac{(\epsilon_1+\epsilon_2)y^{1/3}}{1+3 y^{1/3}}\r)\,, \quad 
g(s) = \frac{\sqrt{s_2^2 s_1^2-4 s_3  s_1^3+18 s_2 s_3 s_1-4 s_2^3-27 s_3^2}}{s_3^3}\,, 
\]
\[
h = \l(\bary{ccc} 0 & 3 & 0 \\ 3 & 0 & 0 \\
 0 & 0 & -\frac{1+ 6 y^{1/3}+9 y^{2/3}}{(\epsilon _1+\epsilon _2)y^{2/3}}
\eary \r)\,, \quad h^{mn} \coloneqq (h^{-1})_{mn}\,, \quad
\cW^{(1)}_{\rm ng} = s_1^3\,, 
\]
\bea
\cW^{(2)}_{\rm ng} &=& 
s_3
   \left(\frac{1}{y^{1/3}}+3\right)+\left(\epsilon _1+\epsilon _2\right) \log \left(s_3\right) - \left(\epsilon _1+\epsilon _2\right) \log \left(-\frac{y^{1/3} \left(\epsilon
   _1+\epsilon _2\right)}{3 y^{1/3}+1}-1\right) \nn \\ &+& \frac{\left(s_3 \left(3 y^{1/3}+1\right)+y^{1/3} \left(\epsilon _1+\epsilon _2\right)\right){}^2}{2 y^{2/3} \left(\epsilon _1+\epsilon _2\right)}
\,. \nn
\eea
 
%
%
%

%
%
The differential operators obtained from the terms $(m,n)=(1,2)$, $(2,1)$ and $(3,3)$ in the Feynman sum \eqref{eq:asympRp2} commute, so we will apply them separately in turn. For $m+n=3$, we have that 
\[
\l[\re^{\frac{z\de^2_{s_1 s_2}}{3}} \re^{-\cW^{(1)}_{\rm ng}/z} g(s) \r]\Bigg|_{s_1=s_2=0} = 
\l[\re^{-\cW^{(1)}_{\rm ng}/z} \re^{\frac{z\de^2_{s_1 s_2}}{3}} g(s) \r]\Bigg|_{s_1=s_2=0} = \sqrt{3}\ri\sum_{l \geq 0} \frac{r_l}{s_3^{l+2}}\,,
\]
where $r_0=-3$, $r_1=1$, $r_2=\frac{4}{9}$, $\dots\,$. Then, 
\[
\mathsf{Asym}_{s^{\rm cr}}[\re^{-\cW/z} \dd \Omega] = \sqrt \frac{-1}{\det h}  \sum_{l \geq 0} r_l 
\l[\re^{ \frac{z\de^2_{s_3}}{2}}\re^{-\cW^{(2)}_{\rm ng}/z} s_3^{-2-l} \r]
\Bigg|_{s_3=\frac{y^{1/3}(\epsilon_1+\epsilon_2) }{1+3 y^{1/3}}}\,.
\]
The prefactor on the r.h.s.~has a simple pole at $\xi=0$, while the summands in the second factor have zeroes of order $l+2$. Therefore the  r.h.s.~vanishes linearly as $\xi\to 0$, with leading coefficient given by the term $l=0$. Notice now that
\[ 
\sqrt \frac{3}{\det h} 
\l[\re^{ \frac{z\de^2_{s_3}}{2}}\re^{-\cW^{(2)}_{\rm ng}/z} s_3^{-2} \r]
\Bigg|_{s_3=-\frac{y^{1/3}(\epsilon_1+\epsilon_2) }{1+3 y^{1/3}}}
= \mathsf{Asym}_{\frac{b}{a}}
[\re^{-ax/z} x^b \dd x]\,,
\] 
where $a=3+y^{-1/3}$ and $b=-(\epsilon_1+\epsilon_2)/z-2$. The r.h.s.~is then the asymptotic tail of the Gamma integral 
\bea
\int_0^\infty \re^{-(3+y^{-1/3})t/z} t^{-(\epsilon_1+\epsilon_2)/z-2} \dd t &=&
\l(\frac{3+y^{-1/3}}{z}\r)^{1+(\epsilon_1+\epsilon_2)/z}
\Gamma\l(-\frac{\epsilon_1+\epsilon_2}{z}-1\r)\,,\nn 
\eea
so that
\[ 
\mathsf{Asym}_{s^{\rm cr}}[\re^{-\cW/z} \dd \Omega] =
-\frac{\xi}{\left(3(\epsilon _1+\epsilon _2)\right)^{1/2}}\frac{\Theta\l(-\frac{z}{\epsilon_1+\epsilon_2} \r)}{\epsilon_1+\epsilon_2-z} + \cO(\xi)^2\,.
\] 
Similarly, we get
\bea 
\mathsf{Asym}_{s^{\rm cr}}[\re^{-\cW/z} s_3 \dd \Omega] &=&
-\sqrt{\frac{1}{3\left(\epsilon _1+\epsilon _2\right)}} \Theta\l(-\frac{z}{\epsilon_1+\epsilon_2} \r) + \cO(\xi) \,, \nn \\
\mathsf{Asym}_{s^{\rm cr}}[\re^{-\cW/z} s_3^2 \dd \Omega] &=&
-\sqrt{\frac{\epsilon _1+\epsilon _2}{3}} \Theta\l(-\frac{z}{\epsilon_1+\epsilon_2} \r)~\frac{1}{\xi} + \cO(1) \,. \nn
\eea 
\subsubsection{Givental's reconstruction at the conifold point}
From the asymptotic analysis of the previous Section, we deduce the following
\begin{thm}[Conifold asymptotics]
For $g>1$, the conifold refined Gromov--Witten potentials of $K_{\bbP^2}$ satisfy
    \[ 
\GW^{\rm cf}_{g} = \frac{(-3)^{g-1}\beta_{2g}(\epsilon_1, \epsilon_2)}{t_{\rm cf}^{2g-2}} +\cO\l(\frac{1}{t_{\rm cf}^{2g-3}}\r)\,.
    \] 
In particular,
\[
\mathscr{F}^{\rm cf}_{0,g} = \frac{(-3)^{g-1} B_{2g}}{2g(2g-2) t_{\rm cf}^{2g-2}} + \cO\l(\frac{1}{t_{\rm cf}^{2g-3}}\r)\,, \quad \mathscr{F}^{\rm cf}_{g,0} = \frac{(-3)^{g-1} B_{2g+2}(1/2)}{2g(2g+1)(2g+2) t_{\rm cf}^{2g-2}}+ \cO\l(\frac{1}{t_{\rm cf}^{2g-3}}\r)\,.
\]
\label{thm:cfasym}
\end{thm}
\begin{proof}
From \eqref{eq:RHb}, \eqref{eq:normNLR}, and $\mathsf{R}^{\rm cf}_{ai} = \mathsf{Cf} (\mathsf{R}_{ai})$, we have
\bea 
\mathsf{Asym}_{s^{\rm cr}}[\re^{-\cW/z} \dd \Omega]
&=& \mathsf{R}^{\rm cf}_{00} N_0\,, \nn \\
\mathsf{Asym}_{s^{\rm cr}}[s_3\re^{-\cW/z} \dd \Omega] &=&
A_{\rm cf}(y) \mathsf{R}^{\rm cf}_{10} N_0\,, \nn \\
\mathsf{Asym}_{s^{\rm cr}}[s_3^2 \re^{-\cW/z} \dd \Omega] &=& 
D_{\rm cf}(y) (\mathsf{R}^{\rm cf}_{20} + E_{\rm cf}(y) \mathsf{R}^{\rm cf}_{10}) N_0\,. \nn 
\eea 
Since
\[
A_{\rm cf}(y) = -1 + \cO(\xi)\,, \quad 
D_{\rm cf}(y) = \frac{1}{\xi} + \cO(1)\,, \quad 
E_{\rm cf}(y) = -\frac{2}{9}(\epsilon_1+\epsilon_2) + \cO(\xi)\,, 
\] 
we then get
\beq
\begin{split}
\mathsf{R}^{\rm cf}_{00} &=-\frac{\xi}{\left(3\epsilon_1 \epsilon_2(\epsilon _1+\epsilon _2)\right)^{1/2}}\frac{\Upsilon\l(\frac{z}{\epsilon_1},\frac{z}{\epsilon_2} \r)}{\epsilon_1+\epsilon_2-z} + \cO(\xi)^2\,, \nn \\
\mathsf{R}^{\rm cf}_{10} &=
-\sqrt{\frac{1}{3\epsilon_1 \epsilon_2\left(\epsilon _1+\epsilon _2\right)}} \Upsilon\l(\frac{z}{\epsilon_1},\frac{z}{\epsilon_2} \r) + \cO(\xi) \,, \nn \\
\mathsf{R}^{\rm cf}_{20} &=
\sqrt{\frac{\epsilon _1+\epsilon _2}{3\epsilon_1 \epsilon_2}}
~\Upsilon\l(\frac{z}{\epsilon_1},\frac{z}{\epsilon_2} \r) + \cO(\xi) \,.
\end{split}
\label{eq:Rcf}
\eeq 
Moreover, we have
\[
Y^{\rm cf}_0 = \theta_y \mathsf{u}_0 = \frac{\epsilon _1+\epsilon _2}{3 \left((1-\xi)^{1/3}-1\right)}=\left(\epsilon _1+\epsilon _2\right)\l(-\frac{1}{\xi }+\frac{1}{3} +\cO\left(\xi \right)\r)\,,
\]
so that
\beq 
(\Delta_0^{\rm cf})^{1/2} \coloneqq
\mathsf{Cf}(\Delta_0^{1/2}) = \sqrt{3 \epsilon_1 \epsilon_2} \mathsf{Cf}(W_0^{1/2}) = -\left(3\epsilon _1\epsilon _2(\epsilon _1+\epsilon _2)^3\right)^{1/2}\l(\frac{1}{\xi }-\frac{2}{3}+\cO\left(\xi\right)\r)\,.
\label{eq:deltacf}
\eeq

We are now ready to compute the dominant contribution to the graph sum \eqref{eq:gtcf} as $\xi \to 0$. Since 
$\mathsf{R}^{\rm cf}_{ai} = \cO(\xi)$ for $i\neq 0$, at the leading order in $t_{\rm cf}$ 
there are no contributions from stable graphs with decoration other than $i=0$.
We have:

\[ 
\GW^{\rm cf}_{g} = \frac{1}{t_{\rm cf}^{2g-2}}\sum_{\Gamma \in \cG^{[0]}_{g,0}} \frac{1}{|\mathrm{Aut(\Gamma)}|} \prod_{e\in \cE_\Gamma} \mathrm{Cont}_{\rm cf}(e)\prod_{d\in \DD_\Gamma} \mathrm{Cont}_{\rm cf}(d) \prod_{v\in \cV_\Gamma} \mathrm{Cont}_{\rm cf}(v) + \cO\l(\frac{1}{t_{\rm cf}^{2g-3}}\r)\,,
\] 
where, from \eqref{eq:Rcf}--\eqref{eq:deltacf},
\begin{align}
\mathrm{Cont}_{\rm cf}(e)  \coloneqq &  [z^{\mathsf{a}(h_1)} w^{\mathsf{a}(h_2)}]\lim_{\xi\to 0}\l( \frac{1- \mathsf{R}^{\rm cf}_{a0}(-z)\eta^{ab}\mathsf{R}^{\rm cf}_{b0}(-w)}{z+w}\r) 
=[z^{\mathsf{a}(h_1)} w^{\mathsf{a}(h_2)}]
\l( \frac{1- \Upsilon\l(-\frac{z}{\epsilon_1},-\frac{z}{\epsilon_2}\r)\Upsilon\l(-\frac{w}{\epsilon_1},-\frac{w}{\epsilon_2}\r)}{z+w}\r)\,, \nn \\
\mathrm{Cont}_{\rm cf}(d)  \coloneqq  & [z^{\mathsf{a}(d)-1} ]\lim_{\xi\to 0}  \l(- (\Delta^{\rm cf}_0)^{1/2} \mathsf{R}^{\rm cf}_{00}(-z)\r) 
= [z^{\mathsf{a}(d)-1} ]
\l(1- \Upsilon\l(-\frac{z}{\epsilon_1},-\frac{z}{\epsilon_2}\r)\l(1-\frac{z}{\epsilon_1+\epsilon_2}\r)^{-1}\r)\,, \nn \\
\mathrm{Cont}_{\rm cf}(v) \coloneqq &  (-\sqrt{3(\epsilon_1+\epsilon_2)^3\epsilon_1 \epsilon_2})^{2\mathsf{g}(v)-2+\mathsf{|H(v)|}} 
\int_{[\Mbar_{\mathsf{g}(v), |\HH(v)|}]} \prod_{k=1}^{|\HH(v)|} \psi_k^{\mathsf{a}(h_k)}\,, \nn
\end{align}
where in the first two equalities we have used \eqref{eq:Rcf}.
    The first equality in the Theorem is then immediate from \eqref{eq:contbetagiv}--\eqref{eq:betagiv}. The last two equalities in the statement of the Theorem are respectively the unrefined ($\epsilon_1 \to - \epsilon_2$) and NS ($\epsilon_1 \to 0$) limits of the first, and follow from \eqref{eq:betalimits}.
\end{proof}

\begin{rmk}
The two special limits in   \cref{thm:cfasym} prove the leading order of the conifold asymptotics conjecture in \cite{HK10:OmegaBG,Coates:2018hms} and \cite[Conj.~1.9]{BFGW21:HAE}. A proof of the full refined conifold gap conjecture,
    \[ 
\GW^{\rm cf}_{g} = \frac{(-3)^{g-1}\beta_{2g}(\epsilon_1, \epsilon_2)}{t_{\rm cf}^{2g-2}} +\cO\l(1\r)\,,
    \] 
appears to require a significant leap in complexity. We have however verified that this holds for low genus: for example, from \eqref{eq:F2ref}, one computes
    \[ 
\GW_2^{\rm cf} =  -\frac{7 \epsilon _1^4+10 \epsilon _2^2 \epsilon _1^2+7 \epsilon _2^4}{1920 \epsilon _1 \epsilon _2}\frac{1}{t_{\rm cf}^{2}} + \cO\l(1\r)\,.
    \] 
\label{rmk:cfgap}
\end{rmk}

\section{Comparison with refined GV/PT theory}
\label{sec:further implications}

\subsection{Refined BPS integrality}
\label{sec: refined GV}
The Gopakumar--Vafa formulation of the refined BPS partition function of a Calabi--Yau threefold $X$ predicts that the invariants $N^{\mathrlap{\beta}}{\mathstrut}_{j_-,j_+}$ in \eqref{eq:GVref} are positive integers \cite{GV98:MthTopStrI,GV98:MthTopStrII}. In this section we will formulate this expectation as precise conjectures in increasing levels of strength.

\subsubsection{Rationality}
\label{sec: rationality conjecture}
Recall the setup and notation from \Cref{sec: refined threefold setup}. We take $Z=X\times \Aaff{2}$ and assume that $\mgp{T}$ acts on the two affine factors via characters $\mchar{q}_1$, $\mchar{q}_2$; assuming that the action on $Z$ is Calabi--Yau, we have
\begin{equation*}
    \kappa = \mchar{q}_1 \mchar{q}_2\,.
\end{equation*}
In the following we will often assume that this $\mgp{T}$-character admits a square root, which we will denote by $\mchar{q}_+ \in K_{\mgp{T}}(\mr{pt})$, and we shall additionally define $\mchar{q}_- \coloneqq  \mchar{q}_+ \mchar{q}_2^{-1} = \mchar{q}_+^{-1} \mchar{q}_1$. The notation is chosen so that
\begin{equation*}
    \epsilon_i = c_1^{\mgp{T}} (\mchar{q}_i)\,,\qquad i \in \{1,2,+,-\}\,.
\end{equation*}
which means that the Chern character $\ch[\mgp{T}] : K_{\mgp{T}} (\mr{pt}) \rightarrow \widehat{R}_{\mgp{T}}$ acts via
\begin{equation*}
    \ch[\mgp{T}] \, \mchar{q}_i = \re^{\epsilon_i}\,.
\end{equation*}

We start with the weakest and most general of our three conjectures, which states that the GW generating series admits a rational lift under this last map.

\begin{conj}
	\label{conj: rational lift}
	Suppose the $\mgp{G}$-action on $Z$ is Calabi--Yau and $\kappa$ admits a square root. Then, for every non-zero effective curve class $\beta$, there exists a rational function
    \begin{equation*}
        \refBPSbarb{Z}{\mgp{G}}{\beta}  \in \mr{Frac}\big( K_{\mgp{G}} (\mr{pt}) \big)
    \end{equation*}
    satisfying
	\begin{equation*}
		\refGWb{Z}{\mgp{G}}{\beta} = \ch[\mgp{G}] \, \refBPSbarb{Z}{\mgp{G}}{\beta}\,.
	\end{equation*}
\end{conj}


\begin{rmk}
    Suppose the above conjecture holds. Then, as is the case for the cohomological GW generating series in \cref{lem: GW contravariance}, the K-theoretic lift $\refBPSbarb{Z}{\mgp{G}}{\beta}$ is contravariant in $\mgp{G}$. Moreover, since $\refGWb{Z}{\mgp{G}}{\beta}$ is an element of the evenly graded part of $\widehat{R}^{\,\mr{loc}}_{\mgp{G}}$ (\Cref{defn: refined GW invariant}), its K-theoretic lift would necessarily be invariant under the dualisation map $\mchar{t} \mapsto \mchar{t}^{\vee} = \mchar{t}^{-1}$.
\end{rmk}

	We formulated the above conjecture for targets $Z$ of the form $X\times \Aaff{2}$, where we can provide solid evidence in its favour. An adventurous reader might however view the statement of \Cref{conj: rational lift} as a speculation about the case of general CY 5-folds.

\subsubsection{Evidence for \cref{conj: rational lift}}

	We will start by the running example of the resolved conifold $X=\Tot \l(\mc{O}_{\bbP^1}(-1) \oplus \mc{O}_{\bbP^1}(-1)\r)$. For (a covering of) a torus action satisfying the conditions of \Cref{prop:ResConifold}, we can reformulate the statement of the \namecref{prop:ResConifold} as
	\begin{equation*}
		\refGWb{Z}{\mgp{T}}{d [\bbP^1]} = \ch[\mgp{T}] \left(\frac{1}{d} \frac{\mchar{q}_+^d}{(1-\mchar{q}_1^d)(1-\mchar{q}_2^d)} \right)\,,
	\end{equation*}
 verifying \cref{conj: rational lift}. Analogously, the Conjecture is seen to be true for the target $X=\Tot(\mc{O}_{\bbP^1}(-2) \oplus \mc{O}_{\bbP^1})$ under all torus actions considered in \Cref{prop: shifted local P1}.

We can further corroborate the evidence for \Cref{conj: rational lift} in the context of local surfaces. Let $S$ be a del Pezzo surface and consider the fibrewise scaling action by $\mgp{T}_{\mr{F}}=\Gm[2]$ on $Z = K_S \times \Aaff{2}$ from \Cref{sec: local surface prelim}. One checks that this action satisfies all conditions listed in \Cref{lem: shape of refGW series 2}, hence the all-genus Gromov--Witten generating series is regular up to a leading factor $(\epsilon_1 \epsilon_2)^{-1}$. Thus, if \Cref{conj: rational lift} holds, the following limit must exist:
\begin{equation*}
	\epsilon_2 \, \refGWb{K_S \times \Aaff{2}}{\mgp{T}_{\mr{F}}}{\beta} ~\Big\lvert_{\epsilon_2 =0} = \ch[\Gm] \, (1-\mchar{q}_2)\cdot \refBPSbarb{K_S \times \Aaff{2}}{\mgp{T}_{\mr{F}}}{\beta}~\Big\lvert_{\mchar{q}_2 =1}\,,
\end{equation*}
where $\ch[\Gm] \mchar{q}_1 = \re^{\epsilon_1}$; in particular, the left-hand side lifts to a rational function in $\mchar{q}_1$. We can show that this is indeed the case.

\begin{prop}
	\label{prop: rationality NS limit}
	Let $X=K_S$ be a local del Pezzo surface together with (a cover of the) $\mgp{T}_{\mr{F}}$-action considered in \Cref{sec: local surface prelim}. Then \Cref{conj: rational lift} holds in the Nekrasov--Shatashvili limit ($\epsilon_2=0$).
\end{prop}
\begin{proof}
	By \cite[Lem.~8.15]{Bou20:QuTropVert}, the generating series of relative Gromov--Witten invariants of $S$ with maximum tangency along a smooth anticanonical curve $D$ has a rational lift: there exists a rational function $\overline{\BPS}_{\beta}(S/ D)\in \bbQ(\mchar{q}_1^{1/2})$ satisfying
	\begin{equation*}
		\GW_{\beta}(S/ D) = \ch[\Gm] \, \overline{\BPS}_{\beta}(S/ D) \,,
	\end{equation*}
	where on the l.h.s.~we use the notation introduced in \Cref{sec: notation relative GW}. Thus, by \Cref{thm: NS vs log num}, we have
	\begin{equation*}
		 \epsilon_2 \, \refGWb{Z}{\mgp{T}_{\mr{F}}}{\beta} ~\Big\lvert_{\epsilon_2 =0} = \ch[\Gm] \, \overline{\BPS}_{\beta}(S/ D)\,. \qedhere
	\end{equation*}
\end{proof}

\begin{rmk}
    Note that the r.h.s.~of the last equation has generally half-integer powers of $\mchar{q}_1$, so one needs to transition to a cover of $\mgp{T}_{\mr{F}}$ on which the square root is defined for $\overline{\Omega}_{\beta}(S/ D)$ to admit an interpretation as an element in K-theory.
\end{rmk}

Further evidence for \Cref{conj: rational lift} is discussed below in the context of stronger conjectures.

\subsubsection{Integrality}
\label{sec: integrality conjecture}
In \cref{sec: rationality conjecture} we did not impose a rigidity assumption on curve classes: adding that strongly constrains the values and pole structure of $\refBPSbarb{Z}{\mgp{G}}{\beta}$ in terms of a set of integer-valued Laurent polynomials. For convenience, we fix a saturated additive subset $\hhh \subseteq \hhh_2(X,\bbZ)^+$ of rigid curve classes in $X$: for instance, when $X$ is a local del Pezzo surface, we may take $\hhh = \hhh_2(X,\bbZ)^+$. Let us assume that \Cref{conj: rational lift} holds for all curve classes in $\hhh$. Then we can recursively define
\begin{equation*}
	\refBPSb{Z}{\mgp{G}}{\beta}(\mathbf{t}) \in \mr{Frac}\big( K_{\mgp{G}} (\mr{pt}) \big)
\end{equation*}
for all $\beta \in \hhh$ by demanding 
\begin{equation*}
	\refBPSbarb{Z}{\mgp{G}}{\beta}(\mathbf{t}) = \sum_{k|\beta} \frac{1}{k[\mchar{q}_1^{k}][\mchar{q}_2^{k}]} \, \refBPSb{Z}{\mgp{G}}{\beta/k}(\mathbf{t}^k) \,,
\end{equation*}
i.e.
\begin{equation*}
	\refBPSb{Z}{\mgp{G}}{\beta}(\mathbf{t}) = \frac{1}{[\mchar{q}_1][\mchar{q}_2]} \sum_{k|\beta} \frac{\mu(k)}{k} \, \refBPSbarb{Z}{\mgp{G}}{\beta/k} (\mathbf{t}^k)\,,
\end{equation*}
where we write $[\mchar{q}]\coloneqq \mchar{q}^{1/2}-\mchar{q}^{-1/2}$ and $\mu$ is the M\"{o}bius function. In the above two equations, we denoted by $\Omega(\mathbf{t}^k)$ the image of the $k^{\rm th}$ Adams operation on an element $\Omega \in K_{\mgp{G}}(\mr{pt})$. The image of a one-dimensional $\mgp{T}$-representation $\mchar{t}$ under this operation is $\mchar{t}^k$.

\begin{conj}[Refined BPS integrality I]
	\label{conj: refined GV integrality}
	Suppose the $\mgp{G}$-action on $Z$ is Calabi--Yau and $\kappa$ admits a square root. Then for all $\beta\in\hhh$ the rational function $\refBPSb{Z}{\mgp{G}}{\beta}$ is an integral Laurent polynomial in $\mchar{q}_+,\mchar{q}_-$:
	\begin{equation*}
		\refBPSb{Z}{\mgp{G}}{\beta}\in \bbZ[\mchar{q}_+^{\pm 1},\mchar{q}_-^{\pm 1}]  \subsetneq \mr{Frac}\big( K_{\mgp{G}} (\mr{pt}) \big) \,.
	\end{equation*}
\end{conj}

\begin{rmk}
    \label{rmk: Kth lift}
	We remark that it is important to assume the curve classes to be rigid in order to guarantee that $\refBPSbarb{Z}{\mgp{G}}{\beta}$ only features denominators of the form
	\begin{equation*}
		\big(1 - \mchar{q}_1^{k_1} \big) \, \big(1 - \mchar{q}_2^{k_2}\big)
	\end{equation*}
	for $k_i\in \bbZ$: a counter-example is once again the local $\bbP^1$ geometry of \Cref{sec: local P1 shifted}. But even when the curve class is non-rigid, the denominators of $\refBPSbarb{Z}{\mgp{G}}{\beta}$ appear to have a very constrained form: in all our examples, they are given by products of factors $(1-\mchar{t})$ for $\mchar{t} \in K_{\mgp{T}}(\mr{pt}) = \Rep \mgp{T}$ a one-dimensional representation. In other words one may wonder if in the non-rigid setup
	\begin{equation*}
		\refBPSbarb{Z}{\mgp{G}}{\beta} \in K_{\mgp{G}}(\mr{pt})^{\mr{loc}} \otimes \bbQ
	\end{equation*}
	where $K_{\mgp{G}}(\mr{pt})^{\mr{loc}}$ is the localisation of $K_{\mgp{G}}(\mr{pt})$ at the augmentation ideal, where the numerators have an integral structure similar to \Cref{conj: refined GV integrality}: cf \Cref{sec: local P1 shifted}.
\end{rmk}

Let us reformulate \Cref{conj: refined GV integrality} in a form closer to conjectures appearing in the physics literature. 
For this, note that, as formulated, \Cref{conj: refined GV integrality} implies the Rigidity \cref{conj: rigidity}
\begin{equation*}
	\epsilon_1\epsilon_2 ~ \refGWb{Z}{\mgp{G}}{\beta} \in \bbQ \llbracket \epsilon_+^2,\epsilon_-^2 \rrbracket\,,
\end{equation*}
and in particular $\refBPSb{Z}{\mgp{G}}{\beta}$ has to be a palindromic polynomial in $\mchar{q}_\pm$. We may then expand it as
\begin{equation*}
	\refBPSb{Z}{\mgp{G}}{\beta} \eqqcolon \sum_{i,j \geq 0}  n_{i,j}^{\beta} \,(-1)^{i+j} [\mchar{q}_+]^{2i} \, [\mchar{q}_-]^{2j}
\end{equation*}
for some $n^{\mathrlap{\beta}}{\mathstrut}_{i,j} \in \bbZ$, which are non-zero for only finitely many $i,j\in \bbZ_{\geq 0}$. Using that
\begin{equation*}
	\ch[\mgp{T}] [\mchar{q}_i^k] = 2 \sinh \tfrac{k \epsilon_i}{2}\,.
\end{equation*}
we find that \Cref{conj: refined GV integrality} is equivalent to the following

\begin{conj}[Refined BPS integrality II]
	\label{conj: refined GV integrality II}
	Suppose that the $\mgp{G}$-action on $Z$ is Calabi--Yau and that $\kappa$ admits a square root. Then, there exist $n^{\mathrlap{\beta}}{\mathstrut}_{i,j}\in \bbZ$  $(i,j\in\bbZ_{\geq 0}, \,\beta\in \hhh)$ such that
	\begin{equation}
		\label{eq: ref BPS expansion}
		\refGWb{Z}{\mgp{G}}{\beta} = \sum_{k|\beta} \sum_{i,j\geq 0} \frac{n_{i,j}^{\beta/k}}{k} \, (-1)^{i+j} \frac{\left(2 \sinh \tfrac{k \epsilon_+}{2}\right)^{2i}\left(2 \sinh \tfrac{k \epsilon_-}{2}\right)^{2j}}{2 \sinh \tfrac{k \epsilon_1}{2} ~~ 2 \sinh \tfrac{k \epsilon_2}{2}}\,.
	\end{equation}
	Moreover, for fixed $\beta$ we have $n^{\mathrlap{\beta}}{\mathstrut}_{i,j} \neq 0$ only for finitely many $i,j\in\bbZ_{\geq 0}$.
\end{conj}

\subsubsection{Evidence for \cref{conj: refined GV integrality,conj: refined GV integrality II}}
\label{sec:refGVevidence}

The first sanity test is again $X=\Tot(\mc{O}_{\bbP^1}(-1) \oplus \mc{O}_{\bbP^1}(-1))$. \Cref{prop:ResConifold} shows that for all $\mgp{T}$-actions specified in the \namecref{prop:ResConifold} there is a single non-zero BPS invariant
	\begin{equation*}
		n_{0,0}^{[\bbP^1]}=1\,.
	\end{equation*}
	In particular, the conjectures hold for these torus actions.

A second test comes from the unrefined limit. Note that from the expansion \eqref{eq: ref BPS expansion} one can recover the original integrality conjecture of \cite{GV98:MthTopStrI,GV98:MthTopStrII}: taking $\epsilon_1= -\epsilon_2 = \epsilon_-$ on both sides of the equation gives
	\begin{equation}
		\label{eq: notation original GV}
		\sum_{g\geq 0} (-\epsilon_-^2)^{g-1} \int_{[\Mbar_g(X,\beta)]^{\vir}}1 = \sum_{k|\beta} \sum_{g\geq 0}\frac{n_{0,g}^{\beta/k}}{k} \, (-1)^{g-1}\left(2 \sinh \tfrac{k \epsilon_-}{2}\right)^{2g-2}
	\end{equation}
	where in order to get the right-hand side we applied \Cref{prop: unref limit}. Gopakumar and Vafa predicted that the coefficients $n^{\mathrlap{\beta}}{\mathstrut}_{0,g}$ appearing in \eqref{eq: notation original GV} are all integers, and furthermore vanish for $g\gg 0$ : these two statements were proven in \cite{IP18:GV} and \cite{DIW21:GVfinite} respectively. Therefore \Cref{conj: refined GV integrality II} (and thus also \Cref{conj: refined GV integrality}) are known to hold in the unrefined limit.
	
	\label{rmk: direct integration BPS}
Further evidence for \Cref{conj: refined GV integrality II} comes from refined mirror symmetry. The expansion \eqref{eq: ref BPS expansion} had previously appeared in the physics literature on refined topological strings \cite[Equation~(3.21)]{HK10:OmegaBG} (see also \cite{CKK14:refBPS,HKK13:OmegaBmodelRigidN2}). 
 In \cite{HKK13:OmegaBmodelRigidN2}, the authors use the refined direct integration method of \Cref{def:dirint}
 to verify
 \Cref{conj: refined GV integrality II} for curve classes of degree up to nine and up to genus $84$. 
Since in \Cref{thm:haeref,thm:crc,thm:cfasym}
we proved 
the refined direct integration constraints up to the  conjectural gap in the conifold expansion, we deduce the following statement.

\begin{prop}[\cite{HKK13:OmegaBmodelRigidN2}]
	Suppose that {(\bf Conifold asymptotics)} holds as formulated in \cref{def:dirint}. Then, \Cref{conj: refined GV integrality II} holds for $K_{\bbP^2}\times \Aaff{2}$ with the fibrewise $\mgp{T}_{\mr{F}}$-action up to curve class degree nine, and modulo contributions from stable maps with $g>84$.
 \label{prop:KP2refintcheck}
\end{prop}

\subsubsection{Geometric interpretation}
\label{sec: geometric interpretation conjecture}
We now restrict the discussion to local del Pezzo surfaces together with the fibrewise $\mgp{T}_{\mr{F}}$-action considered in \Cref{sec: local surface prelim}. We present a conjecture for a geometric description of the coefficients of $\refBPSb{K_S\times \Aaff{2}}{\mgp{T}_{\mr{F}}}{\beta}$.

Let $S$ be a smooth del Pezzo surface, and denote by $M_{\beta}$ the moduli space of one dimensional Gieseker semi-stable sheaves on $S$ (with respect to some fixed polarisation) with support $\beta$ and Euler characteristic one. The Hilbert--Chow morphism
\begin{equation*}
	\pi : M_{\beta} \longrightarrow \mr{Chow}_{\beta}(S)
\end{equation*}
induces a perverse filtration on the cohomology groups of $M_{\beta}$:
\begin{equation*}
	\text{$\hhh$}^{i,j} \coloneqq \mathbb{H}^i \!\left({}^{\mf{p}}\text{$\mathcal{H}$}^j \big({\mathbf{R}}\pi_*\bbQ_{M_{\beta}} [\dim M_{\beta}]\big)\right)\,.
\end{equation*}
We consider a refinement of Maulik and Toda's proposal \cite[Def.~1.1]{MT18:GV} for Gopakumar--Vafa invariants by forming a generating series recording the dimension of these cohomology groups
\begin{equation}
	\refMTb{K_S}{\beta} (\mchar{q}_+,\mchar{q}_-)\coloneqq \sum_{i,j\in \bbZ} \dim \hhh^{i,j} ~ (-\mchar{q}_+)^i \, (-\mchar{q}_-)^j  \quad\in \bbZ[\mchar{q}_+^{\pm 1},\mchar{q}_-^{\pm 1}]\,.
\label{eq:refMTb}
\end{equation}
This refinement has also been investigated
in \cite{KPS23:refBPSP2,KLMP24:CohOneDimSheavesP2} in connection with K-theoretic Pandharipande--Thomas theory. We expect the following relation to Gromov--Witten theory.

\begin{conj}
	\label{conj: MT type BPS local surfaces}
	Suppose $S$ is a del Pezzo surface. Then \Cref{conj: refined GV integrality} holds for $Z=K_S \times \Aaff{2}$  with the fibrewise $\mgp{T}_{\mr{F}}$-action. Moreover, for all $\beta\in \hhh_2 (S,\bbZ)^+$, we have
	\begin{equation*}
		\refBPSb{Z}{\mgp{T}_{\mr{F}}}{\beta} = \refMTb{K_S}{\beta}\,.
	\end{equation*}
\end{conj}
\begin{rmk}
	In our definition of $\refMTb{K_S}{\beta}$ the variables $\mchar{q}_+$, $\mchar{q}_-$ do not enjoy a natural geometric interpretation, which makes the definition rather ad hoc. It therefore seems desirable to find an interpretation of $\refMTb{K_S}{\beta}$ as an honest equivariant index. Such an interpretation also appears necessary in order to generalise \Cref{conj: MT type BPS local surfaces} to geometries other than local surfaces and may provide a lead for the construction of a moduli space of M2 branes \cite{NO14:membranes}.
\end{rmk}

\subsubsection{Evidence for \cref{conj: MT type BPS local surfaces}}
We earlier noticed that, as a consequence of \Cref{lem: shape of refGW series 2}, the function $\refBPSb{Z}{\mgp{T}_{\mr{F}}}{\beta}$ is a palindromic Laurent polynomial in $\mchar{q}_\pm$.
A first check of \cref{conj: MT type BPS local surfaces} is to see if the symmetry $\mchar{q}_\pm \leftrightarrow \mchar{q}_\pm^{-1} $ is  respected by
$\refMTb{K_S}{\beta}$: 
indeed, 
since $M_{\beta}$ is smooth, we have for all $i,j\in \bbZ$ that
	\begin{equation*}
		\hhh^{i,j}\cong \hhh^{-i,j} \qquad \text{and} \qquad \hhh^{i,j}\cong \hhh^{i,-j}
	\end{equation*}
	by the hard Lefschetz theorem for perverse cohomology groups \cite[Thm~2.1.4]{dCMM05:HodgeThAlgMaps}.

A second obvious check comes, again, from the unrefined limit: in this case, the evidence only translates to the compatibility with earlier conjectures.	Restricting \Cref{conj: MT type BPS local surfaces} to the anti-diagonal torus $\Gm_{\mchar{q}_-} \hookrightarrow \mgp{T}_{\mr{F}}$, i.e.~setting $\mchar{q}_+ = 1$, we recover a special instance of a conjecture of Maulik and Toda \cite[Conj.~3.18]{MT18:GV} predicting that the integers $n^{\beta,\MT}_{0,g} \in \bbZ$ defined via
\begin{equation*}
	\sum_{g\geq 0} n_{0,g}^{\beta,\MT} (-1)^{g+1} [\mchar{q}_-]^{2g} \coloneqq \refMTb{K_S}{\beta} (1, \mchar{q}_-)
\end{equation*}
satisfy equation \eqref{eq: notation original GV}. We remark that, for local surfaces, the BPS sheaf agrees with the IC sheaf as proven by Meinhardt in \cite[Thm~1.1]{Mei15:DTinvIC}. By \cite[Lem.~3.11]{MT18:GV} the $\epsilon_1=\epsilon_2=0$ limit of our conjecture is also seen to be compatible with Katz's proposal for genus zero Gopakumar--Vafa invariants \cite{Kat08:GVgenus0}.

A more stringent piece of evidence comes from the case of local $\bbP^2$, by combining a result of Bousseau \cite{Bou19:Takahashi} with the identification of the NS limit with the GW theory of $\bbP^2$ relative to a smooth cubic in \cref{thm: NS vs log num}. For this, observe that by setting $\mchar{q}_+ = \mchar{q}_- = (\mchar{q}_1)^{1/2}$ the generating series $\refMTb{K_S}{\beta}$ specialises to the centred Poincar\'{e} polynomial
\begin{equation*}
	\MT^{\rm NS}_{\beta}(K_S)(\mchar{q}_1) \coloneqq \refMTb{K_S}{\beta}(\mchar{q}_1^{1/2},\mchar{q}_1^{1/2}) = (-\mchar{q}^{-1/2}_1)^{\dim(M_{\beta})} \sum_{j=0}^{2\dim(M_{\beta})} b_{j}(M_{\beta}) ~ (-\mchar{q}_1^{1/2})^{j}\,.
\end{equation*}
\vspace{0.5em}
\begin{thm}[\cite{Bou19:Takahashi,MS20:CohomChiIndep}]
\label{thm: Bousseaus BPS integrality P2 E}
	 For $D$ a smooth cubic in $\bbP^2$, we have\footnote{We corrected the factor $\ell^{-1}$ to $\ell^{-2}$ in the definition of $\bar{F}^{NS}(y^{1/2},Q)$ in \cite[Thm~0.4.6]{Bou19:Takahashi}.} 
	\begin{equation*}
		\GW_{\beta}(\bbP^2/D) = \sum_{k|\beta} \frac{1
 }
{k^2 [\mchar{q}_1^k]} \, \MT^{\rm NS}_{\beta/k}(K_{\bbP^2})(\mchar{q}_1) ~ \Bigg\lvert_{\mchar{q}_1 \to \re^{\epsilon_1}}\,.
	\end{equation*}
\end{thm}

\begin{cor}
    \label{cor: MT BPS local P2 NS}
	\Cref{conj: MT type BPS local surfaces} holds in the Nekrasov--Shatashvili limit ($\epsilon_2=0$) for $S=\bbP^2$.
\end{cor}
\begin{proof}
	Combining \Cref{thm: NS vs log num} with \Cref{thm: Bousseaus BPS integrality P2 E} we obtain
	\begin{equation*}
		\epsilon_2 \, \refGWb{K_{\bbP^2}}{\mgp{T}_{\mr{F}}}{\beta} ~\Big\lvert_{ \epsilon_2 = 0 } = \GW_{\beta}(\bbP^2/D) = \sum_{k|\beta} \frac{1}
{k^2 [\mchar{q}_1^k]} \, \MT^{\rm NS}_{\beta/k}(K_{\bbP^2})(\mchar{q}_1) ~\Bigg\lvert_{\mchar{q}_1 \to \re^{\epsilon_1}} \,.
	\end{equation*}
	The expression on the right-hand side  is of course nothing but the restriction
	\begin{equation*}
		(\mchar{q}_2-1) \sum_{ k \vert \beta} \frac{1}{k[\mchar{q}_1^k][\mchar{q}_2^k]} \, \MT_{\beta/k}(K_{\bbP^2})(\mchar{q}_+^k, \mchar{q}_-^k) ~\Bigg|_{\mchar{q}_1=\re^{\epsilon_1},\,\mchar{q}_2=1}\,.\qedhere
	\end{equation*}
\end{proof}
Note that by the above arguments \Cref{conj: MT type BPS local surfaces} is compatible with the level of generality of \cite[Conj.~8.16]{Bou20:QuTropVert}, which asserts that \Cref{thm: Bousseaus BPS integrality P2 E} generalises to all del Pezzo surfaces.

In the case of local $\bbP^2$, we can also lend further support to the veracity of \Cref{conj: MT type BPS local surfaces} away from the unrefined and NS limits. The refined Maulik--Toda invariants earlier appeared in the work of 
\cite{KPS23:refBPSP2,KLMP24:CohOneDimSheavesP2}; the authors kindly shared the outcome of their calculation of $\refMTb{K_{\bbP^2}}{\beta}$ up to degree seven with us, and we checked that it agrees with  the refined direct integration method \cite[Table~2]{HKK13:OmegaBmodelRigidN2} (see \Cref{rmk: direct integration BPS}). Using the results of \cref{sec: ref mirror sym KP2,sec: ref KP2 beyond LV} and \cref{prop:KP2refintcheck}, we obtain the following statement.

\begin{prop}
    \label{prop:KP2refMTcheck}
	Suppose {\bf (Conifold asymptotics)} holds as formulated in \cref{def:dirint}. Then \Cref{conj: MT type BPS local surfaces} holds for $S=\bbP^2$ up to degree seven, modulo contributions from stable maps with $g>84$.
\end{prop}


\subsection{Refined GW/PT correspondence}
\label{sec: refined PT-GW}
In this section we  formulate a refinement of the GW/PT correspondence conjecture, where the refinement on the stable pairs side is given by the Nekrasov--Okounkov equivariant K-theoretic PT index \cite{NO14:membranes}. 

\subsubsection{K-theoretic Pandharipande--Thomas theory}
\label{sec:KthPT}
Let $\mgp{G}^{\prime}$ be a connected reductive group acting on a smooth quasi-projective Calabi--Yau threefold $X$, and $\mgp{T}^{\prime}$ a maximal torus of $\mgp{G}^{\prime}$. The action lifts to the moduli space $P_n(X,\beta)$ parametrising stable pairs with support $\beta$ and Euler characteristic $n$. As before, we assume that the character $\kappa$ of the induced $\mgp{G}^{\prime}$-action on $\omega_X^{-1}$ admits a square root, denoted by $\mchar{q}_+$. Under this assumption, the authors of \cite{NO14:membranes} prove that the virtual canonical bundle of this moduli space has a $\mgp{G}^{\prime}$-equivariant square root,
\begin{equation*}
    K^{1/2}_{P_n(X,\beta),\vir}\,,
\end{equation*}
 and K-theoretic PT invariants are defined as the $\mgp{G}^{\prime}$-equivariant Euler characteristic of the virtual structure sheaf twisted by this square root,
\begin{equation*}
    \widehat{\mc{O}}^{\vir}_{P_n(X,\beta)} \coloneqq \mc{O}^{\vir}_{P_n(X,\beta)} \otimes K^{1/2}_{P_n(X,\beta),\vir}\,.
\end{equation*}
When the moduli space itself is not proper, but its $\mgp{T}^{\prime}$-fixed locus is, we define the invariants by localisation:
\begin{equation}
    \label{eq: defn Kth PT invariant}
    \refPTnb{n}{X}{\mgp{G}^{\prime}}{\beta} \coloneqq \chi_{ \mgp{T}^{\prime} } \left( P^{\mgp{T}^{\prime}} , \, \frac{\mc{O}^{\vir}_{P^{\mgp{T}^{\prime}}}}{\sum_i(-1)^i\Lambda^i (N^{\vir})^{\vee}} \otimes K^{1/2}_{P,\vir} \Big\rvert_{P^{\mgp{T}^{\prime}}}\right)\qquad \in K_{\mgp{G}^{\prime}}^{\mr{loc}}(\mr{pt})\,,
\end{equation}
where we used the shorthand $P\coloneqq P_{n}(X,\beta)$. Whenever the moduli space of stable pairs is proper, the above equates to
\begin{equation*}
\refPTnb{n}{X}{\mgp{G}^{\prime}}{\beta}=
    \chi_{ \mgp{G}^{\prime} }  \big( P_{n}(X,\beta) , \, \widehat{\mc{O}}^{\vir}_{P_{n}(X,\beta)}\big)\,.
\end{equation*}
In this scenario, there is a much restricted dependence on the equivariant parameters \cite[Thm.~7.1]{NO14:membranes}.

\begin{thm}[Rigidity, \cite{NO14:membranes} ] \label{thm: PT rigidity} 
    Suppose that $P_n(X,\beta)$ is proper, or the $\mgp{G}^{\prime}$-action on $X$ is Calabi--Yau. Then,
    \begin{equation*}
        \refPTnb{n}{X}{\mgp{T}^{\prime}}{\beta} \quad \in \quad \bbZ[\mchar{q}_+^{\pm 1}]  \quad \subseteq \quad K_{\mgp{G}^{\prime}}(\mr{pt})^{\mr{loc}}\,.
    \end{equation*}
    
\end{thm}

\begin{rmk}
The statement of \cite[Thm.~7.1]{NO14:membranes} technically covers only the properness condition in its assumptions. However, the statement holds also without it, provided the $\mgp{G}^{\prime}$-action is Calabi--Yau: this
can be shown to follow from the same arguments used in \cite[Sec.~7.2 \& 7.3]{NO14:membranes}.
\end{rmk}

A consequence of \cref{thm: PT rigidity} is that, for targets which admit a group action with $\kappa \neq 1$, the K-theoretic PT invariants are, indeed, a refinement of the numerical ones, and recover the latter in the limit $\kappa \rightarrow 1$.

\begin{prop} \label{prop: PT unrefined limit}
    If the $\mgp{G}^{\prime}$-action on $X$ is Calabi--Yau we have\footnote{As usual, we adopt \Cref{notation: defn via localisation} and define the right-hand side of the equation via localisation if $P_n(X,\beta)$ is not proper.}
    \begin{equation*}
        \refPTnb{n}{X}{\mgp{G}^{\prime}}{\beta} = \intEquiv{\mgp{G}^{\prime}}_{[P_n(X,\beta)]^{\vir}_{\mgp{G}^{\prime}}} 1 \in \bbZ\,.
    \end{equation*}
\end{prop}

\begin{proof}
    By \Cref{thm: PT rigidity}, we know that $\refPTnb{n}{X}{\mgp{G}^{\prime}}{\beta}\in \bbZ$. This means that after applying virtual Hirzebruch--Riemann--Roch \cite{FG10:virtRR,RS21:EquivVirtGRR} to \eqref{eq: defn Kth PT invariant} the resulting integral
    \begin{equation*}
        \refPTnb{n}{X}{\mgp{G}^{\prime}}{\beta} = \intEquiv{\mgp{T}^{\prime}}_{[P^{{\mgp{T}^{\prime}}}]^{\vir}_{\mgp{T}^{\prime}}} \frac{\ch[\mgp{T}^{\prime}]\!\left(K^{1/2}_{P,\vir} \right) ~ \mr{Td}_{\mgp{T}^{\prime}}\!\left( T^{\vir}_{P} \big\vert_{P^{\mgp{T}^{\prime}}}\right)}{e^{\mgp{T}^{\prime}} \! (N^{\vir})}
    \end{equation*}
    lies in the zero graded piece of $R^{\mr{loc}}_{\mgp{G}^{\prime}} \supset \bbZ$.
    Then, as in the proof of \cite[Prop.~2.22]{Tho20:refVW}, a quick dimension count shows that
    \begin{equation*}
        \frac{\ch[\mgp{T}^{\prime}]\!\left(K^{1/2}_{P,\vir} \right) ~ \mr{Td}_{\mgp{T}^{\prime}}\!\left( T^{\vir}_{P} \big\vert_{P^{\mgp{T}^{\prime}}}\right)}{e^{\mgp{T}^{\prime}} \! (N^{\vir})} = \frac{1}{e^{\mgp{T}^{\prime}} \! (N^{\vir})} + \ldots
    \end{equation*}
    modulo terms in $\Chow^{>\mr{vd}\, P^{{\mgp{T}^{\prime}}}}_{\mgp{T}^{\prime}} \!\!(P^{{\mgp{T}^{\prime}}})^{\mr{loc}}$, meaning that only the leading order can contribute non-trivially:
    \begin{equation*}
    	\refPTnb{n}{X}{\mgp{G}^{\prime}}{\beta} =  \intEquiv{\mgp{T}^{\prime}}_{[P^{{\mgp{T}^{\prime}}}]^{\vir}_{\mgp{T}^{\prime}}} \frac{1}{e^{\mgp{T}^{\prime}} \! (N^{\vir})}\,.
    \end{equation*} 
    By our convention \eqref{eq: defn via localisation}, this is precisely what we wanted to show.
\end{proof}

\subsubsection{The conjecture}
\label{sec:KthPTconj}
We package all K-theoretic PT invariants of $X$ into a generating series
\begin{equation*}
    \refPTb{X}{\mgp{G}^{\prime}}{\beta} \coloneqq \sum_{n\in \bbZ} (-\mchar{q}_-)^n ~ \refPTnb{n}{X}{\mgp{G}^{\prime}}{\beta} \,. 
\end{equation*}
Such as defined, the above is a Laurent series in the formal box counting parameter $\mchar{q}_-$, which \emph{a priori} does not enjoy any geometric interpretation. A key insight of \cite{NO14:membranes} is that $\mchar{q}_-$ should rather be thought of as a coordinate of the extended group
\begin{equation*}
    \mgp{G} \coloneqq \mgp{G}^{\prime} \times \Gm_{\mchar{q}_-}\,,
\end{equation*}
where we naturally identify $K_{\mgp{G}} (\mr{pt}) \cong K_{\mgp{G}^{\prime}}(\mr{pt}) [\mchar{q}_-^{\pm 1}]$.

\begin{conj}[\cite{NO14:membranes}]
    \label{conj: PT rational lift}
     $\refPTb{X}{\mgp{G}^{\prime}}{\beta}$ lifts to an element in $\mr{Frac}\big(K_{\mgp{G}} (\mr{pt})\big)$.
\end{conj}

\begin{rmk}
    In \cite[Conj.~2.1]{NO14:membranes}, the authors implicitly predict that the PT generating series should lift to an element in $K_{\mgp{G}} (\mr{pt})^{\mr{loc}}\otimes \bbQ$, which strongly confines the location of poles. In what follows we will, however, only consider the above weak version of their conjecture.
\end{rmk}

In order to compare with GW theory, we first extend the action of $\mgp{G}^{\prime}$ on $X$ to an action on $Z=X\times \Aaff{2}$ by declaring that both affine factors are acted upon via $\mchar{q}_+$. Together with the action of the anti-diagonal torus, we therefore obtain an action
\begin{equation*}
    \mgp{G} = \mgp{G}^{\prime} \times \Gm_{\mchar{q}_-} \quad \circlearrowright \quad Z\,.
\end{equation*}
We introduce the all-genus disconnected refined Gromov--Witten generating series:
\begin{equation*}
    \refGWdiscb{Z}{\mgp{G}}{\beta} \coloneqq \sum_{g\in \bbZ} \intEquiv{\mgp{G}}_{[\Mbar^{\bullet}_{g}(Z,\beta)]^{\vir}_{\mgp{G}}} 1\,,
\end{equation*}
where $\Mbar^{\bullet}_{g}(Z,\beta)$ is the moduli stack of stable maps with possibly disconnected domain, but without any connected component being contracted to a point. This is related to our earlier setup via
\begin{equation}
    \label{eq: GW conn to disc}
    1 + \sum_{\beta \neq 0} \mc{Q}^{\beta} \, \refGWdiscb{Z}{\mgp{G}}{\beta} = \exp \l(\sum_{\beta \neq 0} \mc{Q}^{\beta} \, \refGWb{Z}{\mgp{G}}{\beta}\r) 
\end{equation}
with the sums taken over a saturated additive subset of $\hhh_2(X,\bbZ)^+$ for which all involved equivariant integrals are well-defined via localisation.

\begin{conj}[Refined PT--GW correspondence]
\label{conj: refined PT GW}
    Assuming \Cref{conj: PT rational lift}, we have
    \begin{equation}
        \label{eq: refined PT GW} 
        \ch[\mgp{G}] \, \refPTb{X}{\mgp{G}^{\prime}}{\beta} = \refGWdiscb{Z}{\mgp{G}}{\beta}\,.
    \end{equation}
\end{conj}



We stated \Cref{conj: refined PT GW} for $Z=X\times \Aaff{2}$, as this is the only setting in which we can currently provide solid evidence. An optimistic reader might, however, wonder if a modification of the conjecture holds for Calabi--Yau fivefolds of the form $Z = \Tot_X(L_1 \oplus L_2)$ where $L_1$, $L_2$ are line bundles on a smooth threefold $X$. We have numerical evidence suggesting the validity of a version of \Cref{conj: refined PT GW}, where the PT generating series is defined as in \cite[Sec.~3.1.8]{NO14:membranes}. It is noteworthy that for local curves and surfaces, this correspondence predicts that a single GW problem relates to multiple PT problems, depending on the choice of line bundles $L_1$ and $L_2$. This expectation is still speculative at this stage, and certainly requires a more in-depth investigation.

\subsubsection{Evidence for \cref{conj: PT rational lift,conj: refined PT GW}}
\label{sec: PT GW unref}
The first test for our conjecture again comes from the resolved conifold, $X = \Tot (\cO_{\bbP^1}(-1) \oplus \cO_{\bbP^1}(-1))$. As the stable pairs moduli space for this target is proper, the PT generating series can be computed using the refined topological vertex \cite{IKV09:RefVert,NO14:membranes} leading to the well-known expression
\begin{equation}
	\label{eq: Kth PT resolved conifold}
	1+\sum_{d >0}\refPTb{X}{\mgp{G}^{\prime}}{d[\bbP^1]} ~ \cQ^{d} = \exp \sum_{k > 0} \frac{\cQ^{k }}{k [\mchar{q}_1^k][\mchar{q}_2^k]}\,.
\end{equation}
The above formula can as well be deduced from the K-theoretic vertex \cite{KOO21:2legKthVertex} and is related to the motivic generating series
\[
1 + \sum_{d,n=0}^\infty \mchar{q}^n \cQ^d [P_n(X,d)]^{\rm mot} = \prod_{m \geq 1} \prod_{j=0}^{m-1} \l(1-\bbL^{-\frac{m+1}{2} +j} \mchar{q}^m \cQ\r)
\]
determined in \cite[Prop.~4.3]{MMNS12:motivicResCon} via the change of variables $\mchar{q}_1 = \mchar{q} \bbL^{1/2}$ and $\mchar{q}_2 = \mchar{q}^{-1}\bbL^{1/2}$: here, $[P_n(X,\beta)]^{\rm mot}$ $\in \bbZ[\bbL^{\pm 1/2}]$ is the virtual motive of $P_n(X,\beta)$, and $\bbL$ is the Lefschetz motive.
Then, \Cref{eq: refined PT GW} applied to \eqref{eq: Kth PT resolved conifold} returns the refined GW generating series \eqref{eq:resconGWref}, verifying \cref{conj: refined PT GW} in this case\footnote{It is worth reiterating that \eqref{eq:resconGWref} was only proven for torus actions scaling a negative and trivial line bundle direction over a fixed point with opposite weights; hence, even for the resolved conifold \Cref{conj: refined PT GW} is still open for general group actions.}. 

Further evidence for \Cref{conj: refined PT GW} comes from its compatibility with the original GW/PT correspondence conjecture \cite{PT09:StabPairs}. To see this, suppose the $\mgp{G}^{\prime}$-action on $X$ is Calabi--Yau. We observed in \Cref{prop: PT unrefined limit} that under this assumption the generating series on the left-hand side of \eqref{eq: refined PT GW} simplifies to
\begin{equation*}
    \refPTb{X}{\mgp{G}^{\prime}}{\beta} = \sum_{n\in\bbZ} (- \mchar{q}_-)^n ~~ \intEquiv{\mgp{G}^{\prime}}_{[P_n(X,\beta)]^{\vir}_{\mgp{G}^{\prime}}} 1  \eqqcolon \PT^{\mr{un}}_{\beta}(X,\mgp{G}^{\prime})\,.
\end{equation*}
Similarly, by \Cref{cor: GW unref limit}, we have
\begin{equation*}
	\refGWdiscb{Z}{\mgp{G}}{\beta} = \sum_{g\in \bbZ} (-\epsilon_-^2)^{g-1} ~~ \intEquiv{\mgp{G}^{\prime}}_{[\Mbar^{\bullet}_{g}(X,\beta)]^{\vir}_{\mgp{G}^{\prime}}} 1 \eqqcolon \refGWdiscb{X}{\mgp{G}^{\prime}}{\beta}\,.
\end{equation*}
Therefore, if $\mgp{G}$' is a Calabi--Yau action, \Cref{conj: PT rational lift} and \labelcref{conj: refined PT GW} reduce to the original GW/PT correspondence conjecture, reformulated as follows.
\begin{conj}[GW/PT correspondence,  \cite{PT09:StabPairs}]
    \label{conj: unref PT-GW}

        $\PT^{\mr{un}}_{\beta}(X,\mgp{G}^{\prime})$ lifts to a rational function in $\mchar{q}_-$. Moreover, under this lift,
        \begin{equation*}
            \PT^{\mr{un}}_{\beta}(X,\mgp{G}^{\prime}) ~\Big\vert_{\mchar{q}_- = \re^{\epsilon_-}}= \GW^{\bullet}_{\beta}(X,\mgp{G}^{\prime})\,.
        \end{equation*}
\end{conj}

	The conjecture is known to be true for local curves \cite{BP08:GWLocCurves,OP10:localDTcurves} and toric varieties \cite{MPT10:K3modularForms,MOOP11:GWDTtoric}, and many more cases are covered by the methods developed by Pandharipande and Pixton \cite{PP17:GWPTquintic}. The general case is covered by recent work of Pardon \cite{Par23:UnivCurveCount}.

Additional evidence is provided by the refined direct integration for local $\bbP^2$. In \Cref{prop:KP2refMTcheck} we have already explained how refined mirror symmetry agrees with the Poincar\'{e} refinement of Maulik--Toda invariants to fairly high genus, which in turn have been checked numerically in \cite{KPS23:refBPSP2,KLMP24:CohOneDimSheavesP2} to match the K-theoretic PT generating series to high degree. This gives the following stable pairs analogue of \cref{prop:KP2refMTcheck}.

\begin{prop}
	Suppose that {(\bf Conifold asymptotics)} holds as formulated in \cref{def:dirint}. Then, \Cref{conj: refined PT GW} holds for $X=K_{\bbP^2}$ with the $\mgp{G}'=\Gm$-scaling action on the fibre direction to curve class degree seven, and modulo contributions from stable maps with $g>84$.
 \label{prop:KP2refPTcheck}
\end{prop}

Finally, recent work of Oberdieck \cite[Thm~1.3]{Ob24:EnriquesKth} presents a conjectural formula for the K-theoretic PT generating series for $X$ a local Enriques  surface. Conforming to the refined GW/PT correspondence in \cref{conj: refined PT GW}, this formula is found to satisfy the holomorphic anomaly equations of the refined direct integration method  as the equivariantly refined GW generating series of $X \times \Aaff{2}$ \cite[Sec.~4.8]{Ob24:EnriquesKth}.

\subsubsection{Implications for the rigidity conjecture}
\label{sec: rigidity implication}
\Cref{conj: refined PT GW} has two important consequences. The first immediate one is that the refined GW/PT correspondence implies the Rigidity \cref{conj: rigidity} for refined GW invariants: this follows from the analogous statement for K-theoretic PT invariants  (\Cref{thm: PT rigidity}).  Let $\hhh$ be a saturated additive subset of $\hhh_2(X,\bbZ)^+$.
\begin{prop}
    \label{prop: PT-GW implies rigidity}
    Suppose \Cref{conj: refined PT GW} holds and $P_n(X,\beta)$ and $\Mbar_{g}^{\bullet}(X,\beta)$ are proper for all $n,g\in \bbZ$ and curve classes $\beta$ in a saturated additive subset $\hhh \subseteq \hhh_2(X,\bbZ)^+$. Then the Gromov--Witten generating series of $X\times \Aaff{2}$ satisfies \Cref{conj: rigidity} for all $\beta \in \hhh$.
\end{prop}
\begin{proof}
    By \Cref{thm: PT rigidity} we have
    \begin{equation*}
        \refPTb{X}{\mgp{G}^{\prime}}{\beta} \in \bbZ[\mchar{q}_+^{\pm 1}](\!(\mchar{q}_-)\!)\,.
    \end{equation*}
    Thus, \eqref{eq: GW conn to disc} and \eqref{eq: refined PT GW} imply that $\refGWb{X}{\mgp{G}}{\beta} \in \widehat{R}^{\,\mr{loc}}_{\mgp{G}}$ can only depend on $\epsilon_+$ and $\epsilon_-$. In combination with \Cref{lem: shape of refGW series} this gives
    \begin{equation*}
    	\epsilon_1 \epsilon_2 \, \refGWb{Z}{\mgp{G}}{\beta} \in \bbQ\llbracket \epsilon_+^2,\epsilon_-^2\rrbracket\,. \qedhere
    \end{equation*} 
\end{proof}

\subsubsection{The vertex}
\label{sec: refined vertex}
As  mentioned in the discussion of the resolved conifold, the K-theoretic PT generating series of a toric Calabi--Yau threefold $X$ is governed by the K-theoretic vertex \cite{NO14:membranes,KOO21:2legKthVertex}, therefore by \Cref{conj: refined PT GW} the same method should compute the equivariant Gromov--Witten generating series of $X\times\Aaff{2}$. If the stable pair moduli spaces are proper, the K-theoretic vertex simplifies to the refined topological vertex \cite{IKV09:RefVert,NO14:membranes,Arb21:KthDT}: here properness guarantees the independence on the choice of the so called preferred direction.

When $X$ is not equirigid, the outcome of a refined topological vertex calculation will only be related to the GW generating series through its conjectural lift to localised equivariant K-theory (see \Cref{rmk: Kth lift}). For a concrete example, let us consider $X=\Tot \cO_{\bbP^1}(-2) \times \bbC$. In \Cref{sec: local P1 shifted} we presented the conjectural formula
\begin{equation*}
    \refGWb{Z}{\TCY}{d[\bbP^1]} = \ch[\TCY] \frac{1}{d} \frac{[\mchar{q}_0^{-1}\mchar{q}_1^{-1}\mchar{q}_2^{-1}]}{[\mchar{q}_0][\mchar{q}_1][\mchar{q}_2]}\,.
\end{equation*}
If we take the limit $\mchar{q}_0^{\pm 1} \rightarrow \infty$ of the K-theoretic lift while keeping $\mchar{q_+}$ constant, we recover the result of a calculation using the refined topological vertex \cite[Sec.~5.3]{IKV09:RefVert}:
\begin{equation*}
    \frac{1}{d} \frac{\mchar{q_+^{\pm 1}}}{[\mchar{q}_1][\mchar{q}_2]}\,.
\end{equation*}
Notice that in order to perform the limit it was crucial to first pass to an analytic continuation in $\epsilon_0, \epsilon_1, \epsilon_2$, meaning that the latter quantity does not have a direct enumerative meaning in terms of stable maps anymore.


\begin{appendix}

\section{The conifold \texorpdfstring{$R$}{R}-matrix}
\label{sec:resconRmat}

Let $Z=\mathrm{Tot} (\cO_{\bbP^1}^{\oplus 2}(-1)\oplus \cO_{\bbP^1}^{\oplus 2})$ be the product of the resolved conifold by the affine plane. It is a GIT quotient of $\bbC^6$ with coordinates $(x_0, x_1, y_1, y_2, z_1, z_2)$,
\[
Z \coloneqq \bbC^6/\!\!/_{\xi>0}  \Gm = (\bbC^6\setminus\Delta)/\Gm\,,
\]
where
\[
\Delta \coloneqq V\l(\bra x_0,x_1\ket\r).
\]
and the quotient torus acts as
\bea
\Gm \times \bbC^6 & \longrightarrow & \bbC^6\,,  \nn \\
(\lambda; x_0, x_1, x_2, y_0,y_1,y_2) & \to & (\lambda x_0, \lambda x_1, \lambda^{-1} y_1, \lambda^{-1}y_2, z_1, z_2)\,. \nn
\eea
$Z$ carries a five-dimensional torus action ${\Tmax}$ descending from an action on $\bbC^6$ with geometric weights $(0, \alpha_0, \alpha_1, \alpha_2, -\epsilon_1, -\epsilon_2)$. Denote by $P_\pm$ the points at zero and infinity in the base $\bbP^1$, $\varepsilon_\pm$ the lifts to $\Chow^{}_{\Tmax}(Z)$ of the respective cohomology classes, and $H\coloneqq c^{{\Tmax}}_1(\cO_{\bbP^1}(1))$ the equivariant point class on the base. The Atiyah--Bott isomorphism is then
\[ 
H \to - \alpha_0 P_-\,,
\] 
and the normal weights at $P_\pm$ are
\[ 
w^{(j)}_+=(\alpha_0, \alpha_1, \alpha_2, -\epsilon_1,-\epsilon_2)_j\,, \qquad 
w^{(j)}_-=(-\alpha_0, \alpha_0+\alpha_1, \alpha_0+\alpha_2, -\epsilon_1,-\epsilon_2)_j\,.
\] 
The three-pointed genus-zero Gromov--Witten invariants of $Z$ are 
\[ 
\bra \phi, \chi, \psi \ket_{0,d} = 
\begin{cases}
\intEquiv{\Tmax}_{[Z]_{\Tmax}} \phi \cup \chi \cup \psi & d=0\,, \\
\frac{1}{\epsilon_1 \epsilon_2} & d>0,~\phi=\chi=\psi=H\,, \\
0 & \mathrm{otherwise}.\\
\end{cases}
\]
The genus zero potential as a function of $t=t_0 \mathbf{1} + t_1 H$ then reads
\[ 
\GW^{\Tmax}_0(t)= 
\frac{\left(\alpha _0+\alpha _1+\alpha _2\right) t_0^3+3 \alpha _1 \alpha _2 t_1 t_0^2-3 \alpha _0 \alpha _1 \alpha _2
   t_1^2 t_0+\alpha _0^2 \alpha _1 \alpha _2 t_1^3}{6 \alpha _1 \left(\alpha _0+\alpha _1\right) \alpha _2 \left(\alpha _0+\alpha _2\right) \epsilon _1 \epsilon _2}
+\frac{\mathrm{Li}_3\left(\re^{t_1}\right)}{\epsilon _1 \epsilon _2}\,.
\] 
In particular, the inverse Gram matrix of the equivariant intersection pairing in cohomology, in the basis $\{\mathbf{1},H\}$ for $\Chow^{}_{\Tmax}(Z)$, is
\beq 
\eta^{-1}=
\left(
\begin{array}{cc}
 \alpha _0 \alpha _1 \alpha _2 \epsilon _1 \epsilon _2 & \alpha _1 \alpha _2 \epsilon _1 \epsilon _2 \\
 \alpha _1 \alpha _2 \epsilon _1 \epsilon _2 & -\left(\alpha _0+\alpha _1+\alpha _2\right) \epsilon _1 \epsilon _2 \\
\end{array}
\right)\,,
\label{eq:etacf}
\eeq
and flat sections $e^{t_0/z} \phi(t_1,z)$ of the Dubrovin connection satisfy the Picard--Fuchs differential equation
\beq
z^2 \phi ''\left(t_1\right)=\frac{z \left(\alpha _0+\alpha _2 \re^{t_1}\right) \phi '\left(t_1\right)-\alpha _1 \re^{t_1} \left(\alpha _2 \phi \left(t_1\right)-z \phi
   '\left(t_1\right)\right)}{\left(\re^{t_1}-1\right)}\,.
   \label{eq:pfcf}
\eeq   
A basis of solutions is given by the components of the $J$-function of $Z$,
\bea 
J(t,z) & \coloneqq & \re^{t_0/z+t_1 p/z} \sum_{d \geq 0} \re^{t_1 d}\frac{\prod_{m=0}^{d-1} (-H+\alpha_1 z) \prod_{m=0}^{d-1} (-H+\alpha_2 z)}{\prod_{m=1}^d (H+m z) \prod_{m=1}^d (H+\alpha_0+m z)}\,, \nn \\
&=& 
z \re^{\frac{t_0}{z}} \, _2F_1\left(-\frac{\alpha _1}{z},-\frac{\alpha _2}{z};1+\frac{\alpha _0}{z};\re^{t_1}\right) \varepsilon_+
+
z \re^{\frac{t_0-\alpha _0 t_1}{z}} \, _2F_1\left(-\frac{\alpha _0+\alpha _1}{z},-\frac{\alpha _0+\alpha _2}{z};1-\frac{\alpha _0}{z};\re^{t_1}\right)\varepsilon_-
\,.\nn
\eea 
The canonical coordinates $u_\pm = t_0 + f_{\pm}(t_1)$
satisfy
\beq
f_{\pm}'(x)^2-\frac{f_\pm'(x) \left(\alpha _0+\left(\alpha _1+\alpha _2\right) \re^x\right)}{\re^x-1} +\frac{\alpha _1 \alpha _2 \re^x}{\re^x-1}= 0\,,
\label{eq:cancf}
\eeq 
with the two values $f_\pm$ corresponding to the two branches of the solutions. From \eqref{eq:etacf} and \eqref{eq:cancf}, we read off the inverse square lengths of the canonical idempotents as
\[ 
\Delta_\pm = \frac{\alpha _1 \left(\alpha _0+\alpha _1\right) \alpha _2 \left(\alpha _0+\alpha _2\right) \left(f'_+-f'_-\right){}^2 \epsilon _1 \epsilon _2}{\left(\alpha _0+\alpha _1+\alpha _2\right) (f'_\mp)^2-2 \alpha _1 \alpha _2 f'_\mp-\alpha _0 \alpha _1 \alpha
   _2}\,.
\] 
\subsection{A one-dimensional equivariant mirror}
Let 
\[ 
\cW(t_0, t_1, s) \coloneqq t_0+\alpha _1 \log \left(1-\re^{t_1}s\right)-\alpha _2 \log s+\left(\alpha _0+\alpha _2\right) \log \left(1-s\right)
\,.
\] 
We claim that $\cW$ is a mirror superpotential for the quantum cohomology of $Z$, with primitive form given by 
\[
\dd \Omega=((\alpha_0 +\alpha_2)\alpha_1 \epsilon_1 \epsilon_2)^{-1/2} \frac{\dd s}{s}\,.
\]
To see this, it suffices to check that the  values of $\cW$ on the critical set $\{s_\pm^{\rm cr}(t_0, t_1) : \de_s \cW(t_0, t_1, s_\pm^{\rm cr}) = 0\}$ coincide with the canonical coordinates, and that $s^2 \de^2_s \cW(t_0, t_1, s_\pm^{\rm cr}) = (\alpha_0 +\alpha_2)\alpha_1 \epsilon_1 \epsilon_2 \Delta_\pm$. We find:
\[ 
s_\pm^{\rm cr} = 
\frac{ \alpha _0+\alpha _1 \re^{t_1}-\alpha _2 \re^{t_1}\pm\sqrt{\left(\alpha _0+\left(\alpha _1-\alpha _2\right) \re^{t_1}\right){}^2+4 \left(\alpha _0+\alpha
   _1\right) \alpha _2 \re^{t_1}}}
   {2 \left(\alpha _0+\alpha _1\right)\re^{t_1}}
   \,,
\] 
from which we verify that
\[ 
\de_{t_0} \cW =1\,, \quad
\de_{t_1} \cW(t_0, t_1, s_\pm^{\rm cr}) = f_\pm(t_0, t_1)\,, \quad s^2\de^2_{s} \cW(t_0, t_1, s_\pm^{\rm cr}) = (\alpha_0 +\alpha_2)\alpha_1 \epsilon_1 \epsilon_2\Delta_\pm\,.
\] 
Let $\Gamma_+ = [l_{0},l_1]$, $\Gamma_- = [l_0, l_{\re^{-t_1}}]$ where $l_x$ is the homology class in $H^1(\bbC \setminus \{0,1,\re^{-t_1}\}, \bbZ)$ of a simple loop around $x \in \bbC$, oriented counter-clockwise. Then, for fixed $z \neq 0$, the twisted periods
\[
\Pi_{\pm} \coloneqq
\frac{((\alpha_0 +\alpha_2)\alpha_1 \epsilon_1 \epsilon_2)^{-1/2}}{(1-\re^{-2\pi \ri \alpha_2/z})
(1-\re^{-2\pi \ri \alpha_1/z})}\int_{\Gamma_\pm} \re^{\cW/z} \frac{\dd s}{s}
\] 
are a basis of solutions of \eqref{eq:pfcf}. In particular, we have
\bea 
\Pi_+ &=& \frac{\re^{t_0/z} ((\alpha_0 +\alpha_2)\alpha_1 \epsilon_1 \epsilon_2)^{-1/2}}{(1-\re^{-2\pi \ri \alpha_2/z})
(1-\re^{-2\pi \ri \alpha_1/z})}
\int_{[l_0,l_1]} (1-s)^{(\alpha _0+\alpha _2)/z} s^{-\alpha _2/z} \left(1-s e^{t_1}\right){}^{\alpha _1/z} \frac{\dd s}{s} \nn \\
&=& \frac{\mathsf{B}\l(-\frac{\alpha_2}{z}, 1+\frac{\alpha_0+\alpha_2}{z}\r)}{\sqrt{(\alpha_0 +\alpha_2)\alpha_1 \epsilon_1 \epsilon_2}} \re^{t_0/z}\,_2F_1\left(-\frac{\alpha _1}{z},-\frac{\alpha _2}{z};1+\frac{\alpha _0}{z};\re^{t_1}\right)
\,, \nn
\eea
from the Euler--Pochhammer representation of the Gauss hypergeometric function, where $\mathrm{B}(x,y)=\Gamma(x) \Gamma(y)/\Gamma(x+y)$ is the Euler Beta function. Likewise, asymptotically around $z=0$,  the saddle-point expansion of $\Pi_{\pm}$ gives formal solutions of the flatness equation of the Dubrovin connection, meromorphically in $\re^{t_1/2}$ away from $t_1=0$ and to all orders in a formal expansion in $z$:
\[ 
\Pi_{\pm} \simeq \frac{\re^{u_\pm/z}}{\sqrt{-2 \pi z}} \mathsf{Asym}_{x^{\rm cr}_\pm} \l[\re^{\cW/z}\dd \Omega\r] = \frac{\re^{u_\pm/z}}{\sqrt{-2 \pi z}} \sum_{n \geq 0} \psi_n(t_1) z^n\,.
\] 
The columns of the $R$-matrix are therefore given by 
\bea
R_\pm(\mathbf{1}) &=& \mathsf{c}_\pm(z) \mathsf{Asym}_{x^{\rm cr}_\pm} \l[\re^{\cW/z}\dd \Omega\r]\,, \nn \\
R_\pm(H) &=& \mathsf{c}_\pm(z) \bigg(f_\pm + z \de_{t_1}\bigg)\mathsf{Asym}_{x^{\rm cr}_\pm} \l[\re^{\cW/z}\dd \Omega\r] \nn
\eea 
for some $\mathsf{c}_\pm(z) \in \bbC\llbracket z\rrbracket$, which is fixed via the quantum Riemann--Roch theorem by imposing
\[ 
R_\pm(\mathbf{1}) \bigg|_{t_0 \to 0, \mathrm{Re}(t_1) \to -\infty} = \prod_{j=1}^5 \frac{1}{\sqrt{w_\pm^{(j)}}}\Theta\l(\frac{z}{w_\pm^{(j)}}\r).
\] 
Let's consider $x_+^{\rm cr}$ for simplicity. We find, from the saddle-point expansion of the Euler Beta integral,
\bea
R_+(\mathbf{1}) \bigg|_{t_0 \to 0, \mathrm{Re}(t_1) \to -\infty} &=& \frac{\mathsf{c}_+(z)}{\sqrt{(\alpha_0 +\alpha_2)\alpha_1 \epsilon_1 \epsilon_2}} \mathsf{Asym}_{-\frac{\alpha_2}{\alpha_0}}\l[(1-s)^{(\alpha _0+\alpha _2)/z} s^{-\alpha _2/z} \frac{\dd s}{s}\r] \nn \\
&=& 
\frac{\mathsf{c}_+(z)}{\sqrt{\alpha_0\alpha_1\alpha_2 \epsilon_1 \epsilon_2}}\frac{\Theta\l(-\frac{z}{\alpha_2}\r)\Theta\l(\frac{z}{\alpha_2+\alpha_0}\r)}{\Theta\l(\frac{z}{\alpha_0}\r)}\,, \nn
\eea 
so that
\[
\mathsf{c}_+(z)=
\frac{\Theta\l(\frac{z}{\epsilon_1}\r)\Theta\l(\frac{z}{\epsilon_2}\r)}{\Theta\l(\frac{z}{\alpha_2+\alpha_0}\r)\Theta\l(\frac{z}{\alpha_1}\r)} \,.
\] 

\begin{rmk}
When $\alpha_0=0$, $\alpha_{1,2}=\epsilon_{1,2}$, we have $c_+(z)=1$. Therefore $R_\pm(\mathbf{1})$ may be written as a Laplace transform of the spectral curve data \cite{Eynard:2011kk}
\beq
\l\{
\bary{l}
x(z) =  \log(1-\re^{t_1+z})+\beta (z-\log(1-\re^{z})) \\[0.3em]
y(z) = z 
\eary
\r.\,,
\label{eq:refspcurve}
\eeq
where $\beta=-\epsilon_2/\epsilon_1$, and we set $t_0=0$. This confirms a heuristic calculation of Eynard--Koz\c caz \cite[Eq.~(8.35)]{Eynard:2011vs}, who derived  \eqref{eq:refspcurve} from a saddle-point asymptotic ansatz for the sum over plane partitions in the corresponding refined topological vertex calculation.
\label{rmk:refmm}
\end{rmk}
\subsubsection{Asymptotics at the conifold point}

We will be interested in the asymptotic behaviour of $f_\pm'$, $\Delta_\pm$, $x_\pm^{\rm cr}$, and $R_\pm(H^a)$ 
%
%
at the conifold point, $t_1=0$. It is readily verified that $f_-$, $\Delta_-$, $x_-^{\rm cr}$, $R_-(H^a)$ are analytic and non-zero at $t_1=0$. On the other hand, we have
%
\[
f_+' =  \frac{\alpha_0+\alpha_1+\alpha_2}{t_1} + \cO(1)\,, \quad
\Delta_+= 
-\frac{\left(\alpha _0+\alpha _1+\alpha _2\right){}^3 \epsilon _1 \epsilon _2}{t_1^2}+\cO\left(\frac{1}{t_1}\right)\,, \quad 
s_+^{\rm cr} =  1-\frac{\left(\alpha _0+\alpha _2\right) t_1}{\alpha _0+\alpha _1+\alpha _2}+\cO\left(t_1^2\right)\,,
\]
\beq 
\cW^{(n)}(s_+^{\rm cr})=  
-(n-1)! \frac{\left(\alpha _0+\alpha _1+\alpha _2\right){}^n \left(\alpha _1^{n-1}+(-1)^n \left(\alpha _0+\alpha
   _2\right){}^{n-1}\right)}{\left(\alpha _1 \left(\alpha _0+\alpha _2\right)\right)^{n-1} t_1^{n}} + \cO\l(\frac{1}{t_1^{n-1}}\r)\,.
\label{eq:cfasymp}
\eeq
%
As for $R_+(H^a)$, recall that
\bea
R_+(\mathbf{1})  &=& 
\frac{\mathsf{c}_+(z)}{\sqrt{\cW''(s_+^{\rm cr})}} \l[\re^{-\frac{z}{2 \cW''(s_+^{\rm cr})}\de^2_{ s}}\l(\frac{\re^{\cW^{\rm ng}/z}}{s}\r)\r]\Bigg|_{s=s_+^{\rm cr}}
\,, \nn \\
R_+(H)  &=& \mathsf{c}_+(z)\l( \frac{2f_+' \cW''(s_+^{\rm cr})- z\de_{t_1} \cW''(s_+^{\rm cr})}{2\cW''(s_+^{\rm cr})^{3/2}}\r) \l[\re^{-\frac{z}{2 \cW''(s_+^{\rm cr})}\de^2_{s}}\l(\frac{\re^{\cW^{\rm ng}/z}}{s}\r) \r]\Bigg|_{s=x_+^{\rm cr}}\nn \\
&+& \frac{\mathsf{c}_+(z)}{\cW''(s_+^{\rm cr})^{1/2}} \de_{t_1}\l[\re^{-\frac{z}{2 \cW''(s_+^{\rm cr})}\de^2_{s}}\l(\frac{\re^{\cW^{\rm ng}/z}}{s}\r) \r]\Bigg|_{s=x_+^{\rm cr}}\,.
\label{eq:RH}
\eea 
From \eqref{eq:cfasymp}, we find that 
\beq
\l[\re^{-\frac{z}{2 \cW''(s_+^{\rm cr})}\de^2_{s}}\l(\frac{\re^{\cW^{\rm ng}/z}}{s}\r)\r]\Bigg|_{s=s_+^{\rm cr}} = \cO(1)
\label{eq:cfasymint}
\eeq
near $t_1=0$, as is readily checked from a Feynman diagram expansion of the left hand side wherein every $n$-valent vertex contribution (arising from $\cW^{(n)}(s_+^{\rm cr})$) diverges as $\cO(t_1^{-n})$, and every half-edge contribution (arising from the propagator $(-\cW^{(2)}(s_+^{\rm cr}))^{-1/2} $) carries a factor of $\cO(t_1)$. Then,
\[
R_+(\mathbf{1})  = 
\cO(t_1)
\,, \quad
R_+(H)  = \cO(1)\,.
\]
We will need the explicit limit of $R_+(H^a)$ as $t_1 \to 0$: by the foregoing discussion,  this can be determined by just truncating each quantity appearing in the r.h.s.~of \eqref{eq:cfasymint} to its leading coefficient in $t_1$ computed in \eqref{eq:cfasymp}. Note that, for $\cW^{(n)}(s_+^{\rm cr})$, this is easily verified to be equal to 
\[ 
-(n-1)! \frac{ (a+b)^n \left(a^{n-1}+(-1)^n b^{n-1}\right)}{\left(a b\right)^{n-1} } = 
\frac{\de^n \l(a \log x + b \log(1-x)\r) }{\de x^n}\bigg|_{x=\frac{a}{a+b}} \,,
\] 
with $a=\alpha_1$, $b=\alpha_0+\alpha_2$.  Since the saddle-point expansion in \eqref{eq:cfasymint} has a finite limit for $t_1 \to 0$, we find
\bea
\lim_{t_1 \to 0}
\l[\re^{-\frac{z}{2 \cW''(s_+^{\rm cr})}\de^2_{s}}\l(\frac{\re^{\cW^{\rm ng}/z}}{s}\r)\r]\Bigg|_{s=s_+^{\rm cr}} &=& \l[\re^{\frac{z \alpha_1 (\alpha_0+\alpha_2)}{2(\alpha_0+\alpha_1+\alpha_2)^3} \de^2_{s}}\re^{\l(\alpha_1 \log x + (\alpha_0+\alpha_2) \log(1-x)\r)^{\rm ng}/z}\r]\Bigg|_{x=\frac{\alpha_1}{\alpha_0+\alpha_1+\alpha_2}} \nn \\
&=& \sqrt{-\frac{(\alpha_0+\alpha_1+\alpha_2)^3}{\alpha_1(\alpha_0+\alpha_2)}}
\mathsf{Asym}_{x=\frac{\alpha_1}{\alpha_0+\alpha_1+\alpha_2}} \l[x^{\alpha_1/z} (1-x)^{(\alpha_0+\alpha_2)/z} \dd x \r] \nn \\
&=& 
\frac{\alpha_0+\alpha_1+\alpha_2}{z+\alpha_0+\alpha_1+\alpha_2}
{\Upsilon\l(\frac{z}{\alpha_1}, \frac{z}{\alpha_0+\alpha_2}\r)}\,. \nn
\eea 
Applying the Stirling expansion of the Beta function $\mathrm{B}(a,b)$ with $a=\alpha_1/z+1$, $b= (\alpha_0+\alpha_2)/z+1$ to the last line, and using \eqref{eq:cfasymp} and \eqref{eq:RH}, we conclude that
\bea
R_+(\mathbf{1}) &=& \frac{t_1}{\sqrt{-(\alpha_0+\alpha_1+\alpha_2)\epsilon _1 \epsilon _2}}
~\frac{\Theta\l(\frac{z}{\epsilon_1}\r) \Theta\l(\frac{z}{\epsilon_2}\r)}{\Theta\l(\frac{z}{\alpha_0+\alpha_1+\alpha_2}\r)\l(\alpha_0+\alpha_1+\alpha_2+z\r)}  +\cO(t_1^2) \nn \\
R_+(H) &=& \sqrt{\frac{1}{ -(\alpha_0+\alpha_1+\alpha_2)\epsilon _1 \epsilon _2}}
~\frac{\Theta\l(\frac{z}{\epsilon_1}\r) \Theta\l(\frac{z}{\epsilon_2}\r)}{\Theta\l(\frac{z}{\alpha_0+\alpha_1+\alpha_2}\r)} + \cO(t_1)\,. \nn
\eea
The result further simplifies when we restrict to the fibrewise action $\mgp{T}_{\rm F} \hookrightarrow \Tmax$, with weights $\alpha_0=0$,  $\alpha_{1,2}=\epsilon_{1,2}$:
\bea
R_+(\mathbf{1})\Big|_{\alpha_0=0, \alpha_{1,2} =\epsilon_{1,2}} &=&  \frac{t_1}{\sqrt{-(\epsilon_1+\epsilon_2)\epsilon _1 \epsilon _2}}~\Upsilon\l(\frac{z}{\epsilon_1}, \frac{z}{\epsilon_2}\r) 
 \frac{1}{\epsilon_1+\epsilon_2 +z} + \cO(t_1^2)\,, \nn \\
R_+(H)\Big|_{\alpha_0=0, \alpha_{1,2} =\epsilon_{1,2}} &=& \sqrt{ \frac{1}{-(\epsilon_1+\epsilon_2)\epsilon _1 \epsilon _2}}
~\Upsilon\l(\frac{z}{\epsilon_1}, \frac{z}{\epsilon_2}\r) + \cO(t_1)\,. \nn
\eea
\subsection{Givental--Teleman reconstruction at the conifold point}

Let's describe the graph contributions relevant for the Givental reconstruction formula for $Z$. We write, for $i,j \in \{+,-\}$, $l \in \bbZ_{\geq 2}$, $n,m \in \bbZ_{\geq 0}$,
\bea 
D_{i,l} & \coloneqq & [z^{l-1}] \big(\mathbf{1}-R_i (-z)(\mathbf{1})\big)\,, \nn \\
E^{ij}_{n,m} & \coloneqq & [z^{n} w^m]\l( \frac{\delta_{ij}-\sum_{ab} R_{i}(-z)(H^a) \eta^{ab} R_{j}(-w)(H^b)}{z+w}\r)\,.
\label{eq:TLRcf}
\eea
For each dilaton leaf $h \in \DD_\Gamma 
$ and edge $e=\{h_1,h_2\} \in \cE_\Gamma$, we define
\[
\mathrm{Cont}(h) \coloneqq 
    T_{\mathsf{q}(\mathsf{v}(h)),\mathsf{a}(h)}\,,
%
%
\quad \mathrm{Cont}(e) \coloneqq
E^{\mathsf{q}(\mathsf{v}(h_1)), \mathsf{q}(\mathsf{v}(h_2)}_{\mathsf{a}(h_1),\mathsf{a}(h_2)}\,.
\] 
For each vertex $v \in \cV_\Gamma$, fixing a labelling $\{h_1, \dots, h_{|\HH(v)|}\}$ of the half-edges incident to $v$, 
we define 
\[
\mathrm{Cont}(v) = (\Delta_i)^{\frac{2\mathsf{g}(v)-2+\mathsf{|H(v)|}}{2}} 
\int_{[\Mbar_{\mathsf{g}(v), |\HH(v)|}]}\prod_{k=1}^{|\HH(v)|} 
\psi_k^{\mathsf{a}(h_k)}\,.
\]
For $t_1 \neq 0$, Givental's reconstruction formula reads
\beq 
\GW_g(Z,\mgp{T}_{\rm F}) = \sum_{\Gamma \in \cG^{[1]}_{g,0}} \frac{1}{|\mathrm{Aut(\Gamma)}|} \prod_{e\in \cE_\Gamma} \mathrm{Cont}(e) \prod_{h\in \DD_\Gamma} \mathrm{Cont}(h) \prod_{v\in \cV_\Gamma} \mathrm{Cont}(v)\,.
\label{eq:GWgrescon}
\eeq 
By \eqref{eq:cfasymp} and \eqref{eq:RH}, $\GW_g(Z,\mgp{T}_{\rm F})$ has a pole of order ${2g-2}$ at $t_1=0$ arising exclusively from contributions of stable graphs with decoration $i=+$ on all vertices: the analyticity of $\Delta_-$ at $t_1=0$ implies that graphs having any vertex coloured by $i=-$ may only contribute sub-leading terms in $1/t_1$. Also, in the edge term contribution \eqref{eq:TLRcf}, only terms with $a=b=1$ contribute to the coefficient of the leading pole, since $R_+(\mathbf{1})=\cO(t_1)$ gives again a sub-leading contribution. Finally, note that stable graphs with $k$-dilaton leaves have the same divergent behaviour as the graph obtained by removing them, since each  addition of a dilaton leaf (which contributes a factor $\cO(t_1)$, by \eqref{eq:RH}) comes balanced with a corresponding factor of $\Delta_+^{1/2}$ (which contributes a factor $\cO(t_1^{-1})$, by \eqref{eq:cfasymp}). 

As far as the dominant asymptotics at the conifold point is concerned, then, the sum over graphs in \eqref{eq:GWgrescon} reduces to a graph expansion akin to that of the higher genus primary potential of a scalar CohFT, but with a slight deformation of the dilaton shift. Explicitly, the coefficient of the leading term in the asymptotics in $1/t_1$ of $\GW_g(Z,\mgp{T}_F)$ will be given by specialising the vertex, dilaton and edge terms in  \eqref{eq:GWgrescon}  to the leading order in $1/t_1$ of $\Delta_+$, $D_{+,l}$, and $E_{n,m}^{++}$ in \eqref{eq:cfasymp} and \eqref{eq:TLRcf}: we find 
\beq
\GW_g(Z,\mgp{T}_{\rm F})= \frac{1}{t_1^{2g-2}}\sum_{\Gamma \in \cG^{[0]}_{g,0}} \frac{1}{|\mathrm{Aut(\Gamma)}|} \prod_{e\in \cE_\Gamma} \mathrm{Cont}_{\rm cf}(e) 
\prod_{d\in \DD_\Gamma} 
\mathrm{Cont}_{\rm cf}(d)
\prod_{v\in \cV_\Gamma} \mathrm{Cont}_{\rm cf}(v) + \cO\l(\frac{1}{t_1^{2g-3}}\r)\,,
\label{eq:GWresconcf}
\eeq 
where, writing $e=\{h_1, h_2\}$,
\beq
\begin{split}
\mathrm{Cont}_{\rm cf}(d) & \coloneqq   [z^{\mathsf{a}(d)-1}]\l(1- \Upsilon\l(-\frac{z}{\epsilon_1},-\frac{z}{\epsilon_2}\r)\l(1-\frac{z}{\epsilon_1+\epsilon_2}\r)^{-1}\r)\,, \\
\mathrm{Cont}_{\rm cf}(e) & \coloneqq   [z^{\mathsf{a}(h_1)} w^{\mathsf{a}(h_2)}]\l( \frac{1- \Upsilon\l(-\frac{z}{\epsilon_1},-\frac{z}{\epsilon_2}\r)\Upsilon\l(-\frac{w}{\epsilon_1},-\frac{w}{\epsilon_2}\r)}{z+w}\r)\,, \\
\mathrm{Cont}_{\rm cf}(v) & \coloneqq  (-(\epsilon_1+\epsilon_2)^3\epsilon_1 \epsilon_2)^{\frac{2\mathsf{g}(v)-2+\mathsf{|H(v)|}}{2}} 
\int_{[\Mbar_{\mathsf{g}(v), |\HH(v)|}]}
\prod_{k=1}^{|\HH(v)|} 
\psi_k^{\mathsf{a}(h_k)}\,.
\end{split}
\label{eq:GWresconcfcont}
\eeq
The graph factors in \eqref{eq:GWresconcfcont} are almost identical to the ones of the special cubic Hodge CohFT \cite{Okounkov:2003aok}, which describes the equivariant GW theory of $\bbC^3$ with respect to a Calabi--Yau torus acting with weights $(\epsilon_1$, $\epsilon_2$, $-\epsilon_1-\epsilon_2$),
and for which the primary higher genus potential is 
the asymptotic expansion of the logarithm of the McMahon function. 
The sole (but crucial) difference here stems from the factor $(1-z/(\epsilon_1+\epsilon_2))^{-1}$ appearing in the expression of the dilaton leaf contribution, $\mathrm{Cont}_{\rm cf}(d)$.

\begin{prop}
Let $X=\mathrm{Tot}\l(\cO_{\bbP^1}^{\oplus 2}(-1)\r)$ and $\mgp{T}_{\rm F}$ be the fibrewise torus action with weights $\alpha_0=0$, $\alpha_{1,2}=\epsilon_{1,2}$. Then,\[
\mgp{GW}_g(X \times \bbC^2, \mgp{T}_{\rm F})
=\frac{\beta_{2g}(\epsilon_1, \epsilon_2)}{t_{1}^{2g-2}}\big(1+\cO(t_1)\big)  \,,
\]
where $t_1 =\log \cQ$. In particular, with notation as in \eqref{eq:GWresconcfcont}, we have
    \beq 
    \sum_{\Gamma \in \cG^{[0]}_{g,0}} \frac{1}{|\mathrm{Aut(\Gamma)}|} \prod_{e\in \cE_\Gamma} \mathrm{Cont}_{\rm cf}(e)\prod_{d\in \DD_\Gamma} 
\mathrm{Cont}_{\rm cf}(d) \prod_{v\in \cV_\Gamma}  \mathrm{Cont}_{\rm cf}(v) =
\beta_{2g}(\epsilon_1, \epsilon_2)\,.
\label{eq:cfgraphsum} 
    \eeq 
    \label{prop:cfgraphsum}
\end{prop}
\begin{proof}
From \Cref{prop:ResConifold} and \eqref{eq:betapol}, we have
\[
\mgp{GW}_g(X \times \bbC^2, \mgp{T}_{\rm F}) = \sum_{d=1}^\infty  \frac{d^{2g-3} \beta_{2g}(\epsilon_1, \epsilon_2)}{(2g-3)!}\re^{d t_1} = \frac{\beta_{2g}(\epsilon_1, \epsilon_2)}{(2g-3)!}\mathrm{Li}_{3-2g}(\re^{t_1})\,. 
\]
The first part of the claim is elementary since, for $n>1$, \[
\mathrm{Li}_{1-n}(\re^{t_1}) = \de_{t_1}^{n} \mathrm{Li}_{1}(\re^{t_1})= - \de_{t_1}^{n} \log(1-\re^{t_1}) = (n-1)! (-{t_1})^{-n} + \cO({t_1}^{1-n})\,.
\]
The stable graph summation formula \eqref{eq:cfgraphsum} in the second part then follows from \eqref{eq:GWresconcf} and \eqref{eq:GWresconcfcont}.
\end{proof}

\end{appendix}


\bibliography{refs_merged}

\newcommand{\noop}[1]{}
\begin{bibdiv}
\begin{biblist}

\bib{ACFW13:ExpDegenPairs}{article}{
      author={Abramovich, D.},
      author={Cadman, C.},
      author={Fantechi, B.},
      author={Wise, J.},
       title={Expanded degenerations and pairs},
        date={2013},
        ISSN={0092-7872},
     journal={Commun. Algebra},
      volume={41},
      number={6},
       pages={2346\ndash 2386},
}

\bib{Aganagic:2006wq}{article}{
      author={Aganagic, M.},
      author={Bouchard, V.},
      author={Klemm, A.},
       title={{Topological strings and (almost) modular forms}},
        date={2008},
     journal={Commun. Math. Phys.},
      volume={277},
       pages={771\ndash 819},
}

\bib{Aganagic:2011mi}{article}{
      author={Aganagic, M.},
      author={Cheng, M. C.~N.},
      author={Dijkgraaf, R.},
      author={Krefl, D.},
      author={Vafa, C.},
       title={{Quantum geometry of refined topological strings}},
        date={2012},
     journal={J. High Energy Phys.},
      volume={11},
       pages={019},
}

\bib{Aganagic:2011sg}{article}{
      author={Aganagic, M.},
      author={Shakirov, S.},
       title={{Knot homology and refined Chern--Simons index}},
        date={2015},
     journal={Commun. Math. Phys.},
      volume={333},
      number={1},
       pages={187\ndash 228},
}

\bib{Alvarez-Gaume:1983uye}{article}{
      author={Alvarez-Gaume, L.},
      author={Freedman, D.~Z.},
       title={{Potentials for the supersymmetric nonlinear sigma nodel}},
        date={1983},
     journal={Commun. Math. Phys.},
      volume={91},
       pages={87},
}

\bib{AFHNZ13:WorldsheetRefTopStr}{article}{
      author={Antoniadis, I.},
      author={Florakis, I.},
      author={Hohenegger, S.},
      author={Narain, K.~S.},
      author={Zein~Assi, A.},
       title={Worldsheet realization of the refined topological string},
        date={2013},
        ISSN={0550-3213},
     journal={Nucl. Phys., B},
      volume={875},
      number={1},
       pages={101\ndash 133},
}

\bib{Arb21:KthDT}{article}{
      author={Arbesfeld, N.},
       title={{{\(K\)}}-theoretic {Donaldson}--{Thomas} theory and the
  {Hilbert} scheme of points on a surface},
    language={English},
        date={2021},
        ISSN={2313-1691},
     journal={Algebr. Geom.},
      volume={8},
      number={5},
       pages={587\ndash 625},
}

\bib{BNTY23:LocOrbSNC}{article}{
      author={Battistella, L.},
      author={Nabijou, N.},
      author={Tseng, H.-H.},
      author={You, F.},
       title={{The local-orbifold correspondence for simple normal crossings
  pairs}},
        date={2023},
        ISSN={1474-7480},
     journal={J. Inst. Math. Jussieu},
      volume={22},
      number={5},
       pages={2515\ndash 2531},
         url={https://doi.org/10.1017/S1474748022000172},
      review={\MR{4624970}},
}

\bib{Bershadsky:1993cx}{article}{
      author={Bershadsky, M.},
      author={Cecotti, S.},
      author={Ooguri, H.},
      author={Vafa, C.},
       title={{Kodaira--Spencer theory of gravity and exact results for quantum
  string amplitudes}},
        date={1994},
     journal={Commun. Math. Phys.},
      volume={165},
       pages={311\ndash 428},
}

\bib{Bouchard:2007ys}{article}{
      author={Bouchard, V.},
      author={Klemm, A.},
      author={Marino, M.},
      author={Pasquetti, S.},
       title={{Remodeling the B-model}},
        date={2009},
     journal={Commun. Math. Phys.},
      volume={287},
       pages={117\ndash 178},
}

\bib{Bou20:QuTropVert}{article}{
      author={Bousseau, P.},
       title={The quantum tropical vertex},
        date={2020},
        ISSN={1465-3060},
     journal={Geom. Topol.},
      volume={24},
      number={3},
       pages={1297\ndash 1379},
         url={https://doi.org/10.2140/gt.2020.24.1297},
}

\bib{Bou19:Takahashi}{article}{
      author={Bousseau, P.},
       title={A proof of {N}. {T}akahashi's conjecture for {$(\mathbb{P}^2,E)$}
  and a refined sheaves/{G}romov--{W}itten correspondence},
        date={2023},
        ISSN={0012-7094,1547-7398},
     journal={Duke Math. J.},
      volume={172},
      number={15},
       pages={2895\ndash 2955},
         url={https://doi.org/10.1215/00127094-2022-0095},
      review={\MR{4675044}},
}

\bib{BFGW21:HAE}{article}{
      author={Bousseau, P.},
      author={Fan, H.},
      author={Guo, S.},
      author={Wu, L.},
       title={Holomorphic anomaly equation for {$(\mathbb{P}^2,E)$} and the
  {N}ekrasov--{S}hatashvili limit of local {$\mathbb{P}^2$}},
        date={2021},
     journal={Forum Math. Pi},
      volume={9},
       pages={Paper No. e3, 57},
         url={https://doi.org/10.1017/fmp.2021.3},
      review={\MR{4254953}},
}

\bib{Brini:2014fea}{article}{
      author={Brini, A.},
      author={Cavalieri, R.},
       title={{Crepant resolutions and open strings II}},
        date={2018},
     journal={\'Epijournal G\'eom. Alg\'ebrique},
      volume={2},
       pages={4580},
}

\bib{Brini:2013zsa}{article}{
      author={Brini, A.},
      author={Cavalieri, R.},
      author={Ross, D.},
       title={{Crepant resolutions and open strings}},
        date={2019},
     journal={J. Reine Angew. Math.},
      volume={755},
       pages={191\ndash 245},
}

\bib{Brini:2010fc}{article}{
      author={Brini, A.},
      author={Marino, M.},
      author={Stevan, S.},
       title={{The uses of the refined matrix model recursion}},
        date={2011},
     journal={J.Math.Phys.},
      volume={52},
       pages={052305},
}

\bib{BP05:CurvesCY3TQFT}{article}{
      author={Bryan, J.},
      author={Pandharipande, R.},
       title={Curves in {C}alabi--{Y}au threefolds and topological quantum
  field theory},
        date={2005},
        ISSN={0012-7094,1547-7398},
     journal={Duke Math. J.},
      volume={126},
      number={2},
       pages={369\ndash 396},
         url={https://doi.org/10.1215/S0012-7094-04-12626-0},
      review={\MR{2115262}},
}

\bib{BP08:GWLocCurves}{article}{
      author={Bryan, J.},
      author={Pandharipande, R.},
       title={The local {G}romov--{W}itten theory of curves},
        date={2008},
        ISSN={0894-0347,1088-6834},
     journal={J. Amer. Math. Soc.},
      volume={21},
      number={1},
       pages={101\ndash 136},
         url={https://doi.org/10.1090/S0894-0347-06-00545-5},
        note={With an appendix by Bryan, C. Faber, A. Okounkov and
  Pandharipande},
      review={\MR{2350052}},
}

\bib{CKYZ99:locMS}{article}{
      author={Chiang, T.M.},
      author={Klemm, A.},
      author={Yau, S-T.},
      author={Zaslow, E.},
       title={{Local mirror symmetry: Calculations and interpretations}},
        date={1999},
     journal={Adv. Theor. Math. Phys.},
      volume={3},
       pages={495\ndash 565},
}

\bib{CKK14:refBPS}{article}{
      author={Choi, J.},
      author={Katz, S.},
      author={Klemm, A.},
       title={The refined {BPS} index from stable pair invariants},
        date={2014},
        ISSN={0010-3616,1432-0916},
     journal={Commun. Math. Phys.},
      volume={328},
      number={3},
       pages={903\ndash 954},
         url={https://doi.org/10.1007/s00220-014-1978-0},
      review={\MR{3201216}},
}

\bib{MR2486673}{article}{
      author={Coates, T.},
       title={On the crepant resolution conjecture in the local case},
        date={2009},
        ISSN={0010-3616},
     journal={Commun. Math. Phys.},
      volume={287},
      number={3},
       pages={1071\ndash 1108},
         url={http://dx.doi.org/10.1007/s00220-008-0715-y},
      review={\MR{2486673 (2010j:14098)}},
}

\bib{MR2510741}{article}{
      author={Coates, T.},
      author={Corti, A.},
      author={Iritani, H.},
      author={Tseng, H-H},
       title={Computing genus-zero twisted {G}romov--{W}itten invariants},
        date={2009},
        ISSN={0012-7094},
     journal={Duke Math. J.},
      volume={147},
      number={3},
       pages={377\ndash 438},
         url={http://dx.doi.org/10.1215/00127094-2009-015},
      review={\MR{MR2510741}},
}

\bib{MR4047552}{article}{
      author={Coates, T.},
      author={Corti, A.},
      author={Iritani, H.},
      author={Tseng, H-H.},
       title={Hodge-theoretic mirror symmetry for toric stacks},
        date={2020},
        ISSN={0022-040X,1945-743X},
     journal={J. Differential Geom.},
      volume={114},
      number={1},
       pages={41\ndash 115},
         url={https://doi.org/10.4310/jdg/1577502022},
      review={\MR{4047552}},
}

\bib{CG07:QuantumRR}{article}{
      author={Coates, T.},
      author={Givental, A.},
       title={Quantum {R}iemann--{R}och, {L}efschetz and {S}erre},
        date={2007},
        ISSN={0003-486X,1939-8980},
     journal={Ann. of Math. (2)},
      volume={165},
      number={1},
       pages={15\ndash 53},
         url={https://doi.org/10.4007/annals.2007.165.15},
      review={\MR{2276766}},
}

\bib{Coates:2012vs}{article}{
      author={Coates, T.},
      author={Iritani, H.},
       title={{On the convergence of Gromov--Witten potentials and Givental's
  formula}},
        date={2015},
        ISSN={0026-2285,1945-2365},
     journal={Michigan Math. J.},
      volume={64},
      number={3},
       pages={587\ndash 631},
}

\bib{Coates:2018hms}{article}{
      author={Coates, T.},
      author={Iritani, H.},
       title={{Gromov--Witten invariants of local $\mathbb{P}^2$ and modular
  forms}},
        date={2021},
     journal={J. Math. (Kyoto)},
      volume={61},
       pages={543\ndash 706},
}

\bib{CMS24:QuasiMapDegen}{article}{
      author={Cobos~Rabano, A},
      author={Manolache, C.},
      author={Shafi, Q.},
       title={The local/logarithmic correspondence and the degeneration formula
  for quasimaps},
        date={2024},
      eprint={2404.01381},
}

\bib{Costello:2012cy}{article}{
      author={Costello, K.~J.},
      author={Li, S.},
       title={{Quantum BCOV theory on Calabi--Yau manifolds and the higher
  genus B-model}},
        date={20121},
      eprint={1201.4501},
}

\bib{dCMM05:HodgeThAlgMaps}{article}{
      author={de~Cataldo, M.~A.},
      author={Migliorini, L.},
       title={The {Hodge} theory of algebraic maps},
        date={2005},
        ISSN={0012-9593},
     journal={Ann. Sci. {\'E}c. Norm. Sup{\'e}r. (4)},
      volume={38},
      number={5},
       pages={693\ndash 750},
}

\bib{dZNPZ22:PlayingMtheory}{article}{
      author={Del~Zotto, M.},
      author={Nekrasov, N.},
      author={Piazzalunga, N.},
      author={Zabzine, M.},
       title={Playing with the index of {M}-theory},
        date={2022},
        ISSN={0010-3616},
     journal={Commun. Math. Phys.},
      volume={396},
      number={2},
       pages={817\ndash 865},
}

\bib{Dijkgraaf:2009pc}{article}{
      author={Dijkgraaf, R.},
      author={Vafa, C.},
       title={Toda theories, matrix models, topological strings, and {$N=2$}
  gauge systems},
        date={20099},
      eprint={0909.2453},
}

\bib{Dijkgraaf:1990qw}{incollection}{
      author={Dijkgraaf, R.},
      author={Verlinde, H.},
      author={Verlinde, E.},
       title={Notes on topological string theory and {$2$}d quantum gravity},
        date={1991},
   booktitle={String theory and quantum gravity ({T}rieste, 1990)},
   publisher={World Sci. Publ., River Edge, NJ},
       pages={91\ndash 156},
      review={\MR{1164017}},
}

\bib{DIW21:GVfinite}{article}{
      author={Doan, A.},
      author={Ionel, E-N.},
      author={Walpuski, T.},
       title={The {G}opakumar--{V}afa finiteness conjecture},
        date={2021},
      eprint={2103.08221},
}

\bib{Dunin-Barkowski:2012vii}{article}{
      author={Dunin-Barkowski, P.},
      author={Shadrin, S.},
      author={Spitz, L.},
       title={{Givental graphs and inversion symmetry}},
        date={2013},
     journal={Lett. Math. Phys.},
      volume={103},
       pages={533\ndash 557},
}

\bib{EG98:equivIntTh}{article}{
      author={Edidin, D.},
      author={Graham, W.},
       title={Equivariant intersection theory ({With} an appendix by {Angelo}
  {Vistoli}: {The} {Chow} ring of {{\({\mathcal M}_2\)}})},
    language={English},
        date={1998},
        ISSN={0020-9910},
     journal={Invent. Math.},
      volume={131},
      number={3},
       pages={595\ndash 644},
}

\bib{EG00:equivRR}{article}{
      author={Edidin, D.},
      author={Graham, W.},
       title={Riemann--{R}och for equivariant {C}how groups},
        date={2000},
        ISSN={0012-7094,1547-7398},
     journal={Duke Math. J.},
      volume={102},
      number={3},
       pages={567\ndash 594},
         url={https://doi.org/10.1215/S0012-7094-00-10239-6},
      review={\MR{1756110}},
}

\bib{Eynard:2011kk}{article}{
      author={Eynard, B.},
       title={{Intersection numbers of spectral curves}},
        date={2011},
      eprint={1104.0176},
}

\bib{Eynard:2011ga}{article}{
      author={Eynard, B.},
       title={{Invariants of spectral curves and intersection theory of moduli
  spaces of complex curves}},
        date={2014},
     journal={Commun. Num. Theor. Phys.},
      volume={08},
       pages={541\ndash 588},
}

\bib{Eynard:2011vs}{article}{
      author={Eynard, B.},
      author={Kozcaz, C.},
       title={{Mirror of the refined topological vertex from a matrix model}},
        date={20117},
      eprint={1107.5181},
}

\bib{Eynard:2007kz}{article}{
      author={Eynard, B.},
      author={Orantin, N.},
       title={{Invariants of algebraic curves and topological expansion}},
        date={2007},
     journal={Commun. Number Theory Phys.},
      volume={1},
      number={2},
       pages={347\ndash 452},
}

\bib{MR1728879}{article}{
      author={Faber, C.},
      author={Pandharipande, R.},
       title={Hodge integrals and {G}romov--{W}itten theory},
        date={2000},
        ISSN={0020-9910},
     journal={Invent. Math.},
      volume={139},
      number={1},
       pages={173\ndash 199},
         url={http://dx.doi.org/10.1007/s002229900028},
      review={\MR{1728879}},
}

\bib{Fang:2016svw}{article}{
      author={Fang, B.},
      author={Liu, C-C.~M.},
      author={Zong, Z.},
       title={{On the remodeling conjecture for toric {C}alabi--{Y}au
  3-orbifolds}},
        date={2020},
     journal={J. Am. Math. Soc.},
      volume={33},
      number={1},
       pages={135\ndash 222},
}

\bib{FG10:virtRR}{article}{
      author={Fantechi, B.},
      author={G\"{o}ttsche, L.},
       title={Riemann--{R}och theorems and elliptic genus for virtually smooth
  schemes},
        date={2010},
        ISSN={1465-3060,1364-0380},
     journal={Geom. Topol.},
      volume={14},
      number={1},
       pages={83\ndash 115},
         url={https://doi.org/10.2140/gt.2010.14.83},
      review={\MR{2578301}},
}

\bib{FMR21:HigherRankDT}{article}{
      author={Fasola, N.},
      author={Monavari, S.},
      author={Ricolfi, A.~T.},
       title={Higher rank {{\(K\)}}-theoretic {Donaldson}--{Thomas} theory of
  points},
        date={2021},
        ISSN={2050-5094},
     journal={Forum Math. Sigma},
      volume={9},
       pages={51},
        note={Id/No e15},
}

\bib{FT23:rankrrankoneDT}{article}{
      author={Feyzbakhsh, S.},
      author={Thomas, R.~P.},
       title={Rank {{\(r\)}} {DT} theory from rank 1},
    language={English},
        date={2023},
        ISSN={0894-0347},
     journal={J. Am. Math. Soc.},
      volume={36},
      number={3},
       pages={795\ndash 826},
}

\bib{MR4562996}{article}{
      author={Giacchetto, A.},
      author={Lewanski, D.},
      author={Norbury, P.},
       title={An intersection-theoretic proof of the {H}arer--{Z}agier
  formula},
        date={2023},
        ISSN={2313-1691,2214-2584},
     journal={Algebr. Geom.},
      volume={10},
      number={2},
       pages={130\ndash 147},
      review={\MR{4562996}},
}

\bib{Givental:1997hm}{incollection}{
      author={Givental, A.},
       title={Elliptic {G}romov--{W}itten invariants and the generalized mirror
  conjecture},
        date={1998},
   booktitle={Integrable systems and algebraic geometry},
   publisher={World Sci. Publ., River Edge, NJ},
       pages={107\ndash 155},
      review={\MR{1672116}},
}

\bib{MR1653024}{article}{
      author={Givental, A.},
       title={A mirror theorem for toric complete intersections},
        date={1998},
     journal={Progr. Math.},
      volume={160},
       pages={141\ndash 175},
      review={\MR{1653024 (2000a:14063)}},
}

\bib{MR1901075}{article}{
      author={Givental, A.},
       title={Gromov--{W}itten invariants and quantization of quadratic
  {H}amiltonians},
        date={2001},
        ISSN={1609-3321},
     journal={Mosc. Math. J.},
      volume={1},
      number={4},
       pages={551\ndash 568, 645},
        note={Dedicated to the memory of I. G. Petrovskii on the occasion of
  his 100th anniversary},
      review={\MR{MR1901075 (2003j:53138)}},
}

\bib{MR1866444}{article}{
      author={Givental, A.},
       title={Semisimple {F}robenius structures at higher genus},
        date={2001},
        ISSN={1073-7928},
     journal={Internat. Math. Res. Notices},
      number={23},
       pages={1265\ndash 1286},
         url={https://doi.org/10.1155/S1073792801000605},
      review={\MR{1866444}},
}

\bib{GV98:MthTopStrI}{article}{
      author={Gopakumar, R.},
      author={Vafa, C.},
       title={M-theory and topological strings--{I}},
        date={1998},
      eprint={hep-th/9809187},
}

\bib{GV98:MthTopStrII}{article}{
      author={Gopakumar, R.},
      author={Vafa, C.},
       title={M-theory and topological strings--{II}},
        date={1998},
      eprint={hep-th/9812127},
}

\bib{GP97:virtloc}{article}{
      author={Graber, T.},
      author={Pandharipande, R.},
       title={Localization of virtual classes},
        date={1999},
        ISSN={0020-9910,1432-1297},
     journal={Invent. Math.},
      volume={135},
      number={2},
       pages={487\ndash 518},
         url={https://doi.org/10.1007/s002220050293},
      review={\MR{1666787}},
}

\bib{Grimm:2007tm}{article}{
      author={Grimm, T.~W.},
      author={Klemm, A.},
      author={Marino, M.},
      author={Weiss, M.},
       title={{Direct integration of the topological string}},
        date={2007},
     journal={J. High Energy Phys.},
      volume={08},
       pages={058},
}

\bib{HKR:HAE}{article}{
      author={Haghighat, B.},
      author={Klemm, A.},
      author={Rauch, M.},
       title={Integrability of the holomorphic anomaly equations},
        date={2008},
        ISSN={1126-6708},
     journal={J. High Energy Phys.},
      volume={2008},
      number={10},
       pages={37},
        note={Id/No 097},
}

\bib{MR848681}{article}{
      author={Harer, J.},
      author={Zagier, D.},
       title={The {E}uler characteristic of the moduli space of curves},
        date={1986},
        ISSN={0020-9910},
     journal={Invent. Math.},
      volume={85},
      number={3},
       pages={457\ndash 485},
         url={http://dx.doi.org/10.1007/BF01390325},
      review={\MR{848681}},
}

\bib{Hori:2000ck}{article}{
      author={Hori, K.},
      author={Iqbal, A.},
      author={Vafa, C.},
       title={{D-branes and mirror symmetry}},
        date={2000},
      eprint={hep-th/0005247},
}

\bib{Hori:2000kt}{article}{
      author={Hori, K.},
      author={Vafa, C.},
       title={{Mirror symmetry}},
        date={2000},
      eprint={hep-th/0002222},
}

\bib{HKK13:OmegaBmodelRigidN2}{article}{
      author={Huang, M-x.},
      author={Kashani-Poor, A-K.},
      author={Klemm, A.},
       title={The {{\(\Omega \)}}-deformed {{\(B\)}}-model for rigid
  {{\(N=2\)}} theories},
        date={2013},
        ISSN={1424-0637},
     journal={Ann. Henri Poincar{\'e}},
      volume={14},
      number={3},
       pages={425\ndash 497},
}

\bib{HK10:OmegaBG}{article}{
      author={Huang, M-x.},
      author={Klemm, A.},
       title={Direct integration for general {$\Omega$} backgrounds},
        date={2012},
        ISSN={1095-0761,1095-0753},
     journal={Adv. Theor. Math. Phys.},
      volume={16},
      number={3},
       pages={805\ndash 849},
      review={\MR{3024275}},
}

\bib{IP18:GV}{article}{
      author={Ionel, E-N.},
      author={Parker, T.~H.},
       title={The {G}opakumar--{V}afa formula for symplectic manifolds},
        date={2018},
        ISSN={0003-486X,1939-8980},
     journal={Ann. of Math. (2)},
      volume={187},
      number={1},
       pages={1\ndash 64},
         url={https://doi.org/10.4007/annals.2018.187.1.1},
      review={\MR{3739228}},
}

\bib{IKV09:RefVert}{article}{
      author={Iqbal, A.},
      author={Kozcaz, C.},
      author={Vafa, C.},
       title={The refined topological vertex},
        date={2009},
        ISSN={1126-6708,1029-8479},
     journal={J. High Energy Phys.},
      number={10},
       pages={069, 58},
         url={https://doi.org/10.1088/1126-6708/2009/10/069},
      review={\MR{2607441}},
}

\bib{MR1950944}{article}{
      author={Jarvis, T.},
      author={Kimura, T.},
       title={Orbifold quantum cohomology of the classifying space of a finite
  group},
        date={2002},
     journal={{\rm in} Orbifolds in mathematics and physics ({M}adison, {WI},
  2001), Contemp. Math.},
      volume={310},
       pages={123\ndash 134},
      review={\MR{MR1950944 (2004a:14056)}},
}

\bib{Kat08:GVgenus0}{article}{
      author={Katz, S.},
       title={Genus zero {Gopakumar}--{Vafa} invariants of contractible
  curves},
        date={2008},
        ISSN={0022-040X},
     journal={J. Differ. Geom.},
      volume={79},
      number={2},
       pages={185\ndash 195},
}

\bib{Kim10:LogStabMaps}{incollection}{
      author={Kim, B.},
       title={Logarithmic stable maps},
        date={2010},
   booktitle={New developments in algebraic geometry, integrable systems and
  mirror symmetry ({RIMS}, {K}yoto, 2008)},
      series={Adv. Stud. Pure Math.},
      volume={59},
   publisher={Math. Soc. Japan, Tokyo},
       pages={167\ndash 200},
         url={https://doi.org/10.2969/aspm/05910167},
      review={\MR{2683209}},
}

\bib{KLR21:degen}{article}{
      author={Kim, B.},
      author={Lho, H.},
      author={Ruddat, H.},
       title={The degeneration formula for stable log maps},
        date={2023},
        ISSN={0025-2611},
     journal={Manuscripta Math.},
      volume={170},
      number={1-2},
       pages={63\ndash 107},
         url={https://doi.org/10.1007/s00229-021-01361-z},
      review={\MR{4533479}},
}

\bib{Klemm:1999gm}{incollection}{
      author={Klemm, A.},
      author={Zaslow, E.},
       title={Local mirror symmetry at higher genus},
        date={2001},
   booktitle={Winter {S}chool on {M}irror {S}ymmetry, {V}ector {B}undles and
  {L}agrangian {S}ubmanifolds ({C}ambridge, {MA}, 1999)},
      series={AMS/IP Stud. Adv. Math.},
      volume={23},
   publisher={Amer. Math. Soc., Providence, RI},
       pages={183\ndash 207},
         url={https://doi.org/10.1090/amsip/023/07},
      review={\MR{1876069}},
}

\bib{KLMP24:CohOneDimSheavesP2}{article}{
      author={Kononov, Y.},
      author={Lim, W.},
      author={Moreira, M.},
      author={Pi, W.},
       title={Cohomology rings of the moduli of one-dimensional sheaves on the
  projective plane},
        date={2024},
      eprint={2403.06277},
}

\bib{KOO21:2legKthVertex}{article}{
      author={Kononov, Y.},
      author={Okounkov, A.},
      author={Osinenko, A.},
       title={The 2-leg vertex in {K}-theoretic {DT} theory},
        date={2021},
        ISSN={0010-3616,1432-0916},
     journal={Commun. Math. Phys.},
      volume={382},
      number={3},
       pages={1579\ndash 1599},
         url={https://doi.org/10.1007/s00220-021-03936-z},
      review={\MR{4232774}},
}

\bib{KPS23:refBPSP2}{article}{
      author={Kononov, Y.},
      author={Pi, W.},
      author={Shen, J.},
       title={{Perverse filtrations, Chern filtrations, and refined BPS
  invariants for local P2}},
        date={2023},
     journal={Adv. Math.},
      volume={433},
       pages={109294},
}

\bib{MR2851153}{article}{
      author={Kontsevich, M.},
      author={Soibelman, Y.},
       title={Cohomological {H}all algebra, exponential {H}odge structures and
  motivic {D}onaldson--{T}homas invariants},
        date={2011},
        ISSN={1931-4523},
     journal={Commun. Number Theory Phys.},
      volume={5},
      number={2},
       pages={231\ndash 352},
         url={https://doi.org/10.4310/CNTP.2011.v5.n2.a1},
      review={\MR{2851153}},
}

\bib{Kozcaz:2018ndf}{article}{
      author={Kozcaz, C.},
      author={Shakirov, S.},
      author={Vafa, C.},
      author={Yan, W.},
       title={{Refined topological branes}},
        date={2021},
     journal={Commun. Math. Phys.},
      volume={385},
      number={2},
       pages={937\ndash 961},
}

\bib{KW11:ExtHAE}{article}{
      author={Krefl, D.},
      author={Walcher, J.},
       title={Extended holomorphic anomaly in gauge theory},
        date={2011},
        ISSN={0377-9017},
     journal={Lett. Math. Phys.},
      volume={95},
      number={1},
       pages={67\ndash 88},
}

\bib{Labastida:1991yq}{article}{
      author={Labastida, J. M.~F.},
      author={Llatas, P.~M.},
       title={{Potentials for topological sigma models}},
        date={1991},
     journal={Phys. Lett. B},
      volume={271},
       pages={101\ndash 108},
}

\bib{lho2018note}{article}{
      author={Lho, H.},
       title={Note on equivariant ${I}$-function of local $\mathbb{P}^{n}$},
        date={2018},
      eprint={1807.05501},
}

\bib{LP18:HAE}{article}{
      author={Lho, H.},
      author={Pandharipande, R.},
       title={Stable quotients and the holomorphic anomaly equation},
        date={2018},
        ISSN={0001-8708,1090-2082},
     journal={Adv. Math.},
      volume={332},
       pages={349\ndash 402},
         url={https://doi.org/10.1016/j.aim.2018.05.020},
      review={\MR{3810256}},
}

\bib{MR3960668}{article}{
      author={Lho, H.},
      author={Pandharipande, R.},
       title={Crepant resolution and the holomorphic anomaly equation for
  {$[\Bbb C^3/\Bbb Z_3]$}. {W}ith an appendix by {T}.~{C}oates and
  {H}.~{I}ritani},
        date={2019},
        ISSN={0024-6115,1460-244X},
     journal={Proc. Lond. Math. Soc. (3)},
      volume={119},
      number={3},
       pages={781\ndash 813},
         url={https://doi.org/10.1112/plms.12248},
      review={\MR{3960668}},
}

\bib{Li01:RelStabMaps}{article}{
      author={Li, J.},
       title={Stable morphisms to singular schemes and relative stable
  morphisms},
        date={2001},
        ISSN={0022-040X,1945-743X},
     journal={J. Differential Geom.},
      volume={57},
      number={3},
       pages={509\ndash 578},
         url={http://projecteuclid.org/euclid.jdg/1090348132},
      review={\MR{1882667}},
}

\bib{Li02:Degen}{article}{
      author={Li, J.},
       title={A degeneration formula of {GW}-invariants},
        date={2002},
        ISSN={0022-040X,1945-743X},
     journal={J. Differential Geom.},
      volume={60},
      number={2},
       pages={199\ndash 293},
         url={http://projecteuclid.org/euclid.jdg/1090351102},
      review={\MR{1938113}},
}

\bib{Ma11:VirtPull}{article}{
      author={Manolache, C.},
       title={Virtual pull-backs},
        date={2012},
        ISSN={1056-3911,1534-7486},
     journal={J. Algebraic Geom.},
      volume={21},
      number={2},
       pages={201\ndash 245},
         url={https://doi.org/10.1090/S1056-3911-2011-00606-1},
      review={\MR{2877433}},
}

\bib{MOOP11:GWDTtoric}{article}{
      author={Maulik, D.},
      author={Oblomkov, A.},
      author={Okounkov, A.},
      author={Pandharipande, R.},
       title={Gromov--{W}itten/{D}onaldson--{T}homas correspondence for toric
  3-folds},
        date={2011},
        ISSN={0020-9910,1432-1297},
     journal={Invent. Math.},
      volume={186},
      number={2},
       pages={435\ndash 479},
         url={https://doi.org/10.1007/s00222-011-0322-y},
      review={\MR{2845622}},
}

\bib{MPT10:K3modularForms}{article}{
      author={Maulik, D.},
      author={Pandharipande, R.},
      author={Thomas, R.~P.},
       title={Curves on {{\(K3\)}} surfaces and modular forms},
        date={2010},
        ISSN={1753-8416},
     journal={J. Topol.},
      volume={3},
      number={4},
       pages={937\ndash 996},
}

\bib{MS20:CohomChiIndep}{article}{
      author={Maulik, D.},
      author={Shen, J.},
       title={Cohomological {$\chi$}-independence for moduli of one-dimensional
  sheaves and moduli of {H}iggs bundles},
        date={2023},
        ISSN={1465-3060,1364-0380},
     journal={Geom. Topol.},
      volume={27},
      number={4},
       pages={1539\ndash 1586},
         url={https://doi.org/10.2140/gt.2023.27.1539},
      review={\MR{4602420}},
}

\bib{MT18:GV}{article}{
      author={Maulik, D.},
      author={Toda, Y.},
       title={Gopakumar--{V}afa invariants via vanishing cycles},
        date={2018},
        ISSN={0020-9910,1432-1297},
     journal={Invent. Math.},
      volume={213},
      number={3},
       pages={1017\ndash 1097},
         url={https://doi.org/10.1007/s00222-018-0800-6},
      review={\MR{3842061}},
}

\bib{Mei15:DTinvIC}{article}{
      author={Meinhardt, S.},
       title={Donaldson--{T}homas invariants vs. intersection cohomology for
  categories of homological dimension one},
        date={2015},
      eprint={1512.03343},
}

\bib{MMNS12:motivicResCon}{article}{
      author={Morrison, A.},
      author={Mozgovoy, S.},
      author={Nagao, K.},
      author={Szendroi, B.},
       title={Motivic {D}onaldson--{T}homas invariants of the conifold and the
  refined topological vertex},
        date={2012},
        ISSN={0001-8708,1090-2082},
     journal={Adv. Math.},
      volume={230},
      number={4-6},
       pages={2065\ndash 2093},
         url={https://doi.org/10.1016/j.aim.2012.03.030},
      review={\MR{2927365}},
}

\bib{Mu83}{incollection}{
      author={Mumford, D.},
       title={Towards an enumerative geometry of the moduli space of curves},
        date={1983},
   booktitle={Arithmetic and geometry, {V}ol. {II}},
      series={Progr. Math.},
      volume={36},
   publisher={Birkh\"{a}user Boston, Boston, MA},
       pages={271\ndash 328},
      review={\MR{717614}},
}

\bib{NO12:WorldsheetOmegaBg}{article}{
      author={Nakayama, Y.},
      author={Ooguri, H.},
       title={Comments on worldsheet description of the {Omega} background},
        date={2012},
        ISSN={0550-3213},
     journal={Nucl. Phys., B},
      volume={856},
      number={2},
       pages={342\ndash 359},
         url={authors.library.caltech.edu/29463/},
}

\bib{Nek03:SWprePInstCounting}{article}{
      author={Nekrasov, N.},
       title={Seiberg-{Witten} prepotential from instanton counting},
        date={2003},
        ISSN={1095-0761},
     journal={Adv. Theor. Math. Phys.},
      volume={7},
      number={5},
       pages={831\ndash 864},
}

\bib{Nek05:Ztheory}{article}{
      author={Nekrasov, N.},
       title={{Z-theory: Chasing M-F-Theory}},
        date={2005},
     journal={Comptes Rendus Physique},
      volume={6},
       pages={261\ndash 269},
}

\bib{MR2491283}{article}{
      author={Nekrasov, N.},
       title={Instanton partition functions and {M}-theory},
        date={2009},
        ISSN={0289-2316,1861-3624},
     journal={Jpn. J. Math.},
      volume={4},
      number={1},
       pages={63\ndash 93},
         url={https://doi.org/10.1007/s11537-009-0853-9},
      review={\MR{2491283}},
}

\bib{Nekrasov:2003rj}{incollection}{
      author={Nekrasov, N.},
      author={Okounkov, A.},
       title={Seiberg-{W}itten theory and random partitions},
        date={2006},
   booktitle={The unity of mathematics},
      series={Progr. Math.},
      volume={244},
   publisher={Birkh\"auser Boston, Boston, MA},
       pages={525\ndash 596},
         url={https://doi.org/10.1007/0-8176-4467-9_15},
      review={\MR{2181816}},
}

\bib{NO14:membranes}{article}{
      author={Nekrasov, N.},
      author={Okounkov, A.},
       title={Membranes and sheaves},
        date={2016},
        ISSN={2313-1691,2214-2584},
     journal={Algebr. Geom.},
      volume={3},
      number={3},
       pages={320\ndash 369},
         url={https://doi.org/10.14231/AG-2016-015},
      review={\MR{3504535}},
}

\bib{Nesterov:2024bdv}{article}{
      author={Nesterov, D.},
       title={{Unramified Gromov--Witten and Gopakumar--Vafa invariants}},
        date={20245},
      eprint={2405.18398},
}

\bib{Ob24:EnriquesKth}{article}{
      author={Oberdieck, G.},
       title={Towards refined curve counting on the {E}nriques surface {I}:
  {K}-theoretic refinements},
        date={2024},
      eprint={2407.21318},
}

\bib{Ob24:EnriquesMotivic}{article}{
      author={Oberdieck, G.},
       title={Towards refined curve counting on the {E}nriques surface {II}:
  Motivic refinements},
        date={2024},
      eprint={2408.02616},
}

\bib{Ok15:lectures}{incollection}{
      author={Okounkov, A.},
       title={Lectures on {K}-theoretic computations in enumerative geometry},
        date={2017},
   booktitle={Geometry of moduli spaces and representation theory},
      series={IAS/Park City Math. Ser.},
      volume={24},
   publisher={Amer. Math. Soc., Providence, RI},
       pages={251\ndash 380},
         url={https://doi.org/10.1090/pcms/024},
      review={\MR{3752463}},
}

\bib{Okounkov:2003aok}{article}{
      author={Okounkov, A.},
      author={Pandharipande, R.},
       title={{Hodge integrals and invariants of the unknot}},
        date={2004},
     journal={Geom. Topol.},
      volume={8},
      number={2},
       pages={675\ndash 699},
}

\bib{OP10:localDTcurves}{article}{
      author={Okounkov, A.},
      author={Pandharipande, R.},
       title={The local {Donaldson}--{Thomas} theory of curves},
        date={2010},
        ISSN={1465-3060},
     journal={Geom. Topol.},
      volume={14},
      number={3},
       pages={1503\ndash 1567},
}

\bib{PP17:GWPTquintic}{article}{
      author={Pandharipande, R.},
      author={Pixton, A.},
       title={Gromov--{Witten}/pairs correspondence for the quintic 3-fold},
    language={English},
        date={2017},
        ISSN={0894-0347},
     journal={J. Am. Math. Soc.},
      volume={30},
      number={2},
       pages={389\ndash 449},
         url={hdl.handle.net/1721.1/115975},
}

\bib{PT09:StabPairs}{article}{
      author={Pandharipande, R.},
      author={Thomas, R.~P.},
       title={Curve counting via stable pairs in the derived category},
        date={2009},
        ISSN={0020-9910,1432-1297},
     journal={Invent. Math.},
      volume={178},
      number={2},
       pages={407\ndash 447},
         url={https://doi.org/10.1007/s00222-009-0203-9},
      review={\MR{2545686}},
}

\bib{Par23:UnivCurveCount}{article}{
      author={Pardon, J.},
       title={Universally counting curves in {C}alabi--{Y}au threefolds},
        date={2023},
      eprint={2308.02948},
}

\bib{RS21:EquivVirtGRR}{article}{
      author={Ravi, C.},
      author={Sreedhar, B.},
       title={Virtual equivariant {Grothendieck}--{Riemann}--{Roch} formula},
    language={English},
        date={2021},
        ISSN={1431-0635},
     journal={Doc. Math.},
      volume={26},
       pages={2061\ndash 2094},
}

\bib{MR0201688}{book}{
      author={Slater, L.~J.},
       title={Generalized hypergeometric functions},
   publisher={Cambridge University Press, Cambridge},
        date={1966},
      review={\MR{0201688}},
}

\bib{MR2528052}{article}{
      author={Spreafico, M.},
       title={On the {B}arnes double zeta and {G}amma functions},
        date={2009},
        ISSN={0022-314X,1096-1658},
     journal={J. Number Theory},
      volume={129},
      number={9},
       pages={2035\ndash 2063},
         url={https://doi.org/10.1016/j.jnt.2009.03.005},
      review={\MR{2528052}},
}

\bib{MR2917177}{article}{
      author={Teleman, C.},
       title={The structure of 2{D} semi-simple field theories},
        date={2012},
        ISSN={0020-9910},
     journal={Invent. Math.},
      volume={188},
      number={3},
       pages={525\ndash 588},
         url={https://doi.org/10.1007/s00222-011-0352-5},
      review={\MR{2917177}},
}

\bib{Tho20:refVW}{article}{
      author={Thomas, R.~P.},
       title={Equivariant {$K$}-theory and refined {V}afa--{W}itten
  invariants},
        date={2020},
        ISSN={0010-3616,1432-0916},
     journal={Commun. Math. Phys.},
      volume={378},
      number={2},
       pages={1451\ndash 1500},
         url={https://doi.org/10.1007/s00220-020-03821-1},
      review={\MR{4134951}},
}

\bib{Th24:RefinedK3}{article}{
      author={Thomas, R.~P.},
       title={Refined sheaf counting on local {K3} surfaces},
        date={2024},
      eprint={2403.12741},
         url={https://arxiv.org/abs/2403.12741},
}

\bib{MR2578300}{article}{
      author={Tseng, H-H.},
       title={Orbifold quantum {R}iemann--{R}och, {L}efschetz and {S}erre},
        date={2010},
        ISSN={1465-3060},
     journal={Geom. Topol.},
      volume={14},
      number={1},
       pages={1\ndash 81},
         url={http://dx.doi.org/10.2140/gt.2010.14.1},
      review={\MR{2578300 (2011c:14147)}},
}

\bib{TY20:HigherGenusRelOrb}{article}{
      author={Tseng, H-H.},
      author={You, F.},
       title={Higher genus relative and orbifold {G}romov--{W}itten
  invariants},
        date={2020},
        ISSN={1465-3060,1364-0380},
     journal={Geom. Topol.},
      volume={24},
      number={6},
       pages={2749\ndash 2779},
         url={https://doi.org/10.2140/gt.2020.24.2749},
      review={\MR{4194303}},
}

\bib{vGGR19}{article}{
      author={van Garrel, M.},
      author={Graber, T.},
      author={Ruddat, H.},
       title={Local {G}romov--{W}itten invariants are log invariants},
        date={2019},
        ISSN={0001-8708,1090-2082},
     journal={Adv. Math.},
      volume={350},
       pages={860\ndash 876},
         url={https://doi.org/10.1016/j.aim.2019.04.063},
      review={\MR{3948687}},
}

\bib{Wit88:TopSigMod}{article}{
      author={Witten, E.},
       title={Topological sigma models},
        date={1988},
        ISSN={0010-3616},
     journal={Commun. Math. Phys.},
      volume={118},
      number={3},
       pages={411\ndash 449},
}

\bib{Witten:1991zz}{incollection}{
      author={Witten, E.},
       title={Mirror manifolds and topological field theory},
        date={1992},
   booktitle={Essays on mirror manifolds},
   publisher={Int. Press, Hong Kong},
       pages={120\ndash 158},
      review={\MR{1191422}},
}

\bib{MR2454324}{incollection}{
      author={Zagier, D.},
      author={Zinger, A.},
       title={Some properties of hypergeometric series associated with mirror
  symmetry},
        date={2008},
   booktitle={Modular forms and string duality},
      series={Fields Inst. Commun.},
      volume={54},
   publisher={Amer. Math. Soc., Providence, RI},
       pages={163\ndash 177},
      review={\MR{2454324}},
}

\end{biblist}
\end{bibdiv}


\end{document}